\documentclass[letterpaper,11pt,oneside,reqno]{amsbook}
\usepackage{bm}
\usepackage{mathrsfs}
\usepackage{amsfonts,amsmath, amssymb,amsthm,amscd,stmaryrd,mathtools}
\usepackage[height=9.6in,width=5.95in]{geometry}
\usepackage[mpexclude,DIV13]{typearea}
\usepackage{verbatim}
\usepackage{hyperref}
\usepackage{graphicx}
\usepackage[latin1]{inputenc}
\usepackage{latexsym}
\usepackage{lscape}
\usepackage{epsfig}
\include{bibtex}

\usepackage{setspace}

\usepackage{parskip}
\newcommand{\e}[0]{\epsilon}

\newcommand{\PP}{\ensuremath{\mathbb{P}}}

\newcommand{\N}{\ensuremath{\mathbb{N}}}
\newcommand{\R}{\ensuremath{\mathbb{R}}}

\newcommand{\Z}{\ensuremath{\mathbb{Z}}}

\newcommand{\E}[0]{\mathbb{E}}

\newtheorem{theorem}{Theorem}[chapter]
\numberwithin{figure}{chapter}

\newtheorem{lemma}[theorem]{Lemma}
\newtheorem{proposition}[theorem]{Proposition}
\newtheorem{corollary}[theorem]{Corollary}

\theoremstyle{definition}

\theoremstyle{definition}

\theoremstyle{definition}
\newtheorem{definition}[theorem]{Definition}

\theoremstyle{definition}

\theoremstyle{definition}
\newtheorem{conjecture}[theorem]{Conjecture}

\newcommand{\intint}[1]{\llbracket 1,#1 \rrbracket}

\newcommand{\mc}{\mathcal}

\newcommand{\ceemack}{\nu}

\newcommand{\nwithzero}{\N}
\newcommand{\nwithoutzero}{\N_+}

\newcommand{\macell}{\ell}

\newcommand{\game}{{\rm Trail}}
\newcommand{\standard}{{\rm Standard}}

\newcommand{\centralratio}{{\rm CenRatio}}
\newcommand{\tis}{\mc{S}_0}

\newcommand\pgameplay[3]{\PP^{#3}_{#1,#2}}
\newcommand\egameplay[3]{\E^{#3}_{#1,#2}}

\newcommand{\macsee}{c}
\newcommand{\macbeta}{\beta}
\newcommand{\macomega}{\omega}
\newcommand{\macdee}{d}
\newcommand{\macphi}{\phi}
\newcommand{\macgsee}{\gamma}
\newcommand{\macgdee}{\delta}

\newcommand{\const}{h}
\newcommand{\genseq}{g}

\newcommand{\adef}{a^{\rm def}}
\newcommand{\bdef}{b^{\rm def}}
\newcommand{\mdef}{m^{\rm def}}
\newcommand{\ndef}{n^{\rm def}}
\newcommand{\asta}{a^{\rm st}}
\newcommand{\bst}{b^{\rm st}}
\newcommand{\mst}{m^{\rm st}}
\newcommand{\nst}{n^{\rm st}}

\newcommand{\solution}{\mathsf{StSol}}
\newcommand{\status}{\mathsf{Stat}}
\newcommand{\minamm}{\mathcal{M}}
\newcommand{\transmap}{\chi}
\newcommand{\PhiMac}{\Phi}

\newcommand{\abmnmac}{\textrm{ABMN}}
\newcommand{\mnmac}{\textrm{MN}}
\newcommand{\ammac}{\textrm{AM}}

\newcommand{\dabmnmac}{\textrm{dABMN}}

\newcommand{\abmnmacspace}{\textrm{ABMN} }
\newcommand{\mnmacspace}{\textrm{MN} }
\newcommand{\ammacspace}{\textrm{AM} }
\newcommand{\amacspace}{\textrm{A} }
\newcommand{\dabmnmacspace}{\textrm{dABMN} }

\newcommand{\cross}{W}
\newcommand{\crossnum}{N}
\newcommand{\totalcost}{\textrm{TotalCost}}

\newcommand{\rocket}{\textrm{Rocket}}
\newcommand{\boost}{\textrm{Boost}}
\newcommand{\drag}{\textrm{Drag}}
\newcommand{\lef}{\textrm{Left}}
\newcommand{\righ}{\textrm{Right}}

\newcommand{\ground}{\textrm{ground}}
\newcommand{\indhyp}{\textrm{IH}}
\newcommand{\access}{\mathsf{A}}

\newcommand{\sdev}{S_-^{\textrm{dev}}}
\newcommand{\sdevmaxine}{S_+^{\textrm{dev}}}

\def\fff#1{&{{\pageref{#1}}}\cr}
\def\hfff#1{\label{#1}}

\begin{document}

\title[The Trail of Lost Pennies]{On the Trail of Lost Pennies: \\
 Player-Funded Tug-Of-War on the Integers}

\author[A. Hammond]{Alan Hammond}
\address{A. Hammond\\
  Department of Mathematics and Statistics\\
 U.C. Berkeley \\
  899 Evans Hall \\
  Berkeley, CA, 94720-3840 \\
  U.S.A.}
  \email{alanmh@berkeley.edu}
  \thanks{The author is supported by NSF grant DMS-$2153359$. }
\keywords{Bidding games on graphs, Colonel Blotto games,  infinite games,  multi-turn games,  non-zero-sum two-player games, player-funded stake-governed games, random-turn games,   strategic move evaluation, tug-of-war, Tullock contests.}

\begin{abstract} 
We study a class of
 two-player, random-turn games
in which a valued resource must be spent during the lifetime of the game so that a favourable outcome is likely secured without wasteful expense.
The Trail of Lost Pennies is the game played on $\Z$. 
It is specified by a parameter $x \in (0,\infty)$. A counter begins at $0 \in \Z$ and lies at some integer at the start of any given $\N$-indexed turn. At this turn, Maxine and Mina offer non-negative real-valued stakes. If Maxine stakes $a$ and Mina $b$, the counter moves one unit to the right with probability $\tfrac{a}{a+b}$; otherwise, it moves to the left by one unit. If the counter location tends to minus infinity over the course of the game, Mina receives a terminal payment of one, and Maxine zero; if the counter tends to plus infinity, then these respective receipts are zero and~$x$. Thus the net receipt to a given player is $-A+B$, where $A$ is the sum of the stakes she offers during the game, and $B$ is her terminal receipt. The game was inspired by  random-turn tug-of-war analysed in the mathematical treatment~\cite{PSSW09} from 2009
but in fact closely resembles the original version of tug-of-war, introduced~\cite{HarrisVickers87} in the economics literature in 1987.
 We show that the game has surprising features. For the class of strategies that satisfy a natural time-invariance, Nash equilibria exist precisely when $x$ lies in the interval $[\lambda,\lambda^{-1}]$, for a certain $\lambda \in (0,1)$. We indicate that the value of $\lambda$ is remarkably close to one, proving that $\lambda \leq 0.999904$ and presenting clear numerical evidence that $\lambda \geq 1 - 10^{-4}$. For each $x \in [\lambda,\lambda^{-1}]$, we exhibit countably many Nash equilibria for the game, showing that each is roughly characterised by an integral {\em battlefield} index: when the counter is around this value, both players stake intensely, with rapid but asymmetric decay in stakes as the counter moves away from the value. 
Our results validate and quantify hypotheses, for strategic resource management and the effect on outcome of discouragement of the less incentivised player, which are similar to premises advanced in economics~\cite{HarrisVickers87,Konrad2012}, and which
  plausibly hold more generally among player-funded stake-governed random-turn games.
Alongside a companion treatment~\cite{HP2022} of games with allocated budgets, we thus offer a detailed mathematical treatment of an illustrative class of tug-of-war games governed by limited resources. We also review the separate developments of tug-of-war in economics and mathematics in the hope that mathematicians may direct further attention to tug-of-war games, in their original resource-allocation guise.
\end{abstract}

\medskip

\maketitle


\setcounter{tocdepth}{2}
\tableofcontents

\chapter{Introduction} 

In random-turn tug-of-war games,
two players push a counter on a gameboard, each winning the right to do so at a given turn according to the flip of a coin, until a terminal payment is made when the counter reaches a boundary set. 
Dating from Harris and Vickers' 1987 article~\cite{HarrisVickers87}, economists have studied such games when the bias of the coin is determined by each player's expenditure of a resource, investigating how players choose to allocate the resource over the lifetime of the game. Since Peres, Schramm, Sheffield and Wilson's 2009 contribution~\cite{PSSW09}, mathematicians specializing in probability and PDE have studied tug-of-war games with much simpler decision rules but whose geometric setting is richer: it is a fair coin that decides turn victor, the game played with small step sizes in a Euclidean domain, and a connection emerging to the infinity Laplacian operator. The two research strands, in economics and mathematics, have been vigorous, but remarkably disjoint. 

The principal aim of this monograph is to present a detailed mathematical inquiry into an infinite-turn stake-governed tug-of-war game, the Trail of Lost Pennies. We find a delicate and surprising structure of Nash equilibria in the game; as we will shortly explain, our comprehensive classification and description of the time-invariant (or Markov perfect) equilibria illustrates and elucidates central themes in the economics literature of tug-of-war. 
 We further aim to direct mathematicians' attention to stake-governed tug-of-war games and to advocate these games  as objects for further inquiry.


The first part of the introduction offers a rough conceptual overview of the motivations and form of the player-funded tug-of-war game that we study and of our main results; and we compare these conclusions to those for a class of allocated-budget tug-of-war games that have been treated in the companion article~\cite{HP2022}. In the second part, we trace the course of the two long, unwoven strands of tug-of-war research, in economics and mathematics.

\section{Rough guide: motivations, model, results}

The first of four subsections presents three broad premises for how resources may be allocated strategically in a random multi-turn game. The second conceptually presents tug-of-war games, first in the simpler setting of constant bias, and then with bias governed by stakes drawn from budgets. These budgets may be allocated by the bank or funded by the players themselves. Our principal game of study, the Trail of Lost Pennies, is of the latter type. In the third subsection, we explain several of our main conclusions about this game, finding it to exemplify two of the three broad premises. The remaining premise is instead illustrated by a class of allocated-budget games analysed in~\cite{HP2022}, as the fourth subsection explains.  

\subsection{Allocating resources strategically in a multi-turn game}\label{s.allocate}

When political opponents buy advertising before an election or rival firms rent stalls in a market place, they allocate precious resources in an effort to win a favourable outcome, if need be at the expense of the opponent. Such competition may be modelled by multi-turn two-player games at each of whose turns each player allocates a resource that she values in an effort to improve her strategic position and raise the chance of a favourable outcome when the game ends. 
In seeking to play such a game effectively, a player must address the basic question: how much of the resource should I allocate for play at the impending turn? To find the answer, the player is led to evaluate the strategic importance of the next turn relative to such a measure for the accumulation of all subsequent turns during the likely lifetime of the game.
So broadly specified, this question may elicit a wide range of conceptual responses. We distinguish three:
\begin{enumerate}
\item
{\bf Compromise between the short and long terms.} In facing a choice of level of resource allocation for the upcoming turn, a player is tempted to offer a large amount, in the hope of securing an improved position at the turn. However, such expenditure  depletes the resource, making it more difficult to stake competitively later in the game. In this view, there is a tension between the short and the long term, with the two time-scales offering upward and downward pressure on resource allocation at a given turn. Presumably then under optimal play the present level of expenditure will maximize a unimodal function. 
\item {\bf An inch in incentive may be a mile for outcome.}
Suppose that the two players will receive different rewards for winning the game. The player with the greater incentive may be led to commit greater resources at any given stage of the game. If the difference in incentive is small, then the resulting increase in stage-win probability for the more committed player may be only slight. But if there are many stages in the game, then these advantages will accumulate to become pronounced, just as a random walk with a slight bias to the right will proceed linearly in that direction in the long term. Moreover, the thus disadvantaged player, in becoming aware of her poor prospects, will be discouraged from committing resources, leading her to stake very little, leaving the other player, even if only slightly incentivised, able to win the game without any real opposition.  
\item {\bf The battlefield. Cut your losses/Foot on gas.} In competing for dominance in some form---say, over a network or geographical space--- in a game taking many turns, both players may begin to recognise as gameplay advances that one has secured a certain advantage. Perhaps one player has won a key battle over a significant part of space, or is simply blessed with a favourable terrain. In such a circumstance, the player in the seemingly weaker position may hesitate to commit much of the valued resource at an upcoming turn. If he makes big stakes consistently, and  performs better at several forthcoming turns, he may for example become able to refight a lost battle, but he knows that his determined opponent may well commit again to win that battle; so his route to victory will likely entail high expenses in the resource, and he may say to himself ``{\em cut your losses}", 
accepting a bad terminal game result but consoled by minimal resource loss in the later game.
To the opposing player we might send a different message, namely {\em don't take your foot off the gas}: if the weakly positioned player will stake very little, his opponent hardly needs to stake intensely, but she should make sure to outspend her opponent significantly, because she wishes to convert a promising position into a terminal victory without setbacks. 
\end{enumerate}
In counterpart to the above headings, the three premises will be called the {\em Present-Future Compromise}; {\em ``Incentive Inch, Outcome Mile"}; and the {\em Battlefield Cyl Fog}.

The topic of strategic importance of positions in multi-turn games is of course natural and important: 
 ~\cite{KleinbergTardos,Azar,Kanoria} have explored such themes in the framework of iterative network-bargaining problems.
 To interpret and investigate the broad and perhaps vague premises that we have just outlined, a suitable class of games is clearly needed. 
Such a class is offered by stake-governed tug-of-war games, in which players win the right to move (a counter along the edges of a graph) according to the flip of a coin with a bias determined by their commitment of resources at the turn. 

\subsection{Tug-of-war, constant-bias and stake-governed}\label{s.tugofwar}

We present tug-of-war first in the simple case of constant bias; and then versions governed by stakes from budgets that are either allocated or funded by players. The Trail Of Lost Pennies, a player-funded game on the integers, is then introduced.

\medskip

\subsubsection{Constant-bias tug-of-war}

Let $G = (V,E)$ be a finite graph equipped with a set $B \subset V$ of boundary vertices on which a function $f:B \to \R$ is defined. Further suppose given a bias parameter $\lambda \in (0,\infty)$. The initial counter location $X_0$ of $\lambda$-biased tug-of-war is taken to be a given $v \in V$. At the $i$\textsuperscript{th} turn (where $i \geq 0$), a coin with heads' probability $\tfrac{\lambda}{1 + \lambda}$  is flipped. If it lands heads, Maxine moves the counter to a neighbour of its present location of her choosing, specifying $X_{i+1} \in V$ with $X_{i+1} \sim X_i$. If it lands tails, then Mina does so. The game ends at time $\tau = \min \{ i \in \N: X_i \in B \}$, when the counter first reaches the boundary, with a payment of $f(X_\tau)$ made from Mina to Maxine. The unbiased case $\lambda =1$ is considered in~\cite{PSSW09}, and the biased $\lambda \in (0,\infty) \setminus \{1 \}$ in~\cite{PPS10}; either way,  
 the value of tug-of-war---in essence, the mean terminal payment under optimal play---exists and equals $h(\lambda,v)$, where $h(\lambda,\cdot):V \to \R$ is the $\lambda$-biased infinity harmonic extension of the boundary data $f:B \to \R$; namely, the unique function satisfying the system of equations
\begin{equation}\label{e.lambdainftylap}
 h(\lambda,v) \, = \, \frac{\lambda}{1+\lambda} \,  \max_{u \sim v} h(\lambda,u)  \, + \,  \frac{1}{1+\lambda} \, \min_{u \sim v} h(\lambda,u) 
\end{equation}
for $v \in V \setminus B$ and $h(\lambda,\cdot) \big\vert_B = f$.


\medskip

\subsubsection{Allocated-budget stake-governed tug-of-war}\label{s.allocatedbudget}

Let the data $G=(V,E)$, $B \subset V$, $f:B \to \R$ and $\lambda \in (0,\infty)$ again be given. 
 The two players are now allocated a budget at the outset of the game: $\lambda$ for Maxine; and one for Mina. Take $X_0 = v \in V$ as above.
At the outset of each turn of the new game, each player retains some part of her initial budget. Each stakes some part of her remaining budget, as she sees fit. These staked funds are withdrawn from the players, so that the new budget of a player, for use at the next turn, is given by subtracting her stake from the budget's present value. If Maxine stakes $a$ and Mina $b$,
a biased coin is then flipped, whose probability of landing heads equals $\tfrac{a}{a+b}$. The turn victor moves the counter according to the rules of constant-bias tug-of-war, and the game ends just as it did originally, with a terminal payment from Mina to Maxine equal to the evaluation of the function $f$ at the terminal counter location in the boundary~$B$. Any remaining budget held by the players is taken from them.

\medskip

\subsubsection{Player-funded stake-governed tug-of-war}\label{s.playerfunded}

 In the allocated-budget stake game, the initially awarded budgets are an irredeemable resource whose sole role is to afford capacity  to the respective players to win moves throughout the lifetime of the game. Maxine and Mina's initial budgets are given finite quantities whose values are part of the game design. The finiteness of these values is what makes the resource precious.
 In the variant we now specify, a different means is used to ensure that the players value this resource: each player can spend without constraint, but she must spend {\em her own money}.
 Mina and Maxine are supposed to be wealthy, so each may stake as freely as she sees fit, but the sequence of stakes made by a player constitutes a running cost that must be offset against 
 the potential benefit that higher expenditure will bring the counter to a more favourable terminal location in the boundary set.

The game data now takes the form  $G=(V,E)$, $B \subset V$ and $f,g:B \to \R$.
Thus a second boundary function $g:B \to \R$ is introduced (and the parameter $\lambda$ is absent, reflecting the players' wealth).
 At the start of the $i$\textsuperscript{th} turn, Maxine stakes $a_i \in [0,\infty)$
and Mina $b_i \in [0,\infty)$. A coin is flipped whose heads' probability is $\tfrac{a_i}{a_i+b_i}$. As above, heads means Maxine wins the right to move the counter to an adjacent vertex; tails, and Mina does. The game ends, as above, when the counter reaches $B$ at time $\tau$. The terminal payments are $f(X_\tau)$ to Maxine, and $g(X_\tau)$ to Mina. 
 Thus the net receipt in the game equals
 \begin{equation}\label{e.netreceipt}
  \textrm{$f(X_\tau) - \sum_{i=0}^{\tau -1} a_{X_i}$ for Maxine, and $g(X_\tau) - \sum_{i=0}^{\tau -1} b_{X_i}$ for Mina} \, .
\end{equation}

In the special case that $g$ equals $-f$, the terminal payments coincide with those in the games already discussed. Unlike these games, the player-funded stake game is not zero-sum, even when $g=-f$;
since zero-sum tools will not be used to analyse the game, it is natural to generalize it with  separate terminal payment functions~$f$ and~$g$.

\medskip

\subsubsection{The Trail of Lost Pennies}

 The constant-bias, and allocated-budget stake-governed, tug-of-war games make little sense on certain infinite graphs. 
 Consider fair-coin tug-of-war on the integers~$\Z$ with nearest-neighbour edges, and a payment of one unit from Mina to Maxine if the counter tends to~$\infty$, and of minus one unit if instead it tends to $-\infty$. The players have no choices for strategy and the counter evolves as simple symmetric random walk, so no terminal payment is made (or perhaps a default rule stipulates the payment). And likewise the game on $\Z$ makes no real sense for the allocated-budget stake game: roughly put, since the game will require infinitely many turns, any positive expenditure of the globally finite budget of a given player is unjustified at any given turn; but if both players consistently stake nothing, then (at least if a symmetric rule is adopted for this circumstance) the counter will again evolve as simple random walk.
 
 It may seem that the self-funded stake game would be equally trivial when played on~$\Z$. Since this gameboard is an infinite transitive graph, there seems to be nothing to distinguish one move from any other in the infinite game, and this seems to make any positive stake by either player at any turn unjustified. This consideration in fact permits one in essence to disprove the uniqueness of a Nash equilibrium in the game, but it is quite false insofar as it suggests that the game is trivial. Indeed, 
 self-funded stake-governed tug-of-war has an intricate theory on $\Z$. In this article, we investigate the new game when the underlying graph is either $\Z$
 or a finite interval therein. We call the game on these graphs the Trail of Lost Pennies---either on $\Z$, or on a given finite integer interval. 
 
 The attentive reader will have noticed imprecisions in our specification of stake games. What happens if both players stake nothing at a given turn? Or if the counter fails to reach the boundary in a suitable sense? We defer resolving these points from the present conceptual overview to Chapter~\ref{c.speclost}, where the Trail of Lost Pennies will be precisely specified.
 
 \subsection{An informal summary of the main conclusions}
 Our principal focus will be on the Trail of Lost Pennies on $\Z$, a game that is necessarily of infinite duration. 
 Our main theorems show that this simple game has a remarkably rich and surprising structure of Nash equilibria, which witnesses two of the broad premises, ``Incentive Inch, Outcome Mile" and the Battlefield Cyl Fog.
 
 In Chapter~\ref{c.speclost}, we will 
  carefully state our principal conclusions about the game. In summary now of the game rules and our results about it, we may say that
 our game data is $\Z$ with nearest-neighbour edges; two boundary points, at $-\infty$ and $\infty$; and a parameter $x \in (0,\infty)$. 
  If the counter tends to~$\infty$ over the course of the game, then a
  terminal payment of one is made to Maxine, and of zero to Mina. If the counter tends to $-\infty$, Mina will receive $x$, and Maxine zero. We will prove that Nash equilibria exist in the game precisely when $x \in [\lambda,\lambda^{-1}]$ for a certain value $\lambda \in (0,1)$. 
 The quantity $\lambda$ is canonically associated to $\Z$ via a simple enough game, 
 making its sheer closeness to the value one remarkable: we will prove that $\lambda \leq 0.999904$ and present clear numerical evidence that $\lambda \geq 1 - 10^{-4}$.
 (On long finite integer intervals, the rise of $x$ roughly from $\lambda$ to $\lambda^{-1}$ is a transition from a phase of Maxine's dominance to one of Mina's via a brief phase with many equilibria.
 An incentive gain of $10^{-4}$ leads to a unit-order gain in outcome: the order of magnitude of the ratio of mile to inch is respected.)
    When $x \in [\lambda,\lambda^{-1}]$ is given, we will exhibit countably many Nash equilibria in the game. Each of these equilibria will be computed explicitly and will be found to witness the Battlefield Cyl Fog: each has an integral battlefield index around which both players stake at unit order; as the counter moves to the right from this battlefield, Mina cuts her losses and stakes minutely, while Maxine also rapidly reduces spending but at a much slower rate (so that her foot is, in a relative sense, on the gas); the roles are reversed to the left of the battlefield. 
 Note that the  countably many Nash equilibria  break  the translation symmetry of~$\Z$: the consideration mooted above, that this symmetry would render the game trivial, is thus invalidated, because it is the collection of equilibria that is translation invariant, even as each equilibrium is a non-trivial object, centred on an integral battlefield, and lacking translation invariance. 
 
 By finding the two premises to be exemplified, our study offers a detailed answer to a question concerning the strategic importance of positions and the optimal allocation of valued resources in multi-turn games. The comprehensive classification and description of Nash equilibria in this study may be viewed as a more general framework that develops certain of the economists' deductions made for specific equilibria in related games that enjoy special properties such as symmetry: for example, and as we will review later in the introduction, economists' investigations of tug-of-war have often addressed the premise of a site of intense competition with rapid and asymmetric decay of effort rates away from the site. 
 The essence of the conclusions that we offer may well be valid more generally:
  as we will discuss further in Section~\ref{s.beyond},  player-funded stake-governed games  may on many graphs
have a space of multiple Nash equilibria that exist only in a narrow range of values of the relative terminal reward and  whose elements are roughly characterized by a small `battlefield' zone that roughly divides the gameboard into two regions in which the subsequent victory of one or other player is close to assured.

\subsection{Overview: allocated-budget stake games and the Present-Future Compromise}\label{s.overviewallocated}

The guiding principle governing allocated-budget stake games stands in contrast to those of the player-funded games that are the focus of this article. The difference is stark: a staccato burst of stakes that punctuates near peace in player-funded games; a much smoother flow of expenditure when fixed budgets are allocated.  In order to put the contrast into relief, we here give a brief overview of the sense in which~\cite{HP2022} shows that the allocated-budget games verify the Present-Future Compromise. Suppose that, in such a game, Maxine adopts a strategy that involves placing deterministic stakes, and that Mina learns this strategy. She may mimic Maxine's strategy, by staking the same proportion of her remaining budget at any given turn as Maxine does of hers. The relative fortune of the players is thus maintained at its initial value, and this permits Mina in effect to reduce the game to the $\lambda$-biased tug-of-war. And Maxine can do likewise if the roles are reversed. Under the premise that players stake deterministic amounts, we find that, under optimal play, the two players stake at any given turn the same proportion of their remaining budget. Call this proportion  the `stake function' $S(\lambda,v)$ when the counter is at $v \in V$, Mina's budget is one, and Maxine's is~$\lambda$. 
We may seek to determine the value of  $S = S(\lambda,v)$ by a perturbative argument: under optimal play, Mina will stake $S$ and Maxine $\lambda S$ at the upcoming turn. If Maxine slightly increases her stake, she will gain from an increased probability of winning the next turn, with a resulting improvement in position; and she will lose, because she will have a smaller budget from the start of the next turn. These gain and loss terms cancel to first order under infinitesimal perturbation, since $\lambda S$ is Maxine's optimal stake. The formula that results is
\begin{equation}\label{e.stakefunction}
 S(\lambda,v) \, = \, \frac{\max_{u \sim v} h(\lambda,u) - \min_{u \sim v} h(\lambda,u)}{(\lambda +1)^2 \tfrac{\partial}{\partial \lambda} h(\lambda,v)} \, ,
\end{equation}
where $h(\lambda,v)$ is the $\lambda$-biased infinity harmonic function in~(\ref{e.lambdainftylap}). This formula witnesses the Present-Future Compromise, with the right-hand numerator representing the short-term demand to stake big and win the next turn, and the denominator reflecting the long-term via the need to retain budget for later turns. But the just sketched derivation is not rigorous. Two notable problems are that the proposed formula depends on the differentiability in $\lambda$ of $h(\lambda,v)$; and that the perturbative argument identifies merely a local saddle point, where a global one is needed. Both of these problems are real. The function $h(\lambda,v)$ does fail to be differentiable at certain values of $\lambda$
for some simple graphs. And the local saddle point often fails to be global, essentially because, when the counter is at one step from the boundary $B$, one or other player may be tempted to {\em go for broke}, staking all of her budget in an effort to close out the game at the present turn. In~\cite{HP2022}, these problems are respectively addressed by restricting the class of game data and by modifying the game rules. The class of root-reward trees, in which $G = (V,E)$ is a tree with $B$ its set of leaves, and $f = {\bf 1}_r$ for a distinguished leaf $r$, is considered, because $\lambda \to h(\lambda,v)$ is then differentiable for each $v \in V$. And the game is modified to a leisurely version, in which, after stakes are offered at any given turn, it is decided that a move will take place only with a small positive probability. The leisurely game disables the efficacy of the go-for-broke deviation and renders the local saddle point global. The principal result of \cite{HP2022} operates under these assumptions and proves in essence that the stake function indeed takes the form~(\ref{e.stakefunction}). It is in this sense that \cite{HP2022} shows that a class of allocated-budget stake-governed tug-of-war games satisfies the Present-Future Compromise.

\section{The two strands of tug-of-war}\label{s.gamespde}

In 1987, Harris and Vickers~\cite{HarrisVickers87} 
proposed a model of two firms that compete to secure a patent.
They called their model tug-of-war, and the notion and this name appear to have been introduced by them. The model is quite close to player-funded stake-governed tug-of-war played on a finite integer interval, and will explain it shortly. The model was intended to describe how firms may compete over time when the effect of a technological advantage is significant but uncertain, and when the firms' decisions are influenced by their relative progress. Their work has inspired a wealth of research, mainly in economics but also among computer scientists, on tug-of-war and related bidding games:~\cite{HarrisVickers87} has a total of almost six hundred citations as of early 2023 according to google scholar.

In 2009, Peres, Schramm, Sheffield and Wilson's~\cite{PSSW09} considered constant-bias tug-of-war with a fair coin. The relation they found between this game in a Euclidean setting and infinity harmonic functions sparked a wave of attention among analysts and probabilists, with almost five hundred works citing this article.

These two contributions 
study variants of the same basic concept---the two sets of authors even appear to have arrived independently at the name `tug of war'---and  
each has influenced many researchers, but they seem to have done so with a remarkable disjointness:  none of the 1100 references that cite one or other article cites both. The author of this monograph was until recently quite unaware
of the origin and development of tug-of-war in economics. Mathematicians may benefit from learning of this literature,
and there is surely value in communication between the concerned communities of researchers.  Next then we review developments for tug-of-war. First economics, and then mathematics.

 \subsection{Tug-of-war in economics}\label{s.economics}
 We begin by reviewing~\cite{HarrisVickers87}, whose treatment is accompanied by the proof-supplying mimeograph~\cite{HarrisVickers86} and Vickers' 1985 DPhil thesis~\cite{Vickers85}.
 
 \medskip
 
 \subsubsection{Harris and Vickers' `Racing with uncertainty'}\label{s.harrisvickers}
 In 1980, Lee and Wilde~\cite{LeeWilde} proposed a single-stage model for the relationship between research effort and outcome. Maxine and Mina  (as we call them)
 select respective effort rates $x,y \in (0,\infty)$. Maxine's discovery time is $M_+$, an exponentially distributed random variable of mean $x^{-1}$; Mina's is $M_-$, an independent exponential of mean~$y^{-1}$. The player with the lower discovery time $M = \min \{ M_-,M_+\}$ wins the race, and receives a reward.   Maxine and Mina  have respective cost functions $c_+(x)$ and $c_-(y)$, and the payments they make  are $M c_+(x)$ and $M c_-(y)$. In other words, $c_-$ and $c_+$ describe the rate of payment incurred by the given player as a function of her effort rate, so the payments reflect that the efforts are sustained for time~$M$.
 Since $M$ has mean $(x+y)^{-1}$,  Maxine and Mina's mean costs with joint effort rates $(x,y)$ are $\tfrac{c_+(x)}{x+y}$ and $\tfrac{c_-(y)}{x+y}$.
 
 Lee and Wilde's model is a single-turn player-funded tug-of-war. It fits in the framework of Subsection~\ref{s.playerfunded} if we take $G$ equal to $\{ -1,0,1\}$ with the adjacency relation, $B = \{-1,1\}$, and
 provided that in the player payments~(\ref{e.netreceipt})
 we replace $a_{X_0}$ by $\tfrac{c_+(a_0)}{a_0+b_0}$ and   $b_{X_0}$ by $\tfrac{c_-(b_0)}{a_0+b_0}$, (where note that $X_0 =1$ and $\tau =1$).

 The original tug-of-war model of Harris and Vickers has gameboard\footnote{For $i,j \in \Z$, $i \leq j$, we write \hfff{intint}$\llbracket i,j \rrbracket = \big\{ k \in \Z: i \leq k \leq j \big\}$.} $\llbracket 0, N+1\rrbracket$ and boundary set $B = \{0,N+1\}$.
 It is the player-funded tug-of-war whose single-step rule is Lee and Wilde's: the players offer effort rates at each turn, and the above rule decides the outcome according to given turn-independent
 cost functions $c_-$ and $c_+$. Formally then we may specify the model by taking $G$ to be an integer interval and replacing $a_{X_i} \to \tfrac{c_+(a_{X_i})}{a_{X_i}+b_{X_i}}$ and  $b_{X_0} \to \tfrac{c_-(b_{X_i})}{a_{X_i}+b_{X_i}}$  in~(\ref{e.netreceipt}).

Also formally, player-funded stake-governed tug-of-war is the special case of Harris and Vickers' model if we take $c_+(x) = x(x+y)$ and $c_-(y) = y(x+y)$, though the dependence of these cost functions on the opponent's strategy makes them inadmissible for this framework.

Harris and Vickers use a fixed point theorem prove that, for a fairly general class of cost functions $c_-$ and $c_+$, there exists a Nash equilibrium that is {\em stationary} in the sense that  players' effort rates are determined merely by present counter location. If Maxine and Mina's effort rates are labelled $x_i$ and $y_i$ when the counter is at $i \in \intint{N}$, then the equilibrium may be called symmetric if $x_i = y_{N - i}$.
And the game itself is called symmetric if $c_- = c_+$ and the terminal payments to the two players are equal after interchange of endpoints.
Harris and Vickers prove that the symmetric game has a symmetric equilibrium, and that there is only one symmetric equilibrium when $c_-(z) = c_+(z) = z^\eta$ for some $\eta > 1$ or when the gameboard is short and under weaker hypotheses on the cost functions. The symmetric equilibrium is computed explicitly in the example where $N$ is even and $\eta =2$. 
The gameboard has a central edge $[N/2,N/2+1]$. 
It is found that 
Maxine's effort $x_i$ peaks at $x_i = N/2+1$, the right end of this edge; that the ratio $x_n/x_{n-1}$ equals $2^{1/2 + (N/2 +1) - n}$; and that $y_n/y_{n-1} = 2^{-1} x_n/x_{n-1}$.
Thus Maxine's effort drops sharply to the right of $N/2 +1$, with a ratio of efforts at consecutive vertices having exponential decay in the distance to $N/2$; Mina's effort drops even more rapidly, with the ratio of her effort to Mina's decaying exponentially in the same sense. This computation illustrates the main message of Harris and Vickers: as they write in the abstract, `in our principal model, which is of a one-dimensional race, it is shown that the leader in the race makes greater efforts than the follower, and efforts increase as the gap between competitors decreases'. In summary, we may say that tug-of-war was first introduced principally to illustrate this effect. In the language of Section~\ref{s.allocate}, the Battlefield Cyl Fog is exemplified, with the battlefield at the vertices $N/2$ and $N/2 +1$.

Harris and Vickers briefly comment on the possible existence of asymmetric equilibria, mentioning that, for certain cost functions, computer simulations appear to indicate a total of $2N-1$ equilibria: one symmetric, and remaining asymmetric in $N$ pairs. 

\medskip

\subsubsection{Further treatments of player-funded tug-of-war in economics}\label{s.tullock}

In 1980, Tullock~\cite{Tullock} considered two-player single-stage contests. Player $A$ stakes $x \in [0,\infty)$ and player $B$, $y \in [0,\infty)$, and the contest is won by $A$ with probability $\tfrac{x^\gamma}{x^\gamma + y^\gamma}$, where $\gamma \in (0,\infty)$ is now known as the Tullock exponent. The case $\gamma =1$ is called a lottery; this is the contest used in each step of the stake-governed tug-of-war games specified in Section~\ref{s.tugofwar}. Taking $\gamma=0$ yields constant-bias tug-of-war with a fair coin. And
the limit $\gamma \to \infty$ leads to  all-pay auctions, in which the higher staking player wins (with a tie-break rule needed for equal stakes), and both stakes are forfeit. 

Konrad and Kovenock~\cite{KonradKovenock05}, Agastya and McAfee~\cite{AgastyaMcAfee} and Konrad~\cite{Konrad2012}, considered player-funded tug-of-war on finite integer intervals with the 
all-pay auction rule used to decide turn victor. Konrad and Kovenock found the battlefield effect to be manifest in a strong form with this rule: players offer stakes in essence at just two vertices. Agastya and McAfee reach rather different conclusions, because of use of a different tie-break rule and of discounting of terminal payments. There are some equilibria where no stakes are ever placed. In others, a battlefield exists, with momentum for the stronger player away from this location. But, contrary to the Battlefield Cyl Fog, stakes rise near the boundary: discounting future payments makes the imminence of victory or defeat significant.  Konrad~\cite{Konrad2012} presents a discussion involving player-funded games with Tullock lottery contests, emphasising the discouraging effect of early losses on a player. (Here, we may compare how the premise ``Incentive Inch, Outcome Mile" evokes discouragement that effects the less motivated player.)

H\"afner~\cite{Hafner17} considered a `tug-of-war team contest': player-funded tug-of-war on an integer interval with all-pay auction rules at each turn, but in which each player is in fact a team of countably many people, with each individual on a team responsible only for the payment at some given turn (but receiving nonetheless the full terminal payment for his team). It is shown that efforts diminish as the stronger team approaches victory: this conclusion, shared by Harris and Vickers (and by the Battlefield Cyl Fog), arises in this case in view of free-rider and effort-discouragement effects for team players staking when victory or defeat is close to assured.   
In~\cite{Hafner16}, H\"afner and Konrad consider a modified setup in which, for each team, a member is assigned to each site in open play, and is responsible for all stake payments made for play at that site. It is found that members of the leading team outspend their opponents (as in the Battlefield Cyl Fog) when there is little discounting of the terminal payment, and that the reverse holds when there is much discounting.  

\medskip

\subsubsection{The player-funded game with the majoritarian objective}

In tug-of-war on an integer interval gameboard symmetric about the origin, one player wins when she secures a number of turn victories that exceeds her opponent's by a certain value. Harris and Vickers further considered the game where the game winner is the first to secure a certain number of victories.  This {\em majoritarian objective} model fits in the framework of Subsection~\ref{s.playerfunded} if we consider what may be called {\em directed rectangle graphs}: that is, 
if we take $G = (V,E)$ to be a directed graph  with $V = \llbracket 0,m \rrbracket \times \llbracket 0,n \rrbracket$, and~$E$ the set of north and east pointing nearest-neighbour edges; $f(x,y) = \lambda {\bf 1}_{y=m}$ and  $g(x,y) = {\bf 1}_{x=n}$; and initially $X_0 = (0,0)$. (Harris and Vickers' stake-rules entail further changes as described in Subsection~\ref{s.harrisvickers}.)
In this model, backward induction readily shows the equilibrium is unique under reasonable hypotheses on the Harris-Vickers cost functions $c_-$ and $c_+$.
Numerical work in~\cite{HarrisVickers87} with $z^\eta$-cost functions indicate some expected effects, such as the rapidly emerging dominance of the stronger player (the presence then of a battlefield); but also of surprises, such as stark oscillations in effort rates in the case $m=n$ as a function of $n \in \intint{10}$.

The Democratic and Republican primary contests to determine the party's nominee in a U.S. presidential election are run on a calendar of several months.
In the `New Hampshire' effect, candidates typically spend a significant fraction of their resources in the earliest contests. Klumpp and Polborn~\cite{KlumppPolborn} considered as a model of this phenomenon majoritarian objective tug-of-war with a Tullock contest rule for deciding turns. The early burst of campaign expenditure and momentum effects discussed are similar to the Battlefield Cyl Fog, and contrast with a more uniform resource allocation in simultaneous contests, so that a sequential rather than simultaneous primary season is typically cheaper for the winning candidate.  

Kovenock and Konrad~\cite{KonradKovenock09} considered the majoritarian objective game with all-pay auction rules at the turn, with the introduction of intermediate prizes, or payoffs to players at non-terminal times. Gelder~\cite{Gelder} shows that discounting dissipates the accrual of momentum to the stronger player in the majoritarian game with all-pay auction rules.
Best-of-three contests have been examined by Malueg and Yates~\cite{MaluegYates}, and by Sela~\cite{Sela}.  Fu, Lu and Pan~\cite{FuLuPan} considered a team contest version of the majoritarian game. 

Ewerhart and Teichgr\"aber~\cite{ET19} introduce a player-funded tug-of-war framework that is more general than integer interval or first-to-$N$ contests. Operating under assumptions including that gameplay before termination has a certain exchangeability of move order, they show that  gameplay takes place in an initial phase, in which no state may be revisited, followed by a later tug-of-war phase. They consider symmetric equilibria, showing that there exists a unique such that is Markov perfect (or stationary, in the above usage), for a broad range of contest rules; and they discuss the problem of designing contests that are optimal for the tournament organizer by maximizing combined player expenditure. 
Existence and uniqueness of the symmetric stationary equilibrium in player-funded tug-of-war with a Tullock contest of exponent one-half at each turn were proved and this equilibrium analysed in~\cite{KST}.

\medskip

\subsubsection{The allocated-budget games}
 
 In 1921, Emile Borel considered a game~\cite{Borel1921} that pits a protagonist now commonly called Colonel Blotto against an opponent. Borel gives an example: players $A$ and $B$ each choose an ordered triple of non-negative reals that sum to one. Player $A$ wins the game if at least two of his values exceed the counterparts for $B$. In general, the colonel and his opponent each allocate a given resource across a fixed number of battlegrounds, with a given battle won by the better committed player, and the battle-total determining the winner of the war. 
 The label `Blotto' is sometimes attached to games where players have fixed resources, whether they be committed simultaneously as in Borel's game, or for example sequentially, as in budget-allocated tug-of-war from Subsection~\ref{s.allocatedbudget}. 
 
 In Section~\ref{s.overviewallocated}, we advanced a premise for how to play the Blotto game with the lottery contest function $\tfrac{x}{x+y}$: stake a proportion~(\ref{e.stakefunction}) of your remaining fortune, whether you be Mina or Maxine. But we also indicated that even tug-of-war on $\llbracket 0,3 \rrbracket$ is a counterexample to this premise. Rather than altering the game (as~\cite{HP2022} did), we may ask whether a class of graphs can be found for which the premise is correct. The majoritarian objective games---tug-of-war on directed rectangle graphs---are a good starting point. In this game on  $\llbracket 0,1 \rrbracket \times \llbracket 0,n \rrbracket$, there are $n$ contests, with one player needing to win every one, and the other requiring only one victory. A symmetry consideration points to a stake proportion of $1/n$ being optimal 
 when $X_0 = (0,0)$. This choice in fact validates~(\ref{e.stakefunction}), as~\cite[Section~$2.6$]{HP2022} explains.  Klumpp, Konrad and Solomon~\cite{KlumpKonradSolomon} generalise this example to show that the game on any directed rectangle graph, and for a broad range of turn contest functions, has an equilibrium with an even split of resources across turns. 
 
 Allocated-budget tug-of-war on a finite integer interval has been examined by Klumpp~\cite{Klumpp}. In note~$1$ on page $27$, he notes an example, similarly as does~\cite{HP2022}, in which a formula equivalent to~(\ref{e.stakefunction})
 is invalidated when the weaker player is one step from winning. Where~\cite{HP2022} turned to a leisurely version of the game to regularise this example, Klumpp considers lowering the Tullock exponent from the lottery value $\gamma = 1$. He surmises from a numerical investigation that the formula in question is likely to offer an equilibrium precisely when $\gamma$ at most one-half.

 \medskip
 
 \subsubsection{An experimental study} 
 
 In 2002, Zizzo~\cite{Zizzo} reported outcomes of a study of human players in a majoritarian objective, first-to-$10$ game.
 The experimental setting was designed to replicate as closely as possible the 
  Harris-Vickers first-to-$10$ game 
   with cost functions $c_-(x) = c_+(x) = x^2$. Leaders clearly outspent their opponents
   only when a significant gap had opened, and the Harris-Vickers predictions appeared not to be well supported by the tests. Two differences between theory and setup should be noted, however.
   
   First, although players retained unspent budgets and thus in effect lost spent monies as  in a player-funded game, they also operated under a budget constraint, which equalled one-half of the terminal payment to the victor. Zizzo acknowledges this difference and writes  ``to minimize the possibility of distortions that might occur ... as a result of a binding budget constraint ... instructions stressed the need for subjects to
retain enough money to fund the later rounds of the competition in an unconstrained way." Such an admonition was likely intended to prevent player bankruptcy, the most blatant manifestation of budget constraint. But it also illustrates the role of the future in the Present-Future Compromise seen in allocated-budget games. The setup was a hybrid of player-funded and allocated-budget rules, and it is conceivable that the outcomes felt the effects of both regimes.

Second, recall that Harris and Vickers' cost functions $c_-(z)$ and $c_+(z)$ are rates of charge to players as a function of their effort rates~$z$. Zizzo's implementation sets costs equal to $z^2$ but treats them simply as the payments made by players at the turn. When Mina and Maxine stake $x$ and $y$ at a given turn, their mean payments at the turn  are $\tfrac{x^2}{x+y}$ and  $\tfrac{y^2}{x+y}$ for quadratic Harris-Vickers; for Zizzo, the payments are deterministic quantities, $x^2$ and $y^2$. Zizzo's framework at a given turn is readily seen to be equivalent to a Tullock contest from Subsection~\ref{s.tullock}, with exponent $\gamma$ equal to one-half.
The player-funded majoritarian objective game with this framework may well have qualitatively similar theoretical features to those that Harris and Vickers identified and predicted, but this premise deserves scrutiny in interpreting Zizzo's work as an experimental test of Harris and Vickers' theory.


\subsection{Tug-of-war in mathematics}

Peres, Schramm, Sheffield and Wilson~\cite{PSSW09} considered constant-bias tug-of-war with a fair coin. Under optimal play, Maxine moves to an $h$-maximizing neighbour of the counter's present location, and Mina to an $h$-minimizing one, where $h$ is the discrete infinity harmonic extension of boundary data $f:B \to \R$; namely, $h(\cdot)$ equals $h(\lambda,\cdot)$ from~(\ref{e.lambdainftylap}) with $\lambda =1$. In this way, $h$ is shown to be the game's value. This is perhaps the simplest connection between discrete infinity harmonic functions and tug-of-war games, but it was not the first. 
Richman games were introduced in~\cite{LLPU} and further studied in~\cite{LLPSU}. A first-price Richman game is an allocated-budget tug-of-war game with an auction to determine the turn victor (so that the higher staking player wins), and with the turn victor paying the opponent his stake, the opponent making no payment at the turn.  The threshold ratio is a value~$t$ such that, if the initial relative ratio of Maxine's budget and the total budget for the two players  exceeds $t$, then Maxine will almost surely win the game if she plays correctly, while Mina will win in this sense if this ratio is less than~$t$. In~\cite{LLPU}, it is proved that the threshold ratio exists as a function of initial counter location and equals the discrete infinity harmonic extension of $f:B \to \R$. Several variants of these games---with all-pay rules, or infinite duration, or in the poormen variant, where the turn victor pays the bank, not the opponent---have been studied~\cite{AJZ,AHC,AHI,AIT20} by theoretical computer scientists, who also address approximate algorithms and computation complexity: see~\cite{AH22} for a survey that treats these directions.   

We have recalled constant-bias tug-of-war in the discrete context of  given graphs in part because the framework of our results is discrete, but the principal aim of~\cite{PSSW09} 
was to relate tug-of-war in a low mesh limit to 
 infinity harmonic functions mapping domains in $\R^d$ to $\R$, these being continuum counterparts to~(\ref{e.lambdainftylap}).
In `$\e$~tug-of-war' for given $\e > 0$,  the gameboard is a domain $D$ in $\R^d$, with boundary data $f: B \to \R$, $B = \partial D$. Mina and Maxine (as we name them) move a counter---a point in $D$---a distance at most~$\e$ at each turn, according to the outcome of a fair coin flip. Mina pays Maxine the evaluation of~$f$ at the counter location when the game ends as the counter reaches~$B$. In~\cite{PSSW09},
 the value of $\e$ tug-of-war is shown in a low-$\e$ limit to converge to the infinity harmonic function that extends $f$ to $D$. This is the viscosity solution~$h:D \to \R$ of the infinity Laplace equation 
 $\sum_{i,j} \partial_{x_i} h \,  \partial_{x_i x_j} h \, \partial_{x_j} h = 0$ subject to $h \big\vert_B = f$.
 The infinity Laplace operator is a degenerate second-order differential operator that was introduced and studied by Aronsson~\cite{Aronsson67,Aronsson68} in 1967 and that has 
   a subtle uniqueness~\cite{Jensen93}  and regularity~\cite{Savin05,EvansSavin}  theory. Its viscosity solutions~\cite{CIL92} are absolutely minimizing Lipschitz functions extending the given boundary data~\cite{Jensen93}. 
   
   As \cite{Manfredi2012} surveys, the game theory connection identified in~\cite{PSSW09} sparked much attention from the PDE community.    
In~\cite{ArmstrongSmart2012}, boundary rules in $\e$ tug-of-war were altered from~\cite{PSSW09}
to obtain more regular game value functions. In~\cite{CCF2012}, random move sizes with a heavy tail (judged on scale~$\e$) were introduced, thus connecting to the infinity fractional Laplacian.
Suppose that $\e$ tug-of-war is played with a random displacement of magnitude $c \e$ of the counter at the end of each move. 
For $p \in (2,\infty)$, the $p$-Laplacian~\cite{Lindqvist} is an interpolation of the usual Laplacian operator (for $p=2$) and the infinity version. 
In~\cite{PeresSheffield}, it is shown that the value of this game converges to a $p$-harmonic function for $p$ suitably chosen as a function of $c$: see~\cite{Lewicka}
for a survey centred around this perspective.  A variant of this game has been used to study $p$-Laplacian obstacle problems in~\cite{LewickaMarta2017}.
The book~\cite{BlancRossi} reviews  the abundant connections between tug-of-war and PDE.

It is the fair case $\lambda=1$ of constant-bias tug-of-war that~\cite{PSSW09} considers. The biased case $\lambda \in (0,\infty)$ was addressed in~\cite{PPS10}: $\e$ tug-of-war was modified so that a suitable drift is maintained in the low-$\e$ limit, and a continuum biased infinity harmonic function is then obtained as limiting game value.
In~\cite{PeresSunic}, a rapid algorithm for computing biased discrete infinity harmonic functions~(\ref{e.lambdainftylap}) was found in terms of a path decomposition of the graph $(V,E)$.
As the bias parameter $\lambda$ varies, this decomposition changes at certain values. It is these changes that cause the lack of differentiability in $(0,\infty) \to \R: \lambda \to h(\lambda,v)$
to which we alluded in Subsection~\ref{s.overviewallocated}.

Tug-of-war played on directed graphs has also garnered attention among mathematicians. 
Many combinatorial games, such as Go and Hex, are zero-sum two-player games in which players alternate in making directed moves by placing counters that cannot be moved later.
In~\cite{PSSW07}, random-turn versions of these `selection games' were considered. Hex is a game discovered by Piet Hein in 1942 in which two players alternately place red or blue hexagons on a hexagonal grid, each seeking to forge a connection between a pair of diagonally opposed boundary segments in his own colour. It is a complex game~\cite{Hayward} which has no explicit solution on $11$ by $11$ boards on which it is played competitively. In~\cite{PSSW07}, the solution to random-turn Hex (among other selection games) was  found explicitly  in terms of the maximum pivotal probability for critical percolation in the unplayed region. This work offers a perspective on  the strategic importance of positions in multi-turn games.

To conclude, we first mention two broad 
antecedents to stake-governed random-turn games. Rufus Isaacs' differential games~\cite{Isaacs65}, include examples where a pursuer seeks to capture an evader in Euclidean space, and each instantaneously selects certain
control variables such as velocity, analogously to Mina and Maxine's ongoing selection of stakes. In Shapley's stochastic games~\cite{Shapley1953},  transition probabilities governed step-by-step by decisions of two players encompass the stake game framework in a broad sense.\footnote{Both differential and stochastic games have attracted ongoing attention. Several chapters of a recent handbook~\cite{HandbookDynamicGameTheory} treat dynamic game theory, and~\cite{StochasticGamesRelatedTopics} collects articles about stochastic games.} 
These are fairly loose connections, however. Regarding tug-of-war itself, we have seen how mathematicians have tended to focus on constant-bias models.
They have explored a much wider geometric setting than the  finite integer intervals and (what we have called) the directed rectangle graphs that are common elements in the economics literature, 
 exploring these games' implications in probability and PDE. Economists originated tug-of-war, considering from the outset bidding rules beyond the trivial constant-bias, and have sustained attention on these games for almost four decades, putting the theme of resource allocation at the heart of their investigations throughout. 
This article and~\cite{HP2022} offer detailed mathematical treatments of resource-allocation tug-of-war games: with allocated budgets in~\cite{HP2022}, and with player funding in the present article.
Raising the prospect of uncovering further beautiful mathematical structure, and informing theoretical analysis of economic behaviour, tug-of-war in its original resource-allocation guise would surely repay the further attention of mathematicians.

\subsection{Acknowledgments}
The author thanks G\'abor Pete for many discussions about stake-governed games. He thanks Christian Ewerhart for interesting discussions, and Judit Z\'ador for help with Mathematica and in preparing the article's figures.  He further thanks a referee for careful reading and valuable comments.
He is supported by the National Science Foundation under DMS grants~$1855550$ and~$2153359$ and by the Simons Foundation as a $2021$ Simons Fellow.

\chapter{Main results}\label{c.speclost}

In this chapter, we specify precisely the Trail of Lost Pennies and state our main conclusions. The structure of the remainder of the article is explained at the end of the chapter.


\section{Game setup, strategies and Nash equilibria}\label{s.gamespec}

The Trail of Lost Pennies on $\Z$ will be denoted by $\game(m_{-\infty},m_\infty,n_{-\infty},n_\infty)$. 
The data that specifies the game takes the form of 
\begin{equation}\label{e.quadruple}
\textrm{a quadruple $\big(m_{-\infty},m_\infty,n_{-\infty},n_\infty \big) \in \R^4$ that satisfies  $m_{-\infty} < m_\infty$ and $n_\infty < n_{-\infty}$} \, .
\end{equation}
For any given $k \in \Z$, we will specify the gameplay of $\game(m_{-\infty},m_\infty,n_{-\infty},n_\infty)$ where the initial location of the counter is equal to $k$.
The counter's location~$X$ will evolve from its initial location $X_0 = k$ in discrete time-steps, the result being a stochastic process $X:\nwithzero \to \Z$ whose law is determined by~$k$ (we take $\N$ to include zero).\hfff{withandwithout} 
This random process is specified by the pair of strategies adopted by Mina (who plays to the left) and Maxine (who plays to the right).

For either player, a strategy~\hfff{strategy} is a map $S:\Z \times \nwithoutzero \to [0,\infty)$. A player who follows the strategy~$S$ stakes $S(X_{i-1},i)$ at the turn with index~$i \geq 1$. 
Let $\mc{S}$ denote the set of strategies. An element of $\mc{S}^2$ is called a strategy pair. A generic element of $\mc{S}^2$ will be written $(S_-,S_+)$, where the respective components are the strategies of Mina and Maxine.

We wish then to specify the gameplay process\hfff{counter}  $X:\nwithzero \to \Z$ as a function of a given element $(S_-,S_+) \in \mc{S}^2$ and the initial location $X_0 = k \in \Z$. 
We will write $\pgameplay{S_-}{S_+}{k}$ for a probability measure that specifies this gameplay process and accompanying aspects of the game; the associated expectation operator will be written $\egameplay{S_-}{S_+}{k}[\cdot]$.
At the turn with index $i \in \nwithoutzero$, Mina stakes $S_-(X_{i-1},i)$ and Maxine stakes  $S_+(X_{i-1},i)$. By sampling of independent randomness, the turn victor is declared to be Maxine with probability  $\tfrac{S_+(X_{i-1},i)}{S_-(X_{i-1},i) + S_+(X_{i-1},i)}$; in the other event, it is declared to be Mina. Maxine will elect to move the counter one place to the right if she is the turn victor; Mina, one place to the left. 
(It is intuitive given the rules of the game that we are specifying that the two players will always elect to move the counter in the said directions, and we will not furnish the straightforward details to the effect that permitting other options changes nothing essential about the game.)
Should neither player make a stake at the given turn---that is, if  $S_-(X_{i-1},i) = S_+(X_{i-1},i) = 0$---then a further rule is needed to permit play to continue. We will declare that, in this event, each player wins the right to move with equal probability (with Maxine moving right, and Mina left, as usual).  

Formally, then, our counter evolution satisfies the condition that,
for $(k,i,\ell) \in \Z \times \nwithoutzero \times \Z$, 
$$
\pgameplay{S_-}{S_+}{k} \Big( X_i - X_{i-1} = \ell \, \Big\vert \, X_j, j \in \llbracket 0, i-1 \rrbracket \Big) \, = \, \tfrac{S_-(X_{i-1},i)}{S_-(X_{i-1},i) + S_+(X_{i-1},i)} {\bf 1}_{\ell = - 1} + \tfrac{S_+(X_{i-1},i)}{S_-(X_{i-1},i) + S_+(X_{i-1},i)} {\bf 1}_{\ell = 1} \, ,
$$
where recall that we use the integer-interval notation $\llbracket i,j \rrbracket = \big\{ \ell \in \Z: i \leq \ell \leq j \big\}$, $i,j \in \Z$.
Note that,
in reading the ratios on the right-hand side in the display, we adopt the convention that $0/0 = 1/2$.

We further wish to specify the other pertinent features of the game when the strategy pair $(S_-,S_+)$ is played. These features are the resulting payoffs to Mina and Maxine. 
Mina's payoff $P_-$ is the sum of a negative term given by the total costs incurred to Mina during gameplay, and a further term that is the terminal payment that is made to her. 
Indeed, we may write\hfff{receiptmina}
\begin{equation}\label{e.minapayoff}
 P_- \, = \,  - \sum_{t = 1}^\infty C_-(t) \, \, + \, \, T_-  \, , 
\end{equation}
where $C_-(t)$ denotes the cost incurred to Mina at the turn with index $t \in \nwithoutzero$,
and $T_-$ equals the terminal payment to Mina. We have then that the cost $C_-(t)$ is equal to Mina's stake $S_-(X_{t-1},t)$.

The terminal payment $T_-$ is in essence given by $n_{-\infty}$ if Mina wins the game by eventually bringing the counter infinitely far to the left; and to $n_\infty$ in the opposing event. However, a precise formulation is needed to make sense of this. We define the {\em escape} event $E$ according to 
\begin{equation}\label{e.escape}
E = \big\{ \lim_n \vert X_n \vert = \infty \big\} \, .
\end{equation}
The {\em left} and {\em right}  escape events\hfff{escape} are given by
\begin{equation}\label{e.leftrightescape}
E_- =  \big\{ \limsup_n X_n  = - \infty \big\} \, \, \, \, \textrm{and} \, \, \, \, 
E_+ =  \big\{ \liminf_n X_n  =  \infty \big\}   \, . 
\end{equation}
 Note that $E_-$ and $E_+$ are disjoint events whose union equals $E$. We regard them as victory events for Mina and Maxine, and accordingly set the terminal payment to Mina as follows:
 \begin{equation}\label{e.terminalmina}
  T_- \, = \, \begin{cases}
  \, \, n_{-\infty}  &  \text{when $E_-$ occurs} \, , \\
 \, \, n_\infty  &  \text{when $E_+$ occurs} \, , \\ 
   \, \, n_*  &  \text{when $E^c$ occurs} 
  \, .
\end{cases}
\end{equation}
Here $n_*$ is a given real value that is at most $n_\infty$. By assigning this terminal payment to Mina in the event of non-escape, we ensure that this payment no more generous than that made in the event~$E_+$ of her defeat.

We may specify Maxine's payoff\hfff{receiptmaxine} 
\begin{equation}\label{e.maxinepayoff}
 P_+ \, = \,  - \sum_{t = 1}^\infty C_+(t) \,  \, + \, \, T_+  
\end{equation}
with counterpart interpretations for the right-hand terms:
 the cost~$C_+(t)$ incurred to Maxine at the turn with index $t \in \nwithoutzero$ equals Maxine's stake $S_+(X_{t-1},t)$, while
 the terminal payment  $T_+$ that she receives is given by   
 \begin{equation}\label{e.terminalmaxine}
  T_+ \, = \, \begin{cases}
  \, \,  m_{-\infty}  &  \text{when $E_-$ occurs} \, , \\
 \, \,  m_\infty  &  \text{when $E_+$ occurs} \, , \\ 
   \, \, m_*  &  \text{when $E^c$ occurs}  \, ,
\end{cases}
 \end{equation}
where $m_*$ is a given real value\footnote{Note that $(m_*,n_*) \in \R^2$ and $(m_{-\infty},m_\infty,n_{-\infty},n_\infty) \in \R^4$ are the parameters that specify the Trail of Lost Pennies on~$\Z$. We thus speak imprecisely when we refer to the latter quadruple as the game's data. Given the upper bounds that we impose on them, the values of $m_*$ and $n_*$ will be immaterial for our analysis.} that is at most $m_{-\infty}$.

The quantities labelled $P$, $C$ and $T$ for Mina and Maxine are determined by the gameplay $X:\nwithzero \to \Z$. The gameplay and these other random variables are thus coupled together under the law\hfff{lawexpect} $\PP_{S_-,S_+}^k$. Of course, starting location $X_0 = k$ and strategy pair $(S_-,S_+)$ are fundamental for determining game outcome including the above described quantities. In our notation, this dependence is communicated by the labels of the law $\PP_{S_-,S_+}^k$, rather than in the annotations $P_-$, $T_-$, and so on.


A strategy is {\em time-invariant}\footnote{This property of strategies corresponds to that called stationary by Harris and Vickers: see Subsection~\ref{s.harrisvickers}. And it is sometimes called Markov, since decisions depend on history only through the present state of play.}
if $S(i,j)$ is independent of $j \in \nwithoutzero$ for every $i \in \Z$. 
The set~\hfff{timeinvariant} of time-invariant strategies will be denoted by~$\tis$. A time-invariant strategy pair~\hfff{typicalnotation} $(S_-,S_+) \in \mc{S}_0^2$
may be identified with a pair of sequences $\big\{ a_i: i \in \Z \big\}$ and  $\big\{ b_i: i \in \Z \big\}$ , where 
 $a_i = S_+(i,j)$ and $b_i = S_-(i,j)$ for $(i,j) \in \Z \times \nwithoutzero$. 


A strategy pair $(S_-,S_+) \in \mc{S}^2$ is a Nash equilibrium if 
$$
 \egameplay{S_-}{S_+}{k} [P_+] \geq  \egameplay{S_-}{S}{k} [P_+]  \, \, \, \,
\textrm{and} \, \, \, \,  \egameplay{S_-}{S_+}{k} [P_-]  \geq  \egameplay{S}{S_+}{k} [P_-] 
$$
 for all $S \in \mc{S}$ and $k \in \Z$. Note that this condition takes a strong form\footnote{The strengthened form is closely related to the notions of subgame perfect, and Markov perfect, equilibrium.}, in that it stipulates the displayed bound for any initial condition $X_0 = k \in \Z$ for the counter location.

Let $\mc{N}(m_{-\infty},m_\infty,n_{-\infty},n_\infty) \subset \mc{S}^2$ denote the set\hfff{nash} of Nash equilibria. Consider
a time-invariant Nash equilibrium, namely an element~$(S_-,S_+)$ of $\tis^2$ that satisfies the above condition: when such a strategy pair is played, neither player would gain by altering strategy, even if the proposed alternative strategy is not time-invariant.

In an abuse of notation,
generic elements of $\tis$, for respective use by Maxine and Mina, will be called
 $\big(a_i:i\in \Z)$ and  $\big(b_i:i\in \Z \big)$. In a further abuse,\hfff{notab} the accompanying element of $\tis^2$
 will be denoted\footnote{In the strategy-pair notation $(S_-,S_+) \in \mc{S}^2$, governed by $- < +$, Mina precedes Maxine. Thus the notation $(b,a)$ for strategy pairs will be standard. We will shortly introduce an $(a,b,m,n)$-quadruple notation for stakes and mean payoffs, in which Maxine precedes Mina (in the sense of `$a$ before $b$' and `$m$ before $n$'). As a result, usages of the form `$(a,b,m,n)$ is the quadruple associated to the Nash equilibrium $(b,a)$' will be made.}   
$\big\{ (b_i,a_i): i \in \Z \big\}$. 

 
 \section{Time-invariant Nash equilibria and the \abmnmacspace  equations}

We now begin to present our main results. We first introduce the \abmnmacspace  equations, which will be fundamental to this study. Theorems~\ref{t.positiveabmn}
 and~\ref{t.minamarginvalues} present basic properties of the equations' solution set, and Theorem~\ref{t.nashabmn}
 is the result that bridges between the equations and the trail game (as we will sometimes informally call the Trail of Lost Pennies). In Theorem~\ref{t.nashequil.prelim}, these theorems are leveraged to characterize when the trail game has time-invariant Nash equilibria in terms of a condition on boundary data involving an important basic quantity, the {\em Mina margin}, which is introduced here. The section ends with Theorem~\ref{t.ajbj}, which offers precise asymptotic decay estimates for Nash equilibria as the index varies away from the {\em battlefield index}, at which the players are most likely to decide the ultimate outcome of a given game; and with its consequence Theorem~\ref{t.unanimity}, which describes gameplay at any Nash equilibrium. 
 
 \begin{definition}\label{d.quadruple}
 Let $\big\{ (b_i,a_i): i \in \Z \big\}$ denote\hfff{aiai} a time-invariant strategy pair: namely,  
 $$
 \big\{ (b_i,a_i): i \in \Z \big\}  \in \tis^2 \, .
 $$
Let $S_-,S_+ \in \mc{S}$ be strategies such that $S_-(i,j) = b_i$ and $S_+(i,j) = a_i$  whenever $(i,j) \in \Z \times \nwithoutzero$. 
 
 Set $m_i = \egameplay{S_-}{S_+}{i} [ P_+]$ and  $n_i = \egameplay{S_-}{S_+}{i}  [P_-]$ for $i \in \Z$. By this means,\hfff{mini}
 we have associated to any element  $\big\{ (b_i,a_i): i \in \Z \big\}  \in \tis^2$
 a $\Z$-indexed quadruple  $\big\{ (a_i,b_i,m_i,n_i): i \in \Z \big\}$ of elements taking values in $[0,\infty)^2 \times \big(\R  \cup \{ -\infty\}\big)^2$.
 \end{definition}
 
 \begin{definition}\label{d.abmn}
The ABMN system on $\Z$ is the set of equations in the four real variables $a_i$, $b_i$, $ m_i$ and $n_i$, indexed by~$i \in \Z$, 
\begin{align*}
  (a_i + b_i)(m_i+a_i)  & =  a_i m_{i+1} + b_i m_{i-1} && \qquad \textrm{ABMN}(1) \\
  (a_i + b_i)(n_i+b_i)  & =  a_i n_{i+1} + b_i n_{i-1}  &&\qquad \textrm{ABMN}(2)   \\
  (a_i + b_i)^2  & =  b_i \big( m_{i+1} - m_{i-1} \big) &&\qquad \textrm{ABMN}(3)   \\
      (a_i + b_i)^2  & =   a_i \big( n_{i-1} - n_{i+1} \big) &&\qquad \textrm{ABMN}(4)   \, ,
  \end{align*}
 where $i$ ranges over $\Z$. We will refer to the above equations throughout in the form \abmnmac$(i)$, for $i \in \{1,2,3,4\}$, rather than by a conventional numerical labelling.
 It is always supposed that $a_i$ and~$b_i$ are non-negative for $i \in \Z$. 
 A solution is said\hfff{boundarydata} to have boundary data $(m_{-\infty},m_\infty,n_{-\infty},n_\infty)$ when 
 \begin{equation}\label{e.boundarydata}
 \lim_{k \to \infty} m_{-k} = m_{-\infty} \, \, \, , \, \, \,
 \lim_{k \to \infty} m_k = m_\infty \, \, \, , \, \, \, 
 \lim_{k \to \infty} n_{-k} = n_{-\infty} \,\,\,\,
 \textrm{and}
 \,\,\,\, \lim_{k \to \infty} n_k = n_\infty \, . 
 \end{equation}
 For such a solution, the {\em Mina margin}\hfff{minamargin} is set equal to $\frac{n_{-\infty} - n_\infty}{m_\infty - m_{-\infty}}$.
 A solution is called {\em positive} if $a_i > 0$ and $b_i > 0$ for all $i \in \Z$.\hfff{positiveabmn} It is called {\em strict} if $m_{i+1} > m_i$ and $n_i > n_{i+1}$ for such $i$.
 \end{definition}
 \begin{theorem}\label{t.positiveabmn}
 Let $\big\{ (a_i,b_i,m_i,n_i) \in (0,\infty)^2 \times \R^2: i \in \Z  \big\}$ be a positive \abmnmacspace  solution.
 \begin{enumerate}
 \item The solution is strict.
 \item The solution has boundary conditions $(m_{-\infty},m_\infty,n_{-\infty},n_\infty)$  that satisfy $m_\infty > m_{-\infty}$ and $n_{-\infty} > n_\infty$.
 \item The values $m_{-\infty}$, $m_\infty$, $n_{-\infty}$ and $n_\infty$ are real numbers.  As such, the Mina margin $\frac{n_{-\infty} - n_\infty}{m_\infty - m_{-\infty}}$ exists and is a positive and finite real number.
 \end{enumerate}
 \end{theorem}
 The Mina margin has a fundamental role to play in determining whether the \abmnmacspace  system can be solved, as we now see. 
 \begin{theorem}\label{t.minamarginvalues}
Invoking Theorem~\ref{t.positiveabmn}(3), we may set $I \subset (0,\infty)$ equal to the set of values of the Mina margin  $\frac{n_{-\infty} - n_\infty}{m_\infty - m_{-\infty}}$, where $\big\{ (a_i,b_i,m_i,n_i) \in (0,\infty)^2 \times \R^2: i \in \Z  \big\}$ 
 ranges over the set of positive \abmnmacspace  solutions. 
 \begin{enumerate}
 \item There exists a value $\lambda \in (0,1]$ such that the set $I$ is equal to the interval $[\lambda,\lambda^{-1}]$.
  \item Moreover, a positive \abmnmacspace  solution exists with boundary data $(m_{-\infty},m_\infty,n_{-\infty},n_\infty) \in \R^4$
  if and only if $m_{-\infty} < m_\infty$, $n_\infty < n_{-\infty}$  and   $\frac{n_{-\infty} - n_\infty}{m_\infty - m_{-\infty}} \in [\lambda,\lambda^{-1}]$.
  \item The value of $\lambda$ is at most $0.999904$. 
  \end{enumerate}
 \end{theorem}
\begin{conjecture}\label{c.lambda} 
The value of $\lambda$ is at least $0.999902$.
 \end{conjecture}
 Evidence for this conjecture will be presented in~Section~\ref{s.minamarginmap}.
 It is perhaps surprising that a value so close to one may emerge from the ABMN system. The next result makes it all the more remarkable.
 \begin{theorem}\label{t.nashabmn}
Let $(m_{-\infty},m_\infty,n_{-\infty},n_\infty) \in \R^4$ satisfy $m_{-\infty}<m_\infty$ and  $n_\infty < n_{-\infty}$. 
 \begin{enumerate}
 \item
 Suppose that an element $\big\{ (b_i,a_i): i \in \Z \big\}$  of $\tis^2$ 
 lies in~$\mc{N}(m_{-\infty},m_\infty,n_{-\infty},n_\infty)$. Then 
 $$
 \textrm{the quadruple $\big\{ (a_i,b_i,m_i,n_i): i \in \Z \big\}$ associated to the element by Definition~\ref{d.quadruple}}
 $$
  is a positive \abmnmacspace  solution with boundary data $(m_{-\infty},m_\infty,n_{-\infty},n_\infty)$.
  \item Conversely, suppose that  $\big\{ (a_i,b_i,m_i,n_i) \in (0,\infty)^2 \times \R^2 : i \in \Z \big\}$ is a positive \abmnmacspace solution
  with boundary  data $(m_{-\infty},m_\infty,n_{-\infty},n_\infty)$. Then  $\big\{ (b_i,a_i): i \in \Z \big\} \in \tis^2$ lies in $\mc{N}(m_{-\infty},m_\infty,n_{-\infty},n_\infty)$.
  \end{enumerate}
 \end{theorem}
 In view of this result, 
 Theorem~\ref{t.minamarginvalues}(3) eliminates the possibility that $\lambda =1$, that Nash equilibria exist only when the players have symmetric roles. 
Conjecture~\ref{c.lambda} implies that the bound on $\lambda$ proved in 
 Theorem~\ref{t.minamarginvalues}(3)  is close to sharp: this quantity, which is canonically associated to the Trail of Lost Pennies on $\Z$---a natural and simple enough game---would then be remarkably close to one, differing from it by less than $10^{-4}$. 

(We suspected for a while that $\lambda =1$. An American Mathematical Monthly article~\cite{PS07} poses the question ``When is $0.999 \ldots$ equal to $1$?", discussing an example of an infinite sum involving harmonic numbers where rigorous bounds suggested that the sum might equal one but for which a later obtained exact formula identified the sum as equalling 0.999222\ldots. Our situation is similar, and might even permit an extra `9' in the posed question. There are differences as well: the premise ``Incentive Inch, Outcome Mile" advances a case that $\lambda$ should be close to one.)
 
 \begin{definition}\label{d.standard}
In its {\em standard} form, the trail game has boundary data that satisfies $m_{-\infty} = 0$, $n_\infty = 0$ and $m_\infty = 1$.
For a game in this form, the game's data is thus specified by one parameter, $n_{-\infty} \in (0,\infty)$. This parameter equals the Mina margin~$\frac{n_{-\infty} - n_\infty}{m_\infty - m_{-\infty}}$.
\end{definition}
  
  Let $x \in (0,\infty)$. By $\standard(x)$, we denote the Trail of Lost Pennies on $\Z$  in its standard form, with the Mina margin equal to $x$.
  That is, $\standard(x)$ equals $\game(0,1,x,0)$, as this game has been specified in~Section~\ref{s.gamespec}.
 
 We briefly recapitulate in the present case the rough argument that motivated the premise ``Incentive Inch, Outcome Mile" in Subsection~\ref{s.allocate}. Suppose that $x$ exceeds one. In playing $\standard(x)$, Maxine has more to play for than does Mina. Maxine may be tempted to outstake Mina, perhaps staking a certain constant multiple $f(x) >1$
 of the stake that Mina offers at any given turn. The resulting gameplay is a walk with a constant bias to the right, making Mina's defeat inevitable---she may as well (or better) have staked nothing. If instead it is $x^{-1}$ that exceeds one, then it is of course Mina who may be tempted by such an approach. Perhaps an argument can be fashioned along these lines to the effect that the game is competitive precisely when $x$ lies in an interval of the form $[\mu,\mu^{-1}]$ for some $\mu \in (0,1]$. This heuristic hardly lacks shortcomings, and it is quite unclear what the value of $\mu$ should be. However, the next result, which anyway follows from Theorems~\ref{t.minamarginvalues} and~\ref{t.nashabmn}, validates its conclusion in a certain sense, with the value of $\mu$ equal to $\lambda$.
\begin{theorem}\label{t.nashequil.prelim}
Recall the quantity $\lambda \in (0,1)$, which is specified and described by Theorem~\ref{t.minamarginvalues}.
For $x \in (0,\infty)$, the  game $\standard(x)$ has a time-invariant Nash equilibrium 
precisely when $x$ lies in $[\lambda,\lambda^{-1}]$.
\end{theorem}
Admitting Conjecture~\ref{c.lambda}, it is tempting to claim that this theorem validates the premise ``Incentive Inch, Outcome Mile" in a precise and quantitative way. In fact, we can surmise that the premise holds for long finite intervals more clearly than we can for~$\Z$, as we will discuss further in Section~\ref{s.playing}.

The shift operator on $\Z$ has a basic role to play as we analyse the Trail of Lost Pennies on this set.
\begin{definition}\label{d.shiftone}
Consider two time-invariant strategy pairs $\big\{ (b_i,a_i): i \in \Z \big\}$ and  $\big\{ (b'_i,a'_i): i \in \Z \big\}$. These pairs are called {\em shift equivalent} if there exists $k \in \Z$
for which   $(b_i,a_i) = (b'_{i+k},a'_{i+k})$ for all $i \in \Z$.
It is straightforward to see that an element  $\big\{ (b_i,a_i): i \in \Z \big\}$ of $\mc{S}_0^2$
lies in $\mc{N}(m_{-\infty},m_\infty,n_{-\infty},n_\infty)$ if and only if every shift equivalent element does so.
\end{definition}

Let $Q:(0,\infty) \to \N$ be such that $Q(x)$ is the maximum cardinality of a set of mutually shift inequivalent time-invariant Nash equilibria for the game $\standard(x)$ for $x \in (0,\infty)$.
 The preceding result implies that the set of $x \in (0,\infty)$ for which $Q(x) > 0$ is equal to the interval $[\lambda,\lambda^{-1}]$---which interval is non-degenerate in view of Theorem~\ref{t.minamarginvalues}(3). In the next result, we assert that a pair of shift inequivalent solutions exist when the Mina margin lies in the interval's interior.
\begin{theorem}\label{t.solutions}
For $x \in (\lambda,\lambda^{-1})$, $Q(x) \geq 2$.
\end{theorem}
We conjecture that no further time-invariant Nash equilibria exist.
\begin{conjecture}\label{c.solutions}
We have that $Q(x) = 2$ when $x \in (\lambda,\lambda^{-1})$ and $Q(x) =1$ when $x \in \{ \lambda,\lambda^{-1} \big\}$.
\end{conjecture}
This conjecture will be discussed in Section~\ref{s.conjectureroute}.

The next result describes  precise asymptotic estimates on four sequences associated to any positive \abmnmacspace solution. In light of Theorem~\ref{t.nashabmn}, it also describes decay rates for the pair of sequences given by any time-invariant Nash equilibrium.
  \begin{definition}\label{d.deltai}
Let  $(a,b,m,n)$ be an \abmnmacspace solution. For $i \in \Z$, set\hfff{phi}   
$\macphi_i = \frac{n_{i-1} - n_i}{m_i - m_{i-1}}$.
  \end{definition}
\begin{definition}\label{d.battlefield}
For an \abmnmacspace solution $(a,b,m,n)$,  the {\em battlefield index}\hfff{battlefield} is the unique value $k \in \Z$ such that $\macphi_k \in (1/3,3]$.
\end{definition}
In Lemma~\ref{l.battlefield}, we will prove the existence and uniqueness claims implicit in the last definition, thus showing that the battlefield index is well-defined. 
\begin{theorem}\label{t.ajbj}
Let $\big\{ (a_i,b_i,m_i,n_i): i \in \Z \big\}$ be a positive  \abmnmacspace solution, and let $k \in \Z$ denote its battlefield index.
\begin{enumerate}
\item There exist positive constants $A$ and $F$ such that,
for $j \geq k$, 
\begin{eqnarray*}
a_j & = & (m_k - m_{k-1})\cdot 2F \cdot  2^{2(j-k)} \exp \big\{ - 2 \cdot 2^{j-k}A     \big\}   \big( 1 + e^{-O(1) 2^{j-k}}\big) \, ; \\
b_j & = & (m_k - m_{k-1})\cdot 4F \cdot  2^{2(j-k)} \exp \big\{ - 3 \cdot 2^{j-k}A     \big\}   \big( 1 + e^{-O(1) 2^{j-k}}\big) \, ; \\
 m_j - m_{j-1} & = & (m_k - m_{k-1})\cdot F \cdot  2^{2(j-k)} \exp \big\{ - 2^{j-k}A     \big\}   \big( 1 + e^{-O(1) 2^{j-k}}\big) \, ; \, \, \,
 \textrm{and} \\
 n_{j-1} - n_j & = & (m_k - m_{k-1})\cdot 2F \cdot  2^{2(j-k)} \exp \big\{ - 2^{j-k+1}A      \big\}   \big( 1 + e^{-O(1) 2^{j-k}}\big) \, .
 \end{eqnarray*}
 The constants $A$ and $F$ may be chosen to lie in a compact interval of $(0,\infty)$ that does not depend on the choice of the solution  $\big\{ (a_i,b_i,m_i,n_i): i \in \Z \big\}$.
 The positive constant that is implicit in the $O$-notation in the four displayed expressions may be chosen independently of this solution.
\item There exist positive constants $B$ and $G$ such that,
 for $j \leq k-1$, 
\begin{eqnarray*}
a_j & = & (n_{k-1} - n_k)\cdot 4G \cdot  2^{2(k-j)} \exp \big\{ - 3 \cdot 2^{k-j}B     \big\}   \big( 1 + e^{-O(1) 2^{k-j}}\big) \, ; \\
b_j & = &  (n_{k-1} - n_k)\cdot 2G \cdot  2^{2(k-j)} \exp \big\{ - 2 \cdot 2^{k-j}B     \big\}   \big( 1 + e^{-O(1) 2^{k-j}}\big) \, ; \\
 m_j - m_{j-1} & = &  (n_{k-1} - n_k)\cdot 2G \cdot  2^{2(k-j)} \exp \big\{ - 2^{k-j+1}B     \big\}   \big( 1 + e^{-O(1) 2^{k-j}}\big) \, ; \, \, \,
 \textrm{and} \\
 n_{j-1} - n_j & = &  (n_{k-1} - n_k)\cdot G \cdot  2^{2(k-j)} \exp \big\{ - 2^{k-j}B      \big\}   \big( 1 + e^{-O(1) 2^{k-j}}\big) \, .
 \end{eqnarray*}
 The conditions on $B$ and $G$ satisfy those set out for $A$ and $F$ in the preceding part; the constant implicit in the $O$-notation satisfies the condition recorded in this part.
 \end{enumerate}
 \end{theorem}
 When $X = k$---when the counter is at the battlefield index---both players spend big to try to win the next move. For example, when $m_\infty - m_{-\infty} = n_{-\infty} - n_\infty =1$,
 so that the difference in terminal receipt between victory and defeat is one unit for each player, then values of Maxine's stake $a_k$
 lie on the interval $[0.12,0.20]$ and values of Mina's stake $b_k$ lie in  $[0.025,0.18]$. (We will shortly present explicit solutions to the \abmnmacspace equations, which validate this assertion: we may use Theorem~\ref{t.altstand}(1), for example. Maxine's expense interval is displaced to the right from Mina's, but the situation is reversed if the counter reaches $k-1$, one place to the left.) 
 These are big expenditures in a single turn of a game with infinitely many. The expenditures drop rapidly as the counter moves away from the battlefield, however. Indeed, if we write  $g_i \ll h_i$ to denote that $g_i \leq \exp \big\{ - e^{ci}  \big\} h_i$ for $i \in \nwithoutzero$ (where $c$ is some given positive constant), then Theorem~\ref{t.ajbj} implies that,
  $$
 \textrm{for} \, \, \,  i \in \nwithoutzero \, \, \, , \, \, \,  
 0 < b_{k+i} \ll a_{k+i} \ll 1 \, \, \, \textrm{and} \, \, \, 0 < a_{k-i} \ll b_{k - i} \ll 1 \, \, \, :
 $$
 to the right of the battlefield, both expenditures drop suddenly; but Maxine, eyeing victory, keeps her foot on the gas, making sure to vastly outspend the loss-cutting Mina; while to the left of the battlefield, the roles are reversed. The Battlefield Cyl Fog of Section~\ref{s.allocate} is thus exemplified.
  We also have that
 $0 < n_{k+i} - n_{k+i+1} \ll m_{k+i+1} - m_{k+i} \ll 1$ and $0 < m_{k-i} - m_{k-i-1} \ll n_{k-i-1} - n_{k-i} \ll 1$; and, by extension, 
 $$
 0 < n_{k+i} - n_\infty \ll m_\infty - m_{k+i}  \ll 1 \, \, \, \textrm{and} \, \, \,   0 < m_{k-i} - m_{-\infty} \ll n_{-\infty} - n_{k-i}  \ll 1 \, .
 $$ 
 Indeed, in the left part of the last display, which is to the right of the battlefield, Mina has essentially (but not absolutely!) thrown in the towel, and her expected payoff $n_{k+i}$ is minutely above her defeat terminal receipt of $n_\infty$. Maxine's average payoff $m_{k+i}$ is just slightly below her victory receipt of $m_\infty$, but her need to keep moving the counter rightwards provides some lower bound on the difference.
 In the right part of the display, roles are naturally reversed.
 
 The players may dread the return of the counter to the battlefield index because this is an expensive occasion for both of them. The next result, a consequence of Theorem~\ref{t.ajbj}, shows that they are typically saved from witnessing this event repeatedly when a Nash equilibrium is played.

 Let $(S_-,S_+) \in \mc{S}^2$ and $i \in \Z$. Under $\PP_{S_-,S_+}^i$,
 the {\em unanimity} event $U$ occurs when all but finitely many of the differences $X_{j+1} - X_j \in \{-1,1\},  j \in \nwithzero$,
 of the gameplay process $X:\N \to \Z$, $X_0 = i$,
 adopt a given value. Writing $U_-$ and $U_+$ for the respective events specified when the given value is $-1$ and $1$,
 the occurrence of these events correspond to victories for Mina and Maxine, and $U$ is the disjoint union of $U_-$ and $U_+$.
 \begin{theorem}\label{t.unanimity} 
 Let $(a,b,m,n)$ denote a positive \abmnmacspace solution on $\Z$ with given boundary data of the form~(\ref{e.quadruple}).
 Suppose that the solution has  battlefield index $k \in \Z$, and let $i \in \Z$.
 Let $(S_-,S_+) \in \tis^2$.
 be given by $(b,a)$, with the usual abuse of notation. 
 \begin{enumerate}
 \item We have that $\PP_{S_-,S_+}^i(U) = 1$.
 \end{enumerate}
 There exist positive constants $C$ and $c$ that may be chosen independently of the element $(S_-,S_+)$
 and the index $i \in \Z$ for which the following hold.
 \begin{enumerate}
 \setcounter{enumi}{1}
 \item If $i \geq k$ 
 then 
 $\PP_{S_-,S_+}^i(U_-) \leq C \exp \big\{- c 2^{k-i}\big\}$.
 \item If $i \leq k - 1$ then $\PP_{S_-,S_+}^i(U_+) \leq C \exp \big\{- c 2^{i-k}\big\}$.
  \end{enumerate}  
 Consider the Trail of Lost Pennies on $\Z$ with given boundary data~(\ref{e.quadruple}).
 Redefine $(S_-,S_+)$ to be an element of $\mc{S}_0^2 \cap \mc{N}(m_{-\infty},m_\infty,n_{-\infty},n_\infty)$. Writing $(b,a)$ for $(S_-,S_+)$,
 the data $(a,b,m,n)$ specified by Definition~\ref{d.quadruple}
 determines the  battlefield index. Suppose that this index is $k \in \Z$, and let $i \in \Z$.
 \begin{enumerate}
 \setcounter{enumi}{3}
 \item
 The preceding three parts remain valid in this framework. 
\end{enumerate}
 \end{theorem}   
 We see then that the player who wins a local victory at (or around) the battlefield index typically comes to entirely dominate the later moves of the game. By playing at a time-invariant Nash equilibrium, players thereby forge an implicit consensus to avoid the mutually destructive circumstance of many returns to the battlefield.
 
 We recapitulated the conclusions of~\cite{HP2022} regarding Nash equilibria in allocated-budget stake-governed games in Section~\ref{s.overviewallocated}. The formula~(\ref{e.stakefunction})  indicates that, for suitable graphs 
 and at Nash equilibrium, each budget must be spent via staking in a more-or-less regular flow, so that the concerned player is competitive throughout the lifetime of the game. 
 Game behaviour under the Battlefield Cyl Fog is very different from under the Present-Future Compromise.
 The empirical stake process of a player at any Nash equilibrium is punctuated by a few brief intense periods as the counter passes through the battlefield. In the large, the only concern for outcome is the answer to the question: on which side of the battlefield does the counter lie?

  \section{Explicit \abmnmacspace solutions}\label{s.solvingabmn}

Here we present an explicit form for all positive \abmnmacspace solutions. It is useful to begin by classifying the solutions into classes, where members of a given class differ in simple ways. If one or other player receives, or must pay, some given amount before a game begins, play will be unaffected---or at least the Nash equilibria will not be. If the unit currency is revalued before play, the outcome will be a mere scaling of all quantities. We identify \abmnmacspace solutions that differ according to translations $\transmap_{x,0}$ or $\transmap_{0,y}$ or dilations $\tau_u$ (where $x,y \in \R$ and $u \in (0,\infty)$) that correspond to such operations. If we can describe one element in each equivalence class, we will be able to describe all solutions. Equivalence classes are naturally parametrized by the positive real quantity $\phi_0 = \tfrac{n_{-1}-n_0}{m_0 - m_{-1}}$, which we call the {\em central ratio}, specified by Definition~\ref{d.deltai}. So there is a one-parameter family of essentially different positive \abmnmacspace solutions. In each equivalence class, we will distinguish two special solutions---the default solution, which has a simpler explicit formula; and the standard solution, which corresponds to a convenient choice of boundary data for the trail game. We will set up this structure and then state the explicit form of the default solution in each equivalence class. 

We consider $\Z$-indexed sequences $\genseq = \{ \genseq_i: i \in \Z \}$. A sequence is {\em monotone} if it is non-decreasing or non-increasing. A bounded monotone sequence $\genseq$
has left and right limits 
$$
\genseq_{-\infty} = \lim_{k \to \infty} \genseq_{-k} \, \, \, \,  \textrm{and} \, \, \, \,  \genseq_\infty = \lim_{k \to \infty} \genseq_k
$$ 
that are elements of~$\R$. We will specify certain bounded monotone sequences~$\genseq$ by giving one of the limiting values, $\genseq_{-\infty}$ or $\genseq_\infty$, alongside the difference sequence $\big\{ \genseq_{i+1} - \genseq_i: i \in \Z \big\}$.

Let $u \in (0,\infty)$ and $v \in \R$. For any sequence $g:\Z \to \R$, we write $u \cdot g:\Z \to \R$ for the sequence given by $(u \cdot g)_i = u \cdot g_i$.

Let $\PhiMac$ denote the space of quadruples of sequences; thus, when $(a,b,m,n) \in \PhiMac$, each component $* \in \{ a,b,m,n \}$ has the form $*:\Z \to \R$. For $u \in (0,\infty)$ and $v_1,v_2 \in \R$, define 
$\tau_u,\transmap_{v_1,v_2}:\PhiMac \to \PhiMac$ so that $\tau_u(a,b,m,n) = \big( u \cdot a,  u \cdot b, u \cdot m, u \cdot n \big)$ and $\transmap_{v_1,v_2}(a,b,m,n) = (a,b,m+v_1,n+v_2)$.

Two solutions of the \abmnmacspace equations on $\Z$ are called {\em equivalent} if one is the image of the other under a composition of the form $\tau_u \circ \transmap_{v_1,v_2}$ for such $u$, $v_1$ and $v_2$ as above. The relation of two such solutions will be denoted by~$\sim$; Proposition~\ref{p.abmnclassify} asserts that~$\sim$ is indeed an equivalence relation.

Let $(a,b,m,n)$ be an \abmnmacspace solution on $\Z$. Recall from Definition~\ref{d.abmn} that the solution's {\em Mina margin}\hfff{centralratio} is defined to be $\tfrac{n_{-\infty} - n_\infty}{m_\infty - m_{-\infty}}$. 
The solution's {\em central ratio}\hfff{centralratio} $\centralratio$ is set equal to $\frac{n_{-1} - n_0}{m_0 - m_{-1}}$. The solution\hfff{standarddefault} is called 
 {\em standard} if $m_{-\infty} = 0$, $n_\infty = 0$ and $m_\infty = 1$.  It is called  {\em default} if  $m_{-\infty} = 0$, $n_\infty = 0$ and $m_0 - m_{-1} = 1$. 
 Compatibly with the usage of Definition~\ref{d.standard}, the Mina margin of a standard solution equals $n_{-\infty}$; note further that the central ratio of a default solution equals $n_{-1} - n_0$.

\begin{proposition}\label{p.default}
For any default \abmnmacspace solution, the value of the central ratio  $\centralratio$ lies in $(0,\infty)$. For any $x \in (0,\infty)$, there is exactly one default solution for which $\centralratio$ equals $x$.
\end{proposition}
\begin{proposition}\label{p.abmnclassify}
\leavevmode
\begin{enumerate}
\item The space of \abmnmacspace solutions is partitioned into equivalence classes by the relation~$\sim$.
\item Each equivalence class contains 
a unique standard solution and 
a unique default solution.
\end{enumerate}
\end{proposition}
Propositions~\ref{p.default} and~\ref{p.abmnclassify} provide a natural labelling of \abmnmacspace solution equivalence classes: 
 any given class is labelled by the value of the central ratio of the unique default solution in the class. The labelling parametrizes the equivalence classes by a copy of $(0,\infty)$.

According to the latter assertion of Proposition~\ref{p.default}, there is a unique default solution to the \abmnmacspace equations whose central ratio equals a given value $x \in (0,\infty)$. The next definitions will enable us to record  the form of this solution in Theorem~\ref{t.defaultexplicit}.

\begin{definition}\label{d.acs}\hfff{smap}
Set $\macomega:(0,\infty) \to (1,\infty)$, $\macomega(x) = \sqrt{8x+1}$. Writing $\omega = \omega(x)$, we further set
$$
\macsee(x)   =    \frac{(\macomega + 3)^2}{16} \, \, , \, \,
\macdee(x)  =  \frac{(\macomega + 3)^2}{8(\macomega + 1)} \, \, \, \, \textrm{and} \,\, \, \, s(x)  =  \frac{(\macomega - 1)^2}{4(\macomega + 7)} \,\,\,\, \textrm{for} \, \, \, \, x \in (0,\infty) \, .
$$
\end{definition}
\begin{definition}\label{d.stabc}
Let $s_{-1}:(0,\infty) \to (0,\infty)$ be given by $s_{-1}(x) = 1/s(1/x)$. 
We now define a collection of functions $s_i:(0,\infty) \to (0,\infty)$ indexed by $i \in \Z$. We begin by setting $s_0(x) = x$ for $x \in (0,\infty)$. 
We then iteratively specify that $s_i(x) = s \big( s_{i-1}(x) \big)$ and  $s_{-i}(x) = s_{-1} \big( s_{-(i-1)}(x) \big)$ 
for $i \in \nwithoutzero$ and $x \in (0,\infty)$. Note that $s_1$ equals $s$ and that the two specifications of $s_{-1}$ coincide. 

Set $\macsee_j,\macdee_j:(0,\infty) \to (0,\infty)$, $j \in \Z$, by means of $\macsee_j(x) =  \macsee (s_j(x))$ and  $\macdee_j(x) =  \macdee (s_j(x))$.
\end{definition}
To get a sense of the maps $s_i$, $i \in \Z$, a few points are worth noting. First, as we will see in Proposition~\ref{p.sminusone}, $s_{-1}$ is the inverse of $s$. Second,
as Lemma~\ref{l.acsfacts}(5) attests, $s(x) < x$ for $x \in (0,\infty)$; the orbit $s_i(x)$ thus decreases or increases from $x$ as $i$ grows to the right or the left.
And third, note that $s(3) = 1/3$. In view of the second point, we see that $(0,\infty) = \cup_{k \in \Z} \, s_k [1/3,3)$ is a partition whose interval elements are arranged in decreasing order in the index~$k$.
\begin{definition}\label{d.zdefault}
For a sequence $h$, we may naturally write $\prod_{i=0}^k h_i = h_0 \cdots h_k$ for $k \in \N$. A convenient device extends this notation to cases where $k \in \Z$ is negative: we set 
$$\prod_{i=0}^k h_i \, = \, \begin{cases}
  \, 1  &  \text{for $k=-1$} \\
 \, h_{k+1}^{-1} \cdots h_{-1}^{-1}  &  \text{for $k \leq -2$} \, .
\end{cases}
$$
Let $x \in (0,\infty)$. This parameter will index four real-valued sequences 
$$
\adef(x),\bdef(x),\mdef(x),\ndef(x):\Z \to (0,\infty)
$$ 
which we denote in the form $\big\{ *^{\rm def}_i(x): i \in \Z \big\}$ for $* \in \{a,b,m,n\}$. 

We begin by specifying $\mdef(x):\Z \to \R$. This is the increasing sequence such that 
$$
\mdef_{-\infty}(x) = 0 \, , \, \, \, \, \textrm{and} \, \, \, \, \mdef_{k+1}(x)- \mdef_k(x) \, = \, \prod_{i=0}^k \big( \macsee_i(x) - 1 \big) \, \, \,  \textrm{for $k \in \Z$} \, . 
$$
Note that $\mdef_0(x) - \mdef_{-1}(x) = 1$ in view of the notation for products. 

Next we set  $\ndef(x):\Z \to \R$. This is the decreasing sequence with 
$$
\ndef_\infty(x) = 0 \, , \, \, \, \,   \textrm{and} \, \, \, \, \ndef_k(x)- \ndef_{k+1}(x)  \, = \, x \prod_{i=0}^k \big( \macdee_i(x) - 1 \big) \, \, \,  \textrm{for $k \in \Z$} \, . 
$$
Note that $n_{-1}(x) - n_0(x) = x$.

To specify $\adef(x),\bdef(x):\Z \to (0,\infty)$, we set 
$$
M_i(x) = \mdef_{i+1}(x) - \mdef_{i-1}(x) \, \, \, \, \textrm{and} \, \, \, \, N_i(x) = \ndef_{i-1}(x) - \ndef_{i+1}(x)
$$ 
for $i \in \Z$. For such~$i$, we take 
$$
\adef_i(x) = \frac{M_i(x)^2 N_i(x)}{\big(M_i(x)+N_i(x)\big)^2} \, \, \, \, \textrm{and} \, \, \, \, 
\bdef_i(x) = \frac{M_i(x) N_i(x)^2}{\big(M_i(x)+N_i(x)\big)^2} \, .
$$
\end{definition}

\begin{theorem}\label{t.defaultexplicit}
Let $x \in (0,\infty)$.
The unique default \abmnmacspace solution with $\centralratio = x$ is the quadruple $\big( \adef_i(x),\bdef_i(x),\mdef_i(x),\ndef_i(x): i \in \Z \big)$ specified in Definition~\ref{d.zdefault}.
\end{theorem}

For $x \in (0,\infty)$, we write $\mc{C}(x)$ for the equivalence class of \abmnmacspace solutions that contains the element $\big( \adef_i(x),\bdef_i(x),\mdef_i(x),\ndef_i(x): i \in \Z \big)$. 
Let $\big(\asta_i(x),\bst_i(x),\mst_i(x),\nst_i(x) : i \in \Z \big)$ denote the unique standard solution in $\mc{C}(x)$.


{\em Remark.}
Let $x \in (0,\infty)$. Set $Z(x) = \mdef_\infty(x)$; which is to say, $Z(x) = \sum_{k \in \Z} \prod_{i=0}^k \big( c_i(x) - 1 \big)$. By Theorem~\ref{t.ajbj} for the $m$-differences, $Z(x) < \infty$.
It is straightforward that
\begin{equation}\label{e.remark}
\Big( \asta_i(x),\bst_i(x),\mst_i(x),\nst_i(x)  : i \in \Z \Big) \, \, \, \, \textrm{equals} \, \, \, \, Z(x)^{-1} \cdot \Big( \adef_i(x),\bdef_i(x),\mdef_i(x),\ndef_i(x) : i \in \Z  \Big) \, .
\end{equation}
\section{The Mina margin map}\label{s.mmm}
According to Theorem~\ref{t.defaultexplicit}, the central ratio $\phi_0$ is a convenient parameter for indexing \abmnmacspace solution equivalence classes. And Theorem~\ref{t.nashequil.prelim} tells us that the Mina margin is a fundamental parameter for locating Nash equilibria in the trail game. The map $(0,\infty) \to (0,\infty)$ from equivalence class index to the Mina margin of any member solution is a natural object that we will use to organize and prove results. We call this function the {\em Mina margin map}\hfff{MMM}. Here, we define it, and state basic properties in Theorem~\ref{t.relativereward}. Theorem~\ref{t.nashequil} shows how to solve the trail game with given boundary data by finding time-invariant Nash equilibria indexed by the map's level sets. Theorem~\ref{t.phithetainverse} states that a reparametrization of the Mina margin map's domain leads to a periodic form for the map that commutes with the shift operator on $\Z$.
\begin{definition}\label{d.r}
Let the Mina margin map $\minamm:(0,\infty) \to (0,\infty)$ be given by $\minamm(x) =\nst_{-\infty}(x)$ for  $x \in (0,\infty)$. Namely, $\minamm(x)$ is 
the Mina margin of  $\big(\asta_i(x),\bst_i(x),\mst_i(x),\nst_i(x): i \in \Z \big)$. 
\end{definition}
\begin{theorem}\label{t.relativereward}
\leavevmode
\begin{enumerate}
\item The function $\minamm:(0,\infty) \to (0,\infty)$ satisfies $\minamm(s(x)) = \minamm(x)$ for $x \in (0,\infty)$.
\item The function $\minamm$ is continuous on $(0,\infty)$ and satisfies
$$
 \minamm(x) \, \, = \, \, \Bigg(  \sum_{k \in \Z} \, \, \prod_{i=0}^k \big( c_i(x) - 1 \big)  \Bigg)^{-1} \, \cdot \, x \, \sum_{k \in \Z} \, \, \prod_{i=0}^k \big( \macdee_i(x) - 1 \big) \, .
$$ 
\item The range $\minamm(0,\infty)$ takes the form $[\lambda,\lambda^{-1}]$, where $\lambda \in (0,0.999904]$ is specified in Theorem~\ref{t.minamarginvalues}.
\end{enumerate}
\end{theorem}


\begin{theorem}\label{t.nashequil}
Let $x \in [\lambda,\lambda^{-1}]$.  Set $X = \big\{ z \in (0,\infty): \minamm(z) = x \big\}$, and let $Y = X \cap (1/3,3]$, so that, as noted after Definition~\ref{d.stabc},  $X = \cup_{k \in \Z} s_k(Y)$.
\begin{enumerate}
\item
The collection of time-invariant Nash equilibria in the game $\standard(x)$ is  given by the set of maps
$$
\Z \to (0,\infty)^2: i \to \big(\bst_i(z),\asta_i(z) \big)
$$
indexed by $z$ in $X$.
\item Alternatively, 
this collection is the set of maps
$$
 \Z \to (0,\infty)^2: i \to \big(\bst_{i+j}(x),\asta_{i+j}(x) \big) \, ,
$$
where now the index $(x,j)$ ranges over $Y \times \Z$.
\end{enumerate}
\end{theorem}

We now develop the notation for the symbolic shift map that was mooted in Definition~\ref{d.shiftone}.  
\begin{definition}\label{d.shiftmap}
 We let $\mc{S}_1$ denote the left shift by one place: this is the map that sends the space of quadruples $(\R^4)^\Z$ to itself by the action 
 $$
  \mc{S}_1 \big\{ (a_i,b_i,m_i,n_i): i \in \Z \big\} =  \big\{ (a_{i+1},b_{i+1},m_{i+1},n_{i+1}): i \in \Z \big\} \, .
 $$
 By iterating this map, we specify the left shift $\mc{S}_k$ by $k$ places, for $k \geq 2$;
 and by specifying $\mc{S}_{-1} = \mc{S}_1^{-1}$ and iterating, we specify the right shift $\mc{S}_{-k}$ by $k$ places, for $k \geq 1$.
 \end{definition}

 The function $s:(0,\infty) \to (0,\infty)$ has been specified formulaically in Definition~\ref{d.acs}, but its fundamental property, stated in Proposition~\ref{p.alphagammaess},  is that $s(\phi_0) = \phi_1$ for any positive \abmnmacspace solution. Writing $x = \phi_0$ and iterating, we then have $s_k(x) = \phi_k$. Just as $\mc{S}_k$ is the $k$-left shift on the space of quadruples, so $s$ plays this role in relation to the $\phi$-sequence $\{ \phi_i: i \in \Z \}$. The next result is equivalent to this fact.
  \begin{proposition}\label{p.shift}
  For $x \in (0,\infty)$ and $k \in \Z$, 
  \begin{equation}\label{e.shift}
    \mc{S}_k \big( \asta(x), \bst(x), \mst(x), \nst(x)  \big) \, = \,  \Big( \asta\big(s_k(x)\big), \bst \big(s_k(x)\big), \mst \big(s_k(x)\big), \nst \big(s_k(x)\big)  \Big) \, . 
 \end{equation}
  \end{proposition}  
 Theorem~\ref{t.relativereward}(1) leads directly to $\minamm(s_k(x)) = \minamm(x)$ for $x \in (0,\infty)$ and $k \in \Z$. To understand the map $\minamm$, we see then that the asymptotics in highly positive and negative $k$
 of the orbits $s_k(x)$ are important. As we will see in Lemma~\ref{l.acsfacts}(2,3), $s(x) \sim x^2/2$ for $0 < x \ll 1$ and $s(x) \sim 2^{-1/2} x^{1/2}$ for $x \gg 1$.
 Thus the forward orbit $s_k(x)$, $k \to \infty$, converges rapidly to zero, while the backward orbit, $k \to -\infty$, grows quickly towards infinity. 
 
We now undertake a change of coordinates of the Mina margin map $\minamm:(0,\infty) \to (0,\infty)$. The domain $(0,\infty)$
will be identified with $\R$ by an increasing bijection $\theta^{-1}$. The goal of the coordinate change is to ensure that the original action $(0,\infty) \to (0,\infty): x \to s_1(x)$ becomes the map $\R \to \R: x \to x-1$. The action of the symbolic sequence shift $\mc{S}_1$ on the $x$-variable, as stated in Proposition~\ref{p.shift}, comes to correspond to a left shift by a unit in the new real variable. This leads to an attractive representation of the Mina margin map in the guise $\R \to (0,\infty): x \to \minamm\big( \theta^{-1}(x)\big)$.
\begin{definition}
Let $q:[1/3,3) \to [0,1)$ be an increasing surjection; for definiteness, we may take $q(x) = 3(x -1/3)/8$.
We specify $\theta:(0,\infty) \to \R$ so that, for $x \in (0,\infty)$, $\theta(x) = k+q\big(s_k(x)\big)$, where $k \in \Z$ is the unique integer such that $s_k(x) \in [1/3,3)$.
Since $\theta:(0,\infty) \to \R$ is an increasing surjection, the inverse $\theta^{-1}: \R \to (0,\infty)$ is well defined. We may thus represent the Mina margin map after domain coordinate change by the function $\psi$, where  
$$
\psi: \R \to (0,\infty) \, \, \, \, , \, \, \, \, \psi(x) = \minamm \big( \theta^{-1}(x) \big) \, .
$$
We define the {\em standard solution} map $\solution:\R \to (\R^4)^\Z$,
$$
 \solution(x) \, = \, \Big( \asta\big( \theta^{-1}(x)\big), \bst\big( \theta^{-1}(x)\big), \mst\big( \theta^{-1}(x)\big), \nst\big( \theta^{-1}(x)\big)   \Big) \, \, \, \, \textrm{for $x \in \R$} \, .
$$
\end{definition}
For $u \in (0,\infty)$ and $j \in \Z$, $\theta\big(s_{-j}(u)\big) - \theta(u) = j$.
For $z = \Theta(1)$, 
the value of $\theta^{-1}(z+k)$ thus tracks that of $s_{-k}(z)$
as $k$ rises, either by growing to infinity (if $k$ is positive); or by decaying to zero (if $k$ is negative). To understand the transformation $\theta^{-1}$, it is thus useful to  introduce a simple explicit function  $\PhiMac:\R \to (0,\infty)$ which is designed so that $s_{-k}(z)$ grows or decays roughly as does $\PhiMac(k)$ for~$\vert k \vert$ large; here $z \in [1/3,3)$, say. 
  
Let $\rm Sign:\R \to \{-1,1\}$ equal ${\rm Sign}(x) = {\bf 1}_{x \geq 0} -  {\bf 1}_{x < 0}$.
Then set 
$$
\PhiMac:\R \to (0,\infty) \, \, \, \, , \, \, \, \, \PhiMac(x) = 2^{{\rm Sign}(x)(2^{\vert x \vert} -1)} \, .
$$

We now present our result concerning the Mina margin map after domain coordinate change. The transformed function $\psi$ is periodic, of unit period; symbolic shift by one place corresponds to unit translation of the domain; and coordinate change asymptotics are, crudely at least, described by $\PhiMac$.

\begin{theorem}\label{t.phithetainverse}
\leavevmode
\begin{enumerate}
\item For $x \in \R$, $\psi(x+1) = \psi(x)$. 
\item 
For $x \in (0,\infty)$ and $k \in \Z$, 
$$
\solution(x+k) = \mc{S}_{-k} \circ \solution (x) \, .
$$
\item There exists a positive constant $C$ such that, for $z \geq 0$,
$$
 2^{2^{z-C}} \leq  \theta^{-1}(z) \leq  2^{2^{z+C}} \, ;
$$
and, for $z < 0$,
$$
 2^{-2^{\vert z \vert +C}} \leq  \theta^{-1}(z) \leq  2^{-2^{\vert z \vert -C}} \, .
$$
\end{enumerate}
\end{theorem}
The map $\PhiMac$ is a simple and explicit surrogate for $\theta^{-1}$, and the transformed Mina margin map $\R \to (0,\infty): x \to \minamm\big(\PhiMac(x)\big)$
shares the periodicity property of $\psi$ in Theorem~\ref{t.phithetainverse}(1) up to a domain perturbation that decays rapidly away from zero. And this surrogate has a more practical version, in which the Mina margin map $\minamm$ is replaced by a counterpart for a trail that is a finite interval, rather than all of $\Z$. These counterpart Mina margin maps $\minamm_{j+1,k+1}$ will be presented in the next section. Plots of several of these maps, indexed by different finite trails, appear in Figure~\ref{f.tmmm}.

\section{The Trail of Lost Pennies on a finite interval}\label{s.finite}

 Even if much of our focus lies with the trail game in the infinite setting, with gameboard~$\Z$,
 it is instructive to introduce and discuss the game whose trail is a finite interval.
 This is the setting used in  economists' treatments of tug-of-war that we surveyed in Section~\ref{s.economics}. 
And it is a more practical setting if two people are to play the game, taking decisions turn-by-turn because, at least for short intervals, the game will end (by the token reaching one end of the interval or the other) in a limited number of moves. The theoretical aspects of the game---time-invariant Nash equilibria; \abmnmacspace solutions and their standard solutions; the Mina margin map---share many basic aspects between the infinite and finite settings. The finite setting permits important objects, such as the Mina margin map, to be plotted in Mathematica, and such investigation has informed several of our main results (in the infinite setting). Our goal then in this section is to communicate the principal aspects of the finite setting so that the reader can interpret pertinent Mathematica plots and understand how these suggest some of our principal results and conjectures. We will also present a conjecture concerning the number of time-invariant Nash equilibria in a symmetric version of the finite game; we will seek to explain why we believe it during the section.  The section contains one result, Proposition~\ref{p.rkvalues}, which we will use and whose proof appears in Section~\ref{s.rolereversal}.
Our basic aim is heuristic, however, and at times our presentation will be informal.

\subsection{Gameplay, strategies and Nash equilibria for the finite trail}\label{s.gsn}
 
Let $j,k \in \nwithzero$. The Trail of Lost Pennies with trail (or gameboard) $\llbracket -j-1,k+1\rrbracket$ 
is specified by 
$$
\big( m_{-j-1},m_{k+1},n_{-j-1},n_{k+1} \big) \in \R^4 \, , 
$$
boundary data  on which the conditions $m_{-j-1} < m_{k+1}$ and $n_{k+1} < m_{-j-1}$ are imposed.
Begun from $\ell$, an element in the field of open play $\llbracket -j,k \rrbracket$, gameplay is a stochastic process $X: \llbracket 0, K \rrbracket \to \llbracket -j-1,k+1 \rrbracket$,  $X_0 = \ell$, where 
$$
K \, = \, \inf \, \Big\{ \, i \in \nwithoutzero: X_i \in \{ -j-1,k+1 \} \, \Big\} \, .
$$
Indeed, with Mina and Maxine playing to the left and right,
the game will end with victory to these respective players when the counter arrives at $-j-1$ or at $k+1$. 

The gameplay is specified by a strategy pair, where 
a strategy is a map $S: \llbracket -  j,k \rrbracket \times \nwithoutzero   \to [0,\infty)$. 
The construction of $X$ from a given location $X_0 = \ell \in \llbracket -j,k \rrbracket$ coincides with that explained in Section~\ref{s.gamespec}, where instances of the trail $\Z$ are replaced by $\llbracket -j,k \rrbracket$, it being understood that the construction stops when $X$ arrives in $\{ - j-1,k+1 \}$.

A strategy $S$ for which $S(\ell,i)$ is independent of $i \in \nwithoutzero$ for all $\ell \in \llbracket -j,k \rrbracket$
is said to be {\em time-invariant}. 
Let $\mc{S}[j,k]$ denote the space of strategies.
For a strategy pair $(S_-,S_+) \in \mc{S}[j,k]$, 
we may reuse notation from the $\Z$-indexed trail game, 
and speak of the law $\PP_{S_-,S_+}^i$ of gameplay $X:\N \to \llbracket -j-1,k+1 \rrbracket$, $X_0 = i$, governed by the pair $(S_-,S_+)$, and stopped on arrival in $\{ -j-1,k+1 \}$. 
Counterpart to~(\ref{e.minapayoff}) and~(\ref{e.maxinepayoff}) are the  $\PP_{S_-,S_+}^i$-almost sure payoff identities
\begin{equation}\label{e.finitepayoff}
P^{j,k}_- \, = \,  - \sum_{i = 1}^\infty C^{j,k}_-(i) \, \, + \, \, T^{j,k}_-  \, \, \, \, \textrm{and} \, \, \, \,
P^{j,k}_+ \, = \,  - \sum_{i = 1}^\infty C^{j,k}_+(i) \, \, + \, \, T^{j,k}_+ \, ,
\end{equation}
where 
the cost  $C^{j,k}_\pm(i)$  incurred to each player at the $i$\textsuperscript{th} turn, $i \in \nwithoutzero$, equals $S_{\pm}(X_{i-1},i)$, as in the original case. 
To specify the terminal payments $T^{j,k}_\pm$, we permit $E_-$ to denote the event that $X$ arrives at the vertex $-j-1$ at some positive time, and~$E_+$ to denote the event that this process instead reaches $k+1$
at some such time. We then adopt~(\ref{e.terminalmina}) and~(\ref{e.terminalmaxine}) for $T_\pm^{j,k}$, where $m_*$ and $n_*$ denote given real values that satisfy $m_* \leq m_{-j-1}$
and $n_* \leq n_{k+1}$.

Definitions concerning Nash equilibria continue to be specified as they are at the end of Section~\ref{s.gamespec}. 
A collection of quadruples $\big\{ (a_i,b_i,m_i,n_i):i \in \llbracket -j,k \rrbracket \big\}$ is associated to any element $\big\{ (b_i,a_i): i \in \llbracket -j,k \rrbracket \big\}$
by Definition~\ref{d.quadruple} after evident changes in notation have been made.

\subsection{The \abmnmacspace equations}
Recall Definition~\ref{d.abmn}. Let $j,k \in \nwithzero$.

The \abmnmacspace system on $\llbracket -j,k \rrbracket$ is the set of equations~\abmnmac$(1,2,3,4)$  in the real variables $a_i$, $b_i$, $m_i$ and $n_i$, where the index $i$ varies over 
$\llbracket -j,k \rrbracket$. 
These equations 
refer to the components of the quadruple $\big( m_{-j-1},m_{k+1},n_{-j-1},n_{k+1} \big) \in \R^4$ which acts as boundary data and for which we suppose a fixed value that satisfies $m_{-j-1}< m_{k+1}$ and $n_{-j-1} > n_{k+1}$. Similarly to Definition~\ref{d.abmn}, a solution is {\em positive} 
if $a_i$ and $b_i$ exceed zero for $i \in \llbracket -j,k \rrbracket$.

\subsection{A result and a conjecture for the finite trail}\label{s.resultandconjecture}

The basic relation between time-invariant Nash equilibria $\big\{ (b_i,a_i): i \in \llbracket -j,k \rrbracket \big\}$
and positive ABMN solutions  $\big\{ (a_i,b_i,m_i,n_i): i \in \llbracket -j,k \rrbracket \big\}$ embodied in Theorem~\ref{t.nashabmn} is maintained.

The trail game on $\llbracket -j-1,k+1 \rrbracket$ is in its {\em standard form} when its boundary data satisfies  $m_{-j-1} = n_{k+1} = 0$ and $m_{k+1}=1$. This class of games is thus parametrized by the Mina margin $n_{-j-1} \in (0,\infty)$. If further $n_{-j-1} =1$, then we speak of the  {\em symmetric} standard game.
Likewise a solution of the ABMN equations on $\llbracket -j,k \rrbracket$ is  {\em standard} when   $m_{-j-1} = n_{k+1} = 0$ and $m_{k+1}=1$. The space of standard solutions may be parametrized by
the central ratio $\centralratio = \tfrac{n_{-1}- n_0}{m_0 - m_{-1}} \in (0,\infty)$. The Mina margin map\hfff{MMM.finite}
 $\minamm_{j+1,k+1}:(0,\infty) \to (0,\infty)$ associates to $x \in (0,\infty)$ the value of the Mina margin~$n_{-j-1}$
of the unique standard \abmnmacspace solution on $\llbracket -j,k \rrbracket$ for which $\centralratio =x$. 

Standard solutions may be computed explicitly, similarly as was~(\ref{e.remark})
in the infinite setting. To obtain the standard solution on $\llbracket -j,k \rrbracket$ with $\centralratio = x \in (0,\infty)$, we start with the restriction of the default solution from Theorem~\ref{t.defaultexplicit} to $\llbracket -j-1, k+1 \rrbracket$. By adding a suitable constant to each $m$-term, and another such to each $n$-term, and then multiplying the result by a suitable scaling factor, we obtain a standard solution whose $\centralratio$ remains equal to $x$ because the additions and the scaling leave this value unchanged.
We thus see that, for $x \in (0,\infty)$,
\begin{equation}\label{e.minammfinite}
 \minamm_{j+1,k+1}(x) \, = \, \frac{n_{-j-1} - n_{k+1}}{m_{k+1} - m_{-j-1}} \, , 
\end{equation}
where $\big\{  (a_i,b_i,m_i,n_i): i \in \llbracket -j,k \rrbracket \big\}$ is any \abmnmacspace solution on $\llbracket -j,k \rrbracket$ such that  $\tfrac{n_{-1}- n_0}{m_0 - m_{-1}} = x$.

\begin{figure}[htbp]
\centering
\includegraphics[width=0.75\textwidth]{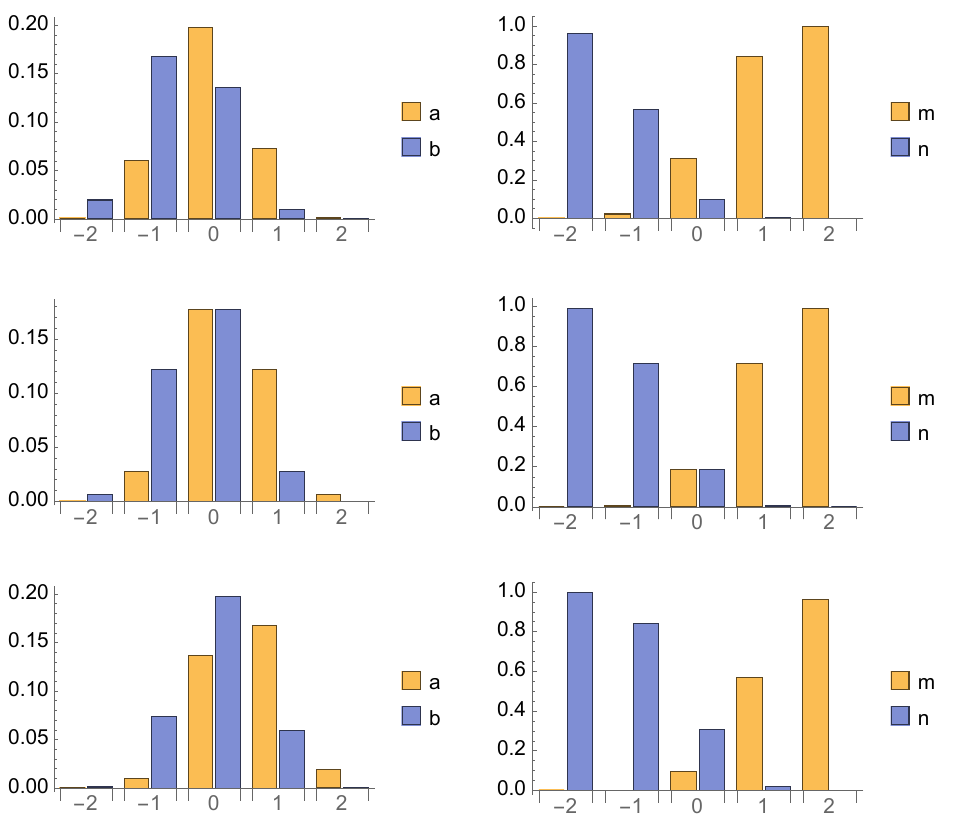}
\caption{The shortest trail with non-unique Nash equilibria for at least some boundary conditions has length six, with five sites in open play. The values 
$x_1 = 1.63$, $x_2 = 3$ and $x_3 = 5.64$ approximate the three solutions of $\minamm_{3,3}(x) = 1$. (There is no error for $x_2$.)
The $(a,b)$ and $(m,n)$ data on $\llbracket - 2,2 \rrbracket$ for the standard solution on $\llbracket -3,3 \rrbracket$ corresponding to $x_1$ appears in the top row; to $x_2$ in the middle; and to $x_3$ in the lower row. The left column thus depicts the three Nash equilibria in standard symmetric game on the shortest trail for which this game may be expected to have several equilibria.
Note that the $x_3$-solution is formed from the $x_1$-solution by role-reversal: that is, by interchanging the roles of $a$ and $b$, and of $m$ and $n$, and by reflecting in the origin.
}\label{f.threenash}
\end{figure}

The trail game on trails $\llbracket -k,k \rrbracket$ of even length differs from that on trails $\llbracket -k-1,k \rrbracket$ of odd length, because the trails in the two classes are reflection symmetric about different objects (the vertex $0$ or the edge $\llbracket -1,0 \rrbracket$). The next result records outcomes of these symmetries for the finite trail Mina margin map.
 \begin{proposition}\label{p.rkvalues}
Let $k \in \nwithoutzero$ and $x \in (0,\infty)$.
\begin{enumerate}
\item We have that $\minamm_{k,k}(x) \cdot \minamm_{k,k} \big( 1/s(x)\big) = 1$.
\item And that $\minamm_{k+1,k}(x) \cdot \minamm_{k+1,k}(x^{-1}) = 1$.
\end{enumerate}
 \end{proposition}
 Here is our conjecture concerning the symmetric form of the finite trail game.
\begin{conjecture}\label{c.tine}
Consider the Trail of Lost Pennies on $\llbracket -j-1,k+1 \rrbracket$ in its symmetric standard form.
The number of time-invariant Nash equilibria equals $\max \big\{ 2(j+k) - 5,1 \big\}$.
\end{conjecture} 
Roughly similar is Harris and Vickers'  numerically based claim (recalled in Subsection~\ref{s.harrisvickers}) of order $2n$ equilibria in their game. 
We mention also that the number of Nash equilibria is odd in almost all finite games~\cite{Wilson1971}.

\begin{figure}[htbp]
\centering
\includegraphics[width=0.75\textwidth]{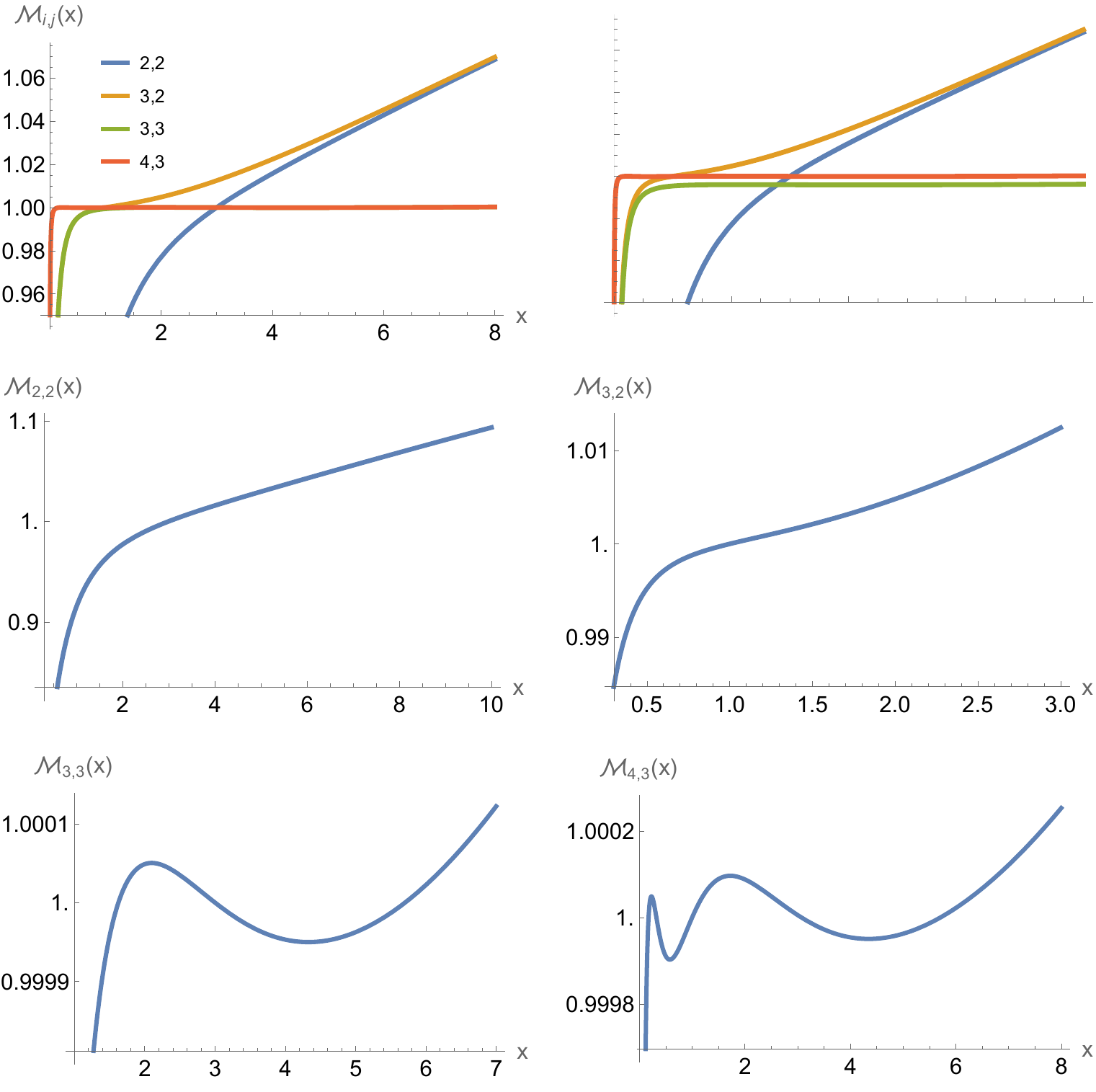}
\caption{Four finite-trail Mina margin maps $\minamm_{j+1,k+1}:(0,\infty) \to (0,\infty)$ are depicted, for values of $j+k$ in $\llbracket 2,5 \rrbracket$.
 {\em Top left.} The four functions are plotted together. {\em Top right.} This is a `Tube Map' of the left-hand graph (a distorted but practical depiction), in which the green curve has been artificially displaced to separate it, so that the viewer may watch the different lines as they run.  \\
The green and red curves seem to suggest that the curves converge to the constant function one as $j+k$ rises, but this impression is false. Indeed, the middle and lower graphs plot the four functions in turn, each on a scale that shows the finer journey of the map as it rises through height one. The maps lose injectivity in the $(j+1,k+1)$-index change $(3,2) \to (3,3)$.  
}\label{f.mmm}
\end{figure}

\begin{figure}[htbp]
\centering
\includegraphics[width=1.05\textwidth]{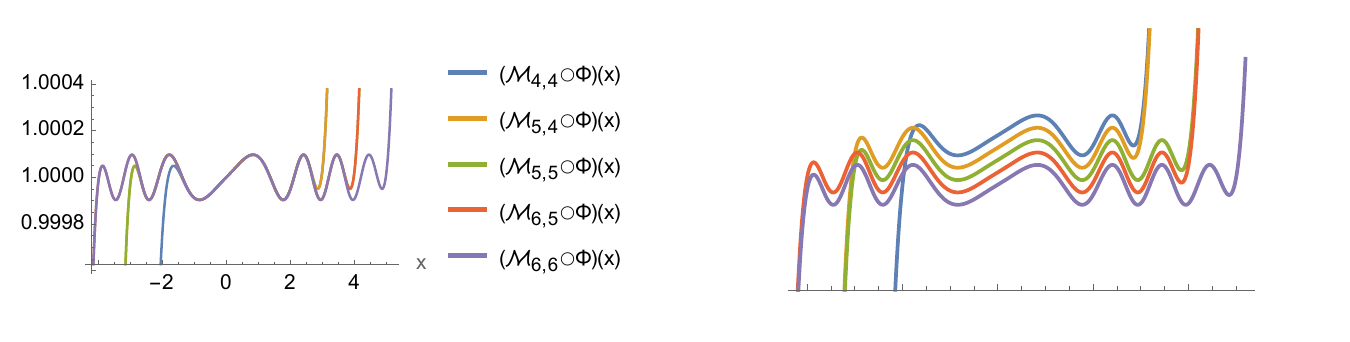}
\caption{{\em Left:} Five $\PhiMac$-transformed finite-trail Mina margin maps $\minamm_{j+1,k+1} \circ \PhiMac:\R \to (0,\infty)$ are depicted, for increasing values of $j+k$ in $\llbracket 6,10 \rrbracket$.
The graphs join and leave a shared highway, which is (up to visually negligible discrepancies)
the graph of the limiting transformed map  $\minamm \circ \PhiMac:\R \to (0,\infty)$. {\em Right:} As in Figure~\ref{f.mmm}(top,right), curves have been artificially displaced so that their routes can be clearly seen.}\label{f.tmmm}
\end{figure}

We offer an explanation of why we believe Conjecture~\ref{c.tine}. By a counterpart to Theorem~\ref{t.nashabmn} (which we have roughly indicated), any time-invariant Nash equilibrium of the symmetric trail game on $\llbracket -j-1,k+1 \rrbracket$ corresponds to a positive \abmnmacspace solution on $\llbracket -j,k \rrbracket$.
This solution must have $m_{-j-1}=n_{k+1}=0$ and $m_{k+1}=1$, as well as $n_{-j-1}=1$. That is, the solution must be standard, and it must satisfy $\minamm_{-j-1,k+1}(x) = 1$, where $x \in (0,\infty)$ is the solution's value of $\centralratio$. We may thus obtain the set of time-invariant Nash equilibria by recording, for each solution $x \in (0,\infty)$ of the equation $\minamm_{-j-1,k+1}(x) = 1$,
the reverse-ordered $(a,b)$-component pair of the unique standard \abmnmacspace solution on $\llbracket -j,k \rrbracket$ whose $\centralratio$-value equals $x$. The case for Conjecture~\ref{c.tine}
thus rests on advancing an argument for the equality
\begin{equation}\label{e.finitenash}
\# \big\{ x \in (0,\infty): \minamm_{j+1,k+1}(x) = 1 \big\} \, = \, \max \big\{ 2(j+k) - 5,1 \big\} \, .
\end{equation}
 Plots of several finite-trail Mina margin maps $(0,\infty) \to (0,\infty): x \to \minamm_{j+1,k+1}(x)$ led to the conjecture. The pattern begins to emerge in the four plots displayed in Figure~\ref{f.mmm}, for which $j+k \in \llbracket 2,5 \rrbracket$. To see the pattern continue, we need higher values of $j+k$. For these, a suitable device is  the finite-trail $\PhiMac$-transformed Mina margin map  $(0,\infty) \to (0,\infty): x \to (\minamm_{j+1,k+1} \circ \PhiMac)(x)$ mentioned at the end of Section~\ref{s.mmm}: see Figure~\ref{f.tmmm} for four depictions, where $j+k \in \llbracket 6,10 \rrbracket$.

\section{Some further formulas}\label{s.formulas}

In this article, we study a new game, presenting conjectures as well as results.
We have derived some formulas of which we do not make use, and we choose to present them as our final results in the introduction because they appear interesting and could be of value in further study of the Trail of Lost Pennies. First we state Theorem~\ref{t.altstand}, an alternative explicit form for standard \abmnmacspace solutions. Then we present the  \amacspace system, which is a closed $\Z$-indexed set of equations that we find in Theorem~\ref{t.symmetric} to describe the $a$- (or $b$-)variables in any time-invariant Nash equilibrium in the special case of the game with a symmetric form of boundary data.

\subsection{Alternative formulas for standard solutions and their Mina margins}

Recall the function $Z:(0,\infty) \to (0,\infty)$, $Z(x) = \mdef_\infty(x) =  \sum_{k \in \Z} \prod_{i=0}^k \big( c_i(x) - 1 \big)$, from the remark that concludes Section~\ref{s.solvingabmn}.
\begin{theorem}\label{t.altstand}
Let $f$, $g$ and $h$ mapping $(0,\infty)$ to itself be specified by 
$$
 f(x) = \frac{x c(x)d(x)}{\big( c(x) + x d(x) \big)^2} \, \,  ; \, \, \, g(x) = Z(x)^{-1}c(x) f(x) \, \,  ; \, \textrm{and} \, \, \, h(x) =  Z(x)^{-1} x d(x) f(x) \, .
$$ 
Let $x \in (0,\infty)$.
\begin{enumerate}
\item For $k \in \Z$, $\asta_k(x) = g \big( s_k(x) \big)$ and $\bst_k(x) = h \big( s_k(x) \big)$.
\item
For $j,k \in \Z$ such that $j < k$,
$$
 \mst_k(x) - \mst_j(x) \, = \, \sum_{i = j+1}^k \frac{1}{Z \big( s_i(x)  \big)} \, \, \, \, \textrm{and} \, \, \, \, 
 \nst_j(x) - \nst_k(x) \, = \, \sum_{i = j+1}^k \frac{s_i(x)}{Z \big( s_i(x)  \big)} \, .
$$
In particular,
$\mst_k(x) - \mst_{k-1}(x) =Z \big( s_k(x)  \big)^{-1}$ and $\nst_{k-1}(x) - \nst_k(x) = s_k(x) Z \big( s_k(x) \big)^{-1}$.
\item For $j,k \in \nwithzero$, 
the finite trail Mina margin map  $\minamm_{j+1,k+1}:(0,\infty) \to (0,\infty)$ satisfies the equation
$$
 \minamm_{j+1,k+1}(x) \, \, = \, \,   \Bigg( \, \sum_{i = -j}^{k+1} \frac{1}{Z \big( s_i(x)  \big)} \, \Bigg)^{-1} \, \cdot \, \sum_{i = -j}^{k+1} \frac{s_i(x)}{Z \big( s_i(x)  \big)} \, \, .
$$
\item The Mina margin map  $\minamm:(0,\infty) \to (0,\infty)$ satisfies
$$
 \minamm(x) \, \, = \, \,   
 \sum_{i \in \Z} \frac{s_i(x)}{Z \big( s_i(x)  \big)} \, \, .
$$
\end{enumerate}
\end{theorem}
\subsection{The game with symmetric boundary data}
As Section~\ref{s.economics} reviewed, economists have often considered  player-funded tug-of-war games in a symmetric form (with players receiving equal terminal rewards), finding and analysing symmetric equilibria (in which the stake of the leading player is determined by the position of the counter relative to the gameboard midpoint rather than by the identity of the player). 
In our case, these symmetries significantly simplify analysis, while leading to equilibria that exemplify generic features. Here we present our conclusions about them. 
The  \amacspace system on $\Z$ is the set of equations\footnote{In Section $III.3$ of Vickers' thesis~\cite{Vickers85}, symmetric equilibria are found in the Harris-Vickers tug-of-war game reviewed in Subsection~\ref{s.harrisvickers}, with cost functions $c_-(x) = c_+(x) = x^2/2$.  Effort rates $x_n$ satisfy $x_{-n} x_{n-1} = x_n^2/2$ for positive and negative $n$. This cousin of the $A$-system is solved  $x_n = x_0 2^{-n(n-1)/2}$. As Subsection~\ref{s.harrisvickers} reviewed, the elements of the Battlefield Cyl Fog are found, but decay rates are $e^{-\Theta(1) n^2}$, rather than doubly exponential in $n$.} in the real variables $A_i$, $i \in \Z$:
\begin{equation}\label{e.a}
 A_{-i-1} (2 A_i + A_{-i}) = A_{i+1}^2 \, ,
\end{equation}
where the index ranges over $\Z$.
We will also speak of the  \amacspace system on $\Z + 1/2$. In this case, the  real variables $A_i$ are indexed by $i$ in the one-half-offset lattice $\Z + 1/2$; the set of equations is given by~(\ref{e.a})
with the index ranging over $\Z + 1/2$.

By the Trail of Lost Pennies on $\Z$ in its symmetric form is meant the game  $\game(m_{-\infty},m_\infty,n_{-\infty},n_\infty)$, where the boundary parameters are supposed to satisfy $m_{-\infty} = n_\infty = 0$ and $m_\infty = n_{-\infty}$. There is thus a one-parameter family of such games, indexed by  $m_\infty \in (0,\infty)$.

\begin{theorem}\label{t.symmetric}
\leavevmode
\begin{enumerate}
\item For $\lambda \in (0,\infty)$, there is a unique solution $\big\{ a_i(\lambda): i \in \Z \big\}$ of  the  \amacspace system on $\Z$ such that $a_0(\lambda) = \lambda$. The solutions satisfy $a_i(\lambda) = \lambda a_i(1)$
for $\lambda \in (0,\infty)$ and $i \in \Z$.
\item For $\lambda \in (0,\infty)$, there is a unique solution $\big\{ A_i(\lambda): i \in \Z + 1/2 \big\}$ of  the  \amacspace system on $\Z+1/2$ such that $A_{1/2}(\lambda) = \lambda$. The solutions satisfy $A_i(\lambda) = \lambda A_i(1)$
for $\lambda \in (0,\infty)$ and $i \in \Z + 1/2$.
 \item In the notation of the first part, let $S_1$ denote the set of strategy pairs $\big( a_{-i+k}(\lambda), a_{i+k}(\lambda) : i \in \Z \big)$ indexed by $k \in \Z$ and $\lambda \in (0,\infty)$.
 In the notation of the second part, let $S_2$ denote the set of the strategy pairs 
  $\big( A_{-i-1/2 +k}(\lambda) ,A_{i+1/2 + k}(\lambda): i \in \Z \big)$
  with the same index set. The elements of $S_1 \cup S_2$ are pairwise distinct time-invariant Nash equilibria for the Trail of Lost Pennies on $\Z$ in its symmetric form.  
  \item Admit Conjecture~\ref{c.solutions} in the special case that $x=1$: namely, suppose that $Q(1)=2$. Then there are no other time-invariant Nash equilibria  for the Trail of Lost Pennies on $\Z$ in its symmetric form than those identified in the preceding part.
\end{enumerate}
\end{theorem}

\section{The article's structure}

There are five further chapters and two appendices. Two basic aspects of later use are treated in Chapter~\ref{c.tools}: a role-reversal symmetry satisfied by the \abmnmacspace system; and the solution of the simplest of the finite trail games, with just one site in open play. The fundamental relationship Theorem~\ref{t.nashabmn} between Nash equilibria and the \abmnmacspace equations is proved in Chapter~\ref{c.nashabmn}.
In Chapter~\ref{c.battlefield}, the explicit solution Theorem~\ref{t.defaultexplicit} 
and asymptotic decay estimates  Theorem~\ref{t.ajbj} are derived, along with the eventual gameplay unanimity  Theorem~\ref{t.unanimity}.
The Mina margin map~$\minamm$ is addressed in Chapter~\ref{c.allminamm}: its approximation by finite-trail counterparts, Theorem~\ref{t.relativereward} and several consequences among our main results; the map's $\PhiMac$-transformed version, and Theorem~\ref{t.phithetainverse}; and an explicitly recorded computation that $\minamm$ evaluated at $0.58$ is bounded away from one above, and thus Theorem~\ref{t.minamarginvalues}(3).
In Chapter~\ref{c.prospects}, we discuss several aspects of our results and proofs and some prospects for further study. 
Appendix~\ref{s.furtherformulas} contains the proofs of the further formulas from Section~\ref{s.formulas}, and Appendix~\ref{s.glossary} provides a glossary of commonly used notation.

\chapter{Some basic tools}\label{c.tools}

Role-reversal symmetry is treated in Section~\ref{s.rolereversal} and the trail game on $\llbracket -1,1 \rrbracket$ in Section~\ref{s.pennyforfeit}.
Later sections introduce some further  basic notation and properties.
\section{Role-reversal symmetry}\label{s.rolereversal}

\begin{definition}\label{d.rolereversal}
The {\em role-reversal map} $\mc{R}$ sends the space of quadruples $\Z \to \R^4$ to itself by mapping
$\big\{ (a_i,b_i,m_i,n_i): i \in \Z \big\}$ to $\big\{ (b_{-i},a_{-i},n_{-i},m_{-i}): i \in \Z \big\}$.
\end{definition}
\begin{proposition}\label{p.rolereversal}
Suppose that $(a,b,m,n) = \big\{ (a_i,b_i,m_i,n_i): i \in \Z \big\}$ is an \abmnmacspace solution. Then so is $\mc{R} (a,b,m,n)$.
\end{proposition}
{\bf Proof.} The result may be verified by inspecting the \abmnmacspace equations. We instead indicate in rough terms a more conceptual, game-theoretic, argument which is available for positive \abmnmacspace solutions if we admit their connection to the trail game via Theorem~\ref{t.nashabmn}. Suppose that a time-invariant Nash equilibrium $(b,a): \Z \to (0,\infty)$  is played in the first instance. If Mina and Maxine swap roles, so that the strategy pair $(a,b)$ is played, each acts in diametric opposition to her interests. But if the gameboard is then reflected through the origin, these interests are reversed, and each plays optimally once more. 
 It is the strategy pair $\Z \to (0,\infty)^2: i \to (a_{-i},b_{-i})$ that is now being played. This pair is a Nash equilibrium (for the game whose boundary data in $\R^4$ is specified by this pair), and the associated quadruple is an \abmnmacspace solution. This quadruple is $\mc{R} (a,b,m,n)$. \qed
  
  We will obtain Proposition~\ref{p.rkvalues} by using the role-reversal map $\mc{R}$ on quadruples whose index set is finite; to do so, we extend our notation to handle this circumstance.
  Let $j,k \in \nwithoutzero$ and let $\big\{ (a_i,b_i,m_i,n_i) \in \R^4: i \in \llbracket -j,k\rrbracket \big\}$. We set $\mc{R}(a,b,m,n)$ equal to  $\big\{ (b_{-i},a_{-i},n_{-i},m_{-i}) \in \R^4: i \in \llbracket -k,j\rrbracket \big\}$.
 Proposition~\ref{p.rolereversal} has a counterpart in the finite case which asserts that
 \begin{eqnarray}
 & & (a,b,m,n): \llbracket -j,k \rrbracket \to \R^4 \, \, \, \textrm{is an \abmnmacspace solution} \label{e.rolereversalfinite} \\
 & \implies & \mc{R}(a,b,m,n): \llbracket -k,j \rrbracket \to \R^4 \, \, \, \textrm{is an \abmnmacspace solution} \nonumber \, .
 \end{eqnarray}
 The \abmnmacspace equations can again be inspected to verify this statement.
 
 We will further consider the left shift $\mc{S}_1$, which sends any quadruple 
 $(a,b,m,n): \llbracket -j,k \rrbracket \to \R^4$ to the quadruple  $\llbracket -j-1,k-1 \rrbracket \to \R^4: i \to (a_{i+1},b_{i+1},m_{i+1},n_{i+1})$.
 
 {\bf Proof of Proposition~\ref{p.rkvalues}(1).} For $x \in (0,\infty)$, let  $\big\{ (a_i,b_i,m_i,n_i): i \in \llbracket -k,k\rrbracket \big\}$ be an \abmnmacspace solution on $\llbracket -k,k \rrbracket$
 such that   $\tfrac{n_{-1} - n_0}{m_0 - m_{-1}}$ equals $x$; that such a solution may be found has been explained in Subsection~\ref{s.resultandconjecture}.
 We {\em claim} that
 \begin{equation}\label{e.rkktwo}
  \minamm_{k,k} \big( \tfrac{m_1 - m_0}{n_0 - n_1} \big) = \minamm_{k,k} \big( \tfrac{n_{-1} - n_0}{m_0 - m_{-1}} \big)^{-1} \, . 
 \end{equation}
 Admitting this claim, we see that $s(x)=  \tfrac{n_0 - n_1}{m_1 - m_0}$ by~(\ref{e.rolereversalfinite}); thus, 
 $1/s(x) =   \tfrac{m_1 - m_0}{n_0 - n_1}$. Using the claim, we confirm Proposition~\ref{p.rkvalues}(1). 
 
To confirm~(\ref{e.rkktwo}), we let $\hat\phi_i$ denote the $\phi_i$-value of $\hat{R}(a,b,m,n)$. The claim follows from
$$
  \minamm_{k,k} \big( \tfrac{m_1 - m_0}{n_0 - n_1} \big)  \, = \, \minamm_{k,k} (\hat\phi_0) \, = \, \frac{\hat{n}_{-k} - \hat{n}_k}{\hat{m}_k - \hat{m}_{-k}} \, = \, \frac{m_k - m_{-k}}{n_{-k} - n_k} 
    \, = \, \minamm_{k,k} \big(  \tfrac{n_{-1} - n_0}{m_0 - m_{-1}} \big)^{-1} \, ,
   $$
   where the second and fourth equalities are due to~(\ref{e.minammfinite}).

 {\bf (2).} We now let  $\big\{ (a_i,b_i,m_i,n_i): 
 i \in \llbracket -k-1,k\rrbracket \big\}$ 
 be an \abmnmacspace solution on $\llbracket -k-1,k \rrbracket$ such that $\tfrac{n_{-1} - n_0}{m_0 - m_{-1}}= x$.
 We consider the operator $\mc{A} = \mc{S}_1 \circ \mc{R}$; note that, directly from~(\ref{e.rolereversalfinite}), $\mc{A}(a,b,m,n)$ is an \abmnmacspace solution, also on the index set $\llbracket -k-1,k\rrbracket$.
 We denote $\mc{A}(a,b,m,n) = \big\{ (\tilde{a}_i,\tilde{b}_i,\tilde{m}_i,\tilde{n}_i): i \in \llbracket -k-1,k\rrbracket \big\}$; and we let $\tilde\phi_i$ denote the $\phi_i$-value of  $\mc{A}(a,b,m,n)$ for $i \in \llbracket -k,k-1 \rrbracket$.
 
 By~(\ref{e.minammfinite}),
 $$
 \minamm_{k+1,k}(\tilde\phi_0) = \frac{\tilde{n}_{-k-1} - \tilde{n}_k}{\tilde{m}_k - \tilde{m}_{-k-1}} = \frac{m_k - m_{-k-1}}{ n_{-k-1} - n_k} \, .
 $$
 And again by~(\ref{e.minammfinite}), $\minamm_{k+1,k}\big( \tfrac{n_{-1} - n_0}{m_0 - m_{-1}} \big) = \tfrac{n_{-k-1} - n_k}{m_k - m_{-k-1}}$. Hence, we obtain
 \begin{equation}\label{e.tworinter}
 \minamm_{k+1,k}(\tilde\phi_0) = \minamm_{k+1,k}\big( \tfrac{n_{-1} - n_0}{m_0 - m_{-1}} \big)^{-1} \, .
 \end{equation}
 Note further that
 $$
 \tilde\phi_0 = \frac{\tilde{n}_{-1} - \tilde{n}_0}{\tilde{m}_0  - \tilde{m}_{-1}} = \frac{m_0 - m_{-1}}{n_{-1} - n_0} \, .
 $$
Since $x = \tfrac{n_{-1} - n_0}{m_0 - m_{-1}}$, we have that $\tilde\phi_0 = 1/x$. 
From~(\ref{e.tworinter}), we thus obtain Proposition~\ref{p.rkvalues}(2). \qed
 
 \begin{corollary}\label{c.rkvalues}
 For $k \in \nwithoutzero$, $\minamm_{k,k}(3) = \minamm_{k+1,k}(1) = 1$.
 \end{corollary}
 {\bf Proof.} By Proposition~\ref{p.rkvalues}(1) and $s(3) = 1/3$, we have that $\minamm_{k,k}(3)^2 = 1$. Since $\minamm_{k,k} > 0$, we obtain $\minamm_{k,k}(3) = 1$.
 By Proposition~\ref{p.rkvalues}(2), $\minamm_{k+1,k}(1)^2 = 1$. Since $\minamm_{k+1,k} > 0$, we confirm that $\minamm_{k+1,k}(1) = 1$. \qed
   
Next we invert the map $s:(0,\infty) \to (0,\infty)$ from Definition~\ref{d.acs}.    
 \begin{proposition}\label{p.sminusone}
 The function $s:(0,\infty) \to (0,\infty)$ from Definition~\ref{d.acs} is invertible, and its inverse is given by
 $$
 s^{-1}(x) \, = \, \frac{1}{s(1/x)} \, \, \, \, , \, \, \, \, \textrm{for $x \in (0,\infty)$} \, .
 $$
 \end{proposition}
 The map $s$ has been defined formulaically, and the proof may be viewed as simply the algebraic task of confirming that $s \big( 1/s(1/x) \big)$ and $1/\big(s(1/s(x)\big)$ equal $x$. 
 The role of $s$ as left shift for the $\phi$-sequence permits the possibly more instructive argument we will give involving role-reversal symmetry. We must make use of the `$s$ as shift' Proposition~\ref{p.shift}, which will be proved early in Chapter~\ref{c.battlefield} when the form of \abmnmacspace solutions is addressed.
 
 But first a quick consequence.
 \begin{corollary}\label{c.sminusone}
  Proposition~\ref{p.sminusone} is also valid when the iterate $s_j$ replaces for $s$, for any $j \in \nwithoutzero$.
 \end{corollary}
 {\bf Proof.} An induction on $j$ that reduces directly to the $j=1$ case.
 
 {\bf Proof of Proposition~\ref{p.sminusone}.}
 It is enough to show that $h:(0,\infty) \to (0,\infty)$ given by $h(x) = 1/s(1/x)$ satisfies
 \begin{equation}\label{e.shhs}
  (s \circ h)(x) = (h \circ s)(x) = x \, .
 \end{equation}
  Set
 $$
 (a_i,b_i,m_i,n_i) = \big( \asta_i(x), \bst_i(x),\mst_i(x),\nst_i(x) \big) \, \, \, , \, \, \, \, \textrm{for $i \in \Z$} \, .
 $$
 We have that $\phi_0 = \tfrac{n_{-1}-n_0}{m_0 - m_{-1}} = x$.
 First note that, by Proposition~\ref{p.shift}, $s(\phi_{-1}) = \phi_0$; or, in other words,
 \begin{equation}\label{e.xformula}
  s \, \bigg( \frac{n_{-2} - n_{-1}}{m_{-1} - m_{-2}} \bigg)  \, = \, \frac{n_{-1} - n_0}{m_0 - m_{-1}} \, = \, x \, .
 \end{equation}
 Let $\hat{\phi}:\Z \to (0,\infty)$ be such that, for $i \in \Z$, $\hat\phi_i$ is the value of $\phi_i$ for the quadruple $\mc{R}(a,b,m,n)$. Note then that
 \begin{equation}\label{e.hatphione}
 \hat\phi_1 \,  =\,  \frac{\hat{n}_0 - \hat{n}_1}{\hat{m}_1 - \hat{m}_0} \, =\,  \frac{m_0 - m_{-1}}{n_{-1} - n_0} \, = \, 1/x \, . 
 \end{equation}
Thus,  note that 
 $$
 s(1/x) =   s(\hat\phi_1) \, = \, \hat\phi_2 \, = \, \frac{\hat{n}_1 - \hat{n}_2}{\hat{m}_2 - \hat{m}_1} \, =\,  \frac{m_{-1} - m_{-2}}{n_{-2} - n_{-1}} \, , 
 $$
where the second equality is justified by 
Propositions~\ref{p.rolereversal} and~\ref{p.shift}. Applying $s$, we find from~(\ref{e.xformula}) that 
$s\big(1/s(1/x) \big) = x$.
We have confirmed that $s\big(h(x)\big) = x$ for $x \in (0,\infty)$.
 
 Next we note that $s(x) = \phi_1 = \tfrac{n_0 - n_1}{m_1 - m_0}$, so that $1/s(x) = \tfrac{m_1 - m_0}{n_0 - n_1} = \hat\phi_0$. But $s(\hat\phi_0) = \hat\phi_1 = 1/x$, by Propositions~\ref{p.shift} and~\ref{p.rolereversal}, and~(\ref{e.hatphione}). Which is to say, $1/s(1/s(x)) = x$, or $h\big(s(x)\big) = x$ for $x \in (0,\infty)$. This completes the derivation of~(\ref{e.shhs}) and thus the proof of Proposition~\ref{p.sminusone}.
 \qed

\section{Penny Forfeit}\label{s.pennyforfeit}
 The simplest case of the finite trail game from Section~\ref{s.finite} has $j=k=0$, when the first move is the last. 
The straightforward solution of this case is already instructive, and we provide it now, calling this game Penny Forfeit.  

Here is an explicit description of this one-turn game.
Maxine and Mina are asked to stake non-negative quantities $a$ and $b$. After these stakes have been submitted, the game victor is declared: this will be Maxine, with probability $\tfrac{a}{a+b}$; otherwise, it will be Mina. If Maxine wins, she receives $m_1$, and Mina $n_1$; if Mina wins, Maxine receives $m_{-1}$ and Mina $n_{-1}$. These four values act as boundary data. They are supposed to be real values that satisfy $m_{-1} < m_1$ and  $n_1 < n_{-1}$.
%
Maxine and Mina's mean winnings in the game are 
\begin{equation}\label{e.maxineminawinnings}
\tfrac{a}{a+b} m_1 + \tfrac{b}{a+b} m_{-1}  -a  \, \, \, \textrm{and} \, \, \,  \tfrac{b}{a+b} n_{-1} + \tfrac{a}{a+b} n_1  - b \, ,
\end{equation}
where in each expression the respective terms are mean terminal receipt in the event of turn victory; such receipt in the event of turn defeat; and the negative contribution from the forfeited stake. 
The pair $(b,a)$ is a Nash equilibrium---a notion that is specified by suitably adapting the definition in Section~\ref{s.gamespec}---
when these last two expressions are both global maxima as the variables $a$ and $b$ are respectively varied over $[0,\infty)$.
\begin{lemma}\label{l.pennyforfeit}
There is a unique solution in $(a,b) \in [0,\infty)^2$ in which the pair of expressions in~(\ref{e.maxineminawinnings}) are both global maxima as the variables $a$ and $b$ are respectively varied over $[0,\infty)$. It is given by 
\begin{equation}\label{e.absolution}
(a,b) \, = \, \bigg(\frac{M^2 N}{(M+N)^2},\frac{M N^2}{(M+N)^2} \bigg) \, , \, \, \, \, \textrm{with} \, \, \, \, M = m_1 - m_{-1} \, \, \, \, \textrm{and} \, \, \, \, N = n_{-1} - n_1 \, ,
\end{equation}
Note that $a$ and $b$ are strictly positive.
\end{lemma}
{\bf Proof.} A critical point $(a,b)$ is given by setting the respective partial derivatives in $a$ and~$b$ of the two expressions equal to zero: the conditions are 
$$
\tfrac{b}{(a+b)^2}(m_1 - m_{-1}) - 1 \, = \, \tfrac{a}{(a+b)^2}(n_{-1} - n_1) - 1 \, = \, 0 \, . 
$$
At least one component in the desired pair $(a,b)$ must be non-zero. 
Indeed, if for example $a$ equals zero,
then an infinitesimal increase of $b$ from zero will increase Mina's expected payoff from $\tfrac{n_{-1}+n_1}{2}$ to $n_{-1}$. 
Restricting then, as we may, to solutions with at least one positive component, 
we see that there exists a unique solution in $(a,b) \in [0,\infty)^2$ of the last displayed equations, and that this solution is given by~(\ref{e.absolution}). This is indeed a global maximum for the pair of expressions in~(\ref{e.maxineminawinnings}) under respective variation of $a$ and~$b$ in~$[0,\infty)$. \qed

{\em Remark.} We see then that Penny Forfeit has a unique Nash equilibrium~$(b,a)$, with $(a,b)$ as just specified. It is straightforward to see that this Nash equilibrium is unique even if we permit the players to offer random stakes. 
  
  \section{The game with a delayed start}\label{s.delayedstart}

We will wish to consider the finite and infinite trail games begun at a turn whose index $\ell \in \nwithoutzero$ is general. For $(i,\ell) \in \Z \times \nwithoutzero$ and $(S_-,S_+) \in \mc{S}$ 
we will write\hfff{delayed} $\PP_{S_-,S_+}^{i,\ell}$ and $\E_{S_-,S_+}^{i,\ell}[\cdot]$  for the law and expectation operator of
gameplay $X: \llbracket \ell,\infty) = \Z \cap [\ell,\infty) \to \Z$, $X_\ell = i$, in the Trail of Lost Pennies on $\Z$, begun at the $(\ell + 1)$\textsuperscript{st} turn at $\ell$. Payoffs, costs and terminal receipts 
$P_\pm$, $C_\pm(u)$ (for $u \in \llbracket \ell,\infty)$) and $T_{\pm}$ remain as specified by Section~\ref{s.gamespec}. Mina and Maxine's payoff identities~(\ref{e.minapayoff}) and~(\ref{e.maxinepayoff}) now take the $\PP_{S_-,S_+}^{i,\ell}$-almost sure form
 \begin{equation}\label{e.delayedpayoff}
 P_\pm \, = \,  - \sum_{j = \ell +1}^\infty C_\pm(j)  \,\, + \,\,  T_\pm \, .
 \end{equation}
Note that  $\PP_{S_-,S_+}^i$ equals $\PP_{S_-,S_+}^{i,0}$.

  \section{Lack of escape entails infinite costs}  
  \begin{lemma}\label{l.dontlookback}
  Let $(S_1,S_2) \in \mc{S} \times \mc{S}_0$ be a strategy pair whose second component is time-invariant.
  Writing $a_i = S_2(i,j)$ for $(i,j) \in \Z \times \nwithoutzero$, suppose that $a_i > 0$ for $i \in \Z$.
For $(i,\ell) \in \Z \times \nwithoutzero$, 
  $\PP_{S_1,S_2}^{i,\ell}(E^c) > 0$ implies that $\E_{S_1,S_2}^{i,\ell} [P_-] = -\infty$.
  \end{lemma}
 {\bf Proof.}
 For $j \in \Z$, let $\cross_-(j)$ denote the event that Mina wins infinitely many turns at which the counter is at $j$. We claim that, up to a $\pgameplay{S_1}{S_2}{i,\ell}$-null set,
 \begin{equation}\label{e.ecomplement}
 E^c \, \subseteq \, \bigcup_{j \in \Z} \cross_-(j) \, .
 \end{equation}
 To see this, set $V_j$ denote the event that the counter visits $j \in \Z$ infinitely often. The occurrence of $E^c$ entails that of $\cup_{j \in \Z} V_j$. If $V_j$ occurs and Mina wins infinitely many of the turns at which $X$ visits $j$,
 then $\cross_-(j)$ occurs. If $V_j$ occurs but Mina does not thus succeed, there are infinitely many occasions on which $X$ leaves $j$ to the right, only to return to $j$ at some later time. Consider the set of turns that occur just before each of these returns. At each, $X$ is at $j+1$ and Mina wins the turn, so that $X$ passes to $j$. Thus, $\cross_-(j+1)$ occurs. We have derived~(\ref{e.ecomplement}).
 
For $j \in \Z$, let $\totalcost_-(j) = \sum_{t=\ell}^\infty {\bf 1}_{X_t =j} C_-(t+1)$ denote Mina's running cost expended at $j$ under $\pgameplay{S_1}{S_2}{i,\ell}$. 
Let  $\crossnum_-(j,j-1) = \sum_{t=\ell}^\infty {\bf 1}_{X_t =j,X_{t+1}=j-1}$ denote the number of turns with index at least $\ell+1$ that are won by Mina and at whose start $X$ visits $j$. 
Since 
$C_-(t+1) = S_1(X_t,t+1)$, we have that
\begin{eqnarray*}
  \egameplay{S_1}{S_2}{i,\ell} \big[ \crossnum_-(j,j-1) \big] 
  & = & \sum_{t=\ell}^\infty \pgameplay{S_1}{S_2}{i,\ell} (X_t =j) \cdot \tfrac{S_1(j,t+1)}{S_1(j,t+1) + a_j}
  \, \leq \, 
   a_j^{-1} \sum_{t=\ell}^\infty  \pgameplay{S_1}{S_2}{i,\ell} (X_t =j) S_1(j,t+1) \\
 & = &
  a_j^{-1} 
 \egameplay{S_1}{S_2}{i,\ell} \big[ \totalcost_-(j) \big]
 \, \leq \, a_j^{-1} 
 \egameplay{S_1}{S_2}{i,\ell} \sum_{t=\ell}^\infty C_-(t) \, .
 \end{eqnarray*}
By~(\ref{e.ecomplement}) and $a_j > 0$ for $j \in \Z$, we see then that, if $\pgameplay{S_1}{S_2}{i,\ell}(E^c) > 0$, then
 $\egameplay{S_1}{S_2}{i,\ell} [\crossnum_-(j,j-1)]$ is infinite for some $j \in \Z$, and thus so is Mina's mean total running cost  
 $\egameplay{S_1}{S_2}{i,\ell} \sum_{t=\ell}^\infty C_-(t)$. Applying $\egameplay{S_1}{S_2}{i,\ell}$  
 to~(\ref{e.delayedpayoff}) with $\pm = -1$ and noting that terminal receipts $T_-$ are almost surely bounded, we find that  Mina's mean payoff $\egameplay{S_1}{S_2}{i,\ell} [P_-]$ equals minus infinity.
 This completes the proof of Lemma~\ref{l.dontlookback}. \qed
 
\section{Relating the finite and infinite trail games}\label{s.relating}
  
  Let $(m_{-\infty},m_\infty,n_{-\infty},n_\infty) \in \R^4$
  satisfy $m_{-\infty} < m_\infty$ and  $n_\infty < n_{-\infty}$. Suppose given a positive \abmnmacspace solution~$\big\{ (a_i,b_i,m_i,n_i): i \in \Z \big\}$ with this boundary data.
  In this setup, we may attach to the infinite trail game an associated game on each finite trail. 
Indeed,  for $j,k \in \N$, this solution provides a quadruple $(m_{-j-1},m_{k+1},n_{-j-1},n_{k+1}) \in \R^4$ that may furnish boundary data for the game on the finite trail  $\llbracket -j-1 ,k+1 \rrbracket$.
The games are then fully specified by setting the non-finishing terminal receipt values $m_*$ and $n_*$ seen (for~$\Z$) in~(\ref{e.terminalmina}) and~(\ref{e.terminalmaxine}). Any given admissible choice of $(m_*,n_*)$ for the infinite game may be used simultaneously for the associated finite trails. Indeed, the requisite conditions $m_* \leq m_{-\infty}$ and $n_* \leq n_\infty$
transmit to counterparts  $m_* \leq m_{-j-1}$ and $n_* \leq n_{k+1}$ for the trail $\llbracket -j-1,k+1 \rrbracket$ in light of Theorem~\ref{t.positiveabmn}(1).

With the infinite trail game pairs to finite counterparts in this way,it is useful to extend the identification to gameplay by specifying a coupling of counter evolution on~$\Z$ and the finite intervals. 
  \begin{definition}\label{d.coupling}
  Let $i \in \Z$ and $(S_-,S_+) \in \mc{S}^2$. Recall that the gameplay $X:\N \to \Z$, $X_0 = i$, of the infinite trail game governed by $(S_-,S_+)$ is specified under the law $\PP_{S_-,S_+}^i$.
 For $j,k \in \nwithoutzero$,  strategy pairs in $\mc{S}[j,k]^2$ for the game with trail $\llbracket -j -1,k+1\rrbracket$ result by restricting the domain of $S_-$ and $S_+$ to $\llbracket -j ,k \rrbracket$. 
   Copies of the gameplay $X^{j,k}:\nwithzero \to \llbracket -j-1,k+1 \rrbracket$, $X^{j,k}_0 = i$, 
   that result from use of these restricted pairs may be coupled under  $\PP_{S_-,S_+}^i$  whenever $j,k \in \nwithoutzero$ are such that $i \in \llbracket -j ,k \rrbracket$.
   To specify these copies, set 
\begin{equation}\label{e.taujk}
\tau^{j,k} \, = \, \inf \big\{ i \in \nwithoutzero: X_i \in \{-j-1,k+1 \} \big\} \, .
\end{equation}
Writing $\wedge$ for minimum, 
we then take $X^{j,k}(u) = X(u \wedge \tau^{j,k})$ for $u \in \nwithzero$. 
  \end{definition} 
  The finite and infinite trail payoffs, costs and terminal receipts $*_\pm^{j,k}$ and $*_\pm$, $* \in \{P,C,T\}$, are coupled under  $\PP_{S_-,S_+}^i$ by this definition, with terminal receipt~$T_\pm^{j,k}$ being one of the components of the quadruple $(m_{-j-1},m_{k+1},n_{-j-1},n_{k+1})$  specified by the given \abmnmacspace solution. We note some basic relationships that result.

\begin{lemma}\label{l.couplingproperties}
Let  $(S_-,S_+) \in \mc{S}^2$. Suppose that $i \in \Z$ and $j,k \in \nwithoutzero$ satisfy $i \in \llbracket -j,k \rrbracket$.
\begin{enumerate}
\item We have that $P_-  - P^{j,k}_- = T_- - T^{j,k}_- - \sum_{t = \tau^{j,k}}^\infty C_-(t)$.
\item And that  $P_-  - P^{j,k}_- \leq T_- - T^{j,k}_-$.
\item Suppose that $\PP_{S_-,S_+}^i \big( \tau^{j,k} < \infty \big) = 1$.   It is  $\PP_{S_-,S_+}^i$-almost certain that  $P_-  - P^{j,k}_- \leq n_{-\infty} - n_\infty$.
\end{enumerate}
  \end{lemma}
{\bf Proof: (1).} This follows from~(\ref{e.minapayoff}),~(\ref{e.finitepayoff}) and $C^{j,k}_-(t) = C_-(t)$ for $u \in \llbracket 0, \tau^{j,k}-1 \rrbracket$. \\
{\bf (2).}  Due to the preceding part and the non-negativity of costs $C_-(t)$. \\
{\bf (3).}  Since  $\tau^{j,k} < \infty$ is certain,
the receipt $T^{j,k}_-$ is a weighted average of $n_{-j-1}$ and $n_{k+1}$. From $\big\{ n_i: i \in \Z \big\}$ being decreasing (this due to Theorem~\ref{t.positiveabmn}(1), because this sequence is the $n$-component of a positive \abmnmacspace solution),
we find that $T^{j,k}_- \geq n_\infty$. Also note that $T_- \leq n_{-\infty}$. The preceding part of the lemma thus implies the stated result.  \qed

\begin{lemma}\label{l.stopping}
Let $(S_-,S_+) \in \mc{S}^2$, $k \in \Z$ and $\ell \in \nwithzero$.
Let $Q \in \nwithzero \cup \{ \infty \}$ be a stopping time with respect to gameplay $X: \nwithzero \to \Z$ under the law $\pgameplay{S_-}{S_+}{k,\ell}$ specified in Section~\ref{s.delayedstart}. Then
\begin{equation}\label{e.stopping}
\egameplay{S_-}{S_+}{k,\ell} [P_-] \, = \,  - \, \egameplay{S_-}{S_+}{k,\ell}  \sum_{t=\ell + 1}^{Q - 1} C_-(t) \,\, + \,\, 
  \egameplay{S_-}{S_+}{k,\ell} \big[ \egameplay{S_-}{S_+}{X(Q)} [P_-] \big] \, . 
\end{equation}
In reading this display in the event that $Q = \infty$, we adopt the conventions that  $\egameplay{S_-}{S_+}{\infty,\ell}[P_-] = n_\infty$ and   $\egameplay{S_-}{S_+}{-\infty,\ell}[P_-] = n_{-\infty}$,
as well as $Q-1 = \infty$. We also have the counterpart identity for Maxine, given by $P_- \to P_+$ and $C_- \to C_+$.
\end{lemma}
{\bf Proof.} The right-hand side of~(\ref{e.delayedpayoff}) with $\pm = -1$ may be written $A_1 + A_2$, where $A_1$ is the sum of costs $C_-(t)$ with $\ell + 1 \leq t < Q$; and $A_2$ is the sum of the higher indexed costs (in the case that $Q$ is finite) and the terminal receipt~$T_-$. 
Since $T_-$ equals $n_{-\infty}$ or $n_\infty$
when the events $E_-$ or $E_+$ occur, we find that, when the mean 
$\egameplay{S_-}{S_+}{k,\ell}$ of~(\ref{e.delayedpayoff}) thus represented is taken, the two right-hand terms in the lemma result. \qed

We have used Theorem~\ref{t.positiveabmn}(1), and we will use it again in a moment. We now give the simple proofs of Theorem~\ref{t.positiveabmn}(1,2).

{\bf Proof of Theorem~\ref{t.positiveabmn}(1).} Since $a_i + b_i > 0$,  \abmnmac$(3)$ implies that $m_{i+1} > m_{i-1}$.  We may rearrange \abmnmac$(1)$ in the form 
$m_i   =  \tfrac{a_i}{a_i + b_i} m_{i+1} + \tfrac{b_i}{a_i + b_i} m_{i-1} - a_i$. Using $m_{i-1} <m_{i+1}$ and $b_i > 0$, we find that $m_i < m_{i+1} - a_i$. Since $a_i > 0$,
$m_i < m_{i+1}$. That $n_{i+1} < n_i$ follows similarly. We have shown that the \abmnmacspace solution $(a,b,m,n)$ is strict.

{ \bf (2).}
The sequences $\big\{ m_i: i \in \Z \big\}$ and  $\big\{ n_i: i \in \Z \big\}$ are increasing and decreasing, by the preceding part. Thus, the limiting values~(\ref{e.boundarydata}) exist, at least as elements of $\R \cup \{ \infty\} \cup \{ - \infty \}$; they satisfy $m_\infty > m_{-\infty}$ and $n_{-\infty} > n_\infty$. \qed

Note that  Theorem~\ref{t.positiveabmn}(1,2) do not exclude the possibilities that $m_\infty$ or $n_{-\infty}$ equals $\infty$ or that  $n_\infty$ or $m_{-\infty}$ equals $-\infty$. We will deduce this when we prove Theorem~\ref{t.positiveabmn}(3). This result will be derived in Section~\ref{s.consequences} as a  consequence of the asymptotic decay estimate Theorem~\ref{t.ajbj}.

The next result interprets the $m$- and $n$-components of a \abmnmacspace solution as mean payoffs. It is couched in the notation of delayed-start games from Section~\ref{s.delayedstart}.
\begin{lemma}\label{l.minipayoff}
Let $\big\{ (a_i,b_i,m_i,n_i) : i \in \Z \big\}$ denote a positive 
solution of the \abmnmacspace equations
with boundary data $(m_{-\infty},m_\infty,n_{-\infty},n_\infty) \in \R^4$.
Let $S_-,S_+ \in \mc{S}$ satisfy
$S_-(i,j) = b_i$ and $S_+(i,j) = a_i$ for $(i,j) \in \Z \times \nwithoutzero$.
\begin{enumerate}
\item 
Let $i \in \Z$ and $\ell \in \nwithzero$. Then $\pgameplay{S_-}{S_+}{i,\ell}(E) =1$.
\item 
 For $i \in \Z$  and $\ell \in \nwithzero$,
$$
m_i = \E_{(S_-,S_+)}^{i,\ell}[P_+]  \, \, \, \, \textrm{and} \, \, \, \, 
n_i =   \E_{(S_-,S_+)}^{i,\ell}  [P_-]  \, .  
$$
\item Let $j,k \in \nwithzero$  and $\ell \in \nwithzero$. For $i \in \llbracket -j,k \rrbracket$,
$$
 m_i =  \E_{(S_-,S_+)}^{i,\ell}  [P^{j,k}_+]   \, \, \, \, \textrm{and} \, \, \, \, 
 n_i =  \E_{(S_-,S_+)}^{i,\ell}  [P^{j,k}_-]  \, .
$$
\end{enumerate}
\end{lemma}
{\bf Proof: (1).}  By Theorem~\ref{t.positiveabmn}(1), $n_i > n_\infty$. But $n_\infty > -\infty$ by assumption. Thus Lemma~\ref{l.dontlookback} implies the sought statement. \\
{\bf (2).} Since $a_i + b_i>0$, $\abmnmac(1)$ may be written in the form
 $m_i =  \tfrac{a_i}{a_i+b_i}m_{i+1} +  \tfrac{b_i}{a_i+b_i}m_{i-1} - a_i$ 
 or equivalently 
 $m_i = \egameplay{S_-}{S_+}{i} [m(X_1)] \, - \, a_i$.
 Iterating, we find that  
 \begin{equation}\label{e.mexpand}
 m_i \, = \, \egameplay{S_-}{S_+}{i,\ell} \, [m(X_{u+1})] \, - \,  \egameplay{S_-}{S_+}{i,\ell} \, \sum_{i=\ell}^u a_{X(i)} 
 \end{equation}
 for any $u \in \nwithzero$, $u \geq \ell$.
 The value of $\lim_{u \to \infty} m(X_u)$ exists on the event $E$, equalling $m_\infty$ or $m_{-\infty}$ according to whether $E_+$ or $E_-$ occurs. By Lemma~\ref{l.minipayoff}(1),
 we see that   $\lim_{u \to \infty} \egameplay{S_-}{S_+}{i,\ell} \, [ m(X_{u+1})]$ equals $m_\infty  \cdot \pgameplay{S_-}{S_+}{i,\ell} (E_+) + m_{-\infty} \cdot  \pgameplay{S_-}{S_+}{i} (E_-)$.
  In the notation of Lemma~\ref{l.stopping}, we find by taking the high-$u$ limit of the preceding display that $m_i$ equals the right-hand side of~(\ref{e.stopping}) with $k=i$ and $Q$ identically equal to infinity. Thus, Lemma~\ref{l.stopping}
  implies that $m_i = \E_{(S_-,S_+)}^{i,\ell} [P_+]$. That   $n_i =   \E_{(S_-,S_+)}^{i,\ell} [ P_-]$ is similarly proved. \\
{\bf (3).} We may obtain~(\ref{e.mexpand}) with $X$ replaced by its stopped version $X^{j,k}$.  By taking the high-$u$ limit, we find that $m_i$ equals  the right-hand side of~(\ref{e.stopping}) with $k=i$ and $Q = \tau^{j,k}$.  From Lemma~\ref{l.stopping} we thus find that
   $m_i = \E_{(S_-,S_+)}^{i,\ell} [P^{j,k}_+]$. That   $n_i =   \E_{(S_-,S_+)}^{i,\ell} [ P^{j,k}_-]$ follows similarly. This completes the proof of Lemma~\ref{l.minipayoff}(3). \qed

Some simple relationships between escape in the finite and infinite trail games are now recorded.
We define the events
$E_-[j,k] = \big\{ X(\tau^{j,k}) = - j-1 \big\}$ and $E_+[j,k] = \big\{ X(\tau^{j,k}) = k+1 \big\}$.
\begin{lemma}\label{l.eminuseplus}
We have that
$$
E_- \, = \, \bigcup_{k=1}^\infty \bigcap_{j=1}^\infty E_-[j,k] \, \, \, \, \textrm{and} \, \, \, \, 
E_+ \, = \, \bigcup_{j=1}^\infty \bigcap_{k=1}^\infty E_+[j,k] \, .
$$
\end{lemma}
{\bf Proof.} These follow from the definitions of the events $E_-$ and $E_+$. \qed
\begin{lemma}\label{l.mn} 
We have that
\begin{equation}\label{e.mn}
\lim_{k \to \infty } \PP_{S_-,S_+}^i \bigg( E_- \setminus \Big\{ \lim_{j \to \infty}  m\big( X_{\tau^{j,k}} \big) = m_{-\infty} \Big\} \bigg) = 0 \, .
\end{equation}
and
$$
\lim_{j \to \infty } \PP_{S_-,S_+}^i \bigg( E_+ \setminus \Big\{ \lim_{k \to \infty}  m\big( X_{\tau^{j,k}} \big) = m_\infty \Big\} \bigg) = 0 \, .
$$
These statements are also valid if we replace all instances of $m$ by $n$. 
\end{lemma}
{\bf Proof}. By Lemma~\ref{l.eminuseplus} for $E_-$, we see that, on this event,
there exists a random value $K \in \nwithoutzero$ such that, for all $j \in \nwithoutzero$,
 $X(\tau^{j,K}) = -j-1$. Since $m_{-i} \to m_{-\infty}$ as $i \to \infty$, we see that, on $E_-$, 
$\lim_j m\big( X(\tau^{j,K})\big) = m_{-\infty}$. Thus, we obtain~(\ref{e.mn}). The other three assertions made by the lemma have the same proof up to evident notational changes. \qed
   
 \chapter{The structure of time-invariant Nash equilibria}\label{c.nashabmn}
The aim of this chapter is to prove Theorem~\ref{t.nashabmn}, our result that relates time-invariant Nash equilibria and positive \abmnmacspace solutions. On the way to this result, we will establish some basic properties of time-invariant Nash equilibria. 
In the first section, we prove Theorem~\ref{t.nashabmn}(1) alongside some simple properties of strategy pairs. 
The second proves Theorem~\ref{t.nashabmn}(2). 

 But first we will address a misgiving which suggests that the structure of Nash equilibria is trivial. If $(S_-,S_+) \in \mc{S}^2$, then $\lambda \cdot (S_-,S_+)$ for $\lambda \in (0,1)$
 is an alternative strategy pair that leads to gameplay and terminal receipts that are identical in law as when $\lambda =1$; but running costs for both players are reduced. The players could effectively eliminate running costs by picking a small value of $\lambda$. Does this mean that $(S_-,S_+) \not\in \mc{N}$ and $\mc{N} = \emptyset$? No: the concept of Nash equilibrium is one of stability under unilateral perturbations, and neither player could alone effect this change in $\lambda$. If the players were to decrease $\lambda$, they would act as airline executives meeting privately to plan a joint price hike: the membership $(S_-,S_+) \in \mc{N}$ is not vulnerable to such illegally cooperative perturbation.

\section{Time-invariant Nash equilibria result in positive \abmnmacspace solutions}

Here, we prove Theorem~\ref{t.nashabmn}(1). Our style of argument is hands on: we build up inferences on the behaviour of a time-invariant Nash equilibrium step-by-step.
With one exception: to close out the proof, we will invoke the unanimity Theorem~\ref{t.unanimity}(2,3), which is argued independently by explicit solution of the $\abmnmacspace$ system in Chapter~\ref{c.battlefield}.

Recall the mean payoff notation~(\ref{e.minapayoff}) and~(\ref{e.maxinepayoff}).   A strategy pair $(S_-,S_+) \in \mc{S}^2$ has {\em finite mean costs} if neither $\E^k_{S_-,S_+}[P_-]$ nor $\E^k_{S_-,S_+}[P_+]$ equals minus infinity, for any $k \in \Z$. 

Let $(S_-,S_+) \in \tis^2$. 
We adopt our standard convention of writing $b_i = S_-(i,j)$ and  $a_i = S_+(i,j)$ for $(i,j) \in \Z \times \nwithoutzero$.
The {\em idle zone}\hfff{idle} $\mc{I} \subset \Z$ is given by $\mc{I} = \big\{ j \in \Z: a_j = b_j = 0 \big\}$. 
\begin{lemma}\label{l.idlezone}
Let  $(S_-,S_+) \in \tis^2$ be such that $\mc{I} \not= \emptyset$. For $k \in \Z$, consider the gameplay $X:\N \to \Z$ under 
$\pgameplay{S_-}{S_+}{k}$. For  $i \in \nwithzero$ given,
condition on the event that $X_i$ is a given element of $\mc{I}$.
(If $i = 0$, suppose that $k \in \mc{I}$.)  Set  $j = \inf \big\{ m  \geq i: X_m \not\in \mc{I} \big\}$. Then the conditional law of $X: \llbracket i,j\rrbracket: \nwithzero \to \Z$ is equal to simple random walk given the value $X_i$ stopped on leaving~$\mc{I}$.
\end{lemma}
{\bf Proof.} At each turn with index in $\llbracket i,j-1 \rrbracket$, neither Mina nor Maxine offers a positive stake, since the $b$ and $a$ values vanish in the idle zone.
The gameplay increments~$X(h+1) - X(h)$ for $h \in \llbracket i,j-1 \rrbracket$ are thus unbiased $\pm 1$ steps as determined by the $0/0 = 1/2$ rule that was specified in Section~\ref{s.gamespec}. \qed

Recall the escape event $E$ from~(\ref{e.escape}). An element of $\mc{S}_0^2$ is non-zero when at least one of its components is not identically zero. 
\begin{proposition}\label{p.fmc}
Let $(S_-,S_+) \in \tis^2$ be non-zero, with finite mean costs.
Then escape occurs almost surely: 
$\pgameplay{S_-}{S_+}{k}(E)=1$
 for $k \in \Z$.
\end{proposition}
{\bf Proof.} Let $k \in \Z$ and suppose that 
$\pgameplay{S_-}{S_+}{k}(E^c)>0$. 
We may find $\ell \in \Z$ such that it is with positive probability that the process $X$, under the law 
$\pgameplay{S_-}{S_+}{k}$, 
visits $\ell$ infinitely often. Since $(S_-,S_+)$ is time-invariant, the strong Markov property implies that
\begin{equation}\label{e.infinitelyoften}
\pgameplay{S_-}{S_+}{k} \Big(  \textrm{$X$ visits $\ell$ infinitely often} \, \Big\vert \,  \textrm{$X$ visits $\ell$ at least once} \, \Big) \, = \, 1 \, .
\end{equation}
Let $i \in \Z \cup \{ - \infty \}$, $j \in \Z \cup \{ \infty\}$, $i \leq j$, be such that at least one of $i$ and $j$ is finite; 
 $\ell \in \llbracket i,j \rrbracket$; $a_m = b_m = 0$ for $m \in \Z \cap (i,j)$; and at least one of $a_m$ and $b_m$ is positive for any endpoint $m \in \{ i,j \}$ that is finite.
 If $\ell \not\in \mc{I}$, then we may take  $i=j=\ell$, with 
  some of these conditions holding vacuously. 
  If $\ell \in \mc{I}$, we may take $i < \ell < j$.

Suppose that $i < \ell < j$. Note that $\llbracket i+1,j-1 \rrbracket \subset \mc{I}$. 
We now consider the conditional law of  $X$ under $\pgameplay{S_-}{S_+}{k}$ given that $X$ visits $\ell$ infinitely often. We invoke~(\ref{e.infinitelyoften}) to note that the conditioning disappears at the first visit of $X$ to $\ell$.
Lemma~\ref{l.idlezone} thus implies that,
on each occasion that $X$ visits~$\ell$, $X$ pursues a simple random walk until it reaches $i$ or $j$. 
 Suppose, without loss of generality, that the index $i$ is finite, and that $a_i > 0$. It is with probability at least $2^{-(\ell-i)}$ that $X$ proceeds from a visit to $\ell$ by means of a string of leftward steps to reach $i$. Later, the conditioned walk $X$ inevitably returns to $\ell$, and a further opportunity to reach $i$ directly ensues. 
 Thus, $X$ will infinitely often visit~$i$, a location to which  $a$ assigns positive value.
(Note that this conclusion also holds trivially in the opposing case, where $i=j=\ell$.)
 The cost $\sum_{t \geq 1} C_+(t)$ incurred by Maxine
(which is specified in Section~\ref{s.gamespec})
is thus seen to be almost surely infinite on the  $\pgameplay{S_-}{S_+}{k}$-positive probability event that $X$ visits $\ell$ infinitely often.
(Were $b_i$ instead supposed positive, then it would be Mina's cost  $\sum_{t \geq 1}C_-(t)$ that is found to be infinite.) This is contrary to our assumption that $(S_-,S_+)$ has finite mean costs. We conclude, as desired, that $\pgameplay{S_-}{S_+}{k}(E) = 1$. \qed

For $S \in \tis$, let $\lef(S) \in \Z \cup \{ -\infty\} \cup \{\infty\}$ denote $\inf \{ i \in \Z : S(i,1) > 0 \}$; and let $\righ(S) \in \Z \cup \{ -\infty\} \cup \{\infty\}$ denote $\sup \{ i \in \Z : S(i,1) > 0 \}$.
We say that $S$ is {\em wide}\hfff{wide} if $\lef(S) = -\infty$ and $\righ(S) = \infty$; if $S$ is not wide, it is {\em narrow}. 

The right rocket $\eta \cdot \rocket^{i\rightarrow}$ at $i \in \Z$ of strength $\eta \in (0,\infty)$ is the element of $\tis$ given 
$$
 \eta \cdot \rocket^{i\rightarrow}_j \, = \, \eta \cdot 2^{-(j-i)-1} {\bf 1}_{j \geq i} \, \, \, , \, \, \, j \in \Z \, .
$$
The counterpart left rocket  $\eta \cdot \rocket^{\leftarrow i} \in \tis$ is 
$$
 \eta \cdot \rocket^{\leftarrow i}_j \, = \, \eta \cdot 2^{-(i-j)-1} {\bf 1}_{j \leq i} \, \, \, , \, \, \, j \in \Z \, .
$$
The right boost at $i \in \Z$ of strength $\eta$ is the map $\boost_\eta^{i\rightarrow}:\tis \to \tis$ that sends $q = (q_i: i \in \Z) \in \tis$
to $q +  \eta \cdot \rocket^{i\rightarrow}$. The corresponding left boost  $\boost_\eta^{i\leftarrow}:\tis \to \tis$ sends $q$
to $q +  \eta \cdot \rocket^{i\leftarrow}$.

The right drag  at $i \in \Z$  is the map $\drag^{i\rightarrow}:\tis \to \tis$ that sends $q  \in \tis$ to the map
$$
\Z \to (0,\infty ): j \to \, \,  \begin{cases}
  \, q_j/2  &  \text{if $j \geq i$} \\
 \, q_j  &  \text{if $j < i$} \, .
\end{cases}
$$
The counterpart left drag  $\drag^{i\leftarrow}:\tis \to \tis$ sends $q  \in \tis$ to 
$$
\Z \to (0,\infty ): j \to \, \,  \begin{cases}
  \, q_j/2  &  \text{if $j \leq i$} \\
 \, q_j  &  \text{if $j > i$} \, .
\end{cases}
$$
\begin{lemma}\label{l.boostdrag}
Let $(S_-,S_+) \in \tis^2$. 
\begin{enumerate}
\item Suppose that the quantities $\righ(S_-)$ and $\righ(S_+)$ are finite. Let $i \in \Z$ exceed their maximum. Then $\egameplay{S_-}{\boost_\eta^{i\rightarrow}(S_+)}{i}[P_+] > \egameplay{S_-}{S_+}{i}[P_+]$ for $\eta \in (0,m_\infty - m_{-\infty})$. 
\item Suppose that $\righ(S_+) = \infty$ and $\righ(S_-) < \infty$. Let $i \in \Z$, $i > \righ(S_-)$, satisfy $S_+(i,1) > 0$. Then  $\egameplay{S_-}{\drag^{i\rightarrow}(S_+)}{i}[P_+] > \egameplay{S_-}{S_+}{i}[P_+]$.
\item If $\lef(S_-)$ and $\lef(S_+)$ exceed $-\infty$ and $i \in \Z$ is less than their minimum, then, provided that $\eta \in (0,n_{-\infty} - n_\infty)$, we have that $\egameplay{\boost_\eta^{\leftarrow i}(S_-)}{S_+}{i}[P_-] > \egameplay{S_-}{S_+}{i}[P_-]$.
\item If $\lef(S_-) = -\infty$ and $\lef(S_+) > - \infty$ and $i \in \Z$, $i < \lef(S_+)$, satisfies $S_-(i,1) > 0$, then  $\egameplay{\drag^{\leftarrow i}(S_-)}{S_+}{i}[P_-] > \egameplay{S_-}{S_+}{i}[P_-]$.
\end{enumerate}
\end{lemma}
{\bf Proof: (1).} The idle zone $\mc{I}$ determined by $(S_-,S_+)$ includes $\llbracket i,\infty)$. By Lemma~\ref{l.idlezone}, $X$ under $\pgameplay{S_-}{S_+}{i}$ thus behaves as a simple random walk when it visits $\llbracket i,\infty)$. Right escape $E_+$ is thus impossible, so Maxine's mean terminal payoff $\egameplay{S_-}{S_+}{i} [T_+]$ is at most $m_{-\infty}$ because it 
is a weighted average of $m_*$ and $m_{-\infty}$. Since $P_+ \leq T_+$ in view of running costs $C_+$ in~(\ref{e.maxinepayoff}) being non-negative, we find that $\egameplay{S_-}{S_+}{i}[P_+] \leq m_{-\infty}$.
Now consider $\pgameplay{S_-}{\boost_\eta^{i\rightarrow}(S_+)}{i}$. Under this law, $X$ proceeds non-randomly by rightward steps, so that $E_+$ occurs almost surely. Since $E_+$ occurs, we have $T_+ = m_\infty$ almost surely. By the non-random rightward movement, we further have that $\sum_{t=1}^\infty C_+(t)$ equals $\sum_{t=1}^\infty \eta \cdot 2^{-t} = \eta$. By~(\ref{e.maxinepayoff}),
we see then that    $\egameplay{S_-}{\boost_\eta^{i\rightarrow}(S_+)}{i}[P_+] = m_\infty - \eta$. This confirms Lemma~\ref{l.boostdrag}(1).

{\bf (2).} Under gameplay governed by the law $\pgameplay{S_-}{S_+}{i}$, Maxine offers a positive stake at $i$, and at infinitely many locations to its right, while Mina offers no stake at or to the right of $i$.
Thus $X:\N \to \Z$, $X_0 = i$, remains always to the right of $X$, and tends to infinity. If Maxine switches from  $S_+$ to $\drag^{i\rightarrow}(S_+)$, the law of gameplay $X$ is unaffected, because the original and altered gameplays may be coupled so that Maxine's altered stake process is one-half of her original one, while Mina's remains identically zero---with the result that Maxine wins exactly the same turns in the altered gameplay as she did in the original one. The value of $T_+$ is thus almost surely equal to $m_\infty$
under 
  $\pgameplay{S_-}{\drag^{i\rightarrow}(S_+)}{i}$ as well under  $\pgameplay{S_-}{S_+}{i}$. But  
  $\egameplay{S_-}{\drag^{i\rightarrow}(S_+)}{i} \sum_{t=1}^\infty C_+(t) 
  = \tfrac{1}{2} \egameplay{S_-}{S_+}{i} \sum_{t=1}^\infty C_+(t) $ and $\egameplay{S_-}{S_+}{i} \sum_{t=1}^\infty C_+(t) \geq  \egameplay{S_-}{S_+}{i} [C_+(1)] > 0$, so that $\egameplay{S_-}{\drag^{i\rightarrow}(S_+)}{i} \sum_{t=1}^\infty C_+(t) <  \egameplay{S_-}{S_+}{i} \sum_{t=1}^\infty C_+(t)$. In summary, the switch to the altered strategy has maintained Maxine's terminal receipt but has reduced her running costs, so that Lemma~\ref{l.boostdrag}(2) holds by~(\ref{e.maxinepayoff}).
  
{\bf (3,4).} The preceding proofs may be readily adapted to prove these statements. \qed

\begin{lemma}\label{l.zeronotnash}
\leavevmode
\begin{enumerate}
\item Any element of $\mc{N}$ has finite mean costs.\footnote{When the value of $(m_{-\infty},m_\infty,n_{-\infty},n_\infty)$ is clear---and it is usually a generic quadruple satisfying~(\ref{e.quadruple})---we will often omit to record this notation when we denote~$\mc{N}$. This includes the present case, where such a generic value is specified by the result, Theorem~\ref{t.nashabmn}(1), that we are seeking to prove.}
\item 
If $(S_-,S_+) \in \tis^2$ is an element of~$\mc{N}$ then  $S_-$ and $S_+$ are wide.
\end{enumerate}
\end{lemma}
In the ensuing proof and later, we denote the identically zero strategy by $0$.

{\bf Proof of Lemma~\ref{l.zeronotnash}(1).} 
 Let $(S_-,S_+) \in \mc{N}$, and let $i \in \Z$. Note that $
 \egameplay{S_-}{S_+}{i} [P_+] \geq  \egameplay{S_-}{0}{i}  [P_+]$.
 In evaluating the latter term, note that no running costs to Maxine have been incurred, so that the quantity is an average of terminal receipts $m_\infty$, $m_{-\infty}$ and $m_*$ to Maxine in the events $E_+$, $E_-$ and $E^c$. 
 We see that  $\egameplay{S_-}{0}{i}  [P_+] \geq \min \{ m_{-\infty},m_\infty,m_* \} = m_* > -\infty$, the latter inequality by assumption.  Likewise,   $\egameplay{S_-}{S_+}{i}[  P_-] > -\infty$. 
 
 {\bf (2).} We argue by contradiction and suppose, without loss of generality---for the other case is similar---that $S_-$ is narrow.
 (Lemma~\ref{l.boostdrag}(4) is not used in the ensuing proof. It is needed for the case whose proof we omit).
  Either $\lef(S_-) > -\infty$ or $\righ(S_-) < \infty$.
 
 Suppose that $\righ(S_-) < \infty$. If $\righ(S_+) < \infty$, then Lemma~\ref{l.boostdrag}(1) provides a strategy $\hat{S}_+$ to Maxine along with a value of $i \in \Z$ such that
 $\egameplay{S_-}{\hat{S}_+}{i}[P_+] > \egameplay{S_-}{S_+}{i}[P_+]$. But this is contrary to $(S_-,S_+) \in \mc{N}$. If $\righ(S_+) = \infty$, then it is Lemma~\ref{l.boostdrag}(2) that provides such $\hat{S}_+ \in \tis$
 and $i \in \Z$. A contradiction has thus been found in the case that  $\righ(S_-) < \infty$.
 
 Suppose now that  $\lef(S_-) > -\infty$. If $\lef(S_+) > -\infty$, then Lemma~\ref{l.boostdrag}(3) furnishes a strategy $\hat{S}_-$ for Mina and an index $i \in \Z$ 
 for which $\egameplay{\hat{S}_-}{S_+}{i}[P_-] > \egameplay{S_-}{S_+}{i}[P_-]$ holds, contrary to $(S_-,S_+) \in \mc{N}$.
 
 The case that  $\lef(S_-) > -\infty$ and $\lef(S_+) = -\infty$ remains. 
 The pair $(S_-,S_+) \in \tis^2 \cap \mc{N}$ is non-zero, because $S_+$ is; it has finite mean costs by Lemma~\ref{l.zeronotnash}(1).
 Thus $\pgameplay{S_-}{S_+}{i}(E^c) = 0$ by Proposition~\ref{p.fmc}. Select $i \in \Z$ for which $S_+(i,1) > 0$ and $S_-(j,1) = 0$ for $j \in (-\infty, i \rrbracket$. Note that 
  $\pgameplay{S_-}{S_+}{i}(E_-) = 0$ because gameplay $X$ is at least $i$ almost surely. Thus,  $\pgameplay{S_-}{S_+}{i}(E_+) = 1$, so that $T_+ = m_\infty$ almost surely. If Maxine drags down her strategy $S_+$ by replacing the stake she offers at $i$ to be one-half of its value, the resulting strategy $\hat{S}_+$
  is such that gameplay $X:\N \to \Z$ is equal under the laws  $\pgameplay{S_-}{S_+}{i}$ and $\pgameplay{S_-}{\hat{S}_+}{i}$; $T_+ = m_\infty$ almost surely under each of them; but $\sum_{t=1}^\infty C_+(t)$ is almost surely less under $\pgameplay{S_-}{\hat{S}_+}{i}$ than it is under $\pgameplay{S_-}{S_+}{i}$, because the value of $C_+(1)$ is lower. 
  Thus~(\ref{e.maxinepayoff}) shows that
   $\egameplay{\hat{S}_-}{S_+}{i}[P_+] > \egameplay{S_-}{S_+}{i}[P_+]$. 
   Again, we have a contradiction to  $(S_-,S_+) \in \mc{N}$. This completes the proof of Lemma~\ref{l.zeronotnash}(2). \qed

 \begin{corollary}\label{c.nashescape}
 For $(S_-,S_+) \in \mc{N} \cap \tis^2$ and $i \in \Z$, $\pgameplay{S_-}{S_+}{i}(E) = 1$. 
 \end{corollary}
{\bf Proof.} Due to Proposition~\ref{p.fmc} and Lemma~\ref{l.zeronotnash}(1,2). \qed

Recall that an element $(S_-,S_+) \in \tis^2$ may be identified as a sequence $\big\{ (b_i,a_i): i \in \Z \big\}$ to which Definition~\ref{d.quadruple}
associates a quadruple  $\big\{ (a_i,b_i,m_i,n_i): i \in \Z \big\}$.
\begin{lemma}\label{l.mnincdec}
Suppose that $(S_-,S_+) \in \mc{N} \cap \tis^2$. Then $m_i \leq m_{i+1}$ and $n_{i+1} \leq n_i$ for $i \in \Z$.
\end{lemma}
{\bf Proof.} Recall that 
$\pgameplay{S_-}{S_+}{i}$ 
denotes the law of gameplay when $X_0 = i$. Let $\sigma_{i+1} \in \nwithoutzero \cup \{ \infty \}$ denote the stopping time $\inf \big\{ \ell \in \nwithoutzero : X_\ell = i+1 \big\}$. 
Noting the non-negativity of running costs $C_+(t)$ in Lemma~\ref{l.stopping} with $k=i$ and $Q = \sigma_{i+1}$, we find that
$$
\egameplay{S_-}{S_+}{i} [P_+] \, \leq \,   \egameplay{S_-}{S_+}{i} \big[ \egameplay{S_-}{S_+}{X(\sigma_{i+1})}[P_+] \big] \, ,
$$
whose left-hand side equals $m_i$ by definition and whose right-hand side takes the form
$$
 m_{i+1} \pgameplay{S_-}{S_+}{i} \big( \sigma_{i+1} < \infty \big) + m_{-\infty} \pgameplay{S_-}{S_+}{i} \big( \sigma_{i+1} = \infty, E \big) 
 + m_* \pgameplay{S_-}{S_+}{i} \big( \sigma_{i+1} = \infty, E^c \big) \, .
 $$
 However, the third term vanishes in view of Corollary~\ref{c.nashescape}. Thus, $m_i$ is seen to be bounded above by a weighted average of $m_{-\infty}$ and $m_{i+1}$.
 To conclude, as we seek to do, that $m_i \leq m_{i+1}$, it is thus enough to argue that $m_{-\infty} \leq m_{i+1}$.
 To obtain this bound, we first {\em claim} that $\egameplay{S_-}{0}{i+1} [P_+] = m_{-\infty}$. To check this, we  invoke Lemma~\ref{l.zeronotnash}(2) to say that $S_-$ is wide.
  Thus, 
 $E_-$, and $T_+ = m_{-\infty}$, are $\pgameplay{S_-}{0}{i+1}$-almost certain. 
The absence of running costs for Maxine means that $P_+ = T_+$ under $\pgameplay{S_-}{0}{i+1}$. This yields the claim. Using it, and $(S_-,S_+) \in \mc{N}$, we find that
$$
 m_{i+1} \, = \, \egameplay{S_-}{S_+}{i+1}[ P_+] \, \geq \, \egameplay{S_-}{0}{i+1} [P_+] = m_{-\infty} \, .
 $$
 We have confirmed that $m_i \leq m_{i+1}$. We omit the similar proof that  $n_{i+1} \leq n_i$.  This completes the proof of Lemma~\ref{l.mnincdec}. \qed

\begin{lemma}\label{l.firstrearranged}
Let  $\big\{ (b_i,a_i): i \in \Z \big\} \in \mc{N} \cap \tis^2$.  
Recall from Definition~\ref{d.quadruple} that $m_i$ equals Maxine's mean receipt when the counter starts at $i \in \Z$.
Suppose that $a_i + b_i > 0$. Then
\begin{equation}\label{e.firstrearranged}
m_i   =  \tfrac{a_i}{a_i + b_i} m_{i+1} + \tfrac{b_i}{a_i + b_i} m_{i-1} - a_i \, .
\end{equation}
\end{lemma}
{\bf Proof.}
 Maxine will spend $a_i$ at the first turn; she will win the turn with probability $\tfrac{a_i}{a_i + b_i}$; if she does so, the counter will reach $i+1$, and her resulting conditional mean receipt will be $m_{i+1}$; if she does not, this receipt will instead be $m_{i-1}$. Note that the two ratios on the right-hand side of~(\ref{e.firstrearranged}) are well defined, because  $a_i + b_i > 0$. \qed

\begin{lemma}\label{l.condpositive}
Let $(S_-,S_+) \in \mc{N} \cap \tis^2$, and let $i \in \Z$. Then $a_i > 0$ implies that $m_{i+1} > m_i$. And $b_i > 0$ implies that $n_{i-1} > n_i$.
\end{lemma}
{\bf Proof.} Lemma~\ref{l.firstrearranged} and $a_i > 0$ imply that $m_i < \max \{ m_{i-1},m_{i+1} \}$. But the maximum is attained by $m_{i+1}$ in view of Lemma~\ref{l.mnincdec}. 
The second assertion in the lemma is similarly obtained.
\qed

\begin{proposition}\label{p.allpositive}
Let $(S_-,S_+) \in \mc{N} \cap \tis^2$. Then $a_i > 0$, $b_i > 0$, $m_{i+1} > m_i$ and $n_i > n_{i+1}$ for all~$i \in \Z$.
\end{proposition}
{\bf Proof.} By Lemma~\ref{l.zeronotnash}(2), 
$S_-$ is wide. 
To show that every $a$-coefficient is positive, it is thus enough to 
argue that $a_i > 0$ implies $a_{i+1} > 0$ for $i \in \Z$,
because every index $i \in \Z$ has a positive $a$-coefficient indexed somewhere to its left.
Suppose to the contrary that $a_i > 0$ but $a_{i+1} = 0$. Applying~(\ref{e.firstrearranged}) at index~$i+1$, we see that $b_{i+1} > 0$ implies that  $m_{i+1} = m_i$. But Lemma~\ref{l.condpositive} and $a_i > 0$ imply that $m_{i+1} > m_i$. Thus, $b_{i+1} = 0$. In view of $a_{i+1} = 0$, we see from~(\ref{e.firstrearranged}) at index $i+1$ (with use of the $0/0 = 1/2$ rule) that $m_{i+1} = \tfrac{m_i + m_{i+2}}{2}$.

However, given that $b_{i+1} = 0$, Maxine can obtain a better outcome when playing at $i+1$ by staking a small positive amount in place of zero. Consider gameplay from~$i+1$ and a deviation~$S'_+$ to Maxine's play in which she nominates a small positive stake~$\zeta$ at the first turn but then adheres to~$S_+$. Because $(S_-,S_+)$ governs gameplay from the second turn, $(m,n)$ accurately reports mean receipts at such times, including as boundary data for the first turn, played at $i+1$. Returning then to~(\ref{e.firstrearranged}), we see that  
a sufficiently small choice of $\zeta$ would yield a value for $\E^{i+1}_{S_-,S_+'}[P_+]$
which is arbitrarily close to $m_{i+2}$. But the latter quantity exceeds  $(m_i + m_{i+2})/2$ because (in view of Lemma~\ref{l.condpositive}, $a_i > 0$ and $a_{i+1} > 0$) we have the bound $m_{i+2} > m_i$.
Thus, $(S_-,S_+) \not\in \mc{N}$, contrary to assumption. We have confirmed that $a_{i+1} > 0$, and thus that every $a$-coefficient is positive.
The argument that $b_i > 0$ for $i \in \Z$ is no different. Lemma~\ref{l.condpositive} then shows that each difference $m_{i+1} - m_i$ and $n_i - n_{i+1}$ is positive. 
This completes the proof of Proposition~\ref{p.allpositive}. \qed

We may now prove the first part of  Theorem~\ref{t.nashabmn}.
 
{\bf Proof of Theorem~\ref{t.nashabmn}(1).}  Suppose  that 
$(S_-,S_+) \in  \tis^2$ is a time-invariant Nash equilibrium for $\game(m_{-\infty},m_\infty,n_{-\infty},n_\infty)$. We abusively identify $(S_-,S_+)$ with the sequence 
 $\big\{ (b_i,a_i): i \in \Z \big\} \in \tis^2$ as usual (and, by doing so, we conform notation with the theorem's statement).
We note at the outset that, in view of Proposition~\ref{p.allpositive}, each $a_i$ and $b_i$, and each difference $m_{i+1} - m_i$ and $n_i - n_{i+1}$, is positive.

Equation \abmnmac$(1)$ results from rearranging the formula in Lemma~\ref{l.firstrearranged}. Equation \abmnmac$(2)$ is similarly derived.

Next we derive \abmnmac$(3,4)$. Recall $S_-(i,j) = b_i$ and $S_+(i,j)=a_i$ for each $(i,j) \in \Z \times \nwithoutzero$. 
We will consider a first-turn deviation to Maxine's strategy, as we did in the proof of Proposition~\ref{p.allpositive}.
Indeed,
for given $i \in \Z$, let  $\hat{S}_+ \in \mc{S}$ denote a perturbed strategy  for Maxine in which only her first-turn stake is altered, and only then if the counter is at $i$. In this way,  $\hat{S}_+(j,k) = a_j$ for $j \in \Z$ and $k \geq 2$; and also for $k=1$ and $j \in \Z$, $j \not= i$. We let $\eta > -a_i$ be small in absolute value, and set $\hat{S}_+(1,i) = a_i + \eta$.

The {\em original} scenario refers to the law $\pgameplay{S_-}{S_+}{i}$, which records counter evolution~$X:\nwithzero \to \Z$ given the initial condition $X_0 = i$ under the strategy pair $(S_-,S_+)$. The {\em altered} scenario refers to the same law, instead governed by the pair $(S_-,\hat{S}_+)$. Let $O_+$ and $A_+$ denote the mean payoff to Maxine in the original and altered scenarios: that is, $O_+ = \egameplay{S_-}{S_+}{i} [P_+]$ and $A_+ = \egameplay{S_-}{\hat{S}_+}{i} [P_+]$.    
Then
$$
 O_+ = \tfrac{a_i}{a_i+b_i} m_{i+1} + \tfrac{b_i}{a_i+b_i} m_{i-1} - a_i \, \, \, \textrm{and} \, \, \, A_+ = \tfrac{a_i+\eta}{a_i+\eta+ b_i} m_{i+1} + \tfrac{b_i}{a_i+\eta+b_i} m_{i-1} - a_i - \eta \, ,
$$
so that
\begin{equation}\label{e.aodifference}
A_+ - O_+ \, = \, \Big( \tfrac{b_i}{(a_i+b_i)^2} (m_{i+1} - m_{i-1})  - 1 \Big) \cdot \eta \cdot \big( 1 + o(1) \big) \, ,
\end{equation}
where the $o(1)$ term is small in the sense of $\vert \eta \vert \to 0$.
Since $(S_-,S_+) \in \mc{N}$, $A_+$ is at most $O_+$, whatever the value of $\eta > - a_i$. The derivative in $\eta$ of $A_+ - O_+$ thus vanishes at zero, so that 
$\tfrac{b_i}{(a_i+b_i)^2} (m_{i+1} - m_{i-1})  - 1 = 0$
or equivalently
\begin{equation}\label{e.bma}
  b_i  (m_{i+1} - m_{i-1})   = (a_i+b_i)^2 \, .
\end{equation}
We now consider the same original scenario alongside a new altered scenario in which it is Mina who adopts a perturbed strategy $\hat{S}_-$ (as a function of a given choice of $i \in \Z$). Analogously to what we have done, we choose $\eta > - b_i$,
and set $\hat{S}_-(j,k) = b_j$ for  $j \in \Z$ and $k \geq 2$  or when $k=1$ and $j \in \Z$, $j\not=i$; and then we set $\hat{S}_-(1,i) = b_i + \eta$. 
We denote by  $O_-$ and $A_-$ Mina's mean payoff in the original and in the newly altered scenarios; to wit,  $O_- = \egameplay{S_-}{S_+}{i} [P_-]$ and $A_- =  \egameplay{\hat{S}_-}{S_+}{i} [P_-]$.
We find then that
$$
 O_- = \tfrac{b_i}{a_i+b_i} n_{i-1} + \tfrac{a_i}{b_i+a_i} n_{i+1} - b_i \, \, \, \textrm{and} \, \, \, A_- = \tfrac{b_i+\eta}{b_i+\eta+ a_i} n_{i-1} + \tfrac{a_i}{a_i+\eta+b_i} n_{i+1} - b_i - \eta \, ;
$$
and, analogously to~(\ref{e.aodifference}),
$$
A_- - O_- \, = \, \Big( \tfrac{a_i}{(a_i+b_i)^2} (n_{i-1} -n_{i+1})  - 1 \Big) \cdot \eta \cdot \big( 1 + o(1) \big) \, .
$$
The condition that $(S_-,S_+) \in \mc{N}$ ensures that $O_- \geq A_-$, whatever the value of $\eta > - b_i$. Thus,
\begin{equation}\label{e.anb}
 a_i  (n_{i-1} - n_{i+1})   = (a_i+b_i)^2 \, .
\end{equation}
The derived equations~(\ref{e.bma}) and~(\ref{e.anb}) are \abmnmac$(3,4)$ with index $i$. 

We have established that $\big\{ (a_i,b_i,m_i,n_i): i \in \Z \big\}$
solves the \abmnmacspace equations. 
To complete the proof of Theorem~\ref{t.nashabmn}(1), it remains to confirm that the boundary values~(\ref{e.boundarydata}) are achieved.
We will argue that $\lim_{i \to \infty} m_{-i} = m_{-\infty}$; the three other limits are similarly shown. 
The sequence $\big\{ m_{-i}: i \in \nwithzero \big\}$ decreases by Proposition~\ref{p.allpositive} to a limiting value that we call $\mathfrak{m}_{-\infty}$.
By Definition~\ref{d.quadruple}, $(S_-,S_+) \in \mc{N}$ and the straightforward `$S_-$ is wide'-consequence~$\pgameplay{S_-}{0}{i} (E_-)=1$, we have that $m_i = \egameplay{S_-}{S_+}{i} [P_+] \geq \egameplay{S_-}{0}{i} [P_+] = m_{-\infty}$, whence $\mathfrak{m}_{-\infty} \geq m_{-\infty}$; we wish to obtain the opposite inequality. Take the mean of~(\ref{e.maxinepayoff}) and
 remove non-negative running costs to find that $m_i \leq \pgameplay{S_-}{S_+}{i}(E_-)\cdot m_{-\infty} +  \pgameplay{S_-}{S_+}{i}(E_+) \cdot m_\infty$  where we invoked Corollary~\ref{c.nashescape}.
 Thus $\mathfrak{m}_{-\infty} \leq m_{-\infty}$ provided that we argue that  $\lim_{i \to -\infty} \pgameplay{S_-}{S_+}{i}(E_+) = 0$: far to the left is the domain of Mina's likely victory. It would be of interest to argue directly; and to do so would be more in keeping with the style of this section. It is quicker however to simply invoke the eventual gameplay unanimity Theorem~\ref{t.unanimity}(3), which will be proved by independent arguments when we find an explicit solution of the \abmnmacspace system in Chapter~\ref{c.battlefield}. 
It should be noted in this regard that the right limit $\mathfrak{m}_\infty$ strictly exceeds $\mathfrak{m}_{-\infty}$ in view of Proposition~\ref{p.allpositive}; this means that the trivial zero  \abmnmacspace solution is out of the picture, enabling the use of  Theorem~\ref{t.unanimity}(3).
(Note also that Theorem~\ref{t.unanimity}(2) is invoked in the corresponding place in two of the three omitted limit derivations.)
We have thus obtained   Theorem~\ref{t.nashabmn}(1).
  \qed


\section{The reverse implication}
Here we prove Theorem~\ref{t.nashabmn}(2). It is here that the infinite-turn nature of the game has to be tamed by comparison with finite-trail counterparts.

We begin by developing definitions and results that will lead to the proof of the desired result at the end of the section.
As such, we now enforce the notation in the hypothesis of Theorem~\ref{t.nashabmn}(2), so that, from now on, 
$$
\text{$\big\{ (a_i,b_i,m_i,n_i) : i \in \Z \big\}$ denotes a positive 
solution of the \abmnmacspace equations}
$$ 
with boundary data $(m_{-\infty},m_\infty,n_{-\infty},n_\infty)$ that satisfies~(\ref{e.quadruple}).
Let $S_-,S_+ \in \mc{S}$ satisfy
\begin{equation}\label{e.ba}
 S_-(i,j) = b_i \, \, \,\,\textrm{and} \, \, \, \, S_+(i,j) = a_i \, \, \, \, \textrm{for $(i,j) \in \Z \times \nwithoutzero$}
 \, .
\end{equation}
Recall that $S_+ = a$ with the usual notational abuse.
\begin{definition}\label{d.mds}
Let $i \in \Z$. The forward play-cone $F_i$ of $i$ is given by 
$$
F_i \, = \, \Big\{ \, (k,\ell) \in \Z \times \nwithoutzero: \vert k - i \vert \leq \ell \, , \, \vert k-i \vert + \ell  \in 2\nwithzero \, \Big\} \, .
$$
This is the set of space-time sites that are in principle accessible for gameplay $X:\nwithzero \to \Z$ under $\pgameplay{S_1}{S_2}{i}$ for some strategy pair $(S_1,S_2) \in \mc{S}^2$.

Let $S \in \mc{S}$. An element $(q,\ell) \in F_i$ such that  $S(q,\ell+1) \not= b_q$ is called a {\em Mina deviation point}. 
The Mina deviation set $\mathsf{D}_-(S,i) \subseteq F_i$ is the collection of Mina deviation points. The strategy $S$ is called {\em deviating for Mina} if  $\mathsf{D}_-(S,i)$ is non-empty.
 A {\em Maxine deviation point} $(q,\ell) \in F_i$ satisfies  $S(q,\ell+1) \not= a_q$. The set $\mathsf{D}_+(S,i)$ of such points is the Maxine deviation set; if  $\mathsf{D}_+(S,i) \not= \emptyset$, then $S$ is deviating for Maxine.  
\end{definition}
When gameplay under $\pgameplay{S}{S_+}{i}$ runs through a Mina deviation point---when $X_\ell = q$ for $(q,\ell) \in \mathsf{D}_-(S,i)$---her stake according to strategy $S$---namely, $S(q,\ell+1)$---may be viewed as a mistake when her opponent plays her element $S_+$
of the putative Nash equilibrium $(S_-,S_+)$. The next result, which is fundamental to proving Theorem~\ref{t.nashabmn}(2), validates this notion. It measures the magnitude of the mistakes that result from a player's deviation in the sense of decrease in mean payoff in finite trail games. It finds the mistakes to be uniformly costly as the finite trails vary.
\begin{proposition}\label{p.jksup}
Let 
$i \in \Z$ be given. 
\begin{enumerate}
\item 
Let $\sdev \in \mc{S}$ be deviating for Mina.  Suppose that $\PP_{\sdev,S_+}^i(E) = 1$.
 Then
$$
\sup \, 
 \E_{\sdev,S_+}^i [P_-^{j,k}]  \, < \, 
 \E_{S_-,S_+}^i [P_-]  \, ,
$$
where the supremum is taken over $j,k \in \nwithoutzero$ such that $i \in \llbracket -j,k \rrbracket$ 
and for which there exists an element $(u,\ell)$ of $\mathsf{D}_-(\sdev,i)$ with $u \in  \llbracket -j,k \rrbracket$.
\item Now suppose that $\sdevmaxine \in \mc{S}$ is deviating for Maxine, and   $\PP_{S_-,\sdevmaxine}^i(E) = 1$.
 Then
$$
\sup \,
 \E_{S_-,\sdevmaxine}^i [P_+^{j,k}]  \, < \, 
 \E_{S_-,S_+}^i [P_+]  \, ,
$$
where now the  supremum is taken over $j,k \in \nwithoutzero$ with $i \in \llbracket -j,k \rrbracket$ and 
for which there exists $(u,\ell) \in \mathsf{D}_+(\sdevmaxine,i)$ such that $u \in  \llbracket -j,k \rrbracket$.
\end{enumerate}
\end{proposition}
It is a short step from the just stated result to the next conclusion, which asserts that a player's deviation will cost her in the infinite trail game. This is in essence what it means for $(S_-,S_+)$
to be a Nash equilibrium. Indeed, we next close out the proof of Theorem~\ref{t.nashabmn} by first deriving Proposition~\ref{p.sminuscomp} from Proposition~\ref{p.jksup}; and second showing how the latter result leads to the desired conclusion. These tasks done, we will turn to the remaining and more substantial one: to prove Proposition~\ref{p.jksup}. 
\begin{proposition}\label{p.sminuscomp}
Let $i \in \Z$. 
\begin{enumerate}
\item Let $\sdev \in \mc{S}$ be deviating for Mina. Then
$$
 \E_{\sdev,S_+}^i [P_-]  <  \E_{S_-,S_+}^i [P_-]  \, .
$$
\item Now  let  $\sdevmaxine \in \mc{S}$ be deviating for Maxine.
 Then
$$
 \E_{S_-,\sdevmaxine}^i [P_+]  <  \E_{S_-,S_+}^i [P_+]  \, .
$$
\end{enumerate}
\end{proposition}
{\bf Proof: (1).}
Suppose first that $\PP_{\sdev,S_+}^i(E^c) > 0$. 
Lemma~\ref{l.dontlookback} implies that 
$\E_{(\sdev,S_+)}^i [P_-]  = -\infty$. But  $ \E_{S_-,S_+}^i [P_-]  = n_i$ by Lemma~\ref{l.minipayoff}(2). We have  that $n_i \geq n_\infty$
since the sequence $\big\{ n_i: i \in \Z \big\}$ decreases for any positive \abmnmacspace solution by Theorem~\ref{t.positiveabmn}(1). And we know that $n_\infty > -\infty$ by hypothesis.
Thus we see that   $ \E_{S_-,S_+}^i [P_-]  > -\infty$, so that  Proposition~\ref{p.sminuscomp}(1) has been established in this case.

Now we suppose instead  that $\PP_{\sdev,S_+}^i(E) = 1$. 
Let $\eta > 0$ denote a parameter whose value we will soon set.
Note that $T_-^{j,k} = n\big(X(\tau^{j,k})  \big)$ for $j,k \in \nwithoutzero$; and that $T_-$ equals $n_{-\infty}$ on $E_-$, and $n_\infty$ on~$E_+$.
By Lemma~\ref{l.mn} and  $\PP_{\sdev,S_+}^i(E) = 1$,
we may thus find   $j_0,k_0 \in \nwithoutzero$ such that, when $j \geq j_0$ and $k \geq k_0$,
$$
\PP_{\sdev,S_+}^i \Big( \big\vert  T_-  - T_-^{j,k} \big\vert \geq \eta \Big) \leq \eta \, .
$$
By Lemma~\ref{l.couplingproperties}(2), we see that
$$
\PP_{\sdev,S_+}^i \Big( P_-  \, \leq \, P^{j,k}_-   + \eta \Big) \geq 1 - \eta  \, .
$$
Finite escape $\tau^{j,k} < \infty$ is almost sure to occur, since  $\PP_{\sdev,S_+}^i(E) = 1$. Thus, by Lemma~\ref{l.couplingproperties}(3),  
$$
\E_{\sdev,S_+}^i [P_-] \leq  \E_{\sdev,S_+}^i  [P^{j,k}_-]  + (1+n_{-\infty}- n_\infty)\eta \, .
$$
Let $\zeta$ denote  the difference of the two sides in the conclusion of Proposition~\ref{p.jksup}(1). Then
$E_{\sdev,S_+}^i [P_-]$ is at most  $\E_{S_-,S_+}^i  [P_-] - \zeta  + (1+n_{-\infty}- n_\infty)\eta$. Choosing $\eta > 0$ so that  $(1+n_{-\infty}- n_\infty)\eta < \zeta$, we obtain  Proposition~\ref{p.sminuscomp}(1).

 By taking $\eta > 0$ to be one-half of the difference of the two sides in the conclusion of Proposition~\ref{p.jksup}(1), the latter result is seen to imply Proposition~\ref{p.sminuscomp}(1). 
 
 {\bf (2).} We omit this similar argument. \qed

 {\bf Proof of Theorem~\ref{t.nashabmn}(2).}
Recall that $\big\{ (a_i,b_i,m_i,n_i): i \in \Z \big\}$ is a positive \abmnmacspace solution
with boundary data $(m_{-\infty},m_\infty,n_{-\infty},n_\infty) \in \R^4$.
 Further recall that $(S_-,S_+) = (b,a)$, with the usual notational abuse. Let $S \in \mc{S}$.
If $S$ is not deviating for Mina, then $\egameplay{S}{S_+}{i} [P_-] = \egameplay{S_-}{S_+}{i} [P_-]$ 
since the laws $\pgameplay{S}{S_+}{i}$ and $\pgameplay{S_-}{S_+}{i}$ are equal.
Otherwise,  $\egameplay{S}{S_+}{i}[P_-] < \egameplay{S_-}{S_+}{i} [P_-]$ by Proposition~\ref{p.sminuscomp}(1). (We recall that implicit in the notation $\pgameplay{S_1}{S_2}{i}$ and $\egameplay{S_1}{S_2}{i}[\cdot]$ are the values $(m_{-\infty},m_\infty,n_{-\infty},n_\infty)$, because these values appear as terminal receipts.)

By Proposition~\ref{p.sminuscomp}(2), it follows similarly that
$\egameplay{S_-}{S}{i} [P_+] < \egameplay{S_-}{S_+}{i} [P_+]$ if $S$ is deviating for Maxine.
Further, $\egameplay{S_-}{S}{i} [P_+] = \egameplay{S_-}{S_+}{i} [P_+]$ if Maxine's $S$ is not deviating. 
We have confirmed that $(S_-,S_+) \in \mc{N}(m_{-\infty},m_\infty,n_{-\infty},n_\infty)$ and thus obtain Theorem~\ref{t.nashabmn}(2). \qed

We now prepare to prove Proposition~\ref{p.jksup}(1). 
(The proof of  Proposition~\ref{p.jksup}(2) is essentially the same.)
 Henceforth,  Proposition~\ref{p.jksup}(1)'s hypotheses are understood to be in force: $S_-$ and $S_+$ are the non-deviating strategies given by~(\ref{e.ba});
  $i \in \Z$ is given; and $\sdev \in \mc{S}$ is deviating for Mina, with $\PP_{\sdev,S_+}^i(E) = 1$.

Let $j,k \in \nwithzero$ satisfy $i \in \llbracket -j,k \rrbracket$. Developing Definition~\ref{d.mds}, we
set
$$
\mathsf{D}_-^{j,k}(S,i) \, = \, \Big\{ (q,\ell) \in \mathsf{D}_-(S,i): q \in \llbracket -j,k \rrbracket \Big\} 
$$
for $S \in \mc{S}$.
It may be that $\mathsf{D}_-^{j,k}(\sdev,i)$ is infinite. 
It serves our purpose to approximate $\sdev$ by strategies for which the counterpart set is finite. We now specify these strategies.
Enumerate  $\mathsf{D}_-^{j,k}(\sdev,i)$ in increasing order of the vertical component, using an arbitrary rule to break the ties that arise when elements share the same height.
For $v \in \nwithoutzero$, let  $\mathsf{D}_{-,v}^{j,k}(\sdev,i)$ denote the set whose elements are the first $v$ members of  $\mathsf{D}_-^{j,k}(\sdev,i)$. 
Let $\sdev[v]$ denote the strategy that equals $\sdev$ on $\mathsf{D}_{-,v}^{j,k}(\sdev,i)$ and $S_-$ otherwise; note that $\mathsf{D}_{-}^{j,k}(\sdev[v],i)$ equals $\mathsf{D}_{-,v}^{j,k}(\sdev,i)$.
We make another basic comparison in terms of the next definition.
\begin{definition}\label{d.ground}
For $S \in \mc{S}$,  let $\ground^{j,k}(S,i) \in \nwithzero \cup \{ \infty\}$ denote  the minimum vertical coordinate assumed by an element of $\mathsf{D}_-^{j,k}(S,i)$, where the value $\infty$ is assumed when  $\mathsf{D}_-^{j,k}(S,i) =  \emptyset$.
\end{definition}
Note then that $\ground^{j,k}(\sdev[v],i)$  is independent of $v \in \nwithoutzero$. 

We wish to argue that Mina's deviant play under the strategies $\sdev[v]$, $v \in \nwithoutzero$, and $\sdev$, is suitably penalized in the trail game on $\llbracket -j-1,k+1 \rrbracket$.
In the notation of the next definition, Lemma~\ref{l.baseconseq} establishes such a conclusion for the finitely deviating strategies $\sdev[v]$: there is a penalty incurred by use of these strategies; and, in a suitable sense, the penalty is uniform among them, and is governed by the limiting strategy~$\sdev$.
After we prove Lemma~\ref{l.baseconseq}, it will remain to address the penalty suffered by using $\sdev$ itself. Definition~\ref{d.merit} speaks of a `strong' penalty as a contrast with a modified definition that will be used to treat the perhaps infinitely deviating $\sdev$, this appearing after the proof of Lemma~\ref{l.baseconseq}.
\begin{definition}\label{d.merit}
Let $S_1,S_2 \in \mc{S}$. Consider the following conditions:
\begin{enumerate}
\item We have that $\E^{u,\ell}_{S_1,S_+} \big[ P^{j,k}_- \big] \leq n_u$ for all $\ell \in \nwithoutzero$ and $u \in \llbracket -j,k \rrbracket$. 
\item Writing $g = \ground^{j,k}(S_1,i)$, consider any $u \in \llbracket -j,k \rrbracket$ for which $(u,g) \in \mathsf{D}_-^{j,k}(S_1,i)$. Then  the value $n_u -  \E^{u,g}_{S_1,S_+} \big[ P^{j,k}_- \big]$ is positive, and indeed is bounded below by a positive quantity that is determined solely by $S_2(u,g+1)$.
\end{enumerate}
If these conditions are met, we say that {\em $S_1$ receives the strong $(i,j,k)$-penalty merited by $S_2$}.
Let $S \in \mc{S}$. If $S$ receives the strong $(i,j,k)$-penalty merited by $S$, we say that {\em $S$  justly receives a  strong $(i,j,k)$-penalty}. 
\end{definition}
(Although it is omitted from the notation of  the strong $(i,j,k)$-penalty, it is the strategy $S_+ = a$ that Mina is facing when she plays $S_1$ or $S_2$. The above definition and the next result are intended to capture the sense of Mina's mistake when she declines to stake at the $b$-level dictated by $S_-$ against Maxine's $a$-stake offered by $S_+$.)
 \begin{lemma}\label{l.baseconseq}
 Let $j,k \in \nwithzero$ satisfy $i \in \llbracket -j,k \rrbracket$.
Let $v$ be at least the number of elements of $\mathsf{D}_-^{j,k}(\sdev,i)$ of minimum height. 
Then $\sdev[v]$ receives the strong $(i,j,k)$-penalty merited by $\sdev$.
 \end{lemma}
 (The value of $g$ implicit in Lemma~\ref{l.baseconseq} does not depend on the value of $v \in \nwithoutzero$ used in $\sdev[v]$, because $\ground^{j,k}(\sdev[v],i)$  is independent of $v \in \nwithoutzero$.)
 
 The finite-error strategies $\sdev[v]$ have been introduced because they may be analysed using the fundamental game-theoretic technique of backwards induction. When Mina uses $\sdev[v]$ for any given $v \in \nwithoutzero$, she never deviates at late enough time. Lemma~\ref{l.minipayoff}(2) then serves to show that she incurs no penalty by doing so. As turn index retreats in backwards induction, Mina will make deviating moves. At the heart of the analysis of the inductive step is the consideration of one turn when Mina deviates. What is being played here is a game of Penny Forfeit, treated in Section~\ref{s.pennyforfeit}. The next result gathers what we need to know about one step in the game.
\begin{lemma}\label{l.onestep}
\leavevmode
\begin{enumerate}
\item
Let $j,k \in \Z$ satisfy $i \in \llbracket -j,k \rrbracket$.
For $\ell \in \nwithoutzero$, let $S_1,S_2 \in \mc{S}$ be such that, if $(u,h) \in \Z \times \nwithoutzero$ satisfies
$S_1(u,h) \not= S_2(u,h)$, then $h \leq \ell$.
Then $\E^{u,h}_{S_1,S_+} \big[ P^{j,k}_- \big] = \E^{u,h}_{S_2,S_+} \big[ P^{j,k}_- \big]$ for any $(u,h) \in \llbracket -j,k \rrbracket \times \llbracket \ell, \infty)$.
\end{enumerate}
Let $S \in \mc{S}$ and $(u,\ell) \in \llbracket -j,k \rrbracket \times \nwithzero$. Suppose that $\E^{v,\ell+1}_{S,S_+} \big[ P^{j,k}_- \big] \leq n_j$ for $v \in \{u-1,u+1\}$.
\begin{enumerate}
 \setcounter{enumi}{1}
\item We have that
$\E^{u,\ell}_{S,S_+} \big[ P^{j,k}_- \big] \leq n_u$. 
\item Suppose further that 
 $(u,\ell) \in \mathsf{D}^{j,k}_-(S,i)$. Then $n_u - \E^{u,\ell}_{S,S_+} \big[ P^{j,k}_- \big]$ is bounded below by a positive quantity that is determined solely by the value of $S(u,\ell+1) \not= b_u$.
\end{enumerate}
\end{lemma}
{\bf Proof: (1).} The laws $\PP^{u,h}_{S_1,S_+}$ and $\PP^{u,h}_{S_2,S_+}$
are identical because $S_1$ and $S_2$ coincide at any point $(u,\ell)$ with $\ell \geq h+1$. \\
{\bf (2).}  Note that
$$
\E^{u,\ell}_{S,S_+} \big[ P^{j,k}_- \big] \, = \, \tfrac{S(u,\ell)}{a_i+S(u,\ell)}  \E^{u-1,\ell+1}_{S,S_+} \big[ P^{j,k}_- \big] 
+ \tfrac{a_i}{a_i+S(u,\ell)} \E^{u+1,\ell+1}_{S,S_+} \big[ P^{j,k}_- \big] - S(u,\ell) \, .
$$
Since $\E^{u-1,\ell+1}_{S,S_+} \big[ P^{j,k}_- \big] \leq n_{u-1}$ and $\E^{u+1,\ell+1}_{S,S_+} \big[ P^{j,k}_- \big] \leq n_{u+1}$, we see that
 $$
\E^{u,\ell}_{S,S_+} \big[ P^{j,k}_- \big] \, \leq \, \tfrac{S(u,\ell)}{a_u+S(u,\ell)}  n_{u-1}
+ \tfrac{a_u}{a_u+S(u,\ell)}  n_{u+1} - S(u,\ell) \, .
$$
 By Lemma~\ref{l.pennyforfeit}, this right-hand side has a unique maximum in $b$ at $b = b_u$, when it assumes the value $n_u$.\\
 {\bf (3).} Since $S(u,\ell+1)$ is not equal to $b_u$, we see that the above right-hand side, and thus $\E^{u,\ell}_{S,S_+} \big[ P^{j,k}_- \big]$, is less than $n_u$. The difference $n_u - \egameplay{S}{S_+}{u,\ell}[P^{j,k}_-]$ is determined solely by $S(u,\ell+1)$. \qed

The next result leads quickly to Lemma~\ref{l.baseconseq}. Indeed, its proof (in a perhaps slightly disguised form) is the backwards inductive argument that underlies Lemma~\ref{l.baseconseq}.
\begin{lemma}\label{l.backwardformal}
Suppose that $S \in \mc{S}$ is such that 
$\mathsf{D}_-^{j,k}(S,i)$ is finite. Then $S$ justly receives a strong $(i,j,k)$-penalty. 
\end{lemma}
{\bf Proof.}
We will induct on the cardinality of $\mathsf{D}_-^{j,k}(S,i)$. 

Let $S \in \mc{S}$. Set $g =  \ground^{j,k}(S,i)$.
 For $\ell \in \nwithzero$, $\ell \not= g$, let $\indhyp(S,\ell)$
denote the assertion that
\begin{equation}\label{e.basicineq}
\E^{u,\ell}_{S,S_+} \big[ P^{j,k}_- \big] \leq n_u \, \, \, \, \textrm{for} \, \, \, \, u \in \llbracket -j,k \rrbracket \, .
\end{equation}
For $\ell =g$, let $\indhyp(S,\ell)$ denote the assertion that the preceding display holds and so does the following.
\begin{eqnarray*}
 & & \textrm{Consider any $u \in \llbracket -j,k \rrbracket$ for which $(u,g) \in \mathsf{D}_-^{j,k}(S,i)$.}\\
 & & \textrm{Then  the value $n_u -  \E^{u,g}_{S,S_+} \big[ P^{j,k}_- \big]$ is positive,} \\
 & & \textrm{and indeed is bounded below by a positive quantity that is determined solely by $S(u,g+1)$.}
 \end{eqnarray*}
We take the inductive hypothesis indexed by $q \in \nwithoutzero$ to be the assertion that the statements  $\indhyp(S,\ell)$, $\ell \in \nwithzero$, are true for each 
$S \in \mc{S}$ such that  $\char"0023 \,  \mathsf{D}_-^{j,k}(S,i) \leq q$.
The base case will be $q = 0$. This is the assertion that~(\ref{e.basicineq}) holds for $\ell \in \nwithoutzero$, when $S \in \mc{S}$
is such that  $\mathsf{D}_-^{j,k}(S,i)$ is empty.

The base case holds by Lemma~\ref{l.minipayoff}(3).


Let $q \in \nwithzero$ and assume the inductive hypothesis indexed by $q$. 
Let $S \in \mc{S}$ be such that  $\char"0023 \, \mathsf{D}_-^{j,k}(i) = q+1$. Again set $g =  \ground^{j,k}(S,i)$.
Let $\hat{S} \in \tis$ be given by
 $$
 \hat{S}(i,\ell) \, = \, \begin{cases}
  \, S_-(i,g+1) = b_i  &  \text{if $\ell = g+1$} \\
 \, S(i,\ell)  &  \text{if $\ell \in \nwithoutzero$, $\ell \not= g+1$} \, .
\end{cases}
$$
for $i \in \Z$. The set  $\mathsf{D}_-^{j,k}(\hat{S},i)$ is formed from
 $\mathsf{D}_-^{j,k}(S,i)$ by the removal of the elements of minimum height---which is height $g$.
Hence,  $\char"0023 \, \mathsf{D}_-^{j,k}(\hat{S},i) < \char"0023 \, \mathsf{D}_-^{j,k}(S,i)$; the hypotheses  $\indhyp(\hat{S},\ell)$, $\ell \in \nwithoutzero$, are thus available.
By Lemma~\ref{l.onestep}(1) with $S_1 = \hat{S}$, $S_2 = S$ and $\ell = g+1$, we find that $\indhyp(S,\ell)$ holds for $\ell \geq g+1$.

 Now consider  $u \in \llbracket -j,k \rrbracket$ such that $(u,g) \in \mathsf{D}_-^{j,k}(S,i)$. Lemma~\ref{l.onestep}(3) (and Lemma~\ref{l.onestep}(2) for other $u \in \llbracket -j,k \rrbracket$) implies $\indhyp(S,g)$.

To complete the inductive step, it remains to verify  $\indhyp(S,\ell)$ for $\ell \in \llbracket 0,g-1 \rrbracket$.
We do so iteratively in decreasing~$\ell$. It is Lemma~\ref{l.onestep}(2) that demonstrates the generic step in this iteration. This completes the proof of
 Lemma~\ref{l.backwardformal}. \qed

 {\bf Proof of Lemma~\ref{l.baseconseq}.} Apply  Lemma~\ref{l.backwardformal}  for $S = \sdev[v]$ for $v \in \nwithoutzero$,
We learn that Definition~\ref{d.merit} holds with $S_1 = S_2 = \sdev[v]$. Thus, the positive quantity in Definition~\ref{d.merit}(2) is determined by $\sdev[v](u,g+1)$. When $v$ satisfies the bound in Lemma~\ref{l.baseconseq}, we have that $\sdev[v](u,g+1)$ equals $\sdev(u,g+1)$. As a result, Definition~\ref{d.merit} holds with $S_1 = \sdev[v]$ and $S_2 = \sdev$. This is what Lemma~\ref{l.baseconseq} asserts. \qed
 
 Lemma~\ref{l.baseconseq} is a stepping stone to a counterpart that describes the penalty incurred by use of the perhaps infinitely deviating strategy $\sdev \in \mc{S}$.
 The counterpart, Lemma~\ref{l.baseconseqtwo}, depends on a variation of Definition~\ref{d.merit}.
\begin{definition}\label{d.just}
Let $S \in \mc{S}$. 
An element $(q,\ell) \in F_i$ is said to be {\em $(S,S_+)$-accessible
from $(i,0)$} 
if $\pgameplay{S}{S_+}{i}(X_\ell = q) > 0$.
Let $\access(S,i)$ denote the set of elements of $F_i$ that are $(S,S_+)$-accessible from $(i,0)$.

Alter Definition~\ref{d.merit} by taking $S_1$ and $S_2$ equal to $S$; the first part to include the condition that the point $(u,\ell)$ belongs to $\access(S,i)$; 
and the second to include the condition that $(u,g) \in \access(S,i)$.
Thus, no requirement is imposed by a given part when 
$(u,\ell)$ or $(u,g)$ is not $(S,S_+)$-accessible from $(i,0)$.
When the altered set of conditions is satisfied, we say that {\em $S$ justly receives a weak $(i,j,k)$-penalty}.
\end{definition}
 \begin{lemma}\label{l.baseconseqtwo}
Let $j,k \in \nwithzero$ satisfy $i \in \llbracket -j,k \rrbracket$. The strategy $\sdev$ justly receives a weak $(i,j,k)$-penalty.
\end{lemma}
 To prove this result, we intend to make use of Proposition~\ref{p.jksup}(1)'s hypothesis  that $\pgameplay{\sdev}{S_+}{i}(E)=1$. Since escape is certain under $\pgameplay{\sdev}{S_+}{i}$, gameplay will exit $\llbracket -j,k \rrbracket$
 in finite time, so that Mina's choice between $S$ and $S[v]$, for high $v$, will typically leave gameplay unaffected. Thus we aim to reduce the proof of the new result to quoting Lemma~\ref{l.baseconseq}.
 To do this, it is useful to state a consequence of   $\pgameplay{\sdev}{S_+}{i}(E)=1$.
\begin{lemma}\label{l.escapepropagate}
Let $(u,\ell) \in \access(\sdev,i)$. Then $\pgameplay{\sdev}{S_+}{u,\ell}(E) = 1$.
\end{lemma}
{\bf Proof.} We have that
$$
 1 \, = \, \pgameplay{\sdev}{S_+}{i}(E) \, = \, \sum_{u \in \Z} \pgameplay{\sdev}{S_+}{i}(X_\ell = u) \cdot \pgameplay{\sdev}{S_+}{u,\ell}(E) \, .
$$
Since  $\pgameplay{\sdev}{S_+}{i}(X_\ell = u) > 0$ if and only if $(u,\ell) \in  \access(\sdev,i)$, we see that $\pgameplay{\sdev}{S_+}{u,\ell}(E)$ equals one when this condition is satisfied. \qed

{\bf Proof of Lemma~\ref{l.baseconseqtwo}.} Let $h(v)$ be the vertical coordinate of the $v$\textsuperscript{th} element of $\mathsf{D}_-^{j,k}(\sdev,i)$. 
Let $(u,\ell) \in \llbracket -j,k \rrbracket \times \nwithzero$.
We will argue that
\begin{equation}\label{e.tau}
\egameplay{\sdev}{S_+}{u,\ell}[P^{j,k}_-] - \egameplay{\sdev[v]}{S_+}{u,\ell}[P^{j,k}_-] \, \leq \, \pgameplay{\sdev}{S_+}{u,\ell} \big(\tau^{j,k} \geq h(v)\big) \cdot (n_{-j-1} - n_{k+1}) \, .
\end{equation}
Given $\tau^{j,k} \leq h(v)$, counter evolution up to time $\tau^{j,k}$ is equal in law under $\pgameplay{\sdev}{S_+}{u,\ell}$ and $\pgameplay{\sdev[v]}{S_+}{u,\ell}$, because there is no discrepancy in the application of the two Maxine strategies during $\llbracket 0, \tau^{j,k} \rrbracket$. This has the upshot that the respective expectations on the left-hand side of~(\ref{e.tau}) take the form 
 of the product $\pgameplay{\sdev}{S_+}{u,\ell} \big(\tau^{j,k} \geq h(v)\big)$ with a further factor, which we label $\alpha$ and $\alpha'$. Now $\alpha$ is a certain average, over locations $v$ in $\llbracket j,k \rrbracket$, of $\pgameplay{\sdev}{S_+}{v}[P^{j,k}_-]$; while $\alpha'$ is another such average of the quantities $\pgameplay{S_-}{S_+}{v}[P^{j,k}_-]$. 
In regard to  $\pgameplay{\sdev}{S_+}{v}[P^{j,k}_-]$, consider (\ref{e.delayedpayoff}) with $\pm = -1$ and $(P,S_-) \to (P^{j,k},\sdev)$; note that  running costs here are non-negative, and that terminal receipt is at most $n_{-j-1}$   by Theorem~\ref{t.positiveabmn}(1): thus these quantities, and so also $\alpha$ itself, are at most $n_{-j-1}$. Turning to~$\alpha'$, note that $\pgameplay{S_-}{S_+}{v}[P^{j,k}_-] = n_v \geq n_{k+1}$, the equality due to Lemma~\ref{l.minipayoff}(2) and the bound to Theorem~\ref{t.positiveabmn}(1). Thus, $\alpha' \geq n_{k+1}$. We have obtained~(\ref{e.tau}).

 Under
$\pgameplay{\sdev[v]}{S_+}{u,\ell}$ given $\tau^{j,k} \geq h(v)$, Mina does not deviate after time $\tau^{j,k}$. 
 By Lemma~\ref{l.minipayoff}(2) and Theorem~\ref{t.positiveabmn}(1),
 the conditional mean of $P^{j,k}_-$
under $\pgameplay{\sdev[v]}{S_+}{u,\ell}$ given that $\tau^{j,k} > h(v)$
 is thus seen to be at least $n_{k+1}$. 
Now  consider (\ref{e.delayedpayoff}) with $\pm = -1$ and $(P,S_-) \to (P^{j,k},\sdev)$; note that  running costs here are non-negative, and that terminal receipt is at most $n_{-j-1}$   by Theorem~\ref{t.positiveabmn}(1). We see then that
 the conditional mean of $P^{j,k}_-$
under $\pgameplay{\sdev}{S_+}{u,\ell}$ given that $\tau^{j,k} > h(v)$ is at most $n_{-j-1}$.


Lemma~\ref{l.baseconseqtwo} will follow from  Lemma~\ref{l.baseconseq} provided
that we show that  the right-hand side of~(\ref{e.tau}) vanishes in high~$v$ whenever $(u,\ell) \in \access(\sdev,i)$.
By Lemma~\ref{l.escapepropagate}, and the hypothesis of Proposition~\ref{p.jksup}(1), we know that  $\pgameplay{\sdev}{S_+}{u,\ell}(E) = 1$.
 Thus, $\tau^{j,k}$ is finite, $\pgameplay{\sdev}{S_+}{u,\ell}$-almost surely.  
The right-hand side of~(\ref{e.tau}) thus indeed tends to zero in the limit of high $v$. Lemma~\ref{l.baseconseq} implies  Lemma~\ref{l.baseconseqtwo}, as we sought to show. \qed

 We are ready for the following proof. 
 
 {\bf Proof of Proposition~\ref{p.jksup}(1).}
 For $j,k \in \nwithoutzero$ such that $i \in \llbracket -j,k \rrbracket$, let $g$ denote the minimum vertical coordinate among elements of $\mathsf{D}_-^{j,k}(\sdev,i)$. 
 Any element $(u,g) \in F_i$ belongs to $\access(\sdev,i)$
 because, under $(\sdev,S_+)$, gameplay is governed before the $g$\textsuperscript{th} turn by the positive-element pair $(S_-,S_+)$. 
  Lemma~\ref{l.baseconseqtwo}
 thus implies that, when $(u,g) \in F_i$, 
 \begin{equation}\label{e.starone}
\E^{u,g}_{\sdev,S_+} \big[ P^{j,k}_- \big]  \leq n_u 
 \end{equation}
 and 
 \begin{equation}\label{e.startwo}
 \E^{u,g}_{\sdev,S_+} \big[ P^{j,k}_- \big] < n_u \, \, \, \, \textrm{if $(u,g) \in \mathsf{D}_-^{j,k}(\sdev,i)$} \, .
 \end{equation}
 Now note that
 \begin{eqnarray*}
 \egameplay{\sdev}{S_+}{i} [P_-^{j,k}] & = & - \,
 \egameplay{\sdev}{S_+}{i} \sum_{t=1}^{g-1} C_-(t)  {\bf 1}_{t < \tau^{j,k}} \\
   & &  \quad  \,\, + \,\, 
 \egameplay{\sdev}{S_+}{i} \big[T_-^{j,k} \cdot {\bf 1}_{\tau^{j,k} \leq g}\big]   \,\, + \,\,  \sum_{\substack{u \in \llbracket -j,k \rrbracket : \\ (u,g) \in F_i}}  \pgameplay{\sdev}{S_+}{i} (X^{j,k}_g = u) \cdot 
 \egameplay{\sdev}{S_+}{u,g} [P_-^{j,k}] \, .
 \end{eqnarray*}
 By way of comparison,
 \begin{eqnarray*}
 \egameplay{S_-}{S_+}{i} [P_-^{j,k}] & = & - \,
 \egameplay{S_-}{S_+}{i} \sum_{t=1}^{g-1} C_-(t)  {\bf 1}_{t < \tau^{j,k}} \\
 & & \quad   \, \, + \, \, 
 \egameplay{\sdev}{S_+}{i} \big[T_-^{j,k} \cdot {\bf 1}_{\tau^{j,k} \leq g}\big] 
  \, \, + \, \,  \sum_{\substack{u \in \llbracket -j,k \rrbracket : \\ (u,g) \in F_i}}  
 \pgameplay{S_-}{S_+}{i} (X^{j,k}_g = u) \cdot 
 \egameplay{S_-}{S_+}{u,g} [P_-^{j,k}] \, ,
 \end{eqnarray*}
 with
 $$
  \egameplay{S_-}{S_+}{u,g} [P_-^{j,k}] = n_u \, \, \, \textrm{for $u \in \llbracket -j,k \rrbracket$}
 $$
 by Lemma~\ref{l.minipayoff}(3).
 Consider a pair $(j,k)$ over which the supremum in Proposition~\ref{p.jksup}(1) is taken. Since $\mathsf{D}_-^{j,k}(\sdev,i)$ is non-empty, we may find
  $q \in \llbracket -j,k \rrbracket$
such that $(q,g) \in \mathsf{D}_-^{j,k}(\sdev,i)$.
 Since $X: \llbracket 0,g \rrbracket \to \Z$ coincides under $\PP_{S_-,S_+}^i$
 and $\PP_{\sdev,S_+}^i$, 
 the terminal payment and running costs terms in the two preceding right-hand sides cancel as we take the difference, so that we obtain
 $$
 \egameplay{S_-}{S_+}{i} [P_-^{j,k}] -
 \egameplay{\sdev}{S_+}{i} [P_-^{j,k}]  \, \, = \, \, \sum_{\substack{u \in \llbracket -j,k \rrbracket : \\ (u,g) \in F_i}} 
 \pgameplay{S_-}{S_+}{i}  (X^{j,k}_g = u) \cdot \Big(n_u -  \egameplay{\sdev}{S_+}{u,g} [P_-^{j,k}] \Big) \, ,
 $$
 where the term in parentheses on the right-hand side is strictly positive if $u = q$ (by~(\ref{e.startwo})), and is non-negative if $u \in \llbracket -j,k \rrbracket$, $u \not= q$ (by~(\ref{e.starone})).
 The bracketed right-hand term is thus positive when we write
 $$
 \egameplay{S_-}{S_+}{i} [P_-^{j,k}] -
 \egameplay{\sdev}{S_+}{i} [P_-^{j,k}]  \, \, \geq \, \, 
 \pgameplay{S_-}{S_+}{i}  (X^{j,k}_g = q) \cdot \Big( \, n_q -  \egameplay{\sdev}{S_+}{q,g} \big[P_-^{j,k} \big] \, \Big)  \, .
 $$
 We claim that $\pgameplay{S_-}{S_+}{i}  (X^{j,k}_g = q) > 0$.
 Indeed, it is enough to find any access route for $X$ from $(i,0)$ to $(q,g)$ that never leaves $\llbracket -j,k \rrbracket \times \llbracket 0,g \rrbracket$, because  the strategies in the pair $(S_-,S_+) = (b,a)$ have positive coefficients; that such a route exists is due to $(q,g) \in F_i$, $i,q \in \llbracket -j,k \rrbracket$ and $k >-j$.
 For example, 
 $\pgameplay{S_-}{S_+}{i}  (X^{j,k}_g = q)$ is at least $\eta^g$, where $\eta = \min \big\{ (a_i \wedge b_i)/(a_i + b_i): i \in \llbracket -j,k \rrbracket \big\}$. That is,
 $$
 \egameplay{S_-}{S_+}{i} [P_-^{j,k}] -
 \egameplay{\sdev}{S_+}{i} [P_-^{j,k}] \, \geq \, \eta^g  \cdot \Big( \, n_q -  \egameplay{\sdev}{S_+}{q,g} \big[P_-^{j,k}\big] \, \Big) \, .
 $$
 Since the positive right-hand side is independent of the choice of the pair $j,k \in \nwithoutzero$ 
 over which the supremum is taken in  Proposition~\ref{p.jksup}(1), we have obtained this result. \qed
 
{\bf Proof of  Proposition~\ref{p.jksup}(2).} The essentially identical argument is omitted. \qed
 
   \chapter{Explicit \abmnmacspace solutions and their consequences}\label{c.battlefield}
  
Here we explicitly solve the \abmnmacspace system, proving Theorem~\ref{t.defaultexplicit}, and its softer cousin Proposition~\ref{p.default}.
Then we analyse the asymptotic decay  in high index values of \abmnmacspace solutions, proving Theorem~\ref{t.ajbj}.
Two consequences of this decay---finiteness of boundary data in Theorem~\ref{t.positiveabmn}(3), and 
  the almost sure eventual unanimity of gameplay in Theorem~\ref{t.unanimity}---are derived.

  \section{Explicit \abmnmacspace solutions}

Fundamental to deriving Theorem~\ref{t.defaultexplicit} is an alternative representation of the \abmnmacspace system that we offer first, in Proposition~\ref{p.abmnsolvesmn}.
    The real-valued variables $\big\{ m_i,n_i: i \in \Z \big\}$ satisfy the \mnmacspace system on $\Z$ if
\begin{align*}
   (m_i - m_{i-1}) (m_{i+1} - m_{i-1} + n_{i-1} - n_{i+1})^2 & \, = \, (m_{i+1} - m_{i-1})^3 && \qquad \textrm{MN}(1) \\
   (n_i - n_{i+1}) (m_{i+1} - m_{i-1} + n_{i-1} - n_{i+1})^2 & \, = \,  (n_{i-1} - n_{i+1})^3 && \qquad \textrm{MN}(2)   \, ,
  \end{align*}
    for $i \in \Z$. As for \abmnmac$(1,2,3,4)$ from Definition~\ref{d.abmn}, we refer to the above equations as  $\textrm{MN}(1)$ and  $\textrm{MN}(2)$
    rather than by the usual convention of numbered equations. 
  \begin{proposition}\label{p.abmnsolvesmn}
  A positive solution of the \abmnmacspace system on $\Z$ solves the \mnmacspace system on $\Z$.
  \end{proposition}
  {\bf Proof.} For $i \in \Z$, set $M_i = m_{i+1} - m_{i-1}$ and $N_i = n_{i-1} - n_{i+1}$. We claim that
  \begin{equation}\label{e.abclaim}
   a_i = \frac{M_i^2 N_i}{(M_i+N_i)^2} \, \, \, , \, \, \,  b_i = \frac{M_i N_i^2}{(M_i+N_i)^2} \, \, \, \, \textrm{and} \, \, \, \, \frac{a_i}{a_i+b_i} =  \frac{M_i}{M_i+N_i} \, .
  \end{equation}
  These follow from \abmnmac$(3,4)$. Expressing \abmnmac$(1)$  in the form~(\ref{e.firstrearranged}),  we find from~(\ref{e.abclaim}) that 
  $$
  m_i \, = \, m_{i-1} +  \frac{M_i^2}{M_i+N_i} -  \frac{M_i^2 N_i}{(M_i+N_i)^2} \, ,
  $$
   whence \mnmac$(1)$ holds. Equation \mnmac$(2)$ is obtained similarly, from \abmnmac$(2)$. \qed
  
  The \mnmacspace framework will permit us to relate \abmnmacspace solutions to the basic functions $\macsee,\macdee,s:(0,\infty) \to (0,\infty)$ from Definition~\ref{d.acs}. 
  Such relations appear in Proposition~\ref{p.alphagammaess} including the key $s(\phi_0)=\phi_1$, which will enable the proof of the `$s$ as shift' Proposition~\ref{p.shift} shortly thereafter. 
    \begin{definition}\label{d.alphagamma}\hfff{alphagamma}
  Let $\macgsee,\macgdee:(0,\infty) \to (0,\infty)$ be given by $\macgsee(x) = \macsee(x)^{-1}$ and $\macgdee(x) = \macdee(x)^{-1}$.
  Set $\macbeta:(0,\infty) \to (0,\infty)$, $\macbeta(x) = \tfrac{\macomega - 1}{4}$, where recall that $\macomega = \sqrt{8x+1}$ for $x \in (0,\infty)$.
  \end{definition}
  
\begin{lemma}\label{l.acsfacts}
\leavevmode
\begin{enumerate}
\item The functions $\macsee,\macdee,s:(0,\infty) \to (0,\infty)$ are increasing.\footnote{Let $* \in \{ c,d,s \}$. By `Lemma~\ref{l.acsfacts}(1:$*$)' will be meant `$*$ is increasing'.}
\item
We have that $s(x) = x^2/2 + O(x^3)$ as $x \searrow 0$.
\item
For $x \in (0,\infty)$, $s(x) = \tfrac{\macbeta(x)^2}{\macbeta(x)+2}$. 
\item  
For $x \in (0,\infty)$, $\macbeta(x) \leq x$.
\item 
For $x \in (0,\infty)$, $s(x) < x$.
\end{enumerate}
\end{lemma}
{\bf Proof: (1).} 
The expressions for $\macsee(x)$, $\macdee(x)$ and $s(x)$ in Definition~\ref{d.acs} are readily seen to be increasing in the variable $\macomega \in (1,\infty)$; since $\macomega = \sqrt{8x +1}$, they are also increasing in $x \in (0,\infty)$.\\
{\bf (2).} We have that $\macomega = \sqrt{8x +1} = 1 + 4x + O(x^2)$, whence 
$$
s(x) = \tfrac{(\omega-1)^2}{4(\omega +7)}= \tfrac{16 x^2 + O(x^3)}{4(8 + O(x))} = x^2/2 + O(x^3)\, . 
$$
{\bf (3).}  This is due to $s(x) = \tfrac{(\macomega-1)^2}{4(\macomega+7)}$ and $\macbeta(x) = (\macomega-1)/4$. \\
{\bf (4).} 
Since $\macomega(x) = \sqrt{8x +1} \leq 4x+1$, $\macbeta(x) \leq x$.  \\
{\bf (5).}
Lemma~\ref{l.acsfacts}(3), $\macbeta > 0$ and Lemma~\ref{l.acsfacts}(4) imply that 
 $$
 s(x) = \tfrac{\macbeta(x)^2}{\macbeta(x) + 2} < \macbeta(x) \leq x
 $$
as desired.
\qed

  Recall Definition~\ref{d.deltai} and suppose given $\big\{ (a_i,b_i,m_i,n_i): i \in \Z \big\}$, a positive \abmnmacspace solution on~$\Z$.
  \begin{proposition}\label{p.alphagammaess}\hfff{gammaess}
  For $i \in \Z$, we have that\footnote{Let $* \in \{\gamma,\delta,s\}$. By `Proposition~\ref{p.alphagammaess}($*$)', we will mean the statement made concerning the labelled quantity.}
  $$
  \macgsee(\macphi_i) = \frac{m_i - m_{i-1}}{m_{i+1} - m_{i-1}} \, \, , \, \,  \macgdee(\macphi_i) = \frac{n_{i-1} - n_i}{n_{i-1} - n_{i+1}} \, \, \, \, \textrm{and} \, \, \, \, s(\macphi_i) = \macphi_{i+1} \, .
$$
  \end{proposition}
  Notation to be used only in the proof of this proposition\footnote{In particular, the temporary usage of $s_i$ introduced in Definition~\ref{d.subscripti} is an abuse, because the denoted quantity is not the function $s_i$; nor is it the value $s_i(x)$ for $x = \phi_0$. Indeed, $s_i(x)$ equals $\phi_i$, while $s_i$ with the temporary usage equals $\phi_{i+1}$.} makes the task to show that $*(\macphi_i)$ equals $*_i$ for $* \in \{\macgsee,\macgdee,s\}$.
  
\begin{figure}[htbp]
\centering
\includegraphics[width=0.75\textwidth]{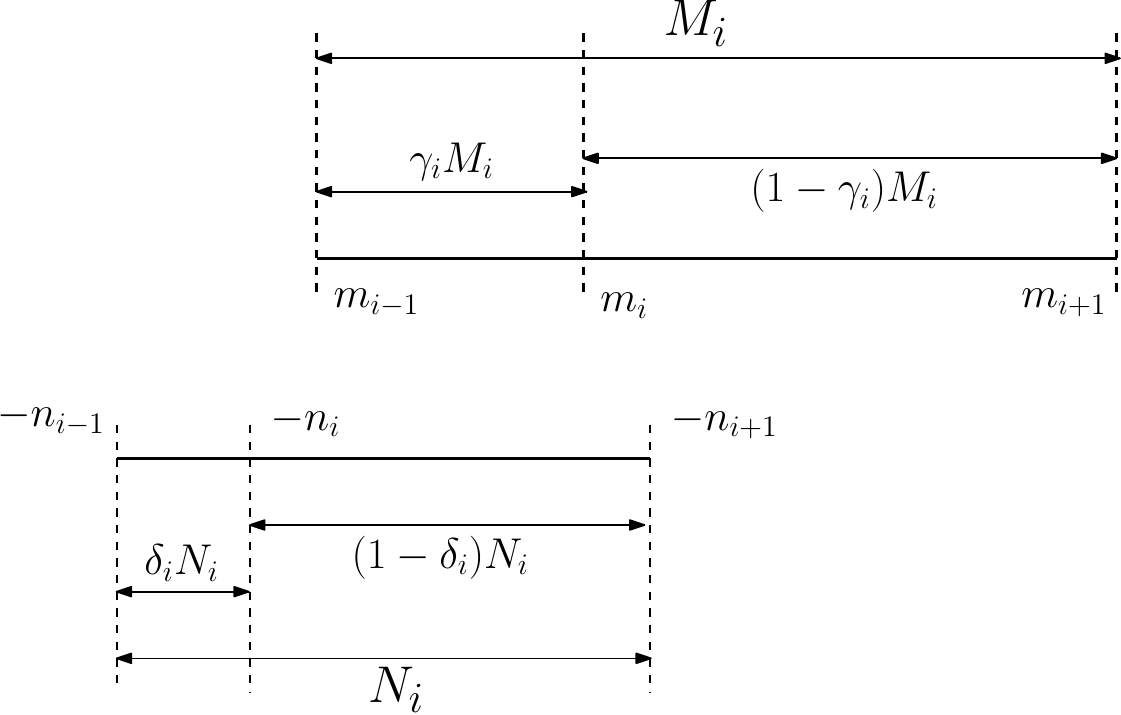}
\caption{Suppose given $\big\{ (a_i,b_i,m_i,n_i): i \in \Z \big\}$, a positive \abmnmacspace solution on~$\Z$. The quantities $M_i$ and $\gamma_i$ are functionals of the triple $(m_{i-1},m_i,m_{i+1})$; likewise for $N_i$ and $\delta_i$ and the triple $(n_{i-1},n_i,n_{i+1})$. The coordinates $-n_{i-1}$, $-n_i$ and $-n_{i+1}$ are depicted in order to maintain the orientation of increasing index on the right. The quantity $\beta_i = N_i/M_i$ specifies the relative lengths of the overall intervals depicted.}\label{f.gammadelta}
\end{figure}

  Figure~\ref{f.gammadelta} illustrates the next definition.
  \begin{definition}\label{d.subscripti}
  For $i \in \Z$, set $\macgsee_i = \tfrac{m_i - m_{i-1}}{m_{i+1} - m_{i-1}}$,  $\macgdee_i = \frac{n_{i-1} - n_i}{n_{i-1} - n_{i+1}}$ and  $s_i = \macphi_{i+1}$. We also set
  $\macbeta_i = \frac{n_{i-1}-n_{i+1}}{m_{i+1} - m_{i-1}}$, and write $\macomega_i =  \macomega(\macphi_i) = \sqrt{8\macphi_i +1}$.  
  \end{definition}
 
\begin{lemma}\label{l.fourfacts}
We have that 
$$
 (1 + \macbeta_i)^2  \macgsee_i = 1 \, \, , \, \, 1 - \macgdee_i = \tfrac{\macbeta_i^2}{(1+\macbeta_i)^2} \, \, ,  \, \, \macphi_i = \macgdee_i\macbeta_i/\macgsee_i \, \, , \, \, \macphi_{i+1} = \tfrac{\macbeta_i(1-\macgdee_i)}{1-\macgsee_i} \, . 
$$
\end{lemma}
{\bf Proof.} Equation~\mnmac(1) implies that  $(1 + \macbeta_i)^2  \macgsee_i = 1$. Equation \mnmac(2) implies  $1 - \macgdee_i = \tfrac{\macbeta_i^2}{(1+\macbeta_i)^2}$. That $\macphi_i = \macgdee_i\macbeta_i/\macgsee_i$ follows by the definitions of the concerned quantities.  
Noting that 
$$
1 - \macgdee_i = \tfrac{n_i - n_{i+1}}{n_{i-1} - n_{i+1}} \, \, \, \, \textrm{and}  \, \, \, \, 1 - \macgsee_i = \tfrac{m_{i+1} - m_i}{m_{i+1} - m_{i-1}} \, ,
$$
we find from the definitions of $\macbeta_i$ and $\macphi_{i+1}$
that $\macphi_{i+1} = \tfrac{\macbeta_i(1-\macgdee_i)}{1-\macgsee_i}$ holds. \qed

\begin{lemma}\label{l.omegai}
For $i \in \Z$,
$$
\macgsee_i^{-1}   =   \tfrac{1}{16} (\macomega_i + 3)^2 \, \, , \, \,
 \macgdee_i^{-1}  =  \frac{(\macomega_i + 3)^2}{8(\macomega + 1)} \, \, , \, \, s_i   =  \frac{(\macomega_i - 1)^2}{4(\macomega_i + 7)}  \, \, \, \textrm{and} \, \, \,
\macbeta_i  =  \tfrac{1}{4} (\macomega_i - 1)
 \, .
 $$
 \end{lemma}
{\bf Proof.} Omitting $i$ subscripts, consider the four equations stated in Lemma~\ref{l.fourfacts} when we take $\macphi \in (0,\infty)$ given. The first and third equations imply that $\macgdee\macbeta(1+\macbeta)^2 = \macphi$. Using the second equation, we find that $(2\macbeta+1)\macbeta = \macphi$; since $\macbeta$ is positive, we confirm that $\macbeta = (\macomega -1)/4$. From the first equation, we then obtain $\macgsee = 16(\macomega +3)^{-2}$. The third equation $\macgdee = \macphi\macgsee/\macbeta$ then yields $\macgdee = \tfrac{16\macphi}{(\macomega +3)^2} \cdot \tfrac{4}{\macomega -1}$ which equals $\tfrac{8(\macomega +1)}{(\macomega +3)^2}$ in view of $\macomega^2 -1 = 8\macphi$. Finally, $s = \macphi_{i+1}$ by definition, so that the fourth equation implies that
$s = \tfrac{\macomega -1}{4} \cdot \tfrac{(\macomega+3)^2 - 8\macomega - 8}{(\macomega+3)^2 - 16}$ whose right-hand side is seen to equal $\tfrac{(\macomega -1)^2}{4(\macomega+7)}$ after cancellation of $\macomega - 1 > 0$
from numerator and denominator.  \qed

{\bf Remark.} Expressed in terms of the variable $\beta = \beta_0 \in (0,\infty)$, the functions $\phi_0$ and $\phi_1$  are increasing bijections $(0,\infty) \to (0,\infty)$ 
that may now be seen to take the forms
$\phi_0 = \beta(2\beta+1)$ and $\phi_1 = \beta^2/(\beta +2)$.
Certainly then, $s:(0.\infty) \to (0,\infty)$, viewed as the function that sends $\phi_0 \in (0,\infty)$ to $\phi_1$, is well defined. Happenstance that $\phi_0$ and $\phi_1$ are both coefficients of quadratic equations in $\beta$. This special feature permits the explicit form for~$s$ in Definition~\ref{d.acs} as a rational function of $\omega = \sqrt{8x+1}$. The Trail of Lost Pennies is akin to an integrable model in this sense. The explicit form means that, in evaluating the Mina margin map~$\mc{M}$, the computation of $s$-iterates in the expression in Theorem~\ref{t.relativereward}(2) does not require numerical approximation. Of course, the expression is a ratio of albeit rapidly convergent infinite sums, so the computation is necessarily approximate in any case. 
 
\begin{lemma}\label{l.omega.asymptotic}
\leavevmode
\begin{enumerate}
\item
We have that
$\macgsee_i^{-1}  -1  = 2\macphi_i + O(\macphi_i^2)$.
\item
And that
$\macbeta_i = \macphi_i + O(\macphi_i^2)$.
\end{enumerate}
\end{lemma}
{\bf Proof: (1).} From Lemma~\ref{l.omegai}, note that $\macgsee_i^{-1} = \tfrac{1}{16} (\omega_i + 3)^2 = \big( 1 + \phi_i + O(\phi_i^2)\big)^2$. \\
{\bf (2).} By the same result, $\macbeta_i = \tfrac{1}{4}(\omega_i - 1) = \macphi_i + O(\macphi_i^2)$.
\qed

{\bf Proof of Proposition~\ref{p.alphagammaess}.}
By Lemma~\ref{l.omegai} and Definitions~\ref{d.acs},~\ref{d.alphagamma} and~\ref{d.subscripti}, 
$$\macgsee_i = 16(\macomega_i + 3)^{-2}  =\macsee(\macphi_i)^{-1} = \macgsee(\macphi_i) \, \, ; \, \, \, \, 
 \macgdee_i = \tfrac{8(\macomega_i +1)}{(\macomega_i + 3)^2}  =\macdee(\macphi_i)^{-1} = \macgdee(\macphi_i) \, \, ;
 $$
  and $s_i = \tfrac{(\macomega_i - 1)^2}{4(\macomega_i +7)} = s(\macphi_i)$. \qed

  {\bf Proof of Proposition~\ref{p.shift}.}
  The symbolic shift map leaves invariant the boundary quadruple of any \abmnmacspace solution. Thus,
  the left-hand quadruple in~(\ref{e.shift}) is a standard solution of the \abmnmacspace system. To identify it as the right-hand quadruple, it is thus enough to show that its  
 $\centralratio$-value equals $s_k(x)$. But  since $x = \phi_0$ this amounts to   $\phi_k = s_k(x)$ which follows from iterating the relation $s(\phi_i) =\phi_{i+1}$ from Proposition~\ref{p.alphagammaess}. \qed

{\bf Proofs of Proposition~\ref{p.default} and Theorem~\ref{t.defaultexplicit}.}
For given $x \in (0,\infty)$, let $(a,b,m,n)$ be an \abmnmacspace solution with $\tfrac{n_{-1} - n_0}{m_0 - m_{-1}}  =x$.
 
Since $c_i(x) = c(s_i(x)) = c(\phi_i)$, Definition~\ref{d.alphagamma} and Proposition~\ref{p.alphagammaess}($\gamma$) imply that
\begin{equation}\label{e.ciformula}
  c_i(x) - 1 = \frac{1 - \gamma(\phi_i)}{\gamma(\phi_i)} = \frac{m_{i+1} - m_i}{m_i - m_{i-1}} \, .
\end{equation}
Adopting the notation in Definition~\ref{d.zdefault}, we find that 
\begin{equation}\label{e.mdifferenceratio}
\frac{m_{j+1} - m_j}{m_0 - m_{-1}} \, = \, \prod_{i=0}^j \big( c_i(x) - 1 \big)
\end{equation}
for any $j \in \Z$. Since a default solution has $m_0 - m_{-1} = 1$ by definition, we deduce that the formula for $\mdef_{k+1} - \mdef_k$ in Definition~\ref{d.zdefault} holds. Similarly to~(\ref{e.ciformula}), we find via  Proposition~\ref{p.alphagammaess}($\delta$) that
$$
  d_i(x) - 1 = \frac{1 - \delta(\phi_i)}{\delta(\phi_i)} = \frac{n_i - n_{i+1}}{n_{i-1} - n_i} \, ,
$$
whence 
$$
\frac{n_j - n_{j+1}}{n_{-1} - n_0} \, = \, \prod_{i=0}^j \big( d_i(x) - 1 \big)
$$
for $j \in \Z$. Since $n_{-1} - n_0 = x(m_0 - m_{-1}) = x$ for any default solution, we find that the formula for $\ndef_k - \ndef_{k+1}$ in Definition~\ref{d.zdefault} is valid. 
Proposition~\ref{p.abmnsolvesmn} implies that the sought formulas for $\adef_i$ and $\bdef_i$ for $i \in \Z$ hold.
The exhibited solution exists and is unique. This completes the proof of Theorem~\ref{t.defaultexplicit}. 
The noted existence and uniqueness also prove Proposition~\ref{p.default}. 
\qed

  \section{Asymptotic decay of solutions}
  Here we prove Theorem~\ref{t.ajbj}.

\begin{lemma}\label{l.deltadecay}
\leavevmode
\begin{enumerate}
\item
 For $\macphi_i \in (0,1)$, we have that
$$
 \macphi_i^2/2 - O (\macphi_i^3) \leq \macphi_{i+1} \leq \macphi_i^2/2 \, ,
$$
where the positive constant implied by the $O$-notation is bounded above in terms of $\const \in (0,1)$, where $\macphi_i \in (0,1-\const)$. 
\item For any $i \in \Z$, $\macphi_{i+1} < \macphi_i$.
\end{enumerate}
\end{lemma}
{\bf Proof: (1).} From Proposition~\ref{p.alphagammaess}(s) 
and Lemma~\ref{l.acsfacts}(3,4), we see  that 
$$
\macphi_{i+1} \leq \macbeta(\macphi_i)^2/2 \leq \macphi_i^2/2 \, .
$$
By Proposition~\ref{p.alphagammaess}(s) and Lemma~\ref{l.acsfacts}(2), 
$\macphi_{i+1} = s(\macphi_i) = \macphi_i^2/2 + O(\macphi_i^3)$. 

{\bf (2).} By Lemma~\ref{l.acsfacts}(5), and 
 Proposition~\ref{p.alphagammaess}(s), 
$\macphi_{i+1} = s(\macphi_i) < \macphi_i$. \qed

We are about to prove Theorem~\ref{t.ajbj}.
Since this result uses the notion of the battlefield index specified in Definition~\ref{d.battlefield}, we now offer a proof that this index is well-defined.
\begin{lemma}\label{l.battlefield}
Let $\big\{ (a_i,b_i,m_i,n_i): i \in \Z \big\}$ be a positive \abmnmacspace solution on $\Z$.
There is a unique value of $k \in \Z$ for which $\macphi_k \in (1/3,3]$.
\end{lemma}
{\bf Proof.}
By Lemma~\ref{l.deltadecay}(2), the sequence $\big\{ \macphi_i: i \in \Z \big\}$ is decreasing. Taking $s(0) = 0$, the value $\lim_{i \to \infty}\macphi_i$
is a fixed point of $s:[0,\infty) \to [0,\infty)$ because $s$ is continuous and $s(\macphi_i) = \macphi_{i+1}$ (the latter by Proposition~\ref{p.alphagammaess}(s)).
But $s(x) < x$ for $x > 0$ by~Lemma~\ref{l.acsfacts}(5). 
Thus, $\macphi_i \searrow 0 $ as $i \to \infty$. The opposite limiting value $\lim_{i \to \infty} \macphi_{-i}$ would also be a fixed point for  $s:[0,\infty) \to [0,\infty)$ 
were it to be finite; we see then that $\lim_{i \to \infty} \macphi_{-i}$ is infinite. 
We may thus set $k \in \Z$ so that $k = \inf \big\{ i \in \Z: \macphi_i \leq 3 \big\}$ and be assured that $k$ is well-defined.
Now, $\macphi_j > 3$ for $j \leq k-1$, while $\macphi_k$, being $s(\macphi_{k-1})$, exceeds $s(3) = 1/3$
by Lemma~\ref{l.acsfacts}(1:$s$).
On the other hand, if $j \geq k+1$, then $\macphi_j \leq \macphi_{k+1} = s(\macphi_k) \leq s(3) = 1/3$. Thus, $k \in \Z$ is the unique index whose $\macphi$-value exceeds one-third and is at most three. \qed 

{\bf Proof of Theorem~\ref{t.ajbj}(1).} For $i \in \nwithzero$, set $\e_i = \macphi_{k+i}/2$ and $g_i = - \log \e_i$. 
By $s(3) = 1/3$ and Lemma~\ref{l.acsfacts}(1:$s$), we have that $s(x) \leq 1/3$ for $x \in (0,3]$. Definition~\ref{d.battlefield} and  $s(\macphi_i) =\macphi_{i+1}$ (from Proposition~\ref{p.alphagammaess}(s)) thus imply that
$\macphi_{k+j} \leq 1/3$ for $j \geq 1$. We may then apply
Lemma~\ref{l.deltadecay}(1) to find that $\e_i^2 \big( 1 - O(\e_i) \big) \leq \e_{i+1} \leq \e_i^2$, 
where the positive constant implicit in the $O$-notation may be chosen independently of the ABMN solution $\big\{ (a_j,b_j,m_j,n_j):j \in \Z \big\}$
and the value of the index $i \geq 1$. (We say that a positive constant is universal, or is bounded universally, if it may be so chosen.) We learn that
\begin{equation}\label{e.twogi}
 2 g_i \, \leq \, g_{i+1} \, \leq \, 2 g_i + O\big(e^{-g_i} \big) \, ,
\end{equation}
where the implicit positive constant is again universal.
Thus, $g_i > \log 6$ for $i \geq 1$, and we may write $g_i = 2^{\macell_i}$ for a real-valued sequence $\{ \macell_i: i \in \nwithoutzero \}$ whose terms are bounded below by $\tfrac{\log \log 6}{\log 2} > 0$. From~(\ref{e.twogi}), we find that
$$
 0 \leq g_{i+1} - 2g_i = \big( 2^{\macell_{i+1} - \macell_i - 1} - 1 \big) 2^{\macell_i +1}   \, = \, O \big( \exp \{ - 2^{\macell_i} \} \big) \, ;
$$
using $\macell_i  > 0$, we readily obtain
$$
  0 \, \leq \, \macell_{i+1} - \macell_i - 1  \, = \, O \big( \exp \{ - 2^{\macell_i} \} \big) \, .
$$
Since $\macell_1 > 0$ and $\macell_{i+1} \geq \macell_i +1$, we have that $\macell_i > i -1$ for $i \geq 1$. Thus,
$$
   0 \, \leq \, \macell_{i+1} - \macell_i - 1  \, = \, O \big( \exp \{ - 2^{i-1}  \} \big) \, .
$$
We may find $B \in \R$ so that $\macell_i = B + i + O \big( \exp \{ - 2^{i-1}  \} \big)$ for $i \in \nwithoutzero$. 
The universal form of $O$ and the fact that $\macell_1$ is bounded (since $\e_1 \in (1/6,3/2]$)
implies that $B$ is bounded in a universal sense.
Set $A = 2^B$ (so that $A$ is bounded away from zero and infinity in a universal sense), and exponentiate with base two to obtain
$$
 g_i = A \cdot 2^{i + O\big(\exp ( -2^{i-1} ) \big)} 
$$
for $i \geq 1$. Since $\phi_{k+i} = 2e^{-g_i}$, we see then that, for $i \geq k+1$,
\begin{equation}\label{e.deltaiformula}
 \macphi_i \, = \, 2 \exp \Big\{ -A \cdot 2^{{i-k} +O\big(\exp ( -2^{i-k-1} ) \big)} \Big\} \, .
\end{equation}
Similarly as we derived~(\ref{e.mdifferenceratio}), we find that
$$
 m_j - m_{j-1} \, = \, (m_k - m_{k-1}) \prod_{i = k}^{j-1} \big( \macgsee_i^{-1} - 1 \big) 
$$
for $j \geq k+1$. By Lemma~\ref{l.omega.asymptotic}(1),
\begin{eqnarray*}
 m_j - m_{j-1} & = & (m_k - m_{k-1}) \prod_{i = k}^{j-1} \Big( \, 4 \exp \Big\{ -A \cdot 2^{i -k+ \kappa_i \exp ( -2^{i-k-1} )   \big)} \Big\} + O(1) e^{-A \cdot 2^{i-k + 1/2 }} \, \Big) \\ 
 & = &  (m_k - m_{k-1}) 4^{j-k} E_{k,j} \prod_{i = k}^{j-1}   \exp \Big\{ -A \cdot 2^{i -k+ \kappa_i \exp ( -2^{i-k-1} )}      \Big\} \, , 
  \end{eqnarray*}
where the values of $\kappa_i$ are bounded above in absolute value (in a universal sense), and 
  where 
 \begin{eqnarray*}
  E_{k,j} & = & \prod_{i=k}^{j-1} \Big( 1 + O(1) \exp \big\{ - A 2^{i-k} \big(2^{1/2} - 2^{\kappa_i \exp \{- 2^{i-k-1} \}} \big) \big\} \Big) \\
  & = & \prod_{i=k}^{j-1} \Big( 1 + \exp \big\{ - O(1)A \cdot 2
  ^{i-k} \big\} \Big)
 \end{eqnarray*}
satisfies $E_{k,j} = E \big( 1 + e^{-O(1)A 2^{j-k}}\big)$ with
$$
E = \prod_{i=k}^\infty \Big( 1 + \exp \big\{ - O(1)A \cdot 2
  ^{i-k} \big\} \Big) \, .
$$
The quantity $E$ is positive and bounded away from zero and infinity universally. 
Note that
\begin{eqnarray*}
  \sum_{i= k}^{j-1} 2^{i -k+ \kappa_i \exp \{ - 2^{i-k-1} \} } & = & 2^{j-k} -  1 + \sum_{i= k}^{j-1} 2^{i-k} \big( 2^{\kappa_i \exp \{- 2^{i-k-1} \} }  - 1 \big) \\
  &  = & 2^{j-k} -  1 + \rho - \sum_{i= j}^\infty 2^{i-k} \big( 2^{\kappa_i \exp \{ - 2^{i-k-1} \} }  - 1 \big) \\
  & = & 2^{j-k} -  1 + \rho  + O(1) e^{-2^{j-k}O(1)}   \, ,
\end{eqnarray*}
where $\rho = \sum_{i= k}^\infty 2^{i-k} \big( 2^{\kappa_i \exp \{ - 2^{i-k -1} \} }  - 1 \big)$.
Thus, $m_j - m_{j-1}$ equals
\begin{eqnarray*}
 &  & (m_k - m_{k-1}) 4^{j-k} \exp \big\{ - 2^{j-k}A     \big\}   E 
 \exp \big\{  A (   1 - \rho ) \big\}  \big( 1 + e^{-O(1)A 2^{j-k}}\big) \\
 & & \qquad \qquad \qquad \qquad \qquad \qquad \times \, \, \,  \Big( 1 + e^{-O(1)A 2^{j-k}}\Big)
 \exp \Big\{ A \cdot O(1) e^{-2^{j-k}O(1)}    \Big\} \\
  & = &  (m_k - m_{k-1}) 4^{j-k} \exp \big\{ - 2^{j-k}A     \big\}   E 
 \exp \big\{  A (   1 - \rho ) \big\}  \big( 1 + e^{-O(1) 2^{j-k}}\big) \, ,
\end{eqnarray*}
where we used that $A = \Theta(1)$---namely, $A$ is bounded away from zero and infinity in a universal sense---for the displayed equality. Set $F$ equal to $E 
 \exp \big\{  A (   1 - \rho ) \big\}$, and note that this positive expression is bounded away from zero and infinity universally. We find that
\begin{equation}\label{e.mjmjminusone}
 m_j - m_{j-1} =  (m_k - m_{k-1})\cdot F \cdot  2^{2(j-k)} \exp \big\{ - 2^{j-k}A     \big\}   \big( 1 + e^{-O(1) 2^{j-k}}\big) \, ,
\end{equation}
which is the inference that Theorem~\ref{t.ajbj} makes for the sequence of $m$-differences.
With $M = m_{j+1} - m_{j-1}$ and $N = n_{j-1} - n_{j+1}$, we have that $a_j = \tfrac{M^2 N}{(M+N)^2}$
and $b_j = \tfrac{M N^2}{(M+N)^2}$ from~(\ref{e.abclaim}). Using the definition of $\beta_j$ in the guise $N = \macbeta_j M$, and Lemma~\ref{l.omega.asymptotic}(2) with $i=j$, we find that
$$
a_j =  ( m_{j+1} - m_{j-1}) \tfrac{\macbeta_j}{(1+\macbeta_j)^2} \, = \, ( m_{j+1} - m_{j-1})  \big(  \macphi_j + O(\macphi_j^2) \big)
$$
and
$$
 b_j =  ( m_{j+1} - m_{j-1}) \tfrac{\macbeta_j^2}{(1+\macbeta_j)^2} \, = \, ( m_{j+1} - m_{j-1})  \big(  \macphi_j^2 + O(\macphi_j^3) \big) \, .
$$
We may use (\ref{e.mjmjminusone}) to replace 
the quantity $m_{j+1} - m_{j-1}$ in these expressions. The expressions in terms of $\macphi_i$ may be bounded by means of~(\ref{e.deltaiformula}):
\begin{eqnarray*}
\macphi_j & = & 2 \exp \big\{ - A \cdot 2^{j-k} \big( 1 + O ( \exp \{ - 2^{j-k-1} \} ) \big) \big\} \\ 
& = &  2 \exp \big\{ - A \cdot 2^{j-k} \big\} \exp \big\{ O (e^{- 2^{j-k}c}) \big\} 
 \, = \, 2 \exp \big\{ - A \cdot 2^{j-k} \big\} \big( 1 + O (e^{- 2^{j-k}c}) \big) \, .
\end{eqnarray*}
Here, the value of $c$ is positive (and universal) in the second line.
We thus obtain the expressions for $a_j$ and $b_j$ in Theorem~\ref{t.ajbj}(1).
It remains to derive the asymptotic expression for the quantity $n_j - n_{j-1}$. Here, we use $n_{j-1} - n_j = \phi_j(m_j - m_{j-1})$,~(\ref{e.mjmjminusone}) and the preceding display. 

{\bf (2).} 
According to Definition~\ref{d.battlefield}, the battlefield index $k \in \Z$ is the unique solution of $\macphi_k \in (1/3,3]$.  Consider the role-reversal transformation that replaces index $i$ by $2k-i$, and $(a,b,m,n)$ by $(b,a,n,m)$. The resulting system is also a solution of the \abmnmacspace system by a minor variation of Proposition~\ref{p.rolereversal}. Write $\hat\phi_i$ for the value of $\phi_i$ in the transformed solution. Then $\hat\phi_i = 1/\phi_{2k+1 -i}$ for $i \in \Z$.
We see then that $\hat\phi_{k+1} \in [1/3,3)$ (so that $k+1$ is the battlefield index of the transformed system except when $\phi_k =1/3$). 
Theorem~\ref{t.ajbj}(2) thus reduces to Theorem~\ref{t.ajbj}(1), because the proof of the latter operates as well as when $\phi_k = 1/3$ as when $\macphi_k \in (1/3,3]$. \qed

  \section{Consequences of asymptotic decay}\label{s.consequences}

 We may now complete the proof of Theorem~\ref{t.positiveabmn}.
 
{\bf Proof of Theorem~\ref{t.positiveabmn}(3).}
By Theorem~\ref{t.positiveabmn}(2), we know that $m_\infty$, $m_{-\infty}$, $n_\infty$ and $n_{-\infty}$ exist as elements of $\R \cup \{ \infty \} \cup \{ - \infty \}$. Since we know that $m_0$ and $n_0$
belong to $\R$, it is enough, in order to exclude the possibility that one of the four quantities is infinite, to argue that
$\lim_{i \to \infty} (m_i - m_0) < \infty$, $\lim_{i \to \infty} (m_{-i} - m_0) > -  \infty$,  $\lim_{i \to \infty} (n_i - n_0) > - \infty$ and  $\lim_{i \to \infty} (n_{-i} - n_0) <  \infty$. 
 These results follow from the asymptotic expressions for $m_j - m_{j-1}$ and $n_{j-1} - n_j$ in Theorem~\ref{t.ajbj}(1,2). \qed
 
 The almost sure occurrence of the unanimity event $U$ is a consequence of Theorem~\ref{t.ajbj}, and we prove it now.
  
 {\bf Proof of Theorem~\ref{t.unanimity}(4).}
 By Theorem~\ref{t.nashabmn}(1), this reduces to Theorem~\ref{t.unanimity}(1,2,3).
 
{\bf (1,2,3).}  We abusively write $(S_-,S_+)= (b,a)$ as usual. Theorem~\ref{t.nashabmn} and Theorem~\ref{t.ajbj}(1)
 imply that, for $i \geq k$, $\tfrac{a_i}{a_i + b_i} 
 = 1 -  2\exp \{ - 2^{i-k}A \}  \big( 1 + e^{-O(1) 2^{i-k}}\big)$.
 Thus, the $\PP_{S_-,S_+}^i$-probability that every move is won by Maxine equals
 $$
  \prod_{j=i}^\infty \tfrac{a_j}{a_j + b_j} \, = \, \prod_{j=i}^\infty \Big( 1 -  2\exp \{ - 2^{j-k}A \}  \big( 1 + e^{-O(1) 2^{j-k}}\big) \Big) \, = \, 1 -  2\exp \{ - 2^{i-k}A \}  \big( 1 + e^{-O(1) 2^{i-k}}\big) \, .
 $$
 This bound proves Theorem~\ref{t.unanimity}(2). The corresponding bound for $i \leq k-1$, and the proof of Theorem~\ref{t.unanimity}(3), are similar. 
 It remains then to derive Theorem~\ref{t.unanimity}(1).
 The displayed and omitted bounds permit us to choose 
 $L \in \N$ such that
 \begin{equation}\label{e.outer}
 \textrm{if $\vert i - k \vert > L$, then $\PP_{S_-,S_+}^ i(U) \geq 1/2$} \, .
 \end{equation}

 The {\em status report} $\status:\N \to \{ I,O,F\}$ is a random process defined under the law  $\PP_{S_-,S_+}^i$ that we will use to prove Theorem~\ref{t.unanimity}(1). This process takes values in a three-point set whose labels denote `inner', `outer' and `final'. 
 To record the status report, 
 we will iteratively specify an increasing sequence $\big\{ \tau_i: i \in \N \big\}$ of times valued in $\nwithzero \cup \{ \infty \}$.
 We set $\tau_0 = 0$. We check whether 
 $\vert X_0 - k \vert \leq L$, where the value of $L$ was specified in the preceding paragraph. If this condition is met then we set $\status(0) = I$.
 If the condition is not met, 
 we  set $\status(0) = O$. 
 
 Let $i \in \nwithoutzero$. Suppose that an initial status report $\status(j) \in \{I,O,F \}$, $j \in \llbracket 0,i-1 \rrbracket$,
 and an increasing sequence $\tau_j \in \N \cup \{ \infty\}$,  $j \in \llbracket 0,i-1 \rrbracket$, has been recorded.
 
 If $\status(i-1) = F$, we set $\tau_i = \infty$ and $\status(i) = F$.
 
 If $\status(i-1)= I$, we set $\tau_i = \tau_{i-1} + L$. We set
 $$
 \status(i) \, = \, \begin{cases}
  \, I  &  \text{if $\vert X_{\tau_i} - k \vert \leq L$} \\
 \, O  &  \text{in the other case} \, .
\end{cases}
$$

 If $\status(i-1) = O$, we begin to view the process $X$ run forward from time~$\tau_{i-1}$.  
 We watch for the first occasion~$F \geq \tau_{i-1} + 2$ at which the sequence of observed differences $X_{j+1}-X_j$, $F-1 \geq j \geq \tau_{i-1}$,
 has assumed both values $-1$ and $1$. If this occasion never occurs, so that $F = \infty$, we set $\status(i) = F$ and $\tau_i = \infty$. If the occasion does occur, we set $\tau_i = F$.
 The last display is used to set $\status(i)$. This completes the description of the iterative scheme for the generic later step indexed by $i \geq 1$.
 
 The status report $\status:\nwithzero \to \{I,O,F\}$ is not a Markov process, but it has simple properties that serve to prove that unanimity~$U$ is an almost sure event under $\PP_{S_-,S_+}^i$ for any $i \in \Z$. Consider then the process $\status$ under the just mentioned law. By construction, $\status$ arrives, and is absorbed, in $F$ precisely when the event $U$ occurs. To prove Theorem~\ref{t.unanimity}(1), our task is thus to show that $\status$ almost surely reaches $F$. Two properties suffice to show this. 
 
 {\em Property~$I$.} Let $j \in \nwithoutzero$.
 Suppose given a status report history $\status_i$, $i \in \llbracket 0,j-1\rrbracket$, for which $\status_{j-1} = I$.
 There exists a constant $c > 0$ that does not depend on this history such that the conditional probability that $\status(j) = O$ is at least $c$.

 {\em Property~$O$.} Let $j \in \nwithoutzero$.
 Suppose given a history $\status_i$, $i \in \llbracket 0,j-1\rrbracket$, for which $\status_{j-1} = O$.
 The conditional probability that $\status(j) = F$ is at least one-half.

Properties $I$ and $O$ show that, whatever the status report history up to a given moment, there is probability at least $c/2$ that one of the next two entries in the report is~$F$. Thus, it is inevitable that the report will eventually contain an entry in~$F$. The proof of Theorem~\ref{t.unanimity}(1) has thus been reduced to the task of deriving the two properties.

The proofs of Properties~$I$ and~$O$ depend on a {\em claim}. This states that all the information in any report history   $\status_i$, $i \in \llbracket 0,j-1\rrbracket$, in which $F$ is not recorded, is contained in the gameplay history 
$X_i$, $i \in \llbracket 0,\tau_{j-1} \rrbracket$. The claim may be proved by induction on $j$. The specifications of $\tau_j$ above are stopping times for the process $X$ that are finite when $\tau_i \in \{ I,O\}$. This proves the claim.

We now prove Property~$I$. The coefficients $a_i$ and $b_i$ are positive by Theorem~\ref{t.nashabmn}; and they are bounded by Theorem~\ref{t.ajbj}. Consider then the event that $X$ makes 
$L$ rightward jumps from time~$\tau_{j-1}$. To find a lower bound on the conditional probability of this event given the circumstance of Property~$I$, note that the claim permits us to further condition on $X$ until time $\tau_{j-1}$. Since $\vert X_{\tau_{j-1}} - k \vert \leq L$, a lower bound is offered by the minimum over $\ell \in \llbracket k-L,k+L \rrbracket$ of the product $\Pi_{i=0}^{L-1} \tfrac{a_{\ell+i}}{a_{\ell + i} + b_{\ell + i}}$. This minimum is positive because the positive and bounded quantities $a$ and $b$ that are involved are finite in number. 

And now we prove Property~$O$. Again, by the claim, we may condition on $X$ until time $\tau_{j-1}$. Since $\vert X_{\tau_{j-1}} - k \vert > L$, we may invoke~(\ref{e.outer}) to show the sought property.

This completes the proof of Theorem~\ref{t.unanimity}(1). \qed

\chapter{The Mina margin map}\label{c.allminamm}


Here we prove our results concerning the Mina margin map in three sections. Finite-trail counterparts to the map are defined and estimated in Section~\ref{s.approxmmm}, and Theorem~\ref{t.relativereward} and several consequences are derived. In Section~\ref{s.mmmtransform}, the $\theta^{-1}$- and $\PhiMac$-transforms of the map are compared, and Theorem~\ref{t.phithetainverse} is proved. In Section~\ref{s.minamarginmap}, the $\lambda \leq 0.999904$ bound
Theorem~\ref{t.minamarginvalues}(3) is derived by a scheme of explicit approximation for a well-chosen value of a suitable finite-trail counterpart for $\minamm$.

\section{Approximating the Mina margin by its finite trail counterpart}\label{s.approxmmm}
Here we prove Theorem~\ref{t.relativereward},  the third part contingent on Theorem~\ref{t.minamarginvalues}(3).
At the end of the section, we prove the consequent  Theorems~\ref{t.minamarginvalues}(1,2);~\ref{t.nashequil.prelim};~\ref{t.solutions}; and~\ref{t.nashequil}(1,2).

\subsection{An explicit form for the finite-trail Mina margin map}
 \begin{lemma}\label{l.ecinvariance}
 Let $x \in (0,\infty)$. For $k,\ell \in \Z \cup \{ \infty\} \cup \{ - \infty\}$, $k < \ell$, 
 the value of $\frac{n_k - n_\ell}{m_\ell - m_k} \in (0,\infty)$
 is a constant function of the element $\big\{ (a_i,b_i,m_i,n_i) \in (0,\infty)^2 \times \R^2 : i \in \Z \big\}$ in the equivalence class~$\mc{C}(x)$.
 \end{lemma}
{\em Remark.} When we write expressions $\frac{n_k - n_\ell}{m_\ell - m_k}$ in this section, we refer to the quantities $\frac{n_k - n_\ell}{m_\ell - m_k}(x)$
that the above lemma identifies; the value of $x \in (0,\infty)$ is often understood. 
 
{ \bf Proof of Lemma~\ref{l.ecinvariance}.} That each expression  $\frac{n_k - n_\ell}{m_\ell - m_k}$ is a finite number follows from Theorem~\ref{t.positiveabmn}(1,3).  Each expression $\frac{n_k - n_\ell}{m_\ell - m_k}$ is invariant under the translation $\transmap_{u,v}$, $u,v \in \R$ and dilation $\tau_x$, $x \in \R$, maps that must be used to interpolate any two elements of~$\mc{C}(x)$. \qed

Recall the functions  $s_j,\macsee_j,\macdee_j:(0,\infty) \to (0,\infty)$, $j \in \Z$, from Definition~\ref{d.stabc}.

Set $P_0 = S_0 = 1$. For $k \in \nwithzero$, we iteratively specify 
\begin{equation}\label{e.prodp}
 P_{k+1}(x) - P_k(x) = \prod_{i=0}^k \big( \macsee_i(x)  - 1 \big) \, ,
\end{equation}
and
\begin{equation}\label{e.prods}
 S_{k+1}(x) - S_k(x) = \prod_{i=0}^k \big( \macdee_i(x) - 1 \big) \, .
\end{equation}
Set $Q_1 = T_1 = 0$. For  $k \in \nwithoutzero$, we then set
\begin{equation}\label{e.prodq}
 Q_{k+1}(x) - Q_k(x) = \prod_{i=1}^k \big( \macsee_{-i}(x) - 1 \big)^{-1} \, ,
\end{equation}
and
\begin{equation}\label{e.prodt}
 T_{k+1}(x) - T_k(x) = \prod_{i=1}^k \big( \macdee_{-i}(x) - 1 \big)^{-1} \, .
\end{equation}
\begin{lemma}\label{l.prodinterpret}
Let $x$ equal $\macphi_0$ from Definition~\ref{d.deltai}. For $k \in \nwithzero$, 
$$
P_k(x) = \tfrac{m_k - m_{-1}}{m_0 - m_{-1}} \, \, \, \, \textrm{and} \, \, \, \,    
S_k(x) = \tfrac{n_{-1} - n_k}{n_{-1} - n_0} \, .
$$
For $\ell \in \nwithoutzero$,
$$
Q_\ell(x) = \tfrac{m_{-1} - m_{-\ell}}{m_0 - m_{-1}}   \, \, \, \, \textrm{and} \, \, \, \,    
T_\ell(x) = \tfrac{n_{-\ell} - n_{-1}}{n_{-1} - n_0}  \, .
$$
\end{lemma}
{\bf Proof.}
 The claimed formula for $P_k(x)$ is trivial when $k=0$. To prove the general formula for $P_k(x)$, it suffices to argue that $\prod_{i=0}^k \big( \macsee_i(x)  - 1 \big)$
 equals $\tfrac{m_{k+1} - m_k}{m_0 - m_{-1}}$ for $k \in \nwithzero$, and we do this by induction on~$k$. The generic step in the induction is enabled by showing that $\macsee_k(x) -1 = \tfrac{m_{k+1} - m_k}{m_k - m_{k-1}}$, which we obtain as follows: 
 $$
 \macsee_k(x) -1 \, = \, \tfrac{1 - \macgsee(s_k(x))}{\macgsee(s_k(x))} \, = \, \tfrac{1 - \macgsee(\macphi_k)}{\macgsee(\macphi_k)} \, = \, \tfrac{m_{k+1} - m_k}{m_k - m_{k-1}} \, , 
 $$
 the respective equalities by Definition~\ref{d.alphagamma}; by iterating Proposition~\ref{p.alphagammaess}(s); and  by Proposition~\ref{p.alphagammaess}($\macgsee$).
 
 Likewise, the claimed formula for $Q_\ell(x)$ is trivial when $\ell=1$.
 Establishing the formula in the general case is a matter of showing that  $\prod_{i=0}^\ell \big( \macsee_{-i}(x)  - 1 \big)^{-1}$ equals  $\tfrac{m_{-\ell} - m_{-\ell-1}}{m_0 - m_{-1}}$ for $\ell \geq 2$.
 The generic inductive step here amounts to showing that $\big( \macsee_{-\ell}(x)  - 1 \big)^{-1} = \tfrac{m_{-\ell} - m_{-\ell-1}}{m_{-\ell + 1} - m_{-\ell}}$ for such $\ell$, and follows, similarly as above, from
 $\big( \macsee_{-\ell}(x)  - 1 \big)^{-1} =  \tfrac{\macgsee(s_{-\ell}(x))}{1-\macgsee(s_{-\ell}(x))}$.
 
 The formulas for $S$ and $T$ follow when the changes 
 $$
 \textrm{$P \to S$, $Q \to T$, $k \to \ell$, $\macsee \to \macdee$, $\macgsee \to \macgdee$ 
and $m_i \to -n_{-i}$} 
$$
are made. \qed

Recall from~(\ref{e.minammfinite}) that the finite-trail Mina margin map  $\minamm_{\ell,k}:(0,\infty) \to (0,\infty)$
satisfies    $\minamm_{\ell,k}(x)
 \, = \, \frac{n_{-\ell} - n_k}{m_k - m_{-\ell}}$
for $k \in \nwithzero$ and $\ell \in \nwithoutzero$, where  $x = \macphi_0$. 
\begin{lemma}\label{l.ratiointerpret}
We have that
$$
   \minamm_{\ell,k}(x)
 \, = \, \frac{x(S_k + T_\ell)}{P_k + Q_\ell} 
 $$
for $k \in \nwithzero$ and $\ell \in \nwithoutzero$. 
\end{lemma}
In reading the proof of this result, recall the notation explained in the remark that follows Lemma~\ref{l.ecinvariance}.

{\bf Proof of Lemma~\ref{l.ratiointerpret}.} By Lemma~\ref{l.prodinterpret}, 
\begin{equation}\label{e.pqstformulas}
m_k - m_{-\ell}  = (P_k+Q_\ell)(m_0 - m_{-1}) \, \, \, \, 
\textrm{and} \, \, \, \,    n_{-\ell} - n_k = (S_k + T_\ell) (n_{-1} - n_0)  \, . 
\end{equation}
But $x = \macphi_0$, which is to say, $x = \tfrac{n_{-1} - n_0}{m_0 - m_{-1}}$. We find then that 
$$
 \minamm_{\ell,k}(x) = \frac{(S_k + T_\ell) (n_{-1} - n_0)}{(P_k+Q_\ell)(m_0 - m_{-1})} =  \frac{x(S_k + T_\ell)}{P_k+Q_\ell} \, ,
$$
as we sought to do. \qed

\subsection{Estimates for the finite trail Mina margin map}
In this subsection, we derive the following compact-uniform Cauchy sequence property of the finite-trail Mina margin maps.
\begin{proposition}\label{p.rkrell}
For $k \geq 0$, $\ell \geq 2$ and $1/3 \leq x \leq 3$,
$$
 \sup_{\substack{i \geq k+1 \\ j \geq \ell +1}} 
  \big\vert \minamm_{i,j}(x) - \minamm_{\ell,k}(x) \big\vert \, \leq \,  3^5 2^{2k-2} 6^{1-2^k} + 3^3 2^{\ell-2}  6^{\ell -2^{\ell-1}} 
\, .
$$
\end{proposition}
The next lemma assembles key elements for the proof of Proposition~\ref{p.rkrell}. We omit to denote the argument `$(x)$'
of $\mc{M}_{\cdot,\cdot}$, $P$, $Q$, $S$ and $T$ as we derive this proposition.
\begin{lemma}\label{l.pqst}
Let $k \in \nwithzero$ and $x \in \R$.
\begin{enumerate}
\item For $k \geq 0$ and $x \leq 3$, $P_{k+1} - P_k \leq 2^{2k} 6^{1-2^k}$.
\item For $k \geq 1$ and $x \geq 1/3$, $Q_{k+1} - Q_k \leq 2^{2k} 6^{1-2^k}$.
\item For $k \geq 0$ and $x \leq 3$, $S_{k+1} - S_k \leq 2^{2k+1} 6^{1-2^{k+1}}$.
\item For $\ell \geq 2$ and $x \geq 1/3$, $T_{\ell+1} - T_\ell \leq 3 (12)^{\ell-1}  6^{1-2^{\ell-1}}$.
\end{enumerate}
\end{lemma}
Two simple lemmas gather estimates needed to prove Lemma~\ref{l.pqst}.
\begin{lemma}\label{l.abounds}
\leavevmode
\begin{enumerate}
 \item For $x \in (0,\infty)$, $s(x) \leq x^2/2$.
 \item For $x \in (0,\infty)$,  $\macsee(x) \leq 1 +  2x$. 
\item  For $x \in (0,\infty)$,  $\macsee(x) \geq 1+x/2$.
 \item For $x \in (0,3]$, $\macdee(x) - 1 \leq 1/3$. 
 \item For $x \in (0,\infty)$, $\macdee(x) \geq 2^{-3/2} x^{1/2}$. \end{enumerate}
 \end{lemma}
{\bf Proof: (1).} Since $\macomega \geq 1$, $\macbeta(x) \geq 0$.
Thus, Lemma~\ref{l.acsfacts}(3) implies that  $s(x) \leq \macbeta(x)^2/2$. So the result reduces to Lemma~\ref{l.acsfacts}(4).  \\
{\bf (2).} By Definition~\ref{d.acs}, $\macsee(x) = \tfrac{(\macomega+3)^2}{16} =  \tfrac{8x+10+6\macomega}{16} \leq 1 + 2x$ where the inequality is due to $\macomega = \sqrt{1+ 8x} \leq 1+4x$ for $x \geq 0$. \\
{\bf (3).} We have that $\macsee(x) =  \tfrac{8x+10+6\macomega}{16} \geq 1 + x/2$ from $\macomega \geq 1$. \\
{\bf (4).} By Lemma~\ref{l.acsfacts}(1:$d$), $\macdee(x) -1 \leq \macdee(3) - 1 =  1/3$.\\
{\bf (5).} Recall that $\macdee(x) = \tfrac{(\macomega + 3)^2}{8(\macomega+1)}$ where $\macomega = \sqrt{8x+1}$. Thus, $\macdee(x) \geq (\macomega +3)/8 \geq 2^{-3/2}x^{1/2}$.  \qed

\begin{lemma}\label{l.stbounds}
Let $j \in \nwithoutzero$. 
\begin{enumerate}
\item For $x \leq 3$, $s_j(x) \leq 2 \cdot 6^{-2^{j-1}}$. 
\item For $x \geq 1/3$, $s_{-j}(x) \geq 2^{-1}  6^{2^{j-1}}$.
\item For $i \geq 1$ and $x \leq 3$, $\macsee_i(x) \leq  1 +  2^2 6^{-2^{i-1}}$. 
\item For $i \in \Z$ and $x \in (0,\infty)$, $\macdee_i(x) - 1 \leq s_i(x)^2$.
\item For $i \geq 2$ and $x \geq 1/3$, $\macdee_{-i}(x)-1  \geq 2^{-1}6^{2^{i-2}-1}$.
\end{enumerate}
\end{lemma}
{\bf Proof: (1).} Note that $s(3) = 1/3$ since $\macomega(3) = 5$. We may thus use Lemma~\ref{l.abounds}(1) to prove the desired statement by induction.  \\
{\bf (2).} Due to the preceding and $s_{-j}(x) = 1/s_j(1/x)$ (which is Corollary~\ref{c.sminusone}). \\
{\bf (3).}
By Lemma~\ref{l.acsfacts}(1:$c$), Lemma~\ref{l.stbounds}(1) and Lemma~\ref{l.abounds}(2), 
$$
\macsee_i(x) = \macsee \big( s_i(x) \big) \leq  \macsee \big( 2 \cdot 6^{-2^{i-1}}  \big) \leq 1 +  2^2 6^{-2^{i-1}}
$$ 
for $i \geq 1$ and $x \leq 3$. \\
{\bf (4).} It is enough to show that $\macdee(x) \leq 1 + x^2$. To see this, note that $\macdee(x) - 1 = \macgdee(x)^{-1} - 1$. From $1 - \macgdee(x) = \big( 1- \tfrac{1}{\macbeta(x) +1} \big)^2$
and Lemma~\ref{l.acsfacts}(4), we find that $\macgdee(x) \geq 1 - \tfrac{x^2}{(1+x)^2}$, so that $\macdee(x) - 1 \leq \tfrac{x^2}{1+2x} \leq x^2$. \\
{\bf (5).} Note that $\macdee_{-i}(x) =  \macdee\big(s_{-i}(x)\big) \geq  \macdee\big(2^{-1}6^{2^{i-1}}\big) \geq 2^{-2}6^{2^{i-2}}$, where the first inequality is due to  Lemma~\ref{l.acsfacts}(1:$d$)  and Lemma~\ref{l.stbounds}(2),
and the second to Lemma~\ref{l.abounds}(5). From this, the sought result follows. \qed

{\bf Proof of Lemma~\ref{l.pqst}: (1).} Note that $\macsee(x) \leq 4$ for $x \in (0,3]$ by Lemma~\ref{l.acsfacts}(1:$c$) and $\macsee(3)=4$. Thus we bound the first term in the product in~(\ref{e.prodp}).
 Bounding the latter terms by Lemma~\ref{l.stbounds}(3), we find that 
$$
 P_{k+1} - P_k \, = \, \prod_{i=0}^k \big( \macsee_i(x) - 1 \big) \, \leq \, 3 \prod_{i=1}^k   2^2 6^{-2^{i-1}}
  \, ,
$$
whence the sought result. 

{\bf (2).}  Note that
$$
\macsee_{-i}(x) -1 = \macsee \big( s_{-i}(x)\big) -1 \geq \macsee \big( 2^{-1} 6^{2^{i-1}} \big) -1 \geq  2^{-2}6^{2^{i-1}}   \, ,
$$
where the first inequality holds when $x \geq 1/3$ in view of Lemma~\ref{l.acsfacts}(1:$c$) and Lemma~\ref{l.stbounds}(2); the second is due to Lemma~\ref{l.abounds}(3).
By~(\ref{e.prodq}), $Q_{k+1} -Q_k \leq \prod_{i=1}^k 2^2 6^{-2^{i-1}}  = 2^{2k} 6^{1 - 2^k}$, whence Lemma~\ref{l.pqst}(2).

{\bf (3).}
Note that
$$
S_{k+1} - S_k \leq  3^{-1}\prod_{i=1}^k 2^2 6^{-2^i} = 2^{2k+1} 6^{1-2^{k+1}} \, ,
$$
where, in the first inequality, the first term in the product expression in~(\ref{e.prods}) is bounded by use of Lemma~\ref{l.abounds}(4), and the later terms are taken care of by the bounds $\macdee_i(x) -1 \leq s_i(x)^2 \leq 2^2 6^{-2^i}$, which are valid for $x \leq 3$ and $i \geq 1$ in view of Lemma~\ref{l.stbounds}(1,4). 
 
 {\bf (4).} 
 Since $s_{-1}(1/3) = 3$, Proposition~\ref{p.sminusone} and Lemma~\ref{l.acsfacts}(1:$s$) imply that $s_{-1}(x) \geq 3$ for $x \geq 1/3$. And since $\macdee(3) = 4/3$, the same result implies that $\big(\macdee_{-1}(x) - 1\big)^{-1} \leq 3$ for such $x$.
 Applying these bounds alongside Lemma~\ref{l.stbounds}(5) to~(\ref{e.prodt}), we see that
 $$
  T_{\ell+1} - T_\ell \, \leq \, 3 \cdot \prod_{i=2}^\ell 2 \cdot 6^{1 - 2^{i-2}}  \, = \, 3 (12)^{\ell-1}  6^{1-2^{\ell-1}}
  $$
 for $x \geq 1/3$ and $\ell \geq 2$. Whence Lemma~\ref{l.pqst}(2). \qed

Two further lemmas will permit the derivation of Proposition~\ref{p.rkrell} from Lemma~\ref{l.pqst}.
\begin{lemma}\label{l.lub}
We have that
$$
x^{-1} \big\vert \minamm_{k+1,\ell} - \minamm_{\ell,k} \big\vert \,
\leq \, \max \, \Big\{ \, S_{k+1} - S_k \, , \,  (S_k + T_\ell)(P_{k+1}-P_k) \, \Big\} \, ,
$$
and
$$
x^{-1} \big\vert \minamm_{\ell,k+1} - \minamm_{\ell,k} \big\vert \,
\leq \, \max \, \Big\{ \, T_{\ell+1} - T_\ell \, , \,  (S_k + T_\ell)(Q_{\ell+1}-Q_\ell) \, \Big\} \, .
$$
\end{lemma}
{\bf Proof.}
Since $P_j$ and $Q_j$ are at least one whenever $j$ is at least one, it is enough to show that
\begin{equation}\label{e.lubone}
x^{-1} \big\vert \minamm_{k+1,\ell} - \minamm_{\ell,k} \big\vert \,
\leq \, \max \, \bigg\{ \, \frac{S_{k+1} - S_k}{P_{k+1}+Q_\ell} \, , \,  \frac{(S_k + T_\ell)(P_{k+1}-P_k)}{(P_{k+1}+Q_\ell)(P_k+Q_\ell)} \, \bigg\} \, ,
\end{equation}
and 
\begin{equation}\label{e.lubtwo}
x^{-1} \big\vert \minamm_{\ell+1,k} - \minamm_{\ell,k} \big\vert \,
\leq \, \max \, \bigg\{ \, \frac{T_{k+1} - T_k}{P_{k+1}+Q_\ell} \, , \,  \frac{(S_k + T_\ell)(Q_{k+1}-Q_k)}{(P_{k+1}+Q_\ell)(P_k+Q_\ell)} \, \bigg\} \, .
\end{equation}
 Note that
\begin{eqnarray}
 x^{-1} \big( \minamm_{k+1,\ell} - \minamm_{\ell,k} \big) & = & \frac{S_{k+1}+T_\ell}{P_{k+1} + Q_\ell} -  \frac{S_k+T_\ell}{ P_k + Q_\ell} \nonumber \\
  & = & \frac{(S_{k+1} + T_\ell)(P_k + Q_\ell)  - (S_k + T_\ell)(P_{k+1} + Q_\ell)}{(P_{k+1} + Q_\ell)(P_k + Q_\ell)} \, . \nonumber
\end{eqnarray}
The numerator in the latter term equals $ (P_k + Q_\ell) (S_{k+1} - S_k)  - (S_k+T_\ell)(P_{k+1} - P_k)$. Since this is 
 a difference of positive terms, and the right-hand denominator above is positive, we obtain~(\ref{e.lubone}).

Note further that
\begin{eqnarray}
 x^{-1} \big( \minamm_{\ell+1,k} - \minamm_{\ell,k} \big) & = & \frac{S_k+T_{\ell+1}}{P_k+ Q_{\ell+1}} -  \frac{S_k+T_\ell}{ P_k + Q_\ell} \nonumber \\
  & = & \frac{(S_k + T_{\ell+1})(P_k + Q_\ell)  - (S_k + T_\ell)(P_k + Q_{\ell+1})}{(P_k + Q_{\ell+1})(P_k + Q_\ell)} \, . \nonumber
\end{eqnarray}
In this case, the numerator in the last line is $(P_k + Q_\ell) (T_{\ell+1} - T_\ell)  - (S_k+T_\ell)(Q_{\ell+1} - Q_\ell)$.
By reasoning as we did above,
 we obtain~(\ref{e.lubtwo}). This completes the proof of Lemma~\ref{l.lub}. \qed

\begin{lemma}\label{l.sup}
\leavevmode
\begin{enumerate}
\item For $x \leq 3$, $\sup_{k \geq 1} S_k \leq 3/2$.
\item For $x \geq 1/3$, $\sup_{k \geq 1} T_k \leq 12$.
\end{enumerate}
\end{lemma}
{\bf Proof: (1).}
By $S_0 =1$ and Lemma~\ref{l.pqst}(2), $\sup_{k \geq 1} S_k \leq 1 + \sum_{k=0}^\infty 2^{2k+1}6^{1-2^{k+1}} = 1 + 3^{-1} + 2^{-4}3^{-7} + 2^{-10}3^{-15} + \cdots = 1.333361\cdots \leq 3/2$.

{\bf (2).} Recall that, by definition, 
 $T_1 = 0$. The quantity $T_2$ equals  $\big(\macdee_{-1}(x) - 1\big)^{-1}$ which is at most $3$ when  $x \geq 1/3$, as we noted 
  in the proof of Lemma~\ref{l.pqst}(4). Using these alongside Lemma~\ref{l.pqst}(4), we find that 
  $\sup_{k \geq 1} T_k \leq 3 + \sum_{\ell = 2}^\infty  3 (12)^{\ell-1}  6^{1-2^{\ell-1}} = 3 + \sum_{k=0}^\infty 2^k 6^{k+3 - 2^{k+1}} \leq 12$. \qed

{\bf Proof of Proposition~\ref{p.rkrell}.}
By Lemma~\ref{l.lub}(1), Lemma~\ref{l.sup}(1,2) and Lemma~\ref{l.pqst}(1,3), we have that
$$
 x^{-1} \big\vert \minamm_{k+1,\ell} - \minamm_{\ell,k} \big\vert \leq \max \big\{ 2^{2k +1} 6^{1 - 2^{k+1}} , 2^{2k-4} 6^{4-2^k} \big\}
$$
for $k \geq 0$ and $\ell \geq 1$, where here we used $\tfrac{27}{2} \cdot 2^{2k} 6^{1 - 2^k}  = 2^{2k-4} 6^{4-2^k}$. And by Lemma~\ref{l.lub}(1), Lemma~\ref{l.sup}(1,2) and Lemma~\ref{l.pqst}(2,4),
$$
 x^{-1} \big\vert \minamm_{\ell+1,k} - \minamm_{\ell,k} \big\vert \leq \max \big\{ 3 (12)^{\ell - 1} 6^{1 - 2^{\ell-1}} , \tfrac{27}{2} 2^{2\ell} 6^{1-2^\ell} \big\}
$$
for $k \geq 0$ and $\ell \geq 2$. In each of the two displayed maximums, it is the first expression which is the greater for the stated ranges of $k$ and $\ell$. 
Set $g(i) = 2^{2i-1} 6^{1-2^i}$ and $h(j) = (12)^{j-1} 6^{1 - 2^{j-1}}$, and note that $g(i+1)/g(i) \leq 1/3$ and $h(j+1)/h(j) \leq 1/3$ provided that
$i \leq k$ and $j \geq \ell$ where $k \geq 0$ and $\ell \geq 2$. What we learn is that
$$
 x^{-1} 
 \sup_{\substack{i \geq k+1 \\ j \geq \ell +1}}  \big\vert \minamm_{i,j} - \minamm_{\ell,k} \big\vert \, \leq \, \tfrac{3}{2} \cdot \Big( 3^3 2^{2k-1} 6^{1-2^k} + 3(12)^{\ell - 1} 6^{1-2^{\ell - 1}} \Big)
$$
for $k \leq 0$, $\ell \geq 2$ and $1/3 \leq x \leq 3$. Using $x \leq 3$, and rewriting the product of three and the right-hand side of this display, we obtain Proposition~\ref{p.rkrell}. \qed

\subsection{Proofs via the finite trail Mina margin map}
$\empty$

{\bf Proof of Theorem~\ref{t.relativereward}(1).}
Note that
\begin{equation}\label{e.rkkconvergence}
 \minamm_{k,k}(x) \, = \, \frac{n_{-k} - n_k}{m_k - m_{-k}} \, \to \, \frac{n_{-\infty} - n_\infty}{m_\infty - m_{-\infty}} \, = \, \minamm(x)
\end{equation}
where the convergence, which is in the limit $k \to \infty$, is explained by the proof of Theorem~\ref{t.positiveabmn}(3); the latter equality here is due to the specification of $\minamm(x)$ in Definition~\ref{d.r} and to Lemma~\ref{l.ecinvariance}. 
Note next that, for $i \in \Z$, the standard element in $\mc{C}\big( s_i(x) \big)$ is equal to the left shift by $i$ places of the standard element in $\mc{C}(x)$. Thus, 
$$
\minamm_{k,k}\big( s_i(x) \big) \, = \, \frac{n_{-k+i} - n_{k+i}}{m_{k+i} - m_{-k+i}} \, .
$$
The left-hand side converges to $\minamm\big( s_i(x) \big)$ in the limit of high $k$, by~(\ref{e.rkkconvergence}). The right-hand side converges to  $\frac{n_{-\infty} - n_\infty}{m_\infty - m_{-\infty}} = \minamm(x)$ 
since $m$- and $n$-differences vanish asymptotically at high values of the index by Theorem~\ref{t.ajbj}. Thus we find that  $\minamm\big( s_i(x) \big) = \minamm(x)$ for $i \in \Z$ and $x \in [1/3,3]$.
Since $\cup_{i \in \Z} s_i[1/3,3] = (0,\infty)$, we see that $\minamm(x)$ exists for all $x \in (0,\infty)$, and that in fact   $\minamm\big( s_i(x) \big) = \minamm(x)$ for $i \in \Z$ and $x \in (0,\infty)$.  This completes the proof of Theorem~\ref{t.relativereward}(1).

{\bf (2).} We first show that $\minamm$ is continuous on $(0,\infty)$. Proposition~\ref{p.rkrell} shows that $\minamm_{k,k}$ converges uniformly as $k \to \infty$ on $[1/3,3]$. By~(\ref{e.rkkconvergence}), the limiting function is the restriction of $\minamm$
to $[1/3,3]$. 
Since the constituent functions $c_i,d_i:(0,\infty) \to (0,\infty)$, $i \in \Z$, are continuous, we see that the map $\minamm_{k,k}:[1/3,3] \to (0,\infty)$ is continuous for any $k \in \nwithoutzero$.
Thus, $\minamm$ is continuous on this interval. But $\minamm(x) = \minamm(s(x))$ for $x \in (0,\infty)$ by Theorem~\ref{t.relativereward}(1), and $\minamm(3) = \minamm(1/3)$ since $s(3) = 1/3$. Since $s:(0,\infty) \to (0,\infty)$ is seen to be continuous from its specification in Definition~\ref{d.acs}, we confirm that $\minamm$ is continuous on~$(0,\infty)$. 

To derive the formula for $\minamm(x)$ claimed in Theorem~\ref{t.relativereward}(2), note, by decoding the notation for products in Definition~\ref{d.zdefault}, that
this formula may expressed in our present notation in the form $\tfrac{n_{-\infty} - n_\infty}{m_\infty - m_{-\infty}} = \tfrac{x(S_\infty +T_\infty)}{P_\infty + Q_\infty}$,
where in fact we have extended this notation to write $*_\infty$ for $\lim_{k \to \infty}*_k$ with $* \in \{P,Q,S,T\}$. Since  $x = \tfrac{m_0 - m_{-1}}{n_{-1} - n_0}$,
the sought formula is a consequence of 
$$
 m_\infty - m_{-\infty} = (P_\infty + Q_\infty) (m_0 - m_{-1}) \, \, \, \, \textrm{and} \, \, \, \, n_{-\infty} -n_\infty = (S_\infty + T_\infty) (n_{-1} - n_0) \, .
$$
To obtain these identities, we take the limit in high $k$ and $\ell$ of the two formulas in~(\ref{e.pqstformulas}), using Theorem~\ref{t.positiveabmn}(3) to justify that the limiting expressions are finite real numbers. This completes the proof of Theorem~\ref{t.relativereward}(2). \qed
 
 \subsection{Some further consequences}

In order to prove Theorem~\ref{t.minamarginvalues} and Theorem~\ref{t.relativereward}(3), we now offer a definition of the quantity $\lambda \in (0,1]$ to which these results refer.
\begin{definition}\label{D.lambda}
We set $\lambda = \inf \{ \minamm(x): x \in [1/3,3] \}$.
\end{definition}
\begin{lemma}\label{l.infimumminamarginmap}
There exists $x_0 \in [1/3,3]$ such that $\minamm(x_0) = \lambda$. We have that 
$$
\lambda \, = \,  \inf \{ \minamm(x): x \in (0,\infty) \} \, .
$$
\end{lemma}
{\bf Proof.}
Since $\minamm:[1/3,3] \to (0,\infty)$ is continuous by Theorem~\ref{t.relativereward}(2), the infimum is attained on $[1/3,3]$, and we may find $x_0 \in [1/3,3]$ so that $\minamm(x_0) = \lambda \in (0,\infty)$.
The proof of Lemma~\ref{l.battlefield} shows that $\cup_{i \in \Z} s_i[1/3,3] = (0,\infty)$. 
By Theorem~\ref{t.relativereward}(3),  we see thus that  $\lambda = \inf \{ \minamm(x): x \in (0,\infty) \}$. \qed
\begin{lemma}\label{l.rangeminamarginmap}
For $x \in (0,\infty)$, $\minamm(x^{-1}) = \minamm(x)^{-1}$. 
\end{lemma}
{\bf Proof.} Recall from Definition~\ref{d.r} that $\minamm(x) =\nst_{-\infty}(x)$ for  $x \in (0,\infty)$ is 
the Mina margin of the standard solution $\big(\asta_i(x),\bst_i(x),\mst_i(x),\nst_i(x): i \in \Z \big)$. The $\phi_0$-value of this solution is equal to  $\tfrac{\nst_{-1}(x) - \nst_0(x)}{\mst_0(x) - \mst_{-1}(x)} = x$. 
By Proposition~\ref{p.rolereversal} and dilation, the quadruple  
$$
 \nst_{-\infty}(x)^{-1} \cdot  \Big( \, \bst_{-i}(x) \, , \, \asta_{-i}(x) \, , \,
 \nst_{-i}(x) \, , \, \mst_{-i}(x): i \in \Z \, \Big)
$$ 
is also a standard \abmnmacspace solution. Its $\phi_1$-value equals  $\tfrac{\mst_0(x) - \mst_{-1}(x)}{\nst_{-1}(x) - \nst_0(x)}  = x^{-1}$.
The left shift by one place of the displayed quadruple is thus a standard \abmnmacspace solution whose $\phi_0$-value equals $x^{-1}.$
The quantity $\minamm(x^{-1})$, which by definition equals $\nst_{-\infty}(x^{-1})$, is thus found to be equal to 
$\nst_{-\infty}(x)^{-1} \cdot  \mst_\infty(x) =\nst_{-\infty}(x)^{-1} = \minamm(x)^{-1}$. Here, we used that  $\mst_\infty(x) = 1$ since the solution $\big(\asta_i(x),\bst_i(x),\mst_i(x),\nst_i(x): i \in \Z \big)$ is standard.
\qed 

\begin{lemma}\label{l.supremumminamarginmap}
 We have that $\lambda^{-1} = \sup \{ \minamm(x): x \in [1,3,3] \} =  \sup \{ \minamm(x): x \in (0,\infty) \}$.
Further, there exists $y_0 \in [1,3,3]$ such that $\minamm(y_0) = \lambda^{-1}$.
\end{lemma}
{\bf Proof.} The ranges $\minamm[1/3,3]$ and $\minamm(0,\infty)$ are invariant under the transformation $z \to z^{-1}$ in view of Lemma~\ref{l.rangeminamarginmap}.
The supremum of the continuous function $\minamm$ is attained on $[1/3,3]$.
 \qed

{\bf Proof of Theorem~\ref{t.relativereward}(3).} 
The range $\minamm[1/3,3]$ has maximum $\lambda^{-1}$ and minimum $\lambda$, by Lemmas~\ref{l.infimumminamarginmap} and~\ref{l.supremumminamarginmap}. By the continuity of $\minamm$ on $[1/3,3]$,
$\minamm[1/3,3]$ is seen to equal $[\lambda,\lambda^{-1}]$. Since  $\cup_{i \in \Z} s_i[1/3,3] = (0,\infty)$, $\minamm(0,\infty)$ equals $\minamm[1/3,3]$. Note that $\lambda \in (0,1]$ since $\minamm[1/3,3] = [\lambda,\lambda^{-1}]$
and $\minamm$ is continuous. The remaining assertion that we need to validate, which is that $\lambda$ is at most $0.999904$, is Theorem~\ref{t.minamarginvalues}(3), whose proof will appear in Section~\ref{s.minamarginmap}. \qed

{\bf Proof of Theorem~\ref{t.minamarginvalues}(1).}
Theorem~\ref{t.relativereward}(3) shows that the set of values of the Mina margins of standard positive \abmnmacspace solutions is equal to $[\lambda,\lambda^{-1}]$.
Now consider an arbitrary positive \abmnmacspace solution. The value of the Mina margin is shared between this solution and the equivalent standard solution. Thus, no new values for the Mina margin emerge as the solution set is enlarged from standard to general. 

{\bf (2).} Consider a positive \abmnmacspace solution $(a,b,m,n)$ with boundary data $(m_{-\infty},m_\infty,n_{-\infty},n_\infty) \in \R^4$. By Theorem~\ref{t.positiveabmn}(1), $(a,b,m,n)$ is strict. Thus, $m_{-\infty} < m_\infty$ and $n_\infty < n_{-\infty}$. By Theorem~\ref{t.minamarginvalues}(1), $\tfrac{n_{-\infty} - n_\infty}{m_\infty - m_{-\infty}} \in [\lambda,\lambda^{-1}]$.

Conversely, suppose that $(m_{-\infty},m_\infty,n_{-\infty},n_\infty) \in \R^4$ satisfies $m_{-\infty} < m_\infty$, $n_\infty < n_{-\infty}$ and $\tfrac{n_{-\infty} - n_\infty}{m_\infty - m_{-\infty}} \in [\lambda,\lambda^{-1}]$. Set $x$ equal to the latter quantity. In the notation of Section~\ref{s.solvingabmn}, the image of the standard solution $\big( \asta_i(x),\bst_i(x),\mst_i(x),\nst_i(x) : i \in \Z \big)$
under the transformation $\transmap_{x,y} \circ \tau_u$,  where $x = m_{-\infty}$, $n = n_{-\infty}$ and $u = m_\infty - m_{-\infty}$,
is a positive \abmnmacspace solution with boundary data  $(m_{-\infty},m_\infty,n_{-\infty},n_\infty)$.  This completes the proof of  Theorem~\ref{t.minamarginvalues}(2). \qed

{\bf Proof of Theorem~\ref{t.nashequil.prelim}.}
Let $x \in (0,\infty)$. By Theorem~\ref{t.nashabmn}, the game $\standard(x)$ has a time-invariant Nash equilibrium if and only if there exists a positive \abmnmacspace solution whose Mina margin equals $x$.
The latter condition is equivalent to $x \in [\lambda,\lambda^{-1}]$ by Theorem~\ref{t.minamarginvalues}(1). \qed

{\bf Proof of Theorem~\ref{t.solutions}.}
Let $y \in [1/3,3]$ with $\minamm(y) = \lambda^{-1}$ be the value~$y_0$ assured by Lemma~\ref{l.supremumminamarginmap}. Let $w \leq y$ be maximal such that $\minamm(w) = \lambda$, where Theorem~\ref{t.relativereward} assures the existence of this quantity. We have that $w < y$
because $\minamm$ is continuous and assigns different values to these two points.
Set $z = s_{-1}(w)$. By Proposition~\ref{p.sminusone} and Lemma~\ref{l.acsfacts}(5), $z > w$. Since $\minamm(z) = \minamm(w) = \lambda$ by Theorem~\ref{t.relativereward}(1), we thus have that $z > y$.

Now let $x \in (\lambda,\lambda^{-1})$. By the continuity of $\minamm$, we may find $u \in (w,y)$ and $v \in (y,z)$ such that $\minamm(u) = \minamm(v) =x$. Note that $w < u < v < z = s_{-1}(w)$.
The quadruples 
$$
\big(\asta_i(u),\bst_i(u),\mst_i(u),\nst_i(u)  : i \in \Z  \big) \, \, \, \, \textrm{and} \, \, \, \,   \big(\asta_i(v),\bst_i(v),\mst_i(v),\nst_i(v)  : i \in \Z  \big)
$$  
are standard \abmnmacspace solutions of Mina margin~$x$. They are shift inequivalent because $u$ is not equal to $s_i(v)$ for any $i \in \Z$. Indeed, the condition $s_i(v) \in [w,s_{-1}(w))$ implies that $i =0$, but $s_0(v) = v \not= u$. This pair of solutions demonstrates that $Q(x) \geq 2$, as required to obtain Theorem~\ref{t.solutions}. \qed

We end this section by proving Theorem~\ref{t.nashequil}, in part as a consequence of Theorem~\ref{t.relativereward}(1). 
This makes now a convenient moment to derive the next result, thus rendering rigorous a verbal argument in the first paragraph of Section~\ref{s.solvingabmn}.
  
{\bf Proof of Proposition~\ref{p.abmnclassify}.}
To prove the two parts of this result, it is enough to argue that there is a unique standard solution, and a unique default solution, to which any positive \abmnmacspace solution is equivalent. Suppose then that $(a,b,m,n) = \big\{ (a_i,b_i,m_i,n_i) \in (0,\infty)^2 \times \R^2 : i \in \Z \big\}$ is a  positive \abmnmacspace solution. The boundary values $m_{-\infty}$ and $n_\infty$ exist as real numbers by Theorem~\ref{t.positiveabmn}(3).
Note then that the translation $(a',b',m',n')  = \tau_{-m_{-\infty},-n_\infty}$ has  $m'_{-\infty} = n'_\infty = 0$. Write $x = \frac{n'_{-1} - n'_0}{m'_0 - m'_{-1}}$
and $y =  m'_\infty - m'_{-\infty}$.
 By applying the dilation $\tau_u$ to  $(a',b',m',n')$, we obtain a default solution if $u = x^{-1}$ and a standard solution if $u = y^{-1}$.
 We have seen that $\tau_{x^{-1}} \circ \transmap_{-m_{-\infty},-n_\infty}(a,b,m,n)$ is a default \abmnmacspace solution. It is clear that any variation of the parameters $(x^{-1},m_{-\infty},n_\infty)$
 will result in an \abmnmacspace solution that fails to be default. Likewise,  $\tau_{y^{-1}} \circ \transmap_{-m_{-\infty},-n_\infty}(a,b,m,n)$ has been shown to be a standard \abmnmacspace solution.
 Any variation of  $(y^{-1},m_{-\infty},n_\infty)$ will result in an \abmnmacspace solution that fails to be standard. Thus we complete the proof of Proposition~\ref{p.abmnclassify}. \qed

{\bf Proof of Theorem~\ref{t.nashequil}(1).} By Theorem~\ref{t.nashabmn}, a time-invariant Nash equilibrium in $\standard(x)$
is the reverse-ordered $(a,b)$-component of a standard \abmnmacspace solution whose Mina margin equals $x$. Since Proposition~\ref{p.abmnclassify}
implies that standard \abmnmacspace solutions are indexed by the value $z \in (0,\infty)$ of their $\centralratio$, we obtain Theorem~\ref{t.nashequil}(1).

{\bf (2).} By Theorem~\ref{t.relativereward}(1) and Proposition~\ref{p.sminusone}, $\minamm(s_k(x))$ equals $\minamm(x)$ for all $x \in (0,\infty)$ and $k \in \Z$.
Thus, the set $X$ is the disjoint union of $s_k(Y)$ as $k$ ranges over~$\Z$. Proposition~\ref{p.shift} then yields Theorem~\ref{t.nashequil}(2). \qed

\section{The Mina margin map after domain coordinate change}\label{s.mmmtransform}
In this section, we prove Theorem~\ref{t.phithetainverse}.
The map  $\PhiMac:\R \to (0,\infty)$ is an increasing surjection, so we may set $\Psi = \PhiMac^{-1}:(0,\infty) \to \R$.
The proof of the theorem will harness the next result.
\begin{lemma}\label{l.thetapsi}
There exists a constant $C > 0$ such that $\vert \theta(x) - \Psi(x) \vert \leq C$ for $x \geq 1/3$.
\end{lemma}
Two further results will serve to prove Lemma~\ref{l.thetapsi}.
\begin{lemma}\label{l.psi}
We have that
$$ \Psi(x)\, = \, \begin{cases}
  \, \, \log_2 \big( \log_2(x) +1 \big)  &  \text{if $x \in [1,\infty)$} \, , \\
 \, \, - \log_2 \big( - \log_2(x) +1 \big)   &  \text{if $x \in (0,1)$} \, .
\end{cases}
$$
\end{lemma}
{\bf Proof.} The formulas follow from the expressions $2^{2^x - 1}$ and $2^{-(2^{-x} -1)}$ for $\PhiMac(x)$ that are respectively valid when $x \geq 0$ and $ x < 0$. \qed

For $x \in [1/3,3]$ and $i \in \nwithzero$, we write $s_{-i}(x)$ in the form $2^{2^i \ceemack_i(x) - 1}$.
\begin{lemma}\label{l.sminusi}
There exists a constant $C > 0$ and a function $\ceemack:[1/3,3] \to (0,\infty)$ such that the bound $\vert \ceemack_i(x) - \ceemack(x) \vert \leq C 2^{-i}$ holds for $x \in [1/3,3]$ and $i \in \nwithzero$. Further, $\inf \big\{ \ceemack(x): x \in [1/3,3] \big\} > 0$.
\end{lemma}
{\bf Proof.}
From the relation 
$$
s_i(x) = 2 s_{-(i-1)}(x)^2 + O \big( s_{-(i-1)}(x) \big)
$$
and the form $s_{-j}(x) = 2^{2^i \ceemack_j - 1}$ (where we write $\ceemack_j = \ceemack_j(x)$), we find that
$$
 2^{2^j \ceemack_j - 1} \, = \, 2^{2^j \ceemack_{j-1}(x) - 1} + O \big( 2^{2^{j-1} \ceemack_{j-1}(x) - 1} \big) \, = \,  2^{2^j \ceemack_{j-1}(x) - 1} \Big( 1+  2^{-2^{j-1} \ceemack_{j-1}(x)} \Big)
$$
so that
$$
2^j \ceemack_j = 2^j \ceemack_{j-1} - 1 + \Theta(1) 2^{-2^{j-1}\ceemack_{j-1}}
$$
and
$$
\ceemack_j = \ceemack_{j-1} - 2^{-j} + \Theta(1) 2^{-j - 2^{j-1}\ceemack_{j-1}} = \ceemack_{j-1} + \Theta(1) 2^{-j} \, .
$$
We thus learn that there exists $\ceemack = \ceemack(x) \in [0,\infty)$ such that $\vert \ceemack_j - \ceemack \vert \leq O(1) 2^{-j}$. We may exclude the possibility that $\ceemack$ equals zero because, in this case, we would have that $\ceemack_j \leq O(1) 2^{-j}$, which would imply the false assertion that $s_{-j}(x) = 2^{2^j \ceemack_j - 1} = O(1)$ is bounded above independently of $j \in \nwithzero$. 

We now argue that $\inf \big\{ \ceemack(x): x \in [1/3,3] \big\} > 0$. This follows from $\ceemack(1/3) > 0$ and the fact that $s_{-i}(x)$ is increasing in $x \in [1/3,3]$ for each $i \in \nwithzero$. This completes the proof of Lemma~\ref{l.sminusi}. \qed

 {\bf Proof of Lemma~\ref{l.thetapsi}.}
 Note that $s_{-i}(x) \geq 3$ for $i \in \nwithoutzero$ and $x \in [1/3,3]$. By Lemma~\ref{l.psi} and $s_{-i}(x) = 2^{2^i \ceemack_i(x) - 1}$, we find that 
 $\Psi\big(s_{-i}(x)\big) = i + \log_2 \ceemack_i$. Indeed, we find that
 $$
  \big\vert \Psi \big( s_{-i}(x) \big) - i \big\vert = \log_2 \big( \ceemack_i/c \big) = \log \big( 1 + \tfrac{\ceemack_i - c}{c} \big) \leq O(1) 2^{-i} \ ,
 $$
 the inequality due to Lemma~\ref{l.sminusi}. 
 
 On the other hand, $\theta\big( s_{-i}(x) \big)$ is equal to the unique value $J \in \Z$ such that $s_J \big( s_{-i}(x) \big) \in [1/3,3)$.
 When $x \in [1/3,3)$, we see then that $J = i$.
 
 We find then that $\vert \Psi(x) - \theta(x) \vert = O(1)$ for $x \geq 1/3$, as we sought to do in proving Lemma~\ref{l.thetapsi}. \qed

{\bf Proof of Theorem~\ref{t.phithetainverse}(1,2).}
We first claim that 
\begin{equation}\label{e.claimsk}
 s_k \big( \theta^{-1}(x+k) \big) = \theta^{-1}(x) \, .
\end{equation}
To check this, note that $\theta \big( s_k(z) \big) = \theta(z) - k$, so that 
 $$
 \theta \Big( s_k \big( \theta^{-1}(x+k) \big) \Big) = \theta \big( \theta^{-1}(x+k)\big) - k = (k+x) - k
 = x \, ,
 $$
 as desired.

Note then that
$$
 \psi(x) =  \minamm \big( \theta^{-1}(x) \big) = \minamm \big( s \big( \theta^{-1}(x+1) \big) =  \minamm \big( \theta^{-1}(x+1) \big)  = \psi(x+1 )
$$ 
where the respective equalities are due to the definition of $\psi$; the above claim with $k=1$; Theorem~\ref{t.relativereward}(1); and the definition of~$\psi$ once more. We have obtained Theorem~\ref{t.phithetainverse}(1).

Note next that $\mc{S}_k \solution(x+k)$ equals
\begin{eqnarray*}
  &  & \mc{S}_k \Big( \asta\big( \theta^{-1}(x+k)  \big),  \bst\big( \theta^{-1}(x+k)  \big),  \mst\big( \theta^{-1}(x+k)  \big),  \nst\big( \theta^{-1}(x+k)  \big) \Big) \\
 & =  & \bigg( \asta \Big( s_k\big(\theta^{-1}(x+k) \big) \Big) ,  \bst \Big( s_k\big(\theta^{-1}(x+k) \big) \Big) , \mst \Big( s_k\big(\theta^{-1}(x+k) \big) \Big) , \nst \Big( s_k\big(\theta^{-1}(x+k) \big) \Big)  \bigg) \, ,
\end{eqnarray*}
the latter equality by Proposition~\ref{p.shift}. Applying~(\ref{e.claimsk}), we find that
 $$
 \mc{S}_k \solution(x+k)  =  \Big( \asta\big( \theta^{-1}(x)  \big),  \bst\big( \theta^{-1}(x)  \big),  \mst\big( \theta^{-1}(x)  \big),  \nst\big( \theta^{-1}(x)  \big) \Big)
 = \solution(x) \, .
 $$
 This implies that  $\solution(x+k) = 
 \mc{S}_{-k} \solution(x)$, which is what Theorem~\ref{t.phithetainverse}(3) asserts.

{\bf  (3).}
Let $z \in \R$ and set $\PhiMac(z) = x$.
Since $\vert \PhiMac(x) - \theta(x) \vert = O(1)$ by Lemma~\ref{l.thetapsi}, and $\PhiMac$ and $\theta$ are increasing, we have that
$$
 \PhiMac\big(z - O(1)\big)
 \leq \theta^{-1}(z) \leq  \PhiMac\big(z + O(1)\big) \, .
$$
Substituting the expressions for $2^{2^z - 1}$ and $2^{-(2^{-z} -1)}$ for $\PhiMac(z)$, valid when $z \geq 0$ and $z < 0$, we obtain Theorem~\ref{t.phithetainverse}(3). \qed

\section{The Mina margin map is not identically equal to one}\label{s.minamarginmap}
Here we prove Theorem~\ref{t.minamarginvalues}(3). We will obtain evidence for Conjecture~\ref{c.lambda} as we do so. 
\begin{proposition}\label{p.thevalueofminamargin}
The value of $\minamm_{5,4}(0.58)$ lies in the interval  $[0.9999032032 , 0.9999032038]$.
\end{proposition}

{\bf Proof of Theorem~\ref{t.minamarginvalues}(3).}  Since $\lambda$ is equal to the infimum of $\minamm(x)$ over $x \in (0,\infty)$, we have that $\lambda \leq \minamm(0.58)$.

Note that Proposition~\ref{p.rkrell} with $(\ell,k) = (5,4)$ implies that the value $\minamm(z) = \lim_j \minamm_{j,j}(z)$ (where $z = 0.58 \in [1/3,3]$)  satisfies
$$
 \big\vert  \minamm(z) - \minamm_{5,4}(z) \big\vert \, \leq \, 3.3 \times 10^{-8} + 5.95 \times 10^{-7} \, \leq \, 6.3 \times 10^{-7} \, .
$$
Applying the upper bound on $\minamm_{5,4}(z)$ in Proposition~\ref{p.thevalueofminamargin}, we find that
$$
 \minamm(0.58) \, \leq \, 0.9999032038 + 6.3 \times 10^{-7} \, = \, 0.9999038338  \, .
$$
We confirm then that $\lambda$, being at most $\minamm(0.58)$, is bounded above by $0.999904$. This completes the proof of Theorem~\ref{t.minamarginvalues}(3). \qed

Numerical work with Mathematica indicates that $\minamm_{5,4}(0.5809)$ equals  $0.999903202726$
to twelve decimal places; that $\minamm_{5,4}(0.5809)$ equals
 $\min \big\{ \minamm_{5,4}(x): x \in [1/3,3] \cap 10^{-4}\Z \big\}$; and that the error between $\inf \big\{ \minamm_{5,4}(x): x \in [1/3,3] \big\}$ and this minimum may jeopardise only the final digit of the twelve. If this evidence is admitted\footnote{Neo Lee and Adam Ousterovitch, U.C. Berkeley undergraduates, are working on putting such evidence on a rigorous footing: a forthcoming article will validate Conjecture~\ref{c.lambda}.}, then the preceding proof yields
 that $\inf \big\{ \minamm(x): x \in [1/3,3] \big\}$ is at least  $\minamm_{5,4}(0.5809) -10^{-11} -  6.3 \times 10^{-7} \geq 0.99990257 \geq 0.999902$; whence
  Conjecture~\ref{c.lambda}. In fact, the conjecture is cautious: $\lambda$ is likely to exceed $0.999903$, as an estimate on a higher indexed $\minamm_{\ell,k}$ might show.
 
The formula for  $(0,\infty) \to \R: x \to \minamm_{5,4}(x)$ in Lemma~\ref{l.ratiointerpret} may be recorded explicitly---it involves several applications of such operations as inverse and square-root---but it is messy, and would occupy several pages of standard print. Arguably a claim that mathematical software evaluates this function at $0.58$ to be within the range claimed by Proposition~\ref{p.thevalueofminamargin} may be admitted as a proof of this result. But a diligent reader who is given this information has no practical way to confirm it. In the following proof, we provide an approximation scheme, from above and below, for computing $\minamm_{5,4}(0.58)$. All quantities in the scheme are values in $10^{-10}\Z$, and the proof is reduced to verifying about fifty explicit statements of the form `if $x=u$, then $f(x)=v$', where $u$ and $v$ are given elements of $10^{-10}\Z$, and $f$ is the application of a function such as $s$, $c$ and $d$ from Definition~\ref{d.acs}
followed by a rounding down or up on to the lattice $10^{-10}\Z$. In this way, the diligent reader has a mundane but manageable task to verify every detail of the derivation of Proposition~\ref{p.thevalueofminamargin}. 

We note that, were $\minamm_{5,4}$ shown to be differentiable, and a suitable bound on its derivative found, then a similarly explicit record of the values of $\minamm_{5,4}$ on a fine enough mesh of points in $[1/3,3]$
would furnish a proof of Conjecture~\ref{c.lambda}. If the number of points in the mesh were large, then a manual check on the explicit bounds would be impracticable, so that in such a case the proof would be at least modestly computer-assisted.

We now turn to introducing and implementing the approximation scheme.
Let  $k \in \nwithzero$, and set 
$$
 \lfloor x \rfloor_k  = 10^{-k} \lfloor 10^k x \rfloor  \, \, \, \, \textrm{and} \, \, \, \,  \ulcorner x \urcorner^k =  10^{-k} \lfloor 10^k x \rfloor  + 10^{-k} \in \R \, .
$$
Namely, the real line is partitioned 
$$
\R \, = \, \bigcup_{j \in \Z} \, 10^{-k} \cdot [j,j+1)
$$
into intervals whose endpoints are consecutive elements in the lattice $10^{-k}\Z\,$; 
 $\big[  \lfloor  x \rfloor_k , 
 \lceil x \rceil^k \big)$ is the unique interval in the partition that contains~$x$.
 From the outset, we set the parameter $k$ equal to ten, and omit to denote it. 
It should thus be understood that  $\lfloor  x \rfloor$ and  
 $\lceil x \rceil$ denote 
 $\lfloor  x \rfloor_{10}$ 
 and $\lceil x \rceil^{10}$, rather than the usual integer roundings of $x$. 
 
 Recall Definition~\ref{d.acs}. 
 We specify $s^\uparrow,s^\downarrow,c^\uparrow,c^\downarrow,d^\uparrow,d^\downarrow:(0,\infty) \to (0,\infty)$
 by setting 
 $$
 *^\uparrow(x) = \lceil *(x) \rceil \, \, \, \textrm{and} \, \, \,   *^\downarrow(x) = \lfloor *(x) \rfloor \, \, \, \textrm{for} \, \, \,  * \in \{s,c,d\} \, .
 $$
 For $x \in (0,\infty)$, we specify  $\big\{ s^\uparrow_i(x): i \in \Z \big\}$ and  $\big\{ s^\downarrow_i(x): i \in \Z \big\}$, the upper and lower $s$-sequences evaluated at $x$.
 Indeed, we set $s_0^\uparrow(x) = s_0^\downarrow(x) = x$.
 For $i \geq 1$, we then iteratively set $s^\uparrow_i(x) = s^\uparrow \big(s^\uparrow_{i-1}(x)\big)$ and  $s^\downarrow_i(x) = s^\downarrow \big(s^\downarrow_{i-1}(x)\big)$.
 We further set $s_{-1}^\uparrow (x) = \lceil s_{-1}(x) \rceil$ and  $s_{-1}^\downarrow (x) = \lfloor s_{-1}(x) \rfloor$ . For $i \leq -2$, we iteratively set $s_{-i}^\uparrow(x) = s_{-1}^\uparrow \big( s_{-i+1}^\uparrow(x) \big)$ and  $s_{-i}^\downarrow(x) = s_{-1}^\downarrow \big( s_{-i+1}^\downarrow(x) \big)$.

 Set $z =0.58$. 
 
 We will write $s_i^\uparrow = s_i^\uparrow(z)$ and  $s_i^\downarrow = s_i^\downarrow(z)$ for $i \in \Z$. In this way, the value of $z = 0.58$ is understood.
 
 We further define the upper and lower $c$- and $d$-sequences, $\big\{c^\uparrow_i,c^\downarrow_i,d^\uparrow_i,d^\downarrow_i: i \in \Z \big\}$, where again the value of $z$ is understood.
  We set $c_i^\uparrow = c^\uparrow(s_i^\uparrow)$,  $c_i^\downarrow = c^\downarrow(s_i^\downarrow)$,  $d_i^\uparrow = d^\uparrow(s_i^\uparrow)$ and  $d_i^\downarrow = d^\downarrow(s_i^\downarrow)$.
  
  The data $\big\{ s^\uparrow_i,s^\downarrow_i,c^\uparrow_i,c^\downarrow_i,d^\uparrow_i,d^\downarrow_i \big\}$, $i \in \llbracket -4,3 \rrbracket$,
  are forty-eight elements of $10^{-10}\Z$. These values are presented in Tables~\ref{t.one} and~\ref{t.two}. Two of the values are known without computation: $s_0^\uparrow = s_0^\downarrow =0.58$.
  The remaining values may be computed, one at a time, where each  step is a computation ${\rm INPUT} \rightarrow {\rm OUTPUT}$
of one element of the lattice $10^{-10}\Z$ from another. Each step takes 
one of the forms $s_i^\uparrow \to s^\uparrow_{i+1}$ for $i \in \llbracket 0,2 \rrbracket$;  
$s_i^\uparrow \to s^\uparrow_{i-1}$ 
for $i \in \llbracket -3,0 \rrbracket$;   
$s_i^\uparrow \to  c_i^\uparrow$ or $s_i^\uparrow \to  d_i^\uparrow$ for $i \in \llbracket -4,3 \rrbracket$;
or it is formed by replacing $\uparrow$ by $\downarrow$ in one of these steps. Forty-six such steps lead to the completion of the two tables, given the two initial entries.
 \begin{table}
 \begin{center}
 \begin{tabular}{| c | c | c |}
 \hline
 $i$ & $s_i^\uparrow$ & $s_i^\downarrow$  \\
 \hline
  -4 & 954911606.03  & 954911605.92    \\ 
 -3 & 21848.5122538904 & 21848.5122525938  \\
 -2 & 102.3071054647 &102.3071054616\\
 -1 & 5.3556473847 & 5.3556473846  \\
0 & 0.5800000000  & 0.5800000000 \\
1 & 0.0504077253 & 0.0504077252  \\
2 & 0.0010408205 & 0.0010408204  \\
3 & 0.0000005392  & 0.0000005391  
 \\
 \hline
 \end{tabular}
 \caption{Values of  $s_i^\uparrow$ and $s_i^\downarrow$ for $i \in \llbracket -4,3 \rrbracket$}\label{t.one}
 \end{center}
 \end{table}

 \begin{table}
 \begin{center}
 \begin{tabular}{| c | c | c | c | c |}
 \hline
 $i$ & $c_i^\uparrow$ & $c_i^\downarrow$ & $d_i^\uparrow$ & $d_i^\downarrow$ \\
 \hline
  -4 & 477488579.78   & 477488579.73  & 10926.0060411432 & 10926.0060404948 \\ 
 -3 & 11081.6603248978 & 11081.6603242447 & 52.8859257466 & 52.8859257450 \\
 -2 & 62.5133614707 & 62.5133614689 & 4.2201465577 & 4.2201465576 \\
 -1 & 5.7859121540 & 5.7859121538 & 1.5182994418 & 1.5182994417 \\
0 & 1.8055756566 & 1.8055756565 & 1.0700124766 & 1.0700124765 \\
1 & 1.0944264319 & 1.0944264316 & 1.0019497202 & 1.0019497201 \\
2 & 1.0020784046 & 1.0020784043 & 1.0000010767 & 1.0000010766 \\
3 & 1.0000010785 & 1.0000010782 & 1.0000000001 & 1.0000000000  \\
 \hline
 \end{tabular}
 \caption{Values of  $c_i^\uparrow$, $c_i^\downarrow$, $d_i^\uparrow$ and $d_i^\downarrow$  for $i \in \llbracket -4,3 \rrbracket$}\label{t.two}
 \end{center}
 \end{table}

According to Lemma~\ref{l.ratiointerpret},
 $$
 \minamm_{5,4}(z) \, = \, \frac{z(S_4+T_5)}{P_4 + Q_5} \, ,
 $$
 where 
 $$
 P_4 = 1 + (c_0 - 1) +  (c_0 - 1)  (c_1 - 1) +  (c_0 - 1)  (c_1 - 1)  (c_2 - 1) +  (c_0 - 1)  (c_1 - 1)  (c_2 - 1)  (c_3 - 1) \, ; 
 $$
 \begin{eqnarray*}
 Q_5 & = & \big( c_{-1} -1 \big)^{-1} +  \big( c_{-1} -1 \big)^{-1}  \big( c_{-2} -1 \big)^{-1} +  \big( c_{-1} -1 \big)^{-1}  \big( c_{-2} -1 \big)^{-1}  \big( c_{-3} -1 \big)^{-1} \\
  & & \qquad \qquad \qquad  + \, \big( c_{-1} -1 \big)^{-1}  \big( c_{-2} -1 \big)^{-1}  \big( c_{-3} -1 \big)^{-1}  \big( c_{-4} -1 \big)^{-1} \, ;
 \end{eqnarray*}
 $$
 S_4 = 1 + (d_0 - 1) +  (d_0 - 1)  (d_1 - 1) +  (d_0 - 1)  (d_1 - 1)  (d_2 - 1) +  (d_0 - 1)  (d_1 - 1)  (d_2 - 1)  (d_3 - 1) \, ; 
 $$
 and
 \begin{eqnarray*}
 T_5 & = & \big( d_{-1} -1 \big)^{-1} +  \big( d_{-1} -1 \big)^{-1}  \big( d_{-2} -1 \big)^{-1} +  \big( d_{-1} -1 \big)^{-1}  \big( d_{-2} -1 \big)^{-1}  \big( d_{-3} -1 \big)^{-1} \\
    & & \qquad \qquad \qquad +  \, \big( d_{-1} -1 \big)^{-1}  \big( d_{-2} -1 \big)^{-1}  \big( d_{-3} -1 \big)^{-1}  \big( d_{-4} -1 \big)^{-1} \, .
 \end{eqnarray*}
 We further specify quantities $*^\uparrow$ and $*^\downarrow$, where $* \in \{ P_4, Q_5, S_4, T_5 \}$.
 To do so, we record variable dependence in the form $P_4 = P_4(c_0,c_1,c_2,c_3)$,  $S_4 = S_4(d_0,d_1,d_2,d_3)$, $Q_5 = Q_5(c_{-1},c_{-2},c_{-3},c_{-4})$
  and $T_5 = T_5(d_{-1},d_{-2},d_{-3},d_{-4})$. We may then set
  $$
  P_4^\uparrow = \lceil P_4(c^\uparrow_0,c^\uparrow_1,c^\uparrow_2,c^\uparrow_3) \rceil \, \, \, \, \textrm{and} \, \, \, \,
  P_4^\downarrow = \lfloor P_4(c^\downarrow_0,c^\downarrow_1,c^\downarrow_2,c^\downarrow_3) \rfloor \, ;
  $$
  $$
  S_4^\uparrow =  \lceil  S_4(d^\uparrow_0,d^\uparrow_1,d^\uparrow_2,d^\uparrow_3)  \rceil \, \, \, \, \textrm{and} \, \, \, \,
  S_4^\downarrow =  \lfloor S_4(d^\downarrow_0,d^\downarrow_1,d^\downarrow_2,d^\downarrow_3) \rfloor \, ;
  $$
  $$
  Q_5^\uparrow =  \lceil  Q_5(c^\downarrow_{-1},c^\downarrow_{-2},c^\downarrow_{-3},c^\downarrow_{-4})  \rceil \, \, \, \, \textrm{and} \, \, \, \,
  Q_5^\downarrow =   \lfloor Q_5(c^\uparrow_{-1},c^\uparrow_{-2},c^\uparrow_{-3},c^\uparrow_{-4}) \rfloor \, ;
  $$
  and
  $$
  T_5^\uparrow =  \lceil  T_5(c^\downarrow_{-1},c^\downarrow_{-2},c^\downarrow_{-3},c^\downarrow_{-4})  \rceil \, \, \, \, \textrm{and} \, \, \, \,
  T_5^\downarrow =  \lfloor T_5(c^\uparrow_{-1},c^\uparrow_{-2},c^\uparrow_{-3},c^\uparrow_{-4}) \rfloor \, .
  $$
 (Note the reversals of the uses of $\downarrow$ and $\uparrow$ in the replaced terms for $Q_5$ and $T_5$.)

The tables then permit us to record the values (all of which are elements of the lattice $10^{-10}\Z$)
\begin{eqnarray}
 S_4^\uparrow & = & 1.0701489815 \, \, \, , \, \, \, S_4^\downarrow =  1.0701489813 \label{e.stpqvalues} \\
 T_5^\uparrow& = & 2.5400964392  \, \, \, , \, \, \,
  T_5^\downarrow = 2.5400964386 \nonumber \\
  P_4^\uparrow & = &  1.8818013910  \, \, \, , \, \, \,
   P_4^\downarrow =1.8818013906  \nonumber \\
   Q_5^\uparrow & = &  0.2123436589  \, \, \, , \, \, \,
   Q_5^\downarrow =   0.2123436587 \, .  \nonumber
\end{eqnarray}

 Next we specify two further elements of $10^{-10}\Z$:
  \begin{equation}\label{e.rupdown}
 \minamm^\uparrow_{5,4}(z) \, = \,  \biggl\lceil \,  \frac{z(S^\uparrow_4+T^\uparrow_5)}{P^\downarrow_4 + Q^\downarrow_5}  \, \biggr\rceil  \, \, \, \, \textrm{and} \, \, \, \, 
 \minamm^\downarrow_{5,4}(z) \, = \,   \biggl\lfloor \, \frac{z(S^\downarrow_4+T^\downarrow_5)}{P^\uparrow_4 + Q^\uparrow_5} \, \biggr\rfloor \, .
 \end{equation}
  
 \begin{lemma}\label{l.fiveshort}
 \begin{enumerate}
 Let $i \in \Z$.
 \item $s_i^\downarrow \leq s_i \leq s_i^\uparrow$.
 \item $c_i^\downarrow \leq c_i \leq c_i^\uparrow$.
 \item $d_i^\downarrow \leq d_i \leq d_i^\uparrow$.
 \item $P_4 \downarrow \leq P_4 \leq P_4^\uparrow$,  $Q_4 \downarrow \leq Q_4 \leq Q_4^\uparrow$, $S_4 \downarrow \leq S_4 \leq S_4^\uparrow$ and $T_4 \downarrow \leq T_4 \leq T_4^\uparrow$.
 \item $\minamm_{5,4}^\downarrow \leq \minamm_{5,4} \leq \minamm_{5,4}^\uparrow$.
 \end{enumerate}
 \end{lemma}
{\bf Proof.} Note that $s^\downarrow(x) \leq s(x) \leq s^\uparrow(x)$ for $x \in (0,\infty)$ by the definitions of $s^\downarrow$ and $s^\uparrow$. By induction on $i \geq 1$,
we will show that $s_i^\uparrow \geq s_i$. Indeed, note that $s_i^\uparrow = s^\uparrow(s_{i-1}^\uparrow) \geq s(s_{i-1}^\uparrow) \geq s(s_{i-1}) = s_i$, where the latter inequality is due to the inductive hypothesis at index $i-1$ and to Lemma~\ref{l.acsfacts}(1:$s$).
We also prove that $s_{-i}^\uparrow \geq s_{-i}$ for $i \geq 1$ by induction on $i$. In this regard, note that $s_{-i-1} = s_{-1}^\uparrow(s_{-i}^\uparrow) \geq s_{-1}(s_{-i}^\uparrow) \geq s_{-1}(s_{-i}) = s_{-1-i}$, where the first bound is due to $s_{-1}^\uparrow(x) \geq s_{-1}(x)$ for $x \in (0,\infty)$, which follows from the definition of $s_{-1}^\uparrow$; the second  is due to the inductive hypothesis at index $i$ and $x \to s_{-1}(x)$ being increasing, which fact follows from  Proposition~\ref{p.sminusone} and Lemma~\ref{l.acsfacts}(1:$s$).

Similar arguments prove that $s_i^\downarrow \leq s_i$ for $i \in \Z$.

{\bf (2).} Note that $c^\uparrow(s_i^\uparrow) \geq c(s_i^\uparrow) \geq c(s_i) = c_i$, where the first bound is due to the definition of $c^\uparrow$ and the second is due to $s_i^\uparrow \geq s_i$ and Lemma~\ref{l.acsfacts}(1:$c$). Similarly may we show that $c_i^\downarrow \leq c_i$.

{\bf (3).} This is similar to the preceding part.

{\bf (4).} Note that $P_4$ is an increasing function of the variables $c_i$, $i \in \llbracket 0,3 \rrbracket$; 
$Q_5$ is decreasing in $c_{-i}$, $i \in \intint{4}$; $S_4$ is increasing in  $d_i$, $i \in \llbracket 0,3 \rrbracket$; and $T_5$ is decreasing in $d_{-i}$, $i \in \intint{4}$.
(The noted properties of $Q_5$ and $T_5$ are valid only insofar as the variables $c_{-i}$ and $d_{-i}$ remain greater than one. But this condition is always met in applications, including the present one.)
Given these monotonicities, Lemma~\ref{l.fiveshort}(2) shows that
$$
 P_4(c^\downarrow_0,c^\downarrow_1,c^\downarrow_2,c^\downarrow_3) \leq  P_4  \leq  P_4(c^\uparrow_0,c^\uparrow_1,c^\uparrow_2,c^\uparrow_3) \, ,
$$
so that the monotonicities of $\lfloor \cdot \rfloor$ and $\lceil \cdot \rceil$ prove the assertions concerning $P_4$.
The derivation for $Q_5$ is similar. So are the others: Lemma~\ref{l.fiveshort}(3) is used in regard to $S_4$; and Lemma~\ref{l.fiveshort}(4) for $T_5$.

{\bf (5).} The expression $\minamm_{5,4}$ is an increasing function of $S_4$ and $T_5$, and it is decreasing in $P_4$ and $Q_5$---thus, we may use the preceding part to reach the desired conclusion. \qed

{\bf Proof of Proposition~\ref{p.thevalueofminamargin}.}
Using the data~(\ref{e.stpqvalues}), note that the expressions~(\ref{e.rupdown}) have the evaluations
$$
 \minamm^\uparrow_{5,4} \, = \,  \biggl\lceil \, 0.58 \times \frac{1.0701489815 + 2.5400964392}{1.8818013906 + 0.2123436587} \,  \biggr\rceil \, = \, 0.9999032038 
 $$
 and
 $$
   \minamm^\downarrow_{5,4} \, = \, \biggl\lfloor \, 0.58 \times \frac{1.0701489813 + 2.5400964386}{1.8818013910 + 0.2123436589} \,  \biggr\rfloor \, = \, 0.9999032032 \, .
   $$
By Lemma~\ref{l.fiveshort}(5), we learn that
$\minamm_{5,4} \in [0.9999032032 , 0.9999032038]$, as Proposition~\ref{p.thevalueofminamargin} states. \qed

\chapter{Trail prospects}\label{c.prospects}

We discuss five topics prompted by the article.

\section{Properties of the Mina margin map and prospective routes to conjectures}\label{s.conjectureroute}

Conjecture~\ref{c.solutions} concerns the level sets of the Mina margin map, and via~(\ref{e.finitenash}), 
Conjecture ~\ref{c.tine} concerns the level sets of the finite trail counterparts to this map. Consider the $\PhiMac$-transformed finite-trail Mina margin maps $\minamm_{j+1,k+1} \circ \PhiMac: \R \to (0,\infty)$
depicted for several pairs~$(j+1,k+1)$ in Figure~\ref{f.tmmm}. In these sketches, there are a total of $2(j+k) - 5$ elements in any level set through which every swerve of the function passes; 
such level sets are indexed by $[\lambda,\lambda^{-1}]$ for $\lambda =0.999903 \cdots$ up to an error that vanishes in high $j$ and $k$; the functions converge to a limit $\minamm \circ \PhiMac: \R \to (0,\infty)$, and this limit has level sets with two elements in each period (such as in $\PhiMac^{-1}(1/3,3]$) for heights in $(\lambda,\lambda^{-1})$. These claims constitute the content of the two conjectures and they can be said to be visually more-or-less evident. But can they be proved? 
In regard to Conjecture~\ref{c.solutions} at least, control on derivatives and explicit evaluation on a suitably fine mesh may be a tractable approach: see the discussion regarding Conjecture~\ref{c.lambda} in Section~\ref{s.minamarginmap}.


\section{The possible existence of further Nash equilibria}

We have studied time-invariant Nash equilibria. It is natural to ask whether further Nash equilibria exist. We discuss two directions.

\subsection{Time-invariant random Nash equilibria}

Our formulation of the notion of Nash equilibrium in Section~\ref{s.gamespec} is deterministic. What if time-invariant random play is permitted? A strategy would consist of a set of laws on the non-negative reals indexed by $\Z$. When such a strategy is played, the stake offered would be sampled from the law indexed by the present counter location, the sampling being independent of other randomness. 
To avoid extra notation, we have not formulated this notion in the main part of the article. We do not believe that non-trivial random time-invariant Nash equilibria exist. Indeed, we remarked after the Penny Forfeit Lemma~\ref{l.pennyforfeit} that random play is suboptimal for the one-step game. By iterating this result and invoking the monotonicity 
in Penny Forfeit argued in the proof of Lemma~\ref{l.onestep}(2),
 the possibility of a non-trivial role for randomness of the form we have discussed may be excluded.

\subsection{Nash equilibria that are not time-invariant}

A deterministic strategy pair that may not be time-invariant takes the form $(b,a):\Z \times \nwithoutzero \to (0,\infty)^2$. We may anticipate that, were such a pair a Nash equilibrium, the naturally associated dynamical quadruple $(a,b,m,n)$, specified by suitably modifying Definition~\ref{d.quadruple}, would satisfy a dynamical form \dabmnmacspace of 
the ABMN system on $\Z$. For simplicity, we describe these equations on a finite trail $\llbracket -K-1,K+1\rrbracket$
and for a finite time interval $\llbracket 0, T \rrbracket$ (so that $K+1,T \in \nwithoutzero$). Boundary data is a quadruple $(m_{-K},m_K,n_K,n_{-K}) \in \R^4$ which equals $(0,1,1,0)$
in the simple symmetric case; and two terminal functions $m_{\rm ter},n_{\rm ter}: \llbracket -K-1,K+1 \rrbracket \to \R$.  
If write $*_i(j) = *(i,j)$ for $* \in \{a,b,m,n\}$, so that, for example, $a_i(j)$ is the stake offered by Maxine at the $j$\textsuperscript{th}
turn in the event that $X_{j-1} =i$, the revised equations are  
\begin{align*}
  \big( a_i(j) + b_i(j) \big)\big(m_i(j) + a_i(j) \big)  & =  a_i(j) m_{i+1}(j+1) + b_i(j) m_{i-1}(j+1) && \qquad \textrm{dABMN}(1) \\
   \big(a_i(j) + b_i(j) \big) \big(n_i(j)+b_i(j) \big)  & =  a_i(j) n_{i+1}(j+1) + b_i(j) n_{i-1}(j+1)  &&\qquad \textrm{dABMN}(2)   \\
   \big(a_i(j) + b_i(j) \big)^2  & =  b_i(j) \big( m_{i+1}(j+1) - m_{i-1}(j+1) \big) &&\qquad \textrm{dABMN}(3)   \\
       \big(a_i(j) + b_i(j) \big)^2  & =   a_i(j) \big( n_{i-1}(j+1) - n_{i+1}(j+1) \big) &&\qquad \textrm{dABMN}(4)   \, ,
  \end{align*}
   where $(i,j)$ ranges over $\llbracket -K,K \rrbracket \times \llbracket 0, T-1 \rrbracket$. Boundary conditions enter via
 \begin{eqnarray*}
   & &  m_{\pm (K+1)}(j) = m_{\pm  (K+1)} \, \, , \, \,   n_{\pm  (K+1)}(j) = n_{\pm  (K+1)} \,  \textrm{for $j \in \intint{K}$; and} \\
   & &  m_i(T) = n_{\rm term}(i) \, \, , \, \, n_i(T) = n_{\rm term}(i)  \, \,  \textrm{for $i \in \llbracket{-K,K \rrbracket}$} \, .
  \end{eqnarray*}
   Of course, any \abmnmacspace solution $(a,b,m,n)$ solves $\dabmnmac$ if we extend notation to set $*_i(j) = *_i$ for $(i,j) \in \Z \times \nwithoutzero$ and $* \in \{a,b,m,n\}$ (and then restrict the domain suitably). Do other solutions of $\dabmnmac$ exist?

Certainly there are some such. Points~$(i,j)$ in $\Z \times \nwithoutzero$ are odd or even according to whether $i+j$ is odd or even. The parity of $j + X_i(j)$ never changes from its initial $j=0$ value in any instance of the trail game. If we select two solutions $(a',b',m',n')$ and $(\hat{a},\hat{b},\hat{m},\hat{n})$ of the ABMN equations, and set
 $$
  (a,b,m,n)(i,j) \, = \, \begin{cases}
  \, \, (a'_i,b_i',m'_i,n'_i)  &  \text{when $i+j$ is even} \, , \\
 \, \,  (\hat{a}_i,\hat{b}_i,\hat{m}_i,\hat{n}_i)  &  \text{when $i+j$ is odd}  
  \, ,
\end{cases}
$$
then $(a,b,m,n)$ solves $\dabmnmac$ (when suitably restricted in the domain) and
 Theorem~\ref{t.nashabmn} directly implies that $(b,a)$ is a Nash equilibrium in the trail game. Conceptually, this is not really a new solution, however. Gameplay resides on the odd or even lattice and use of any of these new Nash equilibria will coincide with that of a time-invariant Nash equilibrium in any given game.

\begin{figure}[htbp]
\centering
\includegraphics[width=0.9\textwidth]{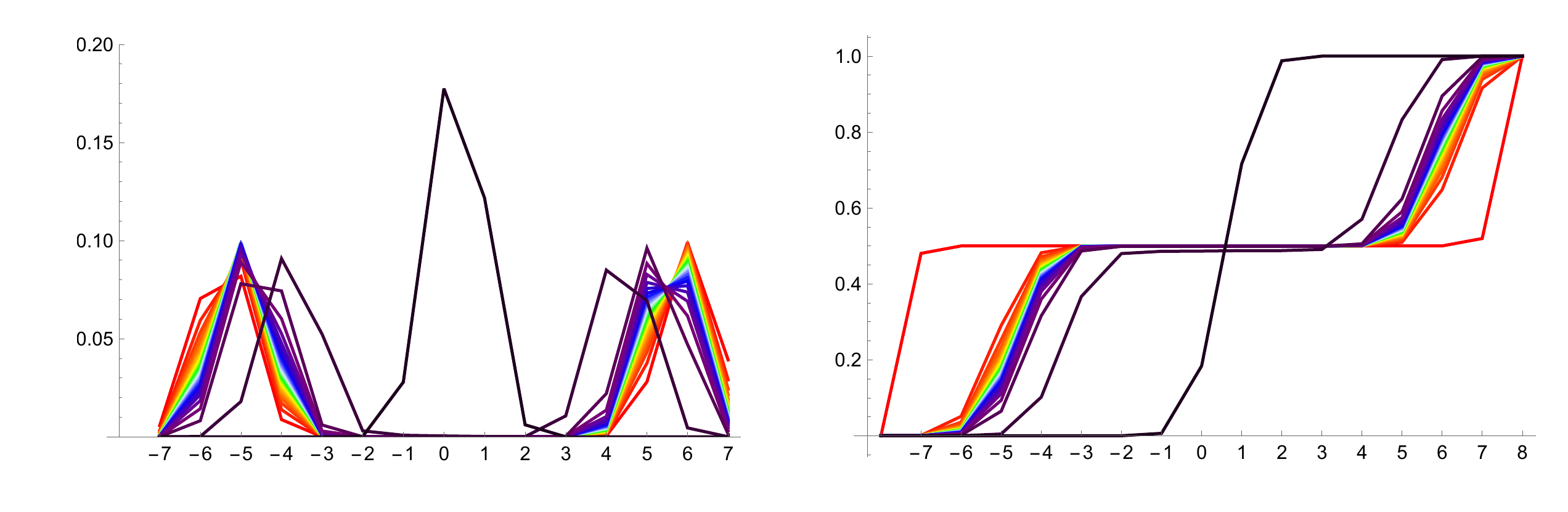}
\caption{The dynamic ABMN equations $\dabmnmac$ on trail $\llbracket -8,8\rrbracket$ and time interval $\llbracket 0,4200 \rrbracket$ are solved with $m_{\rm ter}:\llbracket -8,8 \rrbracket \to [0,1]$, $m_{\rm ter}(-8)=0$, $m_{\rm ter}(8)=1$, rising sharply from zero to run along a rough plateau at height one-half, and ending with a further sharp rise to height one. We set $n_{\rm ter}(\cdot) = m_{\rm ter}(-\cdot)$
and work with a standard symmetric boundary quadruple~$(0,1,1,0)$.  The light grey curve on the right plot, which is most exposed both on the left and the right, is $m_{\rm ter}$. The solutions of the equations are depicted at values of $j$ in $\llbracket 0,4200 \rrbracket$ that are multiples of $140$, so that thirty curves excepting the final condition are depicted in each plot. On the left, the $a$-components of $\dabmnmac$ on the open-play set $\llbracket -7,7 \rrbracket$ for the $j$-values in question are shown (with linear interpolation between integers); on the right, the $m$-components on the trail $\llbracket -8,8\rrbracket$ with such interpolation. The curves are coloured on a bluish spectrum darkening from light grey to black as time passes backwards. These curves make a staccato advance (with this flow of time) from the sides to the centre, with the final black curve in each plot, indexed by $j=0$, representing a single battlefield around the origin. The $b$- and $n$-components are formed by reflecting the $a$- and $m$-components in the vertical axis.}\label{f.dynamicabmn}
\end{figure}

For $k \in \nwithoutzero$, the system $\dabmnmac$ can be solved on $\llbracket -K,K \rrbracket \times \intint{k}$ by choosing a given terminal condition $\big\{ m_i(k),n_i(k) \big): i \in \llbracket -K-1,K+1 \rrbracket  \big\}$
and iteratively solving $\dabmnmac$ for decreasing values of $j$. In searching for a solution that is not invariant in time, we seek a  terminal condition such that,
if this condition is imposed even for a very high value of $k$, the solved solution stabilises for bounded values of $j$ to a form that is not time-invariant. We have tested a few terminal conditions; Figure~\ref{f.dynamicabmn} depicts a solution of $\dabmnmac$ on the trail $\llbracket -8,8 \rrbracket$. The $m$-component of the terminal condition, defined on this trail, rises sharply at both ends, and otherwise has the form of a rough plateau of height one-half; the $n$-component is the reflection of the $m$-component in the vertical axis. The region of each sharp rise for the $m$-component represents a battlefield which may be rather stable as time evolves. As the two plots in Figure~\ref{f.dynamicabmn} show, the two battlefields rapidly approach one another by a short distance, and then remain in a rough stasis in which a gradual movement towards the origin is discernible, before rapidly breaking towards a single central battlefield.
A total of $4200$ time-steps are involved in the simulation, with the $j=140$ black curve in each plot (which is the penultimate depicted in the backwards-in-time evolution) showing an eruption towards the centre, and the final black curve in each plot adopting the central location which is in essence the fixed point of the evolution. Battlefield pairs with greater separation may endure much longer, and may present a metastable effect for $\dabmnmac$ that causes these equations to converge very slowly to fixed points. There is no evidence however of non-convergence: our limited investigation has not produced examples that  support the notion that further time-invariant Nash equilibria exist beyond the simple parity-based class discussed above.

\section{Gameboards beyond $\Z$}\label{s.beyond}
By use of a setup involving directed graphs, self-funded stake-governed random-turn games derived from games such as Hex or Go may be considered. It would be of interest to determine for a suitable class of games whether some of the features of the Trail of Lost Pennies of $\Z$ are present more generally. The central ratio $\tfrac{n_{-1} -n_0}{m_0 - m_{-1}}$ is the ratio of changes in mean payoff for Mina and Maxine arising from Mina's victory at the first turn. This or similar quantities may be considered in suitable infinite games, permitting us to ask whether 
Theorem~\ref{t.nashequil.prelim} generalizes to these games: do Nash equilibria exist precisely when the quantity lies in an interval of the form $[\lambda,\lambda^{-1}]$? Do these game-determined $\lambda$-values differ from one by a notably small but positive quantity, as this value for the trail game on~$\Z$ appears to differ by about $10^{-4}$? Do more general games exemplify as ours does the Battlefield Cyl Fog premise from Section~\ref{s.allocate}, having one (or perhaps several) bounded battlefield regions of the space of configurations on the gameboard, specified by any given Nash equilibrium, in which players concentrate their stake expenditure, with the outcome of turns occurring therein being highly influential in the overall game? 
As we recalled in Section~\ref{s.gamespde}, Peres and \v{S}uni\'c~\cite{PeresSunic} found a path decomposition that leads to rapid computation of $\lambda$-biased infinity harmonic functions. To address the last question, it may be natural at first to work with trees with suitable boundary data, and to seek to iteratively determine Nash equilibria for the player-funded stake-governed tug-of-war along elements in a suitable path decomposition. In this way, solutions for this game in a broader class of graphs
may decompose into components, these being trail games as we have analysed.

\section{The trail game and ``Incentive Inch, Outcome Mile"}\label{s.playing}

We have suggested that the Trail of Lost Pennies advances the premise ``Incentive Inch, Outcome Mile" presented in Section~\ref{s.allocate}. Here we examine this point more closely, for the game on finite intervals and on $\Z$. 
The premise for the trail game on these graphs may be summarised as follows.
\begin{enumerate}
\item When the Mina margin exceeds $\lambda^{-1} + o(1)$, Mina is the dominant player, and is able to secure victory at little cost from any starting position (except perhaps a few near the right boundary).
\item When the Mina margin is less than $\lambda - o(1)$, it is Maxine who is dominant, and is likewise able (except near the left boundary).
\item In a short interval of Mina margin values of length of order~$10^{-4}$ around one, namely $[\lambda + o(1),\lambda^{-1} - o(1)]$, there are many equilibria, and the dominance of one or other player as above arises from any starting vertex (except a few near the boundary) under a suitable equilibrium.
\end{enumerate}
Here, $o(1)$ refers to a positive expression that vanishes in the limit of long finite intervals, and is zero in the case of~$\Z$. The qualifications in parentheses should be omitted in the $\Z$-case.
By `many' in the third part, we mean $2n + \Theta(1)$ for a finite interval of length $n \geq 6$, and countably many for~$\Z$.

\begin{figure}[htbp]
\includegraphics[width=1.0\textwidth]{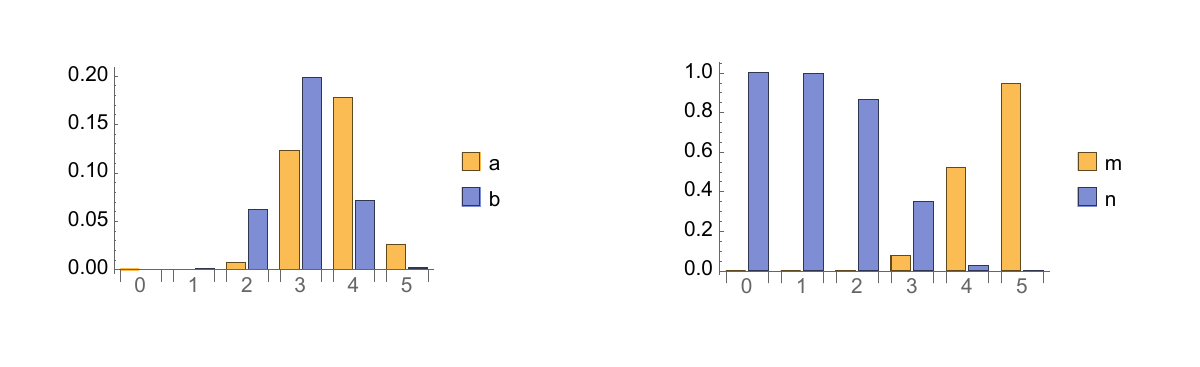}
\caption{The unique standard solution  $(a,b,m,n): \llbracket -5,5 \rrbracket \to (0,\infty)^2 \times \R^2$ for the trail game on $\llbracket -6,6\rrbracket$ with Mina margin equal to $1 + 10^{-4}$ is depicted.
The stakes $(a_i,b_i)$  and expected mean payoffs $(m_i,n_i)$ are displayed for $i \in \llbracket 0,5 \rrbracket$ in the left and right charts. Note that the leftmost displayed data is indexed by the origin: the data to the left, indexed by $\llbracket -5,-1 \rrbracket$, is visually indistinguishable from the zero-indexed data in the two displays.}\label{f.uniquenash}
\end{figure}

To evaluate the premise as just presented in the case of finite trails and~$\Z$, we begin by focusing on a Mina margin value that lies marginally in the regime of Mina's claimed dominance for  the gameboard $\llbracket -6,6 \rrbracket$
(whose set of open play is $\llbracket -5,5 \rrbracket$). The associated $\PhiMac$-transformed Mina margin map $\minamm_{6,6} \circ \PhiMac :\R \to (0,\infty)$ is depicted in Figure~\ref{f.tmmm}(left).
Select $z = 1 + 10^{-4}$, a value that lies slightly above $\lambda^{-1}$ according to Conjecture~\ref{c.lambda} (so that the purple curve in Figure~\ref{f.tmmm}(left) has turned left off the highway to cross this height). Indeed, the equation $\minamm_{6,6}(x) = z$
is found (by some trail-and-error work in Mathematica) to have a unique solution in $x \in \R$, with this solution taking the form $x = 4.04493$ up to five decimal places. 
The corresponding standard solution $(a,b,m,n): \llbracket -5,5 \rrbracket \to (0,\infty)^2 \times \R^2$
is depicted in Figure~\ref{f.uniquenash}.

The solution has central ratio $\tfrac{n_{-1} - n_0}{m_0 - m_{-1}} = \PhiMac(x)$ equal to $46538$ up to rounding error: the origin is comfortably to the left of the battlefield index. Indeed, this index lies at four, since $\phi_4 = 0.719 \cdots \in [1/3,3)$.  Mina's stake at vertex three is the greater, and she dominates staking, and mean payoffs, at vertices two and below. Were we to consider the analogous solution on longer gameboards $\llbracket -\ell,\ell \rrbracket$, $\ell > 6$, we would see that its battlefield index rises  in $\ell$, and that the region around the origin falls progressively further into the territory that Mina controls.
 In this territory, she vastly outbids Maxine, even though all stakes are tiny. The weak limit in high $\ell$ of the gameplay starting at the origin that is governed by the reverse-ordered $(a,b)$-component of the solution is likely to be a deterministic left-moving walk on~$\Z$. The conclusion may seem to be that, when $x$ exceeds $\lambda^{-1}$, the Trail of Lost Pennies on~$\Z$ has become uncompetitive because Mina's position is too strong: she should, it appears, win without expenditure. And as usual likewise for Maxine in the opposing case, when $x$ is less than $\lambda$. But care is needed in this interpretation. After all, the limit in high~$\ell$ of the stakes offered near the origin is zero for both players, and the double-zero strategy will not gratify Mina's ambition to win without cost. Overall, then, the limit from finite gameboards creates a sense of utter dominance for the player with a favourable value of the Mina margin, but our formal results are agnostic: as 
Theorem~\ref{t.nashequil.prelim} shows, there are no time-invariant Nash equilibria in the infinite game; since this is the solution concept we have studied, our results offer no guidance to Mina as she prepares to play the trail game on $\Z$
 with boundary data specifying a value of the Mina margin that exceeds~$\lambda^{-1}$. We expect ``Incentive Inch, Outcome Mile" as presented above to be validated by lengthy finite trails; and we expect the third property, concerning the brief middle phase, to be valid for $\Z$; but the properties concerning the phases of Mina's or Maxine's supposed dominance have no clear mathematical meaning for~$\Z$.


\section{The game of chicken in the Trail of Lost Pennies}
Theorems~\ref{t.nashequil.prelim} and~\ref{t.solutions} show that, for any $x \in (\lambda,\lambda^{-1})$, the trail game with boundary data $(m_{-\infty},m_\infty,n_{-\infty},n_\infty) = (0,1,x,0)$ has at least two distinct time-invariant Nash equilibria of any given integral battlefield index; if $x \in \{- \lambda,\lambda\}$,
then there is at least one such. For any $x \in [-\lambda,\lambda]$, we may thus find an element $\mc{S}_0^2 \cap \mc{N}$ of battlefield index zero. For $k \in \Z$, we denote by 
$(S_-(k),S_+(k))$ the right-shift by $k$ places of $(S_-,S_+)$. This is an element of $\mc{S}_0^2 \cap \mc{N}$ of battlefield index~$k$. Suppose that the counter begins at the origin. 
Game outcomes under the strategy pairs $(S_-(k),S_+(k))$ become more favourable to Mina, and less favourable to Maxine, as the index $k$ increases; as Theorem~\ref{t.unanimity} indicates, the probability of victory for Maxine decays rapidly as $k$ becomes positive. Suppose the game is about to begin, and players must commit to strategies. Mina may consider playing one of the strategies $S_-(k)$ for $k \in \Z$. If her opponent elects to play $S_+(k)$, then Mina would much prefer that the shared value of $k$ be positive; Maxine would naturally prefer a negative choice. But the players must consider the case that the opposing player elects a different value of~$k$. What then?  

For simplicity, consider the symmetric game where $x=1$. Suppose that $S_+ = a$ and $S_-=b$ in the usual notational abuse, where $a_i = b_{-i}$ for $i \in \Z$.
(The choice $(a,b) = (\asta(3),\bst(3))$ meets this condition, as we will see in Proposition~\ref{p.symmetric}(1).)
Let $k \in \nwithoutzero$. Suppose that Mina chooses between  the soft $S_-(-k)$ and the tough  $S_-(k)$, while Maxine elects to play either  the soft $S_+(k)$ or the tough $S_+(-k)$. By this restriction, we consider a two-person game where each player has two alternatives, and in Table~\ref{t.twobytwo}, we depict mean payoffs in a two-by-two array whose rows index Mina's choice, whose columns index Maxine's, and each of whose coordinates contains a list of Mina's and Maxine's mean payoffs when the indexing strategy pair is played. The good outcome $G$ has value $1 - \exp \{- 2^k O(1) \}$.
The medium outcome $M$ takes the form $1/2 - \exp \{- 2^k O(1) \}$.
The bad outcome~$B$ has value  $\exp \{- 2^k O(1) \}$. And the value $C$ of the catastrophic outcome is ... minus infinity!
We will illustrate how to obtain these assertions rather than present formal derivations.

The outcomes of $G$ and $B$ arise in the off-diagonal cases, where Nash equilibria are played, so that the claimed forms for $G$ and $B$ arise from Theorem~\ref{t.ajbj}
in the sense of the paragraph that follows this theorem. Consider the strategy pair $({\rm Soft}=S_-(-k),{\rm Soft}=S_+(k))$. 
At the first turn, Mina is playing $k$ units to the right of her presumed location of the battlefield vertex, as if she has as good as lost already. But Maxine is playing $k$ units to the left of where she is claiming the battlefield index to be, and also in effect nearly admits defeat. Maxine's and Mina's stakes are $a_{-k}$ and $b_k$: both very small, but equal in our special case. So the first turn victor is chosen by the outcome of a fair coin toss. And this early winner will lose even one later move only with probability $\exp \{- 2^k O(1) \}$ as the estimates in  Theorem~\ref{t.ajbj} show, because the victor's stakes rise and her opponent's fall as the counter moves closer to the victor's presumed battlefield location. 
In the case  of  $({\rm Tough}=S_-(k),{\rm Tough}=S_+(-k))$ play, a phenomenon opposite to the eventual unanimity of gameplay in Theorem~\ref{t.unanimity} occurs. The counter location at late time has law approaching an equilibrium which heavily charges the origin and a few nearby sites. When the counter moves slightly to the left of the origin, it comes closer to Maxine's presumed battlefield index at $-k$
than it does to Mina's at $k$, so that Maxine's stake rises far higher than Mina's, and the counter is restored towards the origin. An opposing leftward force naturally acts on the counter when it is in positive territory. The implicit consensus against lengthy play in a bounded region discussed around Theorem~\ref{t.unanimity} has been broken with double-tough play, and the players are trapped in an unending mutually destructive cycle.

 \begin{table}
 \begin{center}
 \begin{tabular}{| c | c | c | }
 \hline
 &  ${\rm Maxine \, \, Soft}: S_+(k)$ & ${\rm Maxine \, \, Tough}: S_+(-k)$  \\
 \hline
  ${\rm Mina \, \,Soft}: S_-(-k)$ &  M, M & B, G \\
  ${\rm Mina \, \,Tough}: S_-(k)$ &  G , B  & C,C    \\ 
 \hline
 \end{tabular}
 \medskip
 \medskip
 \medskip
 \medskip
 \caption{Mina and Maxine choose between their components in two given Nash equilibria with battlefield indices $-k$ and $k$, for some $k \in \nwithoutzero$.
 The respective mean payoffs for Mina and Maxine for each of the four strategy pairs are recorded in each entry of the $2 \times 2$ array. The possible outcomes are $G = {\rm Good}$, ${M = \rm Medium}$,
 $B = {\rm Bad}$ and $C = {\rm Catastrophic}$.}\label{t.twobytwo}
 \end{center}
 \end{table}

In the classic game of chicken~\cite[Chapter~$10$]{Poundstone2011}, two players choose between soft and tough options of swerving or driving straight. When one player drives straight and the other swerves, their payoffs are the pleasure $G$ of winning  and the annoyance $B$ of showing weakness. When both swerve, both receive an intermediate value $M$. When both drive straight, the shared outcome is a highly negative $C$ as the cars crash. We see then that the Trail of Lost Pennies on $\Z$ embeds the game of chicken. The translation symmetry of $\Z$ makes the selection of which Nash equilibrium to play a difficult choice for players who may be infinitely punished for a perhaps unintentionally tough decision. The counterpart embedding of chicken occurs in the finite trail game, where the value of $C$, while often highly negative, is finite.

\section{Play between people and algorithms}

The finite trail $\llbracket -j,k \rrbracket$---perhaps for values of $j$ and $k$ somewhere between one and five---may provide an attractive context for investigating how people or algorithms play the trail game. Given the smallness of $1 - \lambda$ and the multiplicity of Nash equilibria for many of the games with longer trails, it seems fanciful to believe that two people who play the same game repeatedly will typically adhere to such an equilibrium (at least when $j +k$ is high enough). Other strategies may seem natural.

\subsection{Cooperative behaviour}
Trust could be established during iterated play. If two players each stake zero throughout a standard symmetric  trail game on $\llbracket -k,k \rrbracket$, $k \in \nwithoutzero$, whose counter starts at zero, their running costs are zero, and their mean payoffs are one-half (this is because play ends in finite time on a finite trail; we use the $0/0 =1/2$ rule from Section~\ref{s.gamespec}). 
\subsection{Tit-for-tat}
A consistent zero strategy has the flaw of being vulnerable to exploitation, and a player in an iterated game may prefer a tit-for-tat approach: stake zero in every game until the opponent makes a positive stake; in the next game, play more aggressively; revert to playing zero if the opponent reacts modestly to the aggressive play. Of course, there are degrees of aggression that may be adopted. The iterated prisoner's dilemma is a classic example where the Nash equilibrium (which proposes uncooperative play) often predicts wrongly how people will play, and where tit-for-tat and variants thereof are commonly adopted strategies for humans~\cite{DalBoFrechette} that have been found in computer-against-computer tournaments to be effective~\cite{AxelrodHamilton}.

\subsection{The loadsamoney bully}

On a finite trail, the loadsamoney bully chooses $\e > 0$ small and consistently stakes $\e$ against an opponent who stakes zero. He plays aggressively when the opponent makes a positive stake: he may react to a stake~$a > 0$ by staking $2a$ at the next turn, for example. This player wins games against a zero-staking opponent while incurring almost no running cost. He seeks to rein financial terror on the opponent who deviates from a zero strategy, by seeming prepared to win the concerned game no matter what the cost. Hoping to create a sense of formidable financial resources, his long-term plan for the iterated game is to cow the opponent into a submissive zero strategy.

 \appendix
 
\chapter{Proofs of the further formulas}\label{s.furtherformulas}

\section{Deriving  alternative formulas for standard solutions and  Mina margins}

{\bf Proof of Theorem~\ref{t.altstand}(1)}. In the notation of Definition~\ref{d.zdefault},
$\mdef_1(x) - \mdef_0(x) = c(x) -1$ and $\mdef_0(x) - \mdef_{-1}(x) = 1$, so that $M_0(x) = c(x)$. And $\ndef_0(x) - \ndef_1(x) = x(d(x) -1)$, while $\ndef_{-1}(x) - \ndef_0(x) = x$, which leads to   $N_0(x) = xd(x)$.  The formulas for $\adef_0(x)$ and $\bdef_0(x)$ in terms of $M_0(x)$ and $N_0(x)$ yield that
$\adef_0(x) = c(x)f(x)$ and $\bdef_0(x) = xd(x)f(x)$. By~(\ref{e.remark}) and Proposition~\ref{p.shift}, we obtain Theorem~\ref{t.altstand}(1). 

{\bf (2).} By~(\ref{e.remark}), $\mdef_0(x) - \mdef_{-1}(x) = 1$ and $\ndef_{-1}(x) - \ndef_0(x) = x$, we see that 
$\mst_0(x) - \mst_{-1}(x) = Z(x)^{-1}$ and $\nst_{-1}(x) - \nst_0(x) = x Z(x)^{-1}$. Proposition~\ref{p.shift} then implies 
the special cases of
Theorem~\ref{t.altstand}(2) noted after the general cases, which are then derived by summation. 

{\bf (3).}  This follows from~(\ref{e.minammfinite}) in the form
 $\minamm_{j+1,k+1}(x) = \big( \nst_{-j-1} (x) - \nst_{k+1} (x) \big) \big( \mst_{k+1}(x) - \mst_{-j-1}(x) \big)^{-1}$, and the preceding part.

{\bf (4).} By Theorem~\ref{t.positiveabmn}(3), the numerator and denominator in the just obtained expression for $\minamm_{j,k}(x)$ converge in the limits $j,k \to \infty$
to $\nst_{-\infty} (x) - \nst_\infty (x) = \nst_{-\infty} (x)$ and $ \mst_\infty (x) - \mst_{-\infty} (x) = 1$. Since the ratio of the latter quantities equals $\minamm(x)$, 
we obtain  Theorem~\ref{t.altstand}(4) from the limiting numerator. \qed

{\em Remark.} By considering the limiting denominator in the argument for Theorem~\ref{t.altstand}(4), we obtain the identity
$$
 \sum_{i \in \Z} Z\big(s_i(x)\big)^{-1} \, = \, 1 \, \, \, \, \textrm{for $x \in (0,\infty)$} \, .
$$

\section{The symmetric game and the \amacspace equations}

Here we prove Theorem~\ref{t.symmetric}.

\subsection{The \ammacspace equations}

Let $Q \in \Z$ be a parameter; we will consider $Q =0$ and $Q=-1$. The \ammacspace system on $\Z$ with parameter~$Q$ is a set of equations in the real variables $a_i$ and $m_i$ indexed by $i \in \Z$. These equations are
\begin{align*}
  (a_i + a_{Q-i})(m_i+a_i)  & =  a_i m_{i+1} + a_{Q-i} m_{i-1} && \qquad \textrm{AM}(Q,1)  \\
    (a_i + a_{Q-i})^2  & =  a_{Q-i} \big( m_{i+1} - m_{i-1} \big)  && \qquad \textrm{AM}(Q,2) \, 
  \end{align*}
with the index $i$ running over~$\Z$.
\begin{proposition}\label{p.symmetric}
Let $\big\{ (a_i,b_i,m_i,n_i): i \in \Z \big\}$ be an \abmnmacspace solution for which $\centralratio$ is equal to~ a positive parameter $x$, and $m_{-\infty} = n_\infty = 0$.

First let $x$ equal to three. 
\begin{enumerate}
\item For $i \in \Z$, $a_i = b_{-i}$ and $m_i = n_{-i}$.  
\item The $a$- and $m$-variables solve the \ammacspace system on~$\Z$ with parameter $Q= 0$. 
\item The $a$-variables solve the \amacspace system on~$\Z$ (when we take $A_i = a_i$ for $i \in \Z$).
 \end{enumerate}
 Now instead take $x$ equal to one.  
 \begin{enumerate}
 \setcounter{enumi}{3}
\item For $i \in \Z$,  $a_i = b_{-1-i}$ and $m_i = n_{-1-i}$ for $i \in \Z$. 
\item The $a$- and $m$-variables solve the \ammacspace system on~$\Z$ with parameter $Q= -1$.
\item 
With a half-integer offset, the $a$-variables solve the \amacspace system on $\Z + 1/2$: that is, the variables $A_i = a_{i-1/2}$ for $i \in \Z + 1/2$ solve this system.
\end{enumerate}
\end{proposition}
{\bf Proof: (1).} We will argue that $(a,b,m,n)$ is equal to $\mc{R}(a,b,m,n)$, where $\mc{R}$ is the role-reversal map from Definition~\ref{d.rolereversal}.
This will suffice to prove Proposition~\ref{p.symmetric}. For $i \in \Z$, let $\hat\phi_i$ denote the $\phi_i$-value of  $\mc{R}(a,b,m,n)$. 
Note that $\phi_0 = \tfrac{n_{-1} - n_0}{m_0 - m_{-1}} = 3$ and $\hat\phi_1 = \tfrac{m_0 - m_{-1}}{n_{-1} - n_0} = 1/3$. From the latter, and Proposition~\ref{p.sminusone}, $\hat\phi_0 = s_{-1}(1/3) = s^{-1}(1/3) = 3$.
We have learnt that the value of $\centralratio$ is shared between the \abmnmacspace solutions $(a,b,m,n)$ and $\mc{R}(a,b,m,n)$ (the value equals three). 

We write $\mc{R}(a,b,m,n) = (\hat{a},\hat{b},\hat{m},\hat{n})$, and extend the hat-notation to boundary data.
We next argue that this data satisfies
\begin{equation}\label{e.mnzeromn}
 \hat{m}_{-\infty} = \hat{n}_\infty = 0 = m_{-\infty} = n_\infty \, .
\end{equation}
That $m_{-\infty} = n_\infty=0$ holds by assumption. But then we have $\hat{m}_{-\infty} = \lim_{k \to \infty} \hat{m}_{-k} = \lim_{k \to \infty} n_k = n_\infty = 0$. Likewise, $\hat{n}_\infty = 0$.

We further 
argue that
\begin{equation}\label{e.mnzerotwo}
 m_\infty = n_{-\infty} \, \, \, \, \textrm{and} \, \, \, \, \hat{m}_\infty =  \hat{n}_{-\infty}  \, .
\end{equation}
To see these, note that $\minamm(3)=1$, because this is a high $k$ consequence of Corollary~\ref{c.rkvalues}. 
But in the present case $\minamm(3)$ equals   both $n_{-\infty}/m_\infty$ and  $\hat{n}_{-\infty}/\hat{m}_\infty$; whence~(\ref{e.mnzerotwo}). 

The two preceding displays imply that $\hat{m}_\infty$, $n_{-\infty}$,  $\hat{n}_{-\infty}$ and $m_\infty$
have a shared value. Let $u \in (0,\infty)$ denote the reciprocal of this value. The quadruples $\tau_u (a,b,m,n)$ and $\tau_u \, \mc{R}(a,b,m,n)$
are both the unique standard solution with $\centralratio$ equal to three.
 The uniqueness of this solution asserted by Proposition~\ref{p.abmnclassify}(2) implies that these quadruples are equal. But then $(a,b,m,n)$ and $\mc{R}(a,b,m,n)$ have been shown to be equal, as we sought to show in deriving Proposition~\ref{p.symmetric}(1).
 
 {\bf (4).} We will argue that $(a,b,m,n)$ and $\mc{A}(a,b,m,n)$ are equal, where $\mc{A} = \mc{S}_1 \circ \mc{R}$ is the left shift (by one place) of the role reversal map. This suffices to prove
  Proposition~\ref{p.symmetric}(4). With a tilde notation for $\mc{A}(a,b,m,n)$ objects, we have that $\tilde\phi_0 = 1/\phi_0 =1$. Thus,  $(a,b,m,n)$ and $\mc{A}(a,b,m,n)$ are \abmnmacspace solutions, each with $\centralratio = 1$. Since $\minamm(1)=1$ by Corollary~\ref{c.rkvalues} in a high~$k$ limit, we find, similarly to above, that  the values of $\tilde{m}_\infty$, $n_{-\infty}$,  $\tilde{n}_{-\infty}$ and $m_\infty$ coincide. 
  The proof continues to mimic that of  Proposition~\ref{p.symmetric}(1).
  
  {\bf (2).} Given the first part, \abmnmac$(1)$, or \abmnmac$(2)$, implies \ammac$(0,1)$, while \abmnmac$(3)$, or \abmnmac$(4)$, implies \ammac$(0,2)$.
  
  {\bf (5).} Similarly, given the fourth part, and with   \ammac$(-1,1)$ and \ammac$(-1,2)$ being implied.

  {\bf (3,6).}  The \ammacspace system (for general $Q$) satisfies
$$
  (a_i + a_{Q-i})(m_i+a_i)  \, = \, a_i \bigg( \frac{(a_i + a_{Q-i})^2}{a_{Q-i}} + m_{i-1} \bigg) + a_{Q-i} m_{i-1}  \, .
$$
This identity is seen to be equivalent to each of the following equalities in turn:
$$
(a_i + a_{Q-i})(m_i+a_i - m_{i-1}) a_{Q-i} = a_i(a_i + a_{Q-i})^2 \, ;
$$
 $a_i(a_i + a_{Q-i}) = a_{Q-i} (m_i + a_i - m_{i-1})$; and
 $m_i - m_{i-1} = a_i^2 a^{-1}_{Q-i}$. By summing the last of these over two consecutive indices, we find that
$$
 m_{i+1} - m_{i-1} \, = \, \frac{a_i^2}{a_{Q-i}} + \frac{a_{i+1}^2}{a_{Q-i-1}} \, .
$$
Using \ammac$(Q,2)$ here, we obtain    
$$
\frac{(a_i + a_{Q-i})^2}{a_{Q-i}} \, = \, \frac{a_i^2}{a_{Q-i}} + \frac{a_{i+1}^2}{a_{Q-i-1}} \, . 
$$ 
Rearranging, and cancelling $a_{Q - i} \not= 0$, leads to
 $$
 a_{Q-i-1}(2a_i + a_{Q-i}) = a_{i+1}^2
 $$
 for $i \in \Z$.
 When $Q=0$, the $a$-variables are thus seen to satisfy the \amacspace system on~$\Z$, as Proposition~\ref{p.symmetric}(3) asserts. 
 When $Q=-1$, the $A$-variables $\big\{ A_j: j \in \Z +1/2 \big\}$, $A_j = a_{j-1/2}$, solve the system $A_{-i -3/2}(2 A_{i+1/2} + A_{i - 1/2}) = A_{i + 3/2}^2$ as $i$ varies over $\Z$. Setting $k = i + 1/2$,
 we find that 
 $$
 A_{-k-1}(2A_k + A_{-k}) = A_{k+1}^2
 $$ 
 for $k \in \Z + 1/2$. Thus, the variables $\big\{ A_j: j \in \Z +1/2 \big\}$
 solve
  the \amacspace system on~$\Z + 1/2$, as we sought to show in proving Proposition~\ref{p.symmetric}(6).
 This completes the proof of Proposition~\ref{p.symmetric}. \qed

\subsection{Deriving Theorem~\ref{t.symmetric}}

We are ready to derive the remaining result.

{\bf Proof of Theorem~\ref{t.symmetric}(1).} If $\big\{ a_i: i \in \Z \big\}$ solves the \amacspace system on $\Z$ with $a_0 = 1$, then $\lambda a_i$ solves this system with $\lambda a_i = \lambda$. Thus it is without loss of generality that we take $\lambda = 1$.

Given that $a_0 = 1$, the $0$- and $-1$-indexed equations in the \amacspace system on $\Z$ state that $3 a_{-1} = a_1^2$ and $2a_{-1} + a_1 = 1$. These uniquely determine positive $a_{-1}$ and $a_1$.
We may group the remaining \amacspace equations into pairs, indexed by $i=k$ and $i = -k-1$ for $k \geq 1$, with the generic pair taking the form
\begin{equation}\label{e.inputoutput}
 O_- (2I_+ + I_-) = O_+^2
 \, \, \, \, \textrm{and} \, \, \, \, I_+ (2O_- + O_+) = I_-^2
\end{equation}
with input pair $(I-,I_+) = (a_{-k},a_k)$ and output pair $(O_-,O_+) = (a_{-k-1},a_{k+1})$. These equation pairs may be iteratively solved, for $k \in \nwithoutzero$ rising, each step producing a unique pair of positive values in the output. The outcome is the unique solution  $\big\{ a_i: i \in \Z \big\}$ to the \amacspace system subject to $a_0 =1$.

{\bf (2).} Similarly to the first part, we may restrict to the case where $A_{-1/2} = 1$. The $-1/2$-indexed equation then reads $2+ A_{1/2} = A_{1/2}^2$ and produces a unique positive value for $A_{1/2}$.
Inductively we determine the $A_i$-values for $\vert i \vert \leq k -1/2$. The equation pair indexed by $i = -k + 1/2$ and $i = k - 1/2$
takes the form~(\ref{e.inputoutput}) with $(I_-,I_+) = (A_{-k+1/2},A_{k - 1/2})$ and $(O_-,O_+) = A_{-k-1/2,k+1/2}$, with a unique positive pair of outputs given the input pair. In this way, we construct the unique positive solution $\big\{ A_i: i \in \Z + 1/2 \big\}$ of the \amacspace system on $\Z+1/2$ subject to $A_{-1/2} = 1$.

{\bf (3).}
Since $\minamm(3)=1$, Theorem~\ref{t.nashabmn} implies that each strategy pair of the form 
\begin{equation}\label{e.nashcollection}
\lambda \cdot \big( \bst_i(3), \asta_i(3): i \in \Z \big) \, \, \,  \textrm{for} \, \, \, \lambda \in (0,\infty) 
\end{equation}
is a time-invariant Nash equilibrium in the symmetric game.
Note that $\big\{ \big(\bst_i(3), \asta_i(3) \big): i \in \Z \big\}$ equals $\big\{ \big( \asta_{-i}(3), \asta_i(3) \big): i \in \Z \big\}$ by Proposition~\ref{p.symmetric}(1), and that $\big\{ \asta_i(3): i \in \Z \big\}$ solves the \amacspace system on $\Z$ by Proposition~\ref{p.symmetric}(3). By Theorem~\ref{t.symmetric}(1), we see that the collection~(\ref{e.nashcollection}) takes the form 
$$
\Big\{ \lambda \cdot \big( a_{-i}, a_i : i \in \Z \big) : \lambda \in (0,\infty) \Big\} \, ,
$$
where $\big\{ a_i: i \in \Z \big\}$ is the unique solution of the \amacspace system on $\Z$ subject to $a_0 = 1$. The sequence shift operation leaves invariant the space of time-invariant Nash equilibria for a given game, because it preserves boundary data. Since each of the last displayed maps is an equilibrium in the symmetric game, we have found that each element of $S_1$ is also such an equilibrium.

Since $\minamm(1) =1$, similar reasoning, which invokes Proposition~\ref{p.symmetric}(4,6), shows that every element of
$$
\Big\{ \lambda \cdot \big( A_{-i-1/2}, A_{i+1/2}: i \in \Z \big) : \lambda \in (0,\infty) \Big\}
$$
is a time-invariant Nash equilibrium for the symmetric game. All shifts of such equilibria are also equilibria; thus, we see that every element of $S_2$ is a time-invariant Nash equilibrium in the symmetric game. 

It remains to argue that there is no coincidence between two elements specified in the definitions of the sets $S_1$ and $S_2$. Our task is to extract from a strategy pair of the form 
$\lambda \cdot \big( \bst_{i+k}(z), \asta_{i+k}(z) \big)$ that is specified by the data $(\lambda,k,z) \in (0,\infty) \times \Z \times \{ 1,3\}$ the value that this triple takes in the indexing set. Proposition~\ref{p.shift} permits us to write the pair in the form $\lambda \cdot \big( \bst_i(s_k(z)), \asta_i(s_k(z)) \big)$. By Theorem~\ref{t.nashabmn}, we may equally attempt to extract the data 
$(\lambda,k,z)$ from the \abmnmacspace solution  $\lambda \cdot \big( \asta_i(s_k(z)), \bst_i(s_k(z)),  \mst_i(s_k(z)), \nst_i(s_k(z)) \big)$. 
The value of $\centralratio$ for the latter solution lies in the set $\big\{ s_k(z): k \in \Z \, , \, z \in \{ 1,3\} \big\}$. We may extract $k$ since it is the unique integer $i$ for which $s_{-i}(\centralratio) \in(1/3,3]$.
The value of $s_{-k}(\centralratio)$ lies in $\{1,3\}$ and equals $z$.  Since $\mst_\infty - \mst_{-\infty} = 1$, the value of $\lambda$ may be identified as equal to $\mst_\infty(s_k(z)) - \mst_{-\infty}(s_k(z))$.
We have confirmed that there is no coincidence between the elements in the union of $S_1$ and $S_2$.

{\bf (4).} We suppose that $Q(1)=2$. Since $\minamm(1)= \minamm(3) = 1$, there are no further solutions to $\minamm(x) = 1$ among $x \in (1/3,3]$. Theorem~\ref{t.nashequil} implies that
the set of time-invariant Nash equilibria in the game $\standard(1)$ equals the set of pairs of the form $\big( \bst_{i+k}(z), \asta_{i+k}(z): i \in Z \big)$ for $z \in \{ 1,3\}$ and $k \in \Z$.
The game $\standard(1)$ is the symmetric game where $m_\infty =1$. For the general cases, where $m_\infty = \lambda \in (0,\infty)$, the set of Nash equilibria is given by the set of pairs formed by multiplying the indicated expression by $\lambda$. In other words, the set in question is simply the union of $S_1$ and~$S_2$. \qed

\chapter{Glossary of notation}\label{s.glossary}

The article has quite a lot of notation. 
The first table offers principal examples concerning strategy pairs and the associated gameplay. The second concerns situations where an \abmnmacspace solution is supposed given. 
A summarizing phrase is provided for each item, as well as the page number at which the concept is introduced.

\vspace{2mm}

\large{\bf Strategies, gameplay and equilibria}

\vspace{-1mm}

\small{
\bigskip
\def\qq{&}
\begin{center}
\halign{
#\quad\hfill&#\quad\hfill&\quad\hfill#\cr
$\N$, $\nwithoutzero$ \qq the natural numbers, with and without zero \fff{withandwithout}
$\llbracket i,j \rrbracket$ \qq the integer interval $\Z \cap [i,j]$  \fff{intint}
$\mc{S}$ \qq the space of strategies $S:\Z \times \nwithoutzero \to [0,\infty)$ \fff{strategy}
$\tis$ \qq the space of time-invariant strategies \fff{timeinvariant}
$(S_-,S_+) \in \tis^2$ \qq typical notation for a time-invariant strategy pair \fff{typicalnotation}
$S_- = b$ and $S_+ = a$ \qq a notational abuse for $(S_-,S_+) \in \tis^2$; $a_i$ and $b_i$ are stakes at site~$i$ \fff{notab}
$\mc{N}$ \qq the set of Nash equilibria \fff{nash}
$\PP^i_{S_-,S_+}$ and $\E^i_{S_-,S_+}$ \qq  the law and expectation of play from $i \in \Z$ under $(S_-,S_+) \in \mc{S}^2$  \fff{lawexpect}
$\PP^{u,\ell}_{S_-,S_+}$ and $\E^{u,\ell}_{S_-,S_+}$ \qq counterpart objects for play from $u$ with a delayed start-time $\ell \in \nwithoutzero$ \fff{delayed}
$X$ \qq counter evolution $X: \N \to \Z$ \fff{counter}
$P_-$, $T_-$ and $C_-(t)$ \qq Mina's total, and terminal, receipts and time-$t \in \nwithoutzero$ running cost \fff{receiptmina}
$P_+$, $T_+$ and $C_+(t)$ \qq Maxine's counterparts \fff{receiptmaxine}
$m_i$ and $n_i$ \qq Maxine and Mina's mean total receipt for play from $i$ under $(S_-,S_+) \in \mc{N}$ \fff{mini}
$a_i$ and $b_i$ \qq the players' stakes in this scenario  \fff{aiai}
$E$, $E_-$, $E_+$ \qq events of escape $\vert X \vert \to \infty$, left escape $X \to -\infty$ and right escape $X \to \infty$ \fff{escape}
 $\mc{I}$ \qq the idle zone, for $(S_-,S_+) \in \tis^2$, where neither player offers a stake \fff{idle}
 wide \qq an element of $\tis$ of infinite support \fff{wide}
}\end{center}}

\vspace{-4mm}

\large{\bf \abmnmacspace solutions}

\vspace{-1mm}

\small{
\bigskip
\def\qq{&}
\begin{center}
\halign{
#\quad\hfill&#\quad\hfill&\quad\hfill#\cr
$(a,b,m,n):\Z \to (0,\infty)^4$ \qq a positive \abmnmacspace solution \fff{positiveabmn}
$\big( m_{-\infty},m_\infty,n_{-\infty},n_\infty \big)$ \qq the boundary data for such a solution \fff{boundarydata}
standard and default \qq distinguished elements in \abmnmacspace solution equivalence classes  \fff{standarddefault} 
$\phi_i$ \qq the quantity $\frac{n_{i-1} - n_i}{m_i - m_{i-1}}$ \fff{phi}
battlefield index \qq the unique index $i \in \Z$ such that $\phi_i \in (1/3,3]$ \fff{battlefield}
$s:(0,\infty) \to (0,\infty)$ \qq the map that sends $\phi_0$ to $\phi_1$, for any positive \abmnmacspace solution \fff{smap}
central ratio \qq a name for $\phi_0 = \frac{n_{-1} - n_0}{m_0 - m_{-1}}$ \fff{centralratio}
Mina margin \qq the value $\frac{n_{-\infty} - n_\infty}{m_\infty - m_{-\infty}}$ \fff{minamargin}
$\mc{M}:(0,\infty) \to (0,\infty)$ \qq the Mina margin map, sending $\phi_0$ to $\tfrac{n_{-\infty} - n_\infty}{m_\infty - m_{-\infty}}$ \fff{MMM}
$\mc{M}_{j+1,k+1}:(0,\infty) \to (0,\infty)$ \qq the counterpart map  $\phi_0 \to \tfrac{n_{-j-1} - n_{k+1}}{m_{k+1} - m_{-j-1}}$  for trail $\llbracket -j-1,k+1 \rrbracket$ \fff{MMM.finite}
 $M_i$, $N_i$,
$\beta_i$, $\gamma_i$, $c_i$, $\delta_i$, $d_i$ \qq statistics local to site $i$ associated to an \abmnmacspace solution  \fff{gammaess}
$\beta$, $\gamma$, $c$, $\delta$, $d$ \qq functions $(0,\infty) \to (0,\infty)$ such that $*(\phi_i) = *_i$ \fff{alphagamma}
}\end{center}}



\bibliographystyle{plain}

\bibliography{stake}

\end{document}